\documentclass{amsart}
\usepackage{amsmath}
\usepackage{graphicx}
\usepackage{amsfonts}
\usepackage{amssymb}

\newtheorem{theorem}{Theorem}
\theoremstyle{plain}

\newtheorem{corollary}{Corollary}

\newtheorem{definition}{Definition}
\newtheorem{example}{Example}

\newtheorem{lemma}{Lemma}

\newtheorem{proposition}{Proposition}
\newtheorem{remark}{Remark}

\numberwithin{equation}{section}

\begin{document}
\title{Matrads, Biassociahedra, and $A_{\infty }$-Bialgebras}
\author{Samson Saneblidze$^{1}$}
\address{A. Razmadze Mathematical Institute\\
Department of Geometry and Topology \\
M. Aleksidze St. 1 \\
0193 Tbilisi, Georgia}
\email{sane@rmi.acnet.ge}
\author{Ronald Umble $^{2}$}
\address{Department of Mathematics\\
Millersville University of Pennsylvania\\
Millersville, PA. 17551}
\email{ron.umble@millersville.edu}
\thanks{$^{1} $ This research described in this publication was made
possible in part by Award No. GM1-2083 of the U.S. Civilian Research and
Development Foundation for the Independent States of the Former Soviet Union
(CRDF) and by Award No. 99-00817 of INTAS}
\thanks{$^{2}$ This research funded in part by a Millersville University
faculty research grant.}
\date{November 30, 2010}
\subjclass{Primary 55P35, 55P99; Secondary 52B05}
\keywords{$A_{\infty }$-bialgebra, operad, matrad, permutahedron,
biassociahedron}

\begin{abstract}
We introduce the notion of a matrad $M=\left\{ M_{n,m}\right\} $ whose
submodules $M_{\ast ,1}$ and $M_{1,\ast }$ are non-$\Sigma $ operads. We
define the free matrad $\mathcal{H}_{\infty }$ generated by a singleton $%
\theta _{m}^{n}$ in each bidegree $\left( m,n\right) $ and realize $\mathcal{%
H}_{\infty }$ as the cellular chains on a new family of polytopes $\left\{
KK_{n,m}=KK_{m,n}\right\} ,$ called \emph{biassociahedra}, of which $%
KK_{n,1}=KK_{1,n}$ is the associahedron $K_{n}.$ We construct the universal
enveloping functor from matrads to PROPs and define an $A_{\infty }$%
-bialgebra as an algebra over $\mathcal{H}_{\infty }.\vspace*{-0.2in}%
\vspace*{-0.1in}\vspace*{-0.1in}\vspace*{-0.2in}$
\end{abstract}

\maketitle
\tableofcontents

\section{Introduction}

Let $H$ be a DG\ module over a commutative ring $R$ with identity. In
\cite{SU3}, we defined an $A_{\infty}$-bialgebra structure on $H$ in terms of
a square-zero $\circledcirc$-product on the universal PROP $U_{H}=End\left(
TH\right)  $. In this paper we take an alternative point-of-view motivated by
three classical constructions: First, chain maps $Ass\rightarrow U_{H}$ and
$\mathcal{A}_{\infty}\rightarrow U_{H}$ in the category of non-$\Sigma$
operads define strictly (co)associative and $A_{\infty}$-(co)algebra
structures on $H;$ second, there is a minimal resolution of operads
$\mathcal{A}_{\infty}\rightarrow Ass;$ and third, $\mathcal{A}_{\infty}$ is
realized by the cellular chains on the Stasheff associahedra $K=\sqcup K_{n}$
\cite{May}, \cite{MSS}. It is natural, therefore, to envision a category in
which analogs of $Ass$ and $\mathcal{A}_{\infty}$ define strictly
biassociative and $A_{\infty}$-bialgebra structures on $H$.

In this paper we introduce the notion of a \emph{matrad} whose distinguished
objects ${\mathcal{H}}$ and ${\mathcal{H}}_{\infty}$ play the role of $Ass$
and $\mathcal{A}_{\infty}$. But unlike the operadic case, freeness
considerations are subtle since biassociative bialgebras cannot be
simultaneously free and cofree. Although ${\mathcal{H}}$ and ${\mathcal{H}%
}_{\infty}$ are generated by singletons in each bidegree, those in
${\mathcal{H}}$ are module generators while those in ${\mathcal{H}}_{\infty}$
are matrad generators. Indeed, as a non-free matrad, $\mathcal{H}$ has two
matrad generators and ${\mathcal{H}}_{\infty}$ is its minimal resolution. Thus
${\mathcal{H}}$ and ${\mathcal{H}}_{\infty}$ are the smallest possible
constructions that control biassociative bialgebras structures and their
homotopy versions (c.f. \cite{Markl2}, \cite{borya}, \cite{Pirashvili}).

Given a finite sequence $\mathbf{x}$ in $\mathbb{N}$, let $\left\vert
\mathbf{x}\right\vert =\sum x_{i}.$ A matrad $(M,\gamma)$ is a bigraded module
$M=\left\{  M_{n,m}\right\}  _{m,n\geq1}$ together with a family of structure
maps
\[
\gamma=\left\{  \gamma_{\mathbf{x}}^{\mathbf{y}}:\Gamma_{p}^{\mathbf{y}%
}(M)\otimes\Gamma_{\mathbf{x}}^{q}(M)\rightarrow M_{\left\vert \mathbf{y}%
\right\vert ,\left\vert \mathbf{x}\right\vert }\right\}  _{\mathbf{x\times
y}\in\mathbb{N}^{1\times p}\times\mathbb{N}^{q\times1}}%
\]
defined on certain submodules $\Gamma_{p}^{\mathbf{y}}(M)\otimes
\Gamma_{\mathbf{x}}^{q}(M)\subseteq\bigotimes_{j=1}^{q}M_{y_{j},p}%
\otimes\bigotimes_{i=1}^{p}M_{q,x_{i}}$ and generated by certain components of
the S-U diagonal $\Delta_{P}$ on permutahedra \cite{SU2}; its substructures
$\left(  \Gamma_{1}^{\mathbf{y}}(M),\gamma\right)  $ and $\left(
\Gamma_{\mathbf{x}}^{1}(M),\gamma\right)  $ are non-$\Sigma$ operads. We think
of monomials in $\Gamma_{\mathbf{x}}^{q}(M)$ as $p$-fold tensor products of
multilinear operations, each with $q$ outputs, and monomials in $\Gamma
_{p}^{\mathbf{y}}(M)$ as $q$-fold tensor products of multilinear operations,
each with $p$ inputs.

A general $\operatorname*{PROP}$, and $U_{H}$ in particular, admits a
canonical matrad structure and chain maps ${\mathcal{H}}\rightarrow U_{H}$ and
${\mathcal{H}}_{\infty}\rightarrow U_{H}$ in the category of matrads define
biassociative bialgebra and $A_{\infty}$-bialgebra structures on $H.$
Furthermore, $\mathcal{H}_{\infty}$ is realized by the cellular chains on a
new family of polytopes $KK=\bigsqcup\nolimits_{m,n\geq1}KK_{n,m},$ called
$\emph{biassociahedra,}$ of which $KK_{m,n}=KK_{n,m},$ and $KK_{1,n}$ is the
Stasheff associahedron $K_{n}.$ We identify the top dimensional cell of
$KK_{n,m}$ with the indecomposable matrad generator $\theta_{m}^{n}$
represented graphically by a \textquotedblleft double
corolla\textquotedblright\ with data flowing upward through $m$ input and $n$
output channels (see Figure 1). The action of the matrad product $\gamma$ on
the submodule $\mathbf{\Theta}=\left\langle \theta_{m}^{n}\right\rangle
_{m,n\geq1}$ generates ${\mathcal{H}}_{\infty};$ we define a differential
$\partial$ of degree $-1$ on $\theta_{m}^{n}$ and extend it as a derivation 
of $\gamma$ (as in Example 6.11).
\vspace{-0.1in}
\begin{align*}
&
\begin{array}
[c]{c}%
\theta_{m}^{n}\leftrightarrow
\raisebox{-0.0553in}{\includegraphics[
trim=0.000000in 0.140808in 0.000000in 0.254377in,
height=0.1833in,
width=0.1842in
]%
{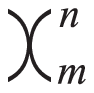}%
}
\text{ }:=
\end{array}%
\begin{array}
[c]{c}%
n\\
{\includegraphics[
height=0.4445in,
width=0.4506in
]%
{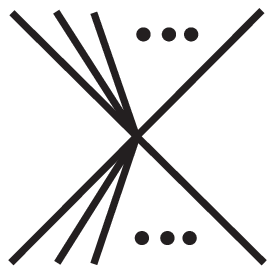}%
}
\\
m
\end{array}
\\
&  \hspace*{0.5in}\text{Figure 1.}%
\end{align*}

Among the various attempts to construct homotopy versions of bialgebras,
recent independent results of Markl and Shoikhet are related to ours through
the theory of PROPs: In characteristic 0, the low-order relations in Markl's
version of a homotopy bialgebra \cite{Markl2} agree with our $A_{\infty}%
$-bialgebra relations and Shoikhet's composition product on the universal
preCROC \cite{borya} agrees with our prematrad operation on $U_{H}$. Thus we
construct a functor from matrads to PROPs called the \emph{universal
enveloping functor}.

The paper is organized as follows: In Section 2 we construct the
biassociahedra $KK_{n,m}$ in the range $m+n\leq6$. These polytopes have a
simple description in terms of the S-U\ diagonal $\Delta_{K}$ on associahedra
\cite{SU2} and demonstrate the general construction while avoiding the
complicating subtleties. In Section 3 we discuss various submodules of $TTM$
(the tensor module of $TM$) that model the geometry of our construction. We
introduce the notion of a prematrad in Section 4, the notion of a
$k$-approximation in Section 5, and a matrad in Section 6. We construct the
posets $\mathcal{PP}$ and $\mathcal{KK}$ in Section 7, introduce the notion of
the combinatorial join of permutahedra in Section 8, and construct $PP$ and
$KK$ as geometric realizations of $\mathcal{PP}$ and $\mathcal{KK}$ in Section
9. We identify the cellular chains of $KK$ with the $A_{\infty}$-bialgebra
matrad ${\mathcal{H}}_{\infty}$ and prove that the restriction of the free
resolution of prematrads $\rho^{^{{\operatorname*{pre}}}}%
:F^{^{{\operatorname*{pre}}}}(\Theta)\rightarrow{\mathcal{H}}$ to
${\mathcal{H}}_{\infty}$ is a free resolution in the category of matrads.

\section{Low Dimensional Biassociahedra}

Our construction of the biassociahedra $\left\{  KK_{n,m}\right\}  $ in the
range $1\leq m,n\leq4$ and $m+n\leq6$ is controlled by the S-U diagonal on
associahedra $K$; the polytope $KK_{n,m}$ is identical to $B_{m}^{n}$
constructed by M. Markl in \cite{Markl2}. In the course of his construction,
Markl makes arbitrary choices, which correspond to choices we made when
constructing the S-U diagonal $\Delta_{K}$. So for us, all choices in our
construction were made a priori once and for all.

\subsection{The Fraction Product}

Let $\Theta=\left\langle \theta_{m}^{n}\neq0\mid\theta_{1}^{1}=\mathbf{1}%
\right\rangle _{m,n\geq1}$ and let $M=\left\{  M_{n,m}\right\}  _{m,n\geq1}$
be the free $\operatorname*{PROP}$ generated by $\Theta$. For simplicity, we
assume that $M$ is a free bigraded $\mathbb{Z}_{2}$-module; sign
considerations that arise over a general ground ring will be addressed in
subsequent sections. For $p,q\geq1,$ let $\mathbf{x\times y}=\left(
x_{1},\ldots,x_{p}\right)  \times(y_{1},\ldots,$ $y_{q})\in\mathbb{N}%
^{p}\times\mathbb{N}^{q}.$ In \cite{Markl2}, M. Markl defined the
submodule\emph{\ }$S$ of \emph{special elements} in $M$ whose additive
generators are monomials $\alpha_{\mathbf{x}}^{\mathbf{y}}$ expressed as
\textquotedblleft elementary fractions\textquotedblright\ of the form%
\begin{equation}
\alpha_{\mathbf{x}}^{\mathbf{y}}=\left(  \alpha_{p}^{y_{1}}\cdots\alpha
_{p}^{y_{q}}\right)  /\left(  \alpha_{x_{1}}^{q}\cdots\alpha_{x_{p}}%
^{q}\right)  ,\label{fraction}%
\end{equation}
where $\alpha_{x_{i}}^{q}$ and $\alpha_{p}^{y_{j}}$ are additive generators of
$S$ and the $j^{th}$ output of $\alpha_{x_{i}}^{q}$ is linked to the $i^{th}$
input of $\alpha_{p}^{y_{j}}$ (juxtaposition in the numerator and denominator
denotes tensor product). Thus $\dim\alpha_{\mathbf{x}}^{\mathbf{y}}=\sum
_{i,j}\dim\alpha_{x_{i}}^{q}+\dim\alpha_{p}^{y_{j}},$ and $\alpha_{\mathbf{x}%
}^{\mathbf{y}}$ is represented graphically by a connected non-planar graph
under the identification $\theta_{m}^{n}\leftrightarrow$
\raisebox{-0.0553in}{\includegraphics[
trim=0.000000in 0.140808in 0.000000in 0.254377in,
height=0.1833in,
width=0.1842in
]%
{dblvee_nm.eps}%
}
\ (see Example \ref{alpha_3,3}). We refer to $\mathbf{x}$ and $\mathbf{y}$ as
the \emph{leaf sequences} of $\alpha_{\mathbf{x}}^{\mathbf{y}}$.

Let $TM$ denote the tensor module of $M.$ All elementary fractions define a
non-associative \emph{fraction product \ }$/:TM\otimes TM\rightarrow S.$ For
example, in the iterated fraction%
\[
A/B/C=\curlywedge/\left(  \curlywedge\text{ }\mathbf{1}\right)  /(\curlywedge
\text{ }\mathbf{1}\text{ }\mathbf{1})
\]
with $\curlywedge\leftrightarrow\theta_{2}^{1},$ we have $\left(  A/B\right)
/C\neq0$ and $A/\left(  B/C\right)  =0.$ Notationally, let $\mathbf{M}%
_{\mathbf{x}}^{q}=M_{q,x_{1}}\otimes\cdots\otimes M_{q,x_{p}}$ and
$\mathbf{M}_{p}^{\mathbf{y}}=M_{y_{1},p}\otimes\cdots\otimes M_{y_{q},p};$
then the fraction product (\ref{fraction}) can be expressed in terms of our
\emph{prematrad product }$\gamma_{\mathbf{x}}^{\mathbf{y}}:\mathbf{M}%
_{p}^{\mathbf{y}}\otimes\mathbf{M}_{\mathbf{x}}^{q}\rightarrow M_{\left\vert
\mathbf{y}\right\vert ,\left\vert \mathbf{x}\right\vert }$ as%
\[
\left(  \alpha_{p}^{y_{1}}\cdots\alpha_{p}^{y_{q}}\right)  /\left(
\alpha_{x_{1}}^{q}\cdots\alpha_{x_{p}}^{q}\right)  =\gamma_{\mathbf{x}%
}^{\mathbf{y}}\left(  \alpha_{p}^{y_{1}}\cdots\alpha_{p}^{y_{q}};\alpha
_{x_{1}}^{q}\cdots\alpha_{x_{p}}^{q}\right)  .
\]
Note that the free action of $\gamma$ on up-rooted $m$-leaf corollas
$\curlywedge_{m}$ (or on down-rooted $n$-leaf corollas $\curlyvee^{n}$)
generates all planar rooted trees with levels (PLTs), and the projection to
planar rooted trees (PRTs) by forgetting levels induces the standard operadic product.

The matrad $\mathcal{H}_{\infty}$ is a proper submodule of $S$ and the matrad
product $\gamma$ is the fraction product restricted to $\mathcal{H}_{\infty}$.
Thus we may regard the $A_{\infty}$-operad $\mathcal{A}_{\infty}$ as either%
\begin{align}
\mathcal{H}_{1,\ast} &  =\left(  \left.  \left\langle \theta_{m}^{1}\mid
m\geq1\right\rangle \right/  \sim,\text{ }\partial\right)  ,\text{ where
}\partial\left(  \theta_{m}^{1}\right)  =\sum_{\substack{\alpha_{m}^{1}%
\in\mathcal{H}_{1,\ast}\cdot\mathcal{H}_{1,\ast} \\\dim\alpha_{m}^{1}=m-3
}}\alpha_{m}^{1}\text{ \ or}\label{operads}\\
\mathcal{H}_{\ast,1} &  =\left(  \left.  \left\langle \theta_{1}^{n}\mid
n\geq1\right\rangle \right/  \sim,\text{ }\partial\right)  ,\text{ where
}\partial\left(  \theta_{1}^{n}\right)  =\sum_{\substack{\alpha_{1}^{n}%
\in\mathcal{H}_{\ast,1}\cdot\mathcal{H}_{\ast,1} \\\dim\alpha_{1}^{n}=n-3
}}\alpha_{1}^{n}.\nonumber
\end{align}

\subsection{Low Dimensional Matrad Products\label{low-products}}

Let us construct $\{\mathcal{H}_{n,m}\}_{m+n\leq6}$ inductively as stage
$\mathcal{F}_{6}$ of the increasing filtration $\mathcal{F}_{k}=\bigoplus
{}_{m+n\leq k}\mathcal{H}_{n,m}$. Our construction is controlled by the S-U
diagonal $\Delta_{K}$ on cellular chains of the associahedra $K$ (see
Subsection \ref{S-U} and \cite{SU2}), which in the range of dimensions
considered here is given by%
\begin{align*}
\Delta_{K}(\text{ }
\raisebox{-0.0415in}{\includegraphics[
height=0.1885in,
width=0.1929in
]%
{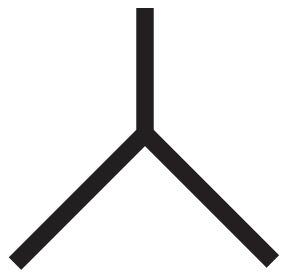}%
}
\text{ })\text{ } &  =\hspace*{0.07in}
\raisebox{-0.0415in}{\includegraphics[
height=0.1885in,
width=0.1929in
]%
{T2.eps}%
}
\text{ }\otimes\text{\ }
\raisebox{-0.0415in}{\includegraphics[
height=0.1885in,
width=0.1929in
]%
{T2.eps}%
}
\text{ ,}\\
\Delta_{K}(\text{ }
\raisebox{-0.0346in}{\includegraphics[
height=0.2101in,
width=0.2145in
]%
{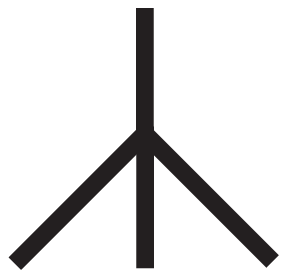}%
}
\text{ })\text{ } &  =\text{ }
\raisebox{-0.0415in}{\includegraphics[
height=0.2101in,
width=0.2145in
]%
{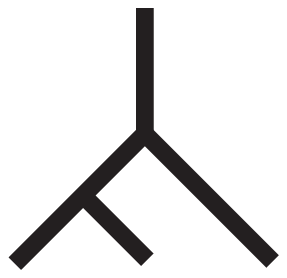}%
}
\text{ }\otimes\text{ }
\raisebox{-0.0346in}{\includegraphics[
height=0.2101in,
width=0.2145in
]%
{T3.eps}%
}
\text{ \ }+\text{ \ }
\raisebox{-0.0346in}{\includegraphics[
height=0.2101in,
width=0.2145in
]%
{T3.eps}%
}
\text{ }\otimes\text{ }
\raisebox{-0.0346in}{\includegraphics[
height=0.2101in,
width=0.2145in
]%
{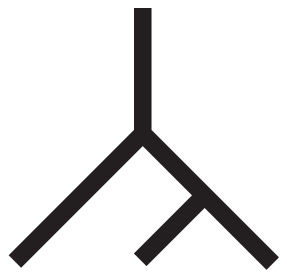}%
}
\text{ , \ and}\\
\Delta_{K}(\,
\raisebox{-0.0692in}{\includegraphics[
height=0.2127in,
width=0.2179in
]%
{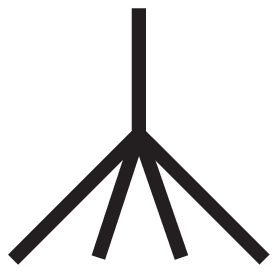}%
}
\text{ })\text{ } &  =\text{\ }
\raisebox{-0.0761in}{\includegraphics[
height=0.2145in,
width=0.2179in
]%
{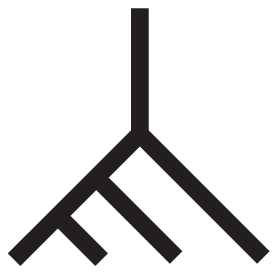}%
}
\text{ }\otimes\text{\ }
\raisebox{-0.0692in}{\includegraphics[
height=0.2127in,
width=0.2179in
]%
{T4.eps}%
}
\text{ \ }+\text{ \ }
\raisebox{-0.0692in}{\includegraphics[
height=0.2127in,
width=0.2179in
]%
{T4.eps}%
}
\text{ }\otimes\text{\ }
\raisebox{-0.0692in}{\includegraphics[
height=0.2145in,
width=0.2179in
]%
{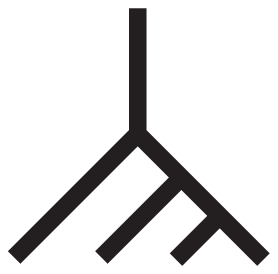}%
}
\text{ \ }+\text{ \
\raisebox{-0.0692in}{\includegraphics[
height=0.2145in,
width=0.2179in
]%
{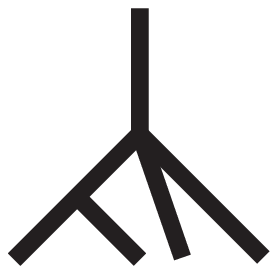}%
}
}\otimes\text{ }
\raisebox{-0.0623in}{\includegraphics[
height=0.2145in,
width=0.2179in
]%
{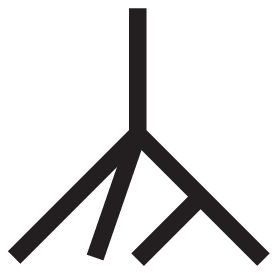}%
}
\\
&  \hspace*{-0.2in}\hspace*{0.04in}+\text{ \
\raisebox{-0.0553in}{\includegraphics[
height=0.2145in,
width=0.2179in
]%
{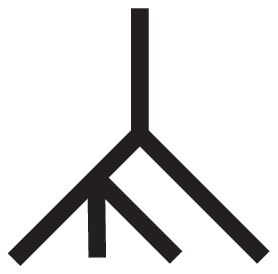}%
}
\ }\otimes\text{
\raisebox{-0.0553in}{\includegraphics[
height=0.2145in,
width=0.2179in
]%
{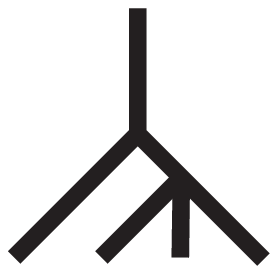}%
}
\ }+\text{ \
\raisebox{-0.0553in}{\includegraphics[
height=0.2145in,
width=0.2179in
]%
{T31-2.eps}%
}
}\otimes\text{\
\raisebox{-0.0553in}{\includegraphics[
height=0.2145in,
width=0.2179in
]%
{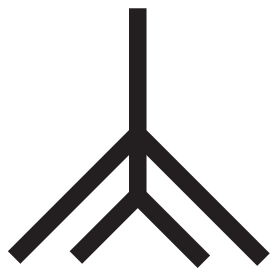}%
}
\ }+\text{ \ }%
\raisebox{-0.0553in}{\includegraphics[
height=0.2145in,
width=0.2179in
]%
{T121-3.eps}%
}
\text{ }\otimes\text{\ }%
\raisebox{-0.0553in}{\includegraphics[
height=0.2145in,
width=0.2179in
]%
{T13-2.eps}%
}
\text{ \ }.
\end{align*}

Note that $\Delta_{K}$ agrees with the Alexander-Whitney diagonal on
$K_{2}=\ast$ and $K_{3}=I.$ Define $\Delta_{K}^{(0)}=${$\mathbf{1}$}$;$ for
each $k\geq1,$ define the \emph{(left) }$k$\emph{-fold} \emph{iterate of
}$\Delta_{K}$ by%
\[
\Delta_{K}^{(k)}=\left(  \Delta_{K}\otimes{\mathbf{1}}^{\otimes k-1}\right)
\Delta_{K}^{(k-1)}%
\]
and view $\Delta_{K}^{(k)}(\curlywedge_{p})$ as a $\left(  p-2\right)
$-dimensional subcomplex of $K_{p}^{\times k+1},$ and dually for $\Delta
_{K}^{(k)}(\curlyvee^{p}).$

Referring to (\ref{operads}) above, define $\mathcal{F}_{3}=\mathcal{H}%
_{1,1}\oplus\mathcal{H}_{1,2}\oplus\mathcal{H}_{2,1}.$ To construct
$\mathcal{F}_{4},$ use the generators of $\mathcal{F}_{3}$ to construct all
possible elementary fractions with two inputs and two outputs. There are
exactly two such elementary fractions, namely,%
\[
\alpha_{2}^{2}\leftrightarrow\curlyvee\diagup\curlywedge\ \ \text{and
}\ \alpha_{11}^{11}\leftrightarrow\left(  \curlywedge\curlywedge\right)
\diagup\left(  \curlyvee\curlyvee\right)  ,\text{ }%
\]
each of dimension zero. Let $\mathcal{H}_{2,2}=\left\langle \theta_{2}%
^{2},\alpha_{2}^{2},\alpha_{11}^{11}\right\rangle $ and define $\partial
\left(  \theta_{2}^{2}\right)  =\alpha_{2}^{2}+\alpha_{11}^{11}.$ Then
$\mathcal{H}_{2,2}$ is a proper submodule of $M_{2,2}$ and $KK_{2,2}$ is an
interval $I$ whose edge is identified with $\theta_{2}^{2}$ and whose vertices
are identified with $\alpha_{2}^{2}$ and $\alpha_{11}^{11}$. Define
$\mathcal{F}_{4}=\mathcal{F}_{3}\oplus\mathcal{H}_{1,3}\oplus\mathcal{H}%
_{2,2}\oplus\mathcal{H}_{3,1}.$

Although all fraction products are used to construct $\mathcal{F}_{4},$ more
fractions than we need appear at the next stage of the construction and
beyond. Note that each numerator or denominator of $\alpha_{2}^{2}$ and
$\alpha_{11}^{11}$ is an iterated S-U diagonal $\Delta_{K}^{\left(  k\right)
}(\curlyvee)$ or $\Delta_{K}^{\left(  k\right)  }(\curlywedge)$ for some
$k=0,1.$ Indeed, the components of $\Delta_{K}^{\left(  k\right)  }\left(
\curlywedge_{m}\right)  $ and $\Delta_{K}^{\left(  k\right)  }\left(
\curlyvee^{n}\right)  $ will determine which fraction products to admit and
which to discard.

Continuing with the construction of $\mathcal{F}_{5},$ use the generators of
$\mathcal{F}_{4}$ to construct all possible fractions with three inputs and
two outputs. There are 18 of these: one in dimension 2, nine in dimension 1
and eight in dimension 0. Since $KK_{2,3}$ is to have a single top dimensional
(indecomposable) 2-cell, we must discard the 2-dimensional (decomposable)
generator%
\[
e\text{ \ = \
\raisebox{-0.1669in}{\includegraphics[
height=0.4255in,
width=0.5725in
]%
{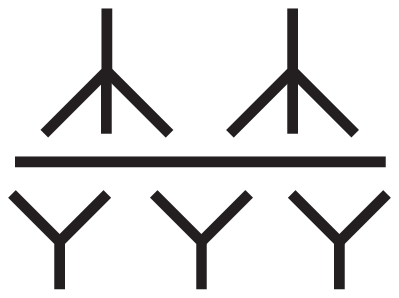}%
}
}%
\]
as well as the appropriate components of its boundary. Note that $e$
represents a square given by the Cartesian product of the three points in the
denominator with the two intervals in the numerator. Thus the boundary of $e$
consists of the four edges%
\begin{equation}
{\includegraphics[
height=0.4298in,
width=3.0649in
]%
{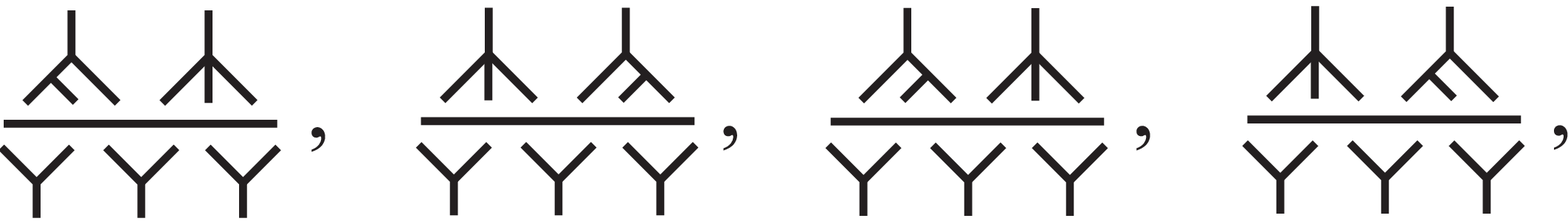}%
}
\label{four edges}%
\end{equation}
the first two of which contain components of $\Delta_{K}^{\left(  1\right)
}(\curlywedge_{3})$ and $\Delta_{K}^{\left(  2\right)  }(\curlyvee)$ in their
numerators and denominators. Our selection rule admits the first two edges and
their vertices.

Express each of the factors $\alpha_{x_{i}}^{q}$ and $\alpha_{p}^{y_{j}}$ in
$\alpha_{\mathbf{x}}^{\mathbf{y}}=\gamma_{\mathbf{x}}^{\mathbf{y}}\left(
\alpha_{p}^{y_{1}}\cdots\alpha_{p}^{y_{q}};\alpha_{x_{1}}^{q}\cdots
\alpha_{x_{p}}^{q}\right)  $ in terms of their respective \emph{leaf
sequences} $\mathbf{x}_{i},$ $\mathbf{q}_{i},$ $\mathbf{y}_{j}$ and
$\mathbf{p}_{j}$ so that%
\[
\alpha_{\mathbf{x}}^{\mathbf{y}}=\frac{\alpha_{\mathbf{p}_{1}}^{\mathbf{y}%
_{1}}\cdots\alpha_{\mathbf{p}_{q}}^{\mathbf{y}_{q}}}{\alpha_{\mathbf{x}_{1}%
}^{\mathbf{q}_{1}}\cdots\alpha_{\mathbf{x}_{p}}^{\mathbf{q}_{p}}}.
\]
Then $\left(  \mathbf{p}_{1},\ldots,\mathbf{p}_{q}\right)  $ and $\left(
\mathbf{q}_{1},\ldots,\mathbf{q}_{p}\right)  $ define the \emph{upper and
lower contact sequences of }$\alpha_{\mathbf{x}}^{\mathbf{y}},$ respectively.

\begin{example}
\label{alpha_3,3}The upper and lower contact sequences of
\begin{equation}
\alpha_{111}^{12}\text{ \ }=\text{ \ }\frac{\theta_{3}^{1}\alpha_{12}^{11}%
}{\theta_{1}^{2}\theta_{1}^{2}\theta_{1}^{2}}\text{ \ }\leftrightarrow\ \
\raisebox{-0.2361in}{\includegraphics[
height=0.5665in,
width=0.6953in
]%
{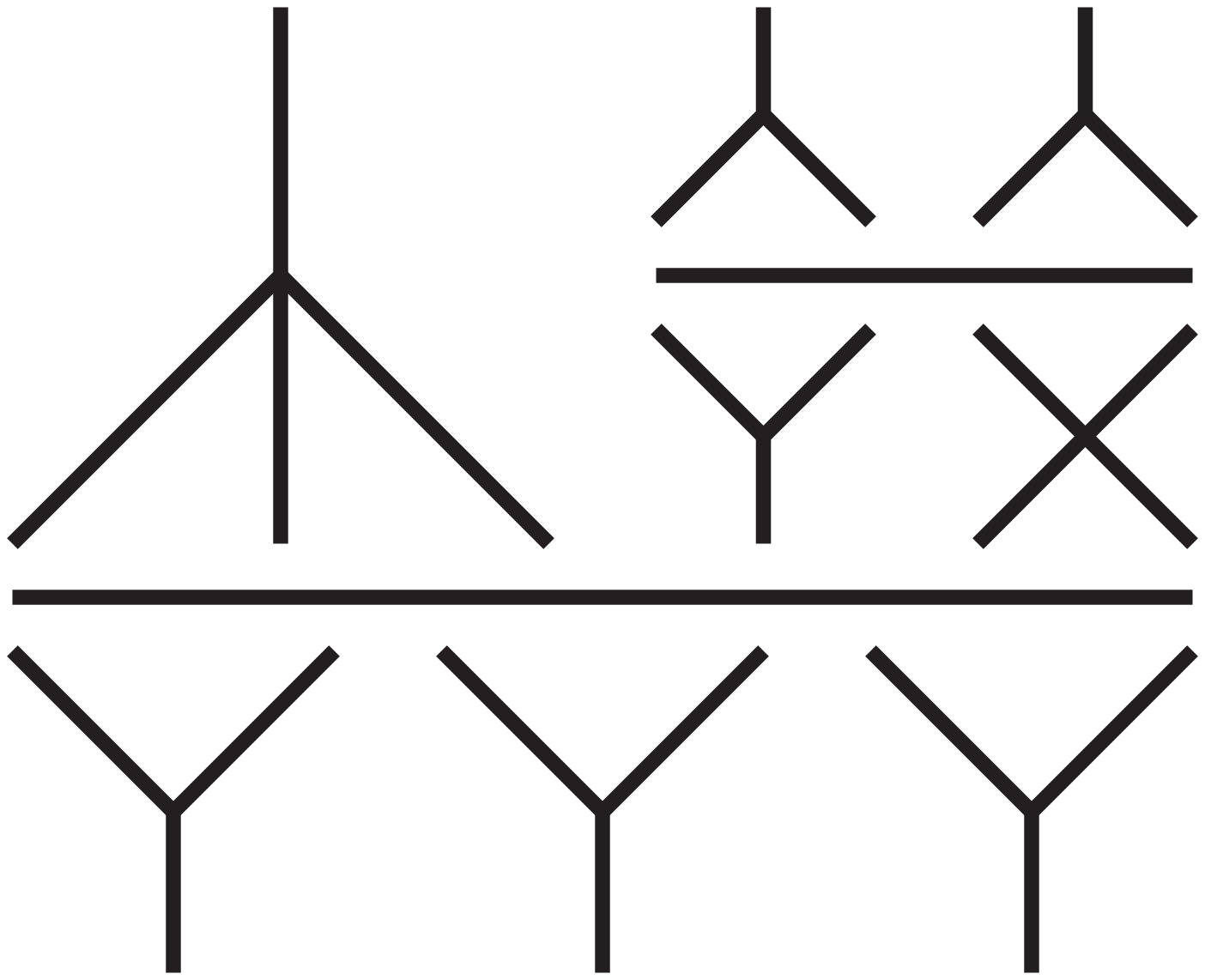}%
}
\text{ \ }=\text{ \
\raisebox{-0.2499in}{\includegraphics[
height=0.5915in,
width=0.4791in
]%
{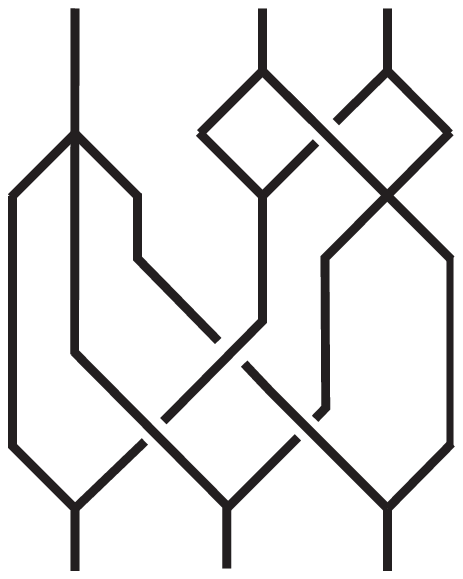}%
}
}\label{alpha3,3}%
\end{equation}
are $\left(  \left(  3\right)  ,\left(  1,2\right)  \right)  $\ and $\left(
\left(  2\right)  ,\left(  2\right)  ,\left(  2\right)  \right)  ,$ respectively.
\end{example}

\noindent A non-vanishing \emph{matrad monomial of codimension }1\emph{\ }%
\[
\alpha_{\mathbf{x}}^{\mathbf{y}}=\frac{\alpha_{\mathbf{p}_{1}}^{\mathbf{y}%
_{1}}\cdots\alpha_{\mathbf{p}_{q}}^{\mathbf{y}_{q}}}{\alpha_{\mathbf{x}_{1}%
}^{\mathbf{q}_{1}}\cdots\alpha_{\mathbf{x}_{p}}^{\mathbf{q}_{p}}}\in
M_{\left\vert \mathbf{y}\right\vert ,\left\vert \mathbf{x}\right\vert }\text{
with }\left\vert \mathbf{x}\right\vert +\left\vert \mathbf{y}\right\vert \leq6
\]
satisfies\emph{\ }the following two conditions:

\begin{enumerate}
\item[\textit{(i)}] The upper contact sequence $\left(  \mathbf{p}_{1}%
,\ldots,\mathbf{p}_{q}\right)  $ is the list of leaf sequences in some
component of $\Delta_{K}^{(q-1)}(\curlywedge_{p})$.

\item[\textit{(ii)}] The lower contact sequence $\left(  \mathbf{q}_{1}%
,\ldots,\mathbf{q}_{p}\right)  $ is the list of leaf sequences in some
component of $\Delta_{K}^{(p-1)}(\curlyvee^{q})$.
\end{enumerate}

\begin{example}
The elementary fraction $\alpha_{111}^{12}=\left(  \theta_{3}^{1}\alpha
_{12}^{11}\right)  /\left(  \theta_{1}^{2}\theta_{1}^{2}\theta_{1}^{2}\right)
$ \ in Example \ref{alpha_3,3} is a non-vanishing 2-dimensional matrad
generator since its upper contact sequence $\left(  \left(  3\right)  ,\left(
1,2\right)  \right)  $ is the list of leaf sequences in the component
\[
\raisebox{-0.0346in}{\includegraphics[
height=0.2101in,
width=0.2145in
]%
{T3.eps}%
}
\text{\ }\otimes\text{\ }
\raisebox{-0.0346in}{\includegraphics[
height=0.2101in,
width=0.2145in
]%
{T12-2.eps}%
}
\ \ \text{of \ }\Delta_{K}^{\left(  1\right)  }\left(  \text{ }
\raisebox{-0.0346in}{\includegraphics[
height=0.2101in,
width=0.2145in
]%
{T3.eps}%
}
\text{\ }\right)
\]
and its lower contact sequence $\left(  \left(  2\right)  ,\left(  2\right)
,\left(  2\right)  \right)  $ is the list of leaf sequences in
\[
\Delta_{K}^{\left(  2\right)  }(\mathsf{Y})=\mathsf{Y}\otimes\mathsf{Y}%
\otimes\mathsf{Y}.
\]

\end{example}

Having discarded the last two fractions in (\ref{four edges}), our selection
rule admits seven 1-dimensional generators labeling the edges of $KK_{2,3}$.
Now linearly extend the boundary map $\partial$ to these matrad generators and
compute the seven admissible 0-dimensional generators labeling the vertices of
$KK_{2,3}$ (see Figure 12). Then in addition to the 2-dimensional generator
$e$ and the last two 1-dimensional generators in (\ref{four edges}), our
selection rule discards the common vertex%
\[
\raisebox{-0.0355in}{\includegraphics[
height=0.4246in,
width=0.5708in
]%
{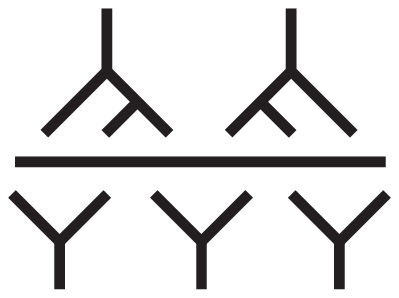}%
}
\text{\ \ .}%
\]
Different elementary fractions may represent the same element. For example,
\begin{equation}
\raisebox{-0.4938in}{\includegraphics[
height=0.5976in,
width=1.388in
]%
{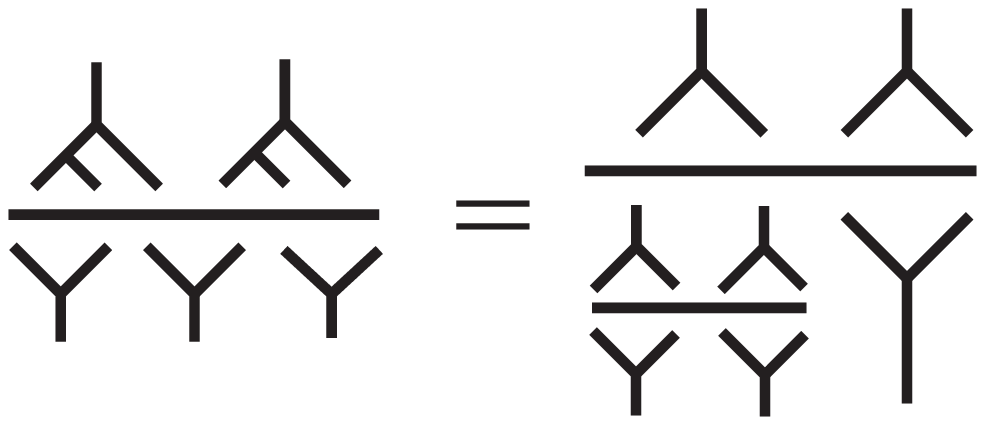}%
}
\text{ \ }.\label{non-unique}%
\end{equation}
The associativity and unit axioms in the definition of a prematrad (see
Definition \ref{prematrad} below) identify various representations such as these.

Finally, $\mathcal{H}_{2,3}$ is the proper submodule of $M_{2,3}$ generated by
$\theta_{3}^{2}$ and the 14 admissible fractions $\alpha_{3}^{2}$ given by the
selection rule above. Define
\[
\partial\left(  \text{
\raisebox{-0.0623in}{\includegraphics[
height=0.2179in,
width=0.2179in
]%
{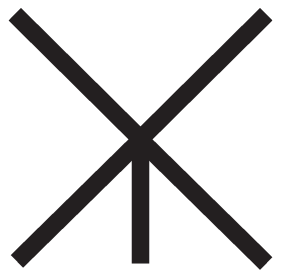}%
}
}\right)  =\
\raisebox{-0.1669in}{\includegraphics[
height=0.4238in,
width=0.1903in
]%
{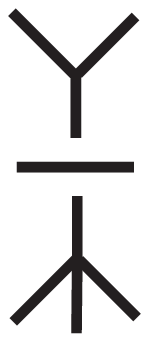}%
}
\ +\
\raisebox{-0.1669in}{\includegraphics[
height=0.4082in,
width=0.3355in
]%
{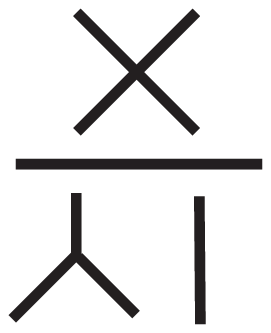}%
}
\text{ }+\
\raisebox{-0.1669in}{\includegraphics[
height=0.4082in,
width=0.3226in
]%
{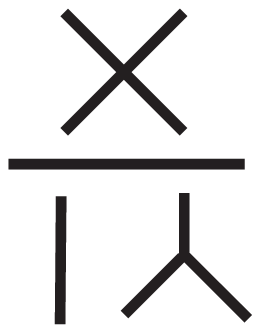}%
}
\ +\
\raisebox{-0.1669in}{\includegraphics[
height=0.4065in,
width=0.4065in
]%
{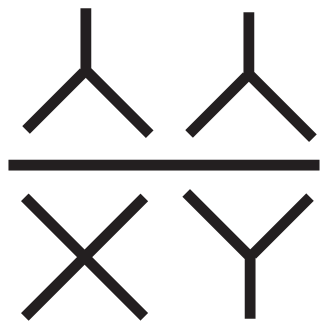}%
}
\ +\
\raisebox{-0.1669in}{\includegraphics[
height=0.4065in,
width=0.4065in
]%
{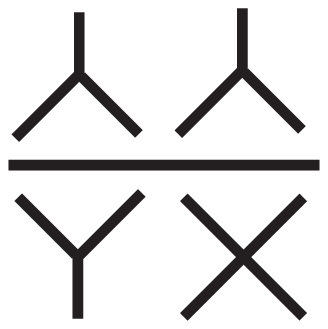}%
}
\ +\
\raisebox{-0.1461in}{\includegraphics[
height=0.416in,
width=0.4938in
]%
{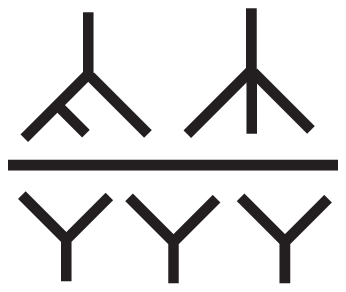}%
}
\ +\
\raisebox{-0.1392in}{\includegraphics[
height=0.416in,
width=0.4938in
]%
{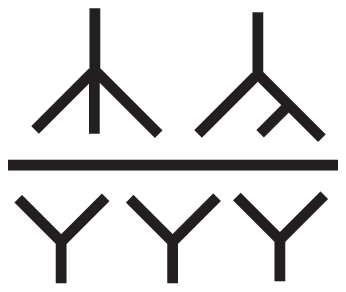}%
}
\text{\ ;}%
\]
then $KK_{2,3}$ is the heptagon pictured in Figure 2. Since $KK_{3,2}$ is
homeomorphic to $KK_{2,3}$ (see Figure 19), we simultaneously obtain
$\mathcal{H}_{3,2}.$ Define $\mathcal{F}_{5}=\mathcal{F}_{4}\oplus
\mathcal{H}_{1,4}\oplus\mathcal{H}_{2,3}\oplus\mathcal{H}_{3,2}\oplus
\mathcal{H}_{4,1}.$

We continue with the construction of $\mathcal{F}_{6}.$ Use the generators of
$\mathcal{F}_{5}$ to construct all possible fractions with two inputs and four
outputs. Using the selection rule defined above, admit all elementary
fractions in dimension 2 whose upper and lower contact sequences agree with
lists of leaf sequences of components in $\Delta_{K}^{(k)}(\curlywedge_{p}) $
or $\Delta_{K}^{(k)}(\curlyvee^{q})$; these represent the 2-faces of
$KK_{4,2}.$ Linearly extend the boundary map $\partial$ and compute the
admissible generators in dimensions 0 and 1. Let $\mathcal{H}_{4,2}$ be the
proper submodule of $M_{4,2}$ generated by $\theta_{2}^{4}$ and all admissible
fractions $\alpha_{2}^{4},$ and define%
\begin{align*}
\partial\left(
\raisebox{-0.0692in}{\includegraphics[
height=0.2171in,
width=0.2179in
]%
{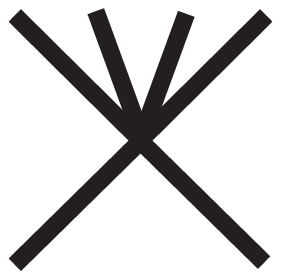}%
}
\right)  \  &  =\
\raisebox{-0.16in}{\includegraphics[
height=0.3952in,
width=0.1799in
]%
{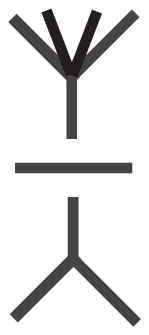}%
}
\ +\
\raisebox{-0.1461in}{\includegraphics[
height=0.3701in,
width=0.3044in
]%
{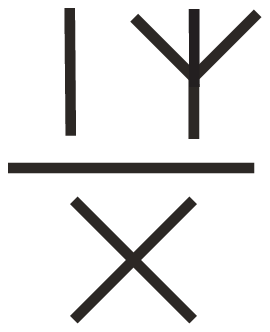}%
}
\ +\
\raisebox{-0.1531in}{\includegraphics[
height=0.3701in,
width=0.2966in
]%
{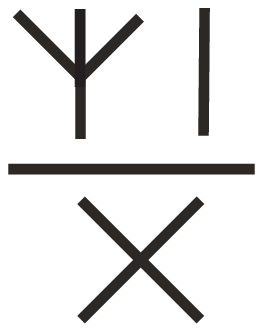}%
}
\ +\
\raisebox{-0.1461in}{\includegraphics[
height=0.3736in,
width=0.3157in
]%
{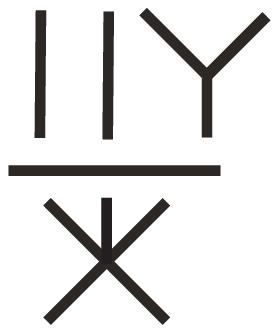}%
}
\ +\
\raisebox{-0.1461in}{\includegraphics[
height=0.3745in,
width=0.2603in
]%
{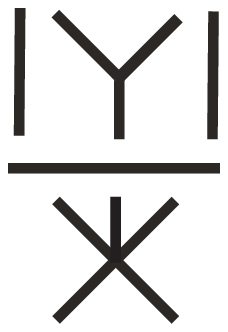}%
}
\ +\
\raisebox{-0.1392in}{\includegraphics[
height=0.3745in,
width=0.2888in
]%
{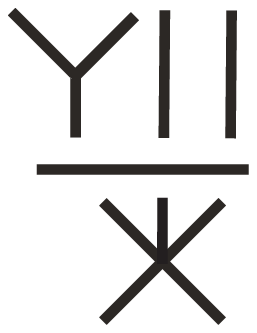}%
}%
\\
&  \hspace{-0.5in}\text{ }+\
\raisebox{-0.1461in}{\includegraphics[
height=0.3632in,
width=0.3667in
]%
{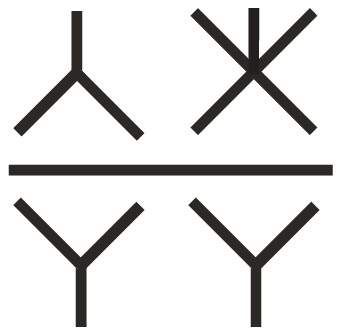}%
}
+\
\raisebox{-0.1392in}{\includegraphics[
height=0.3632in,
width=0.3667in
]%
{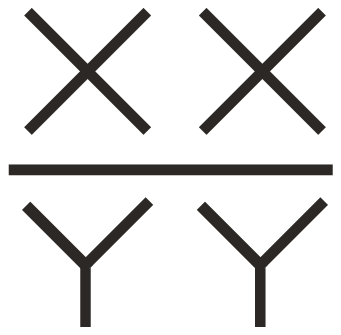}%
}
+\
\raisebox{-0.16in}{\includegraphics[
height=0.4151in,
width=0.4194in
]%
{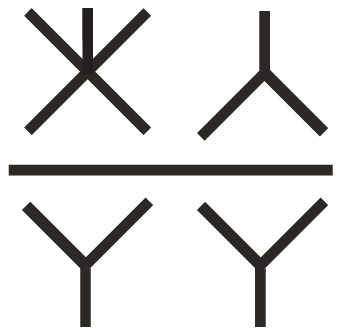}%
}
\ +\
\raisebox{-0.1877in}{\includegraphics[
height=0.4021in,
width=0.4635in
]%
{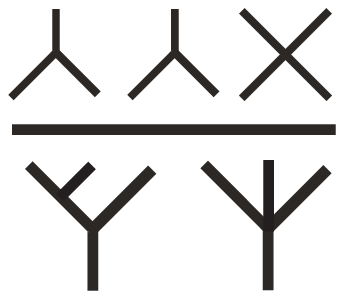}%
}
\ +\
\raisebox{-0.1877in}{\includegraphics[
height=0.4021in,
width=0.4583in
]%
{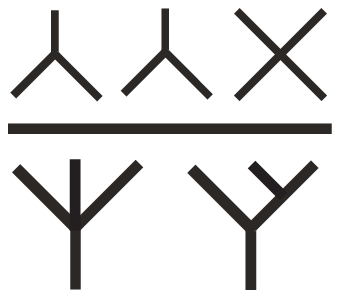}%
}
\text{ }+\
\raisebox{-0.1807in}{\includegraphics[
height=0.3978in,
width=0.4635in
]%
{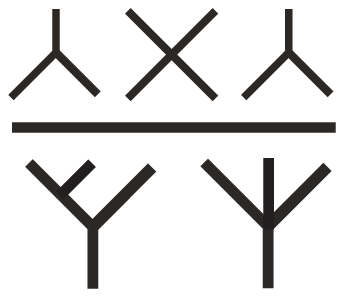}%
}
\\
&  \hspace{-0.5in}\text{ }+\
\raisebox{-0.1877in}{\includegraphics[
height=0.3978in,
width=0.4635in
]%
{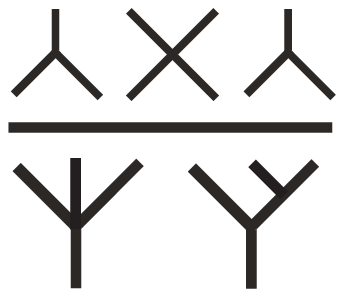}%
}
\text{ }+\text{
\raisebox{-0.1807in}{\includegraphics[
height=0.4013in,
width=0.4583in
]%
{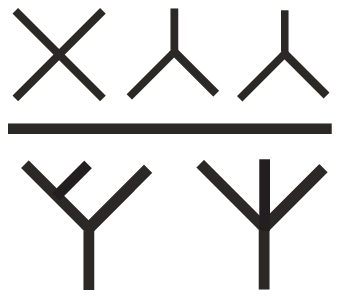}%
}%
}+\ \text{
\raisebox{-0.1807in}{\includegraphics[
height=0.4021in,
width=0.4635in
]%
{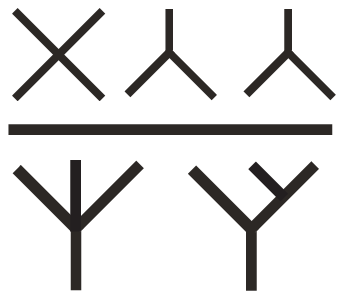}%
}%
}+\
\raisebox{-0.1946in}{\includegraphics[
height=0.4272in,
width=0.5656in
]%
{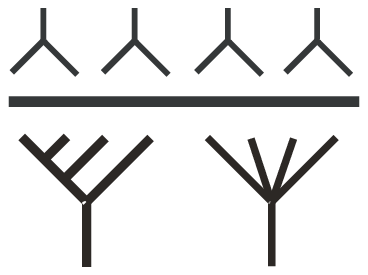}%
}
\ +\
\raisebox{-0.1946in}{\includegraphics[
height=0.4272in,
width=0.5656in
]%
{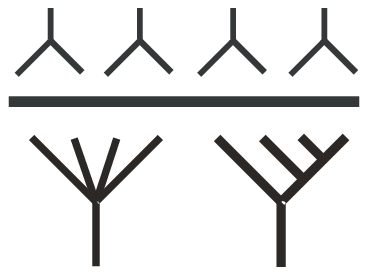}%
}
\\
&  \hspace{-0.5in}\text{ }+\
\raisebox{-0.1877in}{\includegraphics[
height=0.4281in,
width=0.5604in
]%
{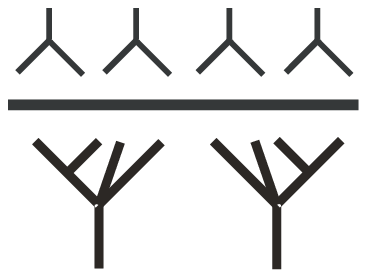}%
}
+\
\raisebox{-0.1738in}{\includegraphics[
height=0.4281in,
width=0.5604in
]%
{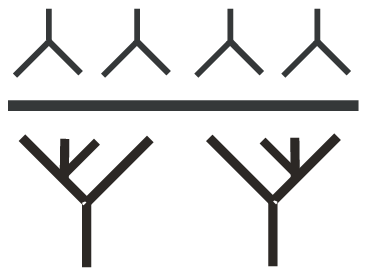}%
}
\ +\
\raisebox{-0.1877in}{\includegraphics[
height=0.4281in,
width=0.5604in
]%
{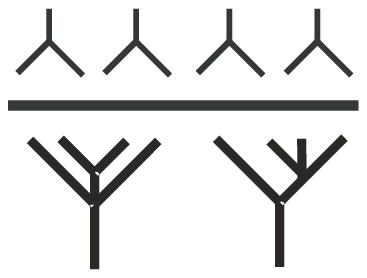}%
}
\text{ }+\
\raisebox{-0.1877in}{\includegraphics[
height=0.4281in,
width=0.5604in
]%
{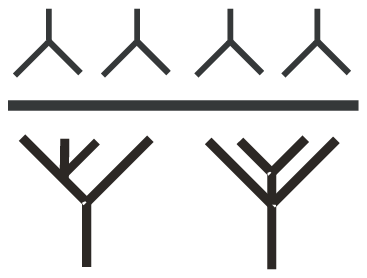}%
}
\text{ };
\end{align*}
then $KK_{4,2}$ is the 3-dimensional polytope pictured in Figures 3 and
21.\smallskip%

\begin{center}
\includegraphics[
height=1.5696in,
width=1.5982in
]%
{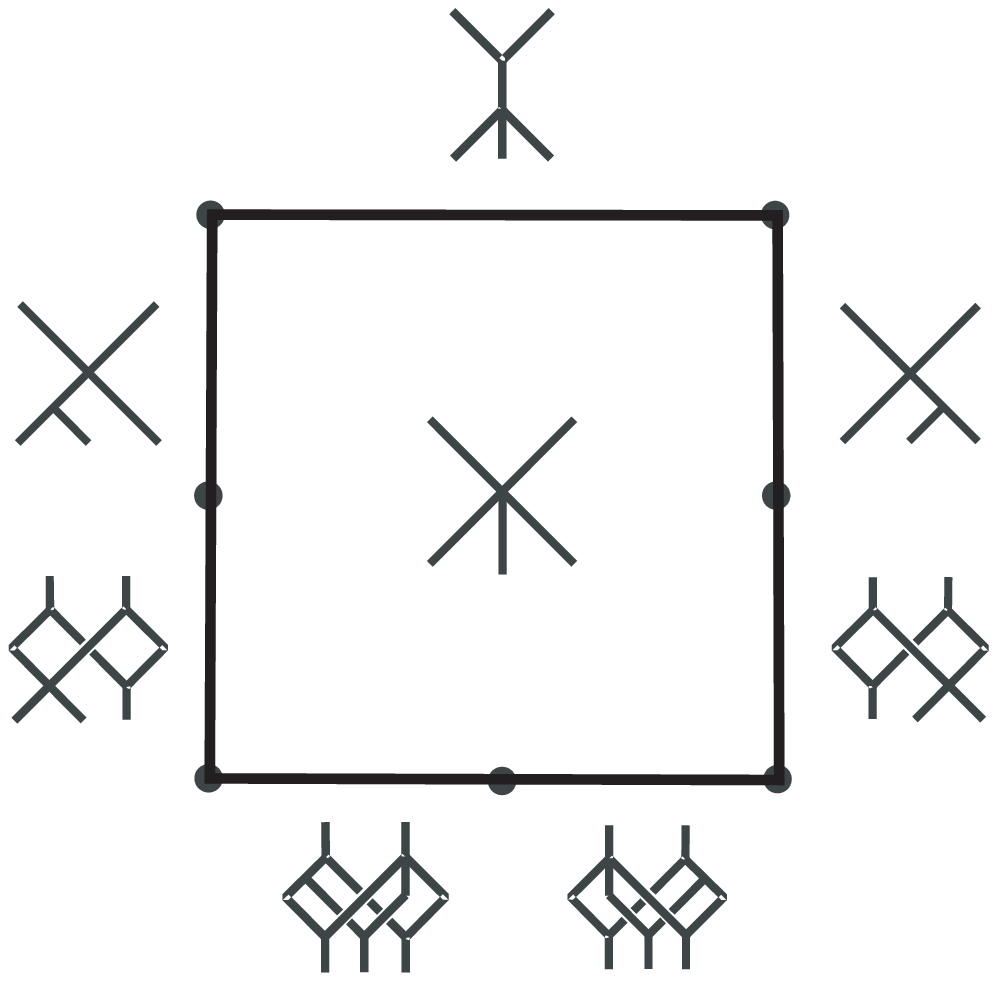}%
\\
Figure 2. The biassociahedron $KK_{2,3}.$%
\end{center}

\begin{center}
\includegraphics[
trim=0.000000in -0.203507in 0.000000in 0.000000in,
height=3.4471in,
width=4.9882in
]%
{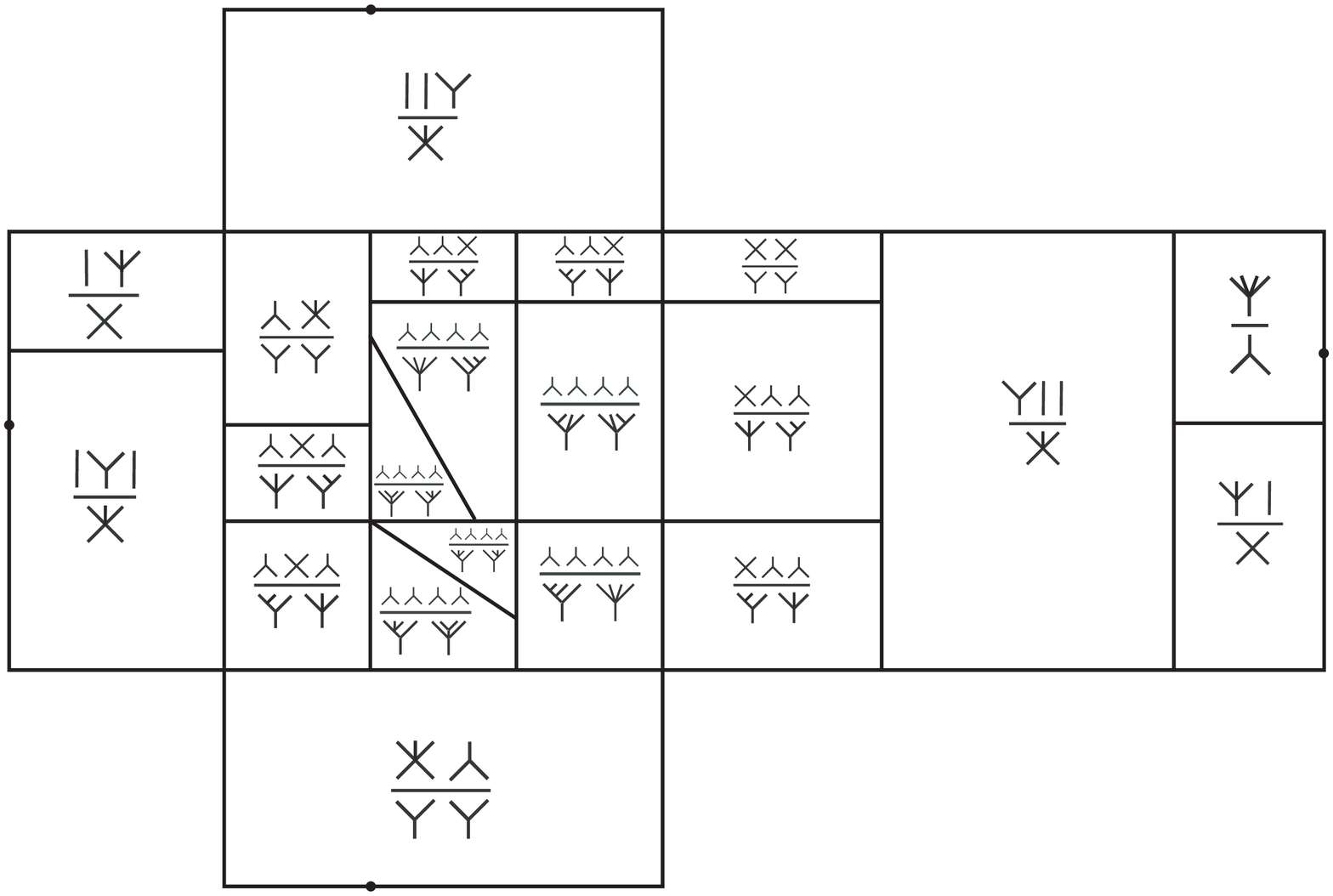}%
\\
Figure 3. The 2-skeleton of $KK_{4,2}.$%
\end{center}

\vspace*{-0.1in}
Again, $KK_{2,4}$ is homeomorphic to $KK_{4,2}$ (see Figure 22) and we
simultaneously obtain $\mathcal{H}_{2,4}$. While similar, the calculations for
$KK_{3,3}$ (see Figure 20) also involve elements such as (\ref{alpha3,3}); we
leave this case to the reader. Define $\mathcal{F}_{6}=\mathcal{F}_{5}%
\oplus\mathcal{H}_{1,5}\oplus\mathcal{H}_{2,4}\oplus\mathcal{H}_{3,3}%
\oplus\mathcal{H}_{4,2}\oplus\mathcal{H}_{5,1}.$

Note that all fractions in $\mathcal{F}_{6}$ are \textquotedblleft
operadic\textquotedblright\ in the sense that each contact sequence is
identified with some component of $\Delta_{K}^{\left(  p\right)  }\left(
\curlywedge_{q}\right)  $ or $\Delta_{K}^{\left(  p\right)  }\left(
\curlyvee^{q}\right)  .$ When $n>6,$ however, $\mathcal{F}_{n}$ contains
\textquotedblleft matradic\textquotedblright\ fractions whose contact
sequences are identified with components of $\Delta_{P}^{\left(  p\right)
}\left(  \curlywedge_{q}\right)  $ or $\Delta_{P}^{\left(  p\right)  }\left(
\curlyvee^{q}\right)  ,$ the iterated diagonal on the permutahedron $P_{q-1}$.
In $\mathcal{F}_{7},$ for example, there is the fraction%
\[
\raisebox{-0.3191in}{\includegraphics[
height=0.7766in,
width=0.9452in
]%
{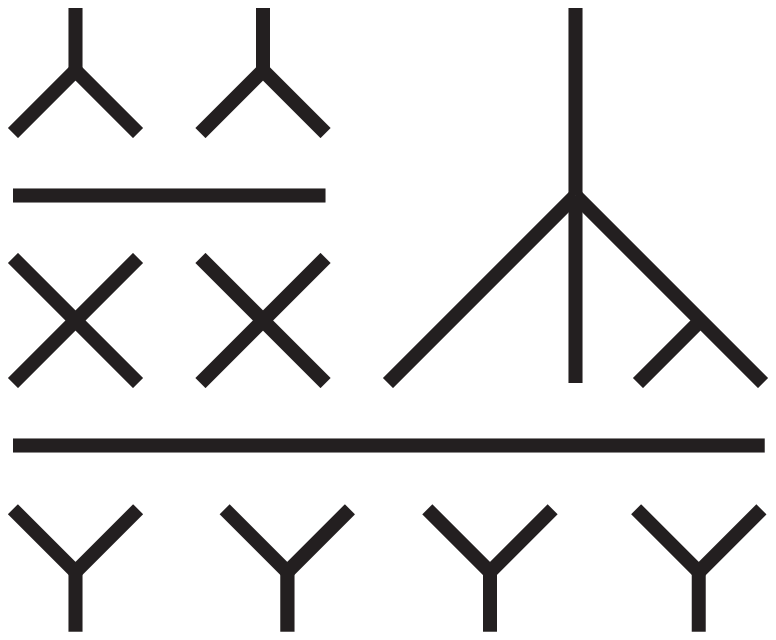}%
}
\text{ .}%
\]
Nevertheless, the low dimensional examples discussed here demonstrate the
general principle, and with this in mind we proceed with the general construction.

\section{Submodules of $TTM$}

Let $M=\left\{  M_{n,m}\right\}  _{m,n\geq1}$ be a bigraded module over a
commutative ring $R$ with identity $1_{R}$. Various submodules of $TTM$ will
be important in our work, the most basic of which is the $q\times p$
\emph{matrix submodule }$(M^{\otimes p})^{\otimes q}$. The name
\textquotedblleft matrix submodule\textquotedblright\ is motivated by the fact
that each pair of $q\times p$ matrices $X=\left(  x_{ij}\right)  ,$ $Y=\left(
y_{ij}\right)  \in\mathbb{N}^{q\times p}$ with $p,q\geq1$ uniquely determines
the submodule%
\[
{M}_{Y,X}=\left(  M_{y_{1,1},x_{1,1}}\otimes\cdots\otimes M_{y_{1,p},x_{1,p}%
}\right)  \otimes\cdots\otimes\left(  M_{y_{q,1},x_{q,1}}\otimes\cdots\otimes
M_{y_{q,p},x_{q,p}}\right)  \subset TTM.
\]
Fix a set of bihomogeneous module generators $G\subset{M}$. A \emph{monomial
in} $TM$ is an element of $G^{\otimes p}$ and a \emph{monomial in} $TTM$ is an
element of $\left(  G^{\otimes p}\right)  ^{\otimes q}.$ Thus $A\in\left(
G^{\otimes p}\right)  ^{\otimes q}$ is naturally represented by the $q\times
p$ matrix
\[
\left[  A\right]  =\left[
\begin{array}
[c]{ccc}%
g_{x_{1,1}}^{y_{1,1}} & \cdots & g_{x_{1,p}}^{y_{1,p}}\\
\vdots &  & \vdots\\
g_{x_{q,1}}^{y_{q,1}} & \cdots & g_{x_{q,p}}^{y_{q,p}}%
\end{array}
\right]
\]
with entries in $G$ and rows identified with elements of $G^{\otimes p}.$ Thus
$A$ is the $q$-fold tensor product of the rows of $\left[  A\right]  ,$ and we
refer to $A$ as a $q\times p$ \emph{monomial }(we use the symbols $A$ and
$\left[  A\right]  $ interchangeably) Consequently,
\[
\left(  M^{\otimes p}\right)  ^{\otimes q}=\bigoplus_{X,Y\in\mathbb{N}%
^{q\times p}}M_{Y,X}%
\]
and we refer to%
\[
\overline{\mathbf{M}}=\bigoplus_{\substack{X,Y\in\mathbb{N}^{q\times p}
\\p,q\geq1}}M_{Y,X}\text{ \ and \ }\overline{\mathbf{V}}=\bigoplus
_{\substack{X,Y\in\mathbb{N}^{1\times p}\cup\mathbb{N}^{q\times1}\text{ }
\\p,q\geq1 }}M_{Y,X}%
\]
as the \emph{matrix submodule }and\emph{\ }the\emph{\ vector submodule,}
respectively. The matrix transpose $A\mapsto A^{T}$ induces the permutation of
tensor factors $\sigma_{p,q}:\left(  M^{\otimes p}\right)  ^{\otimes
q}\overset{\approx}{\rightarrow}\left(  M^{\otimes q}\right)  ^{\otimes p}$
given by%
\begin{equation}
\left(  \alpha_{x_{1,1}}^{y_{1,1}}\otimes\cdots\otimes\alpha_{x_{1,p}%
}^{y_{1,p}}\right)  \otimes\cdots\otimes\left(  \alpha_{x_{q,1}}^{y_{q,1}%
}\otimes\cdots\otimes\alpha_{x_{q,p}}^{y_{q,p}}\right)  \mapsto
\label{permutation}%
\end{equation}%
\[
\left(  \alpha_{x_{1,1}}^{y_{1,1}}\otimes\cdots\otimes\alpha_{x_{q,1}%
}^{y_{q,1}}\right)  \otimes\cdots\otimes\left(  \alpha_{x_{1,p}}^{y_{1,p}%
}\otimes\cdots\otimes\alpha_{x_{q,p}}^{y_{q,p}}\right)  .
\]

Throughout this paper, $\mathbf{x}$ and $\mathbf{y}$ denote matrices in
$\mathbb{N}^{q\times p}$ with constant columns and constant rows,
respectively; $\mathbf{x}$ and $\mathbf{y}$ will often be row and column
matrices. Define%
\[
\overline{\mathbf{M}}_{{\operatorname*{row}}}=\bigoplus_{\mathbf{x}%
,Y\in\mathbb{N}^{q\times p};\text{ }p,q\geq1}M_{Y,\mathbf{x}}\text{ \ \ and
\ \ }\overline{\mathbf{M}}^{{\operatorname*{col}}}=\bigoplus_{X,\mathbf{y}%
\in\mathbb{N}^{q\times p};\text{ }p,q\geq1}M_{\mathbf{y},X}.
\]
Consider $q\times p$ monomials
\[
A=\left[
\begin{array}
[c]{ccc}%
\alpha_{x_{1}}^{y_{1,1}} & \cdots & \alpha_{x_{p}}^{y_{1,p}}\\
\vdots &  & \vdots\\
\alpha_{x_{1}}^{y_{q,1}} & \cdots & \alpha_{x_{p}}^{y_{q,p}}%
\end{array}
\right]  \in M_{Y,\mathbf{x}}\text{ \ and \ }B=\left[
\begin{array}
[c]{ccc}%
\beta_{x_{1,1}}^{y_{1}} & \cdots & \beta_{x_{1,p}}^{y_{1}}\\
\vdots &  & \vdots\\
\beta_{x_{q,1}}^{y_{q}} & \cdots & \beta_{x_{q,p}}^{y_{q}}%
\end{array}
\right]  \in M_{\mathbf{y},X}.
\]
The \emph{row (or coderivation) leaf sequence} \emph{of} $A$ is the $p$-tuple
of lower (input) indices ${\operatorname*{rls}}\left(  A\right)  =\left(
x_{1},\ldots,x_{p}\right)  $ along each row of $A.$ Dually, the \emph{column}
\emph{(or derivation) leaf sequence} \emph{of} $B$ is the $q$-tuple of upper
(output) indices ${\operatorname*{cls}}\left(  B\right)  =\left(  y_{1}%
,\ldots,y_{q}\right)  ^{T}$ along each column of $B.$ Pictorially, each graph
in the $j^{th}$ column of $A$ has $x_{j}$ inputs and each graph in the
$i^{th}$ row of $B$ has $y_{i}$ outputs (see Figure 4).

The \emph{bisequence submodule }of $TTM$ is the intersection
\[
\mathbf{M}=\overline{\mathbf{M}}_{{\operatorname*{row}}}\cap\overline
{\mathbf{M}}^{{\operatorname*{col}}}=\bigoplus_{\mathbf{x},\mathbf{y}%
\in\mathbb{N}^{q\times p};\text{ }p,q\geq1}M_{\mathbf{y},\mathbf{x}}.
\]
A $q\times p$ monomial $A\in\mathbf{M}$ is represented as a \emph{bisequence}
\emph{matrix}
\begin{equation}
A=\left[
\begin{array}
[c]{lll}%
\alpha_{x_{1}}^{y_{1}} & \cdots & \alpha_{x_{p}}^{y_{1}}\\
\multicolumn{1}{c}{\vdots} & \multicolumn{1}{c}{} & \multicolumn{1}{c}{\vdots
}\\
\alpha_{x_{1}}^{y_{q}} & \cdots & \alpha_{x_{p}}^{y_{q}}%
\end{array}
\right]  ;\label{bisequence}%
\end{equation}
in this case ${\operatorname*{rls}}\left(  A\right)  =\left(  x_{1}%
,\ldots,x_{p}\right)  $ and ${\operatorname*{cls}}\left(  A\right)  =\left(
y_{1},\ldots,y_{q}\right)  ^{T}.$ Let
\[
\mathbf{M}_{\mathbf{x}}^{\mathbf{y}}=\left\langle A\in\mathbf{M}\mid
\mathbf{x}={\operatorname*{rls}}\left(  A\right)  \text{ and }\mathbf{y}%
={\operatorname*{cls}}\left(  A\right)  \right\rangle ;
\]
then
\[
\mathbf{M}=\bigoplus_{\substack{\mathbf{x}\times\mathbf{y}\in\mathbb{N}%
^{1\times p}\times\mathbb{N}^{q\times1} \\p,q\geq1}}\mathbf{M}_{\mathbf{x}%
}^{\mathbf{y}}.
\]

Given a finite sequence of positive integers $\mathbf{u}=\left(  u_{1}%
,\ldots,u_{k}\right)  ,$ let $\left\vert \mathbf{u}\right\vert =\Sigma u_{i}$.
By identifying $\left(  H^{\otimes q}\right)  ^{\otimes p}\approx\left(
H^{\otimes p}\right)  ^{\otimes q}$ with $\left(  q,p\right)  \in
\mathbb{N}^{2} $, we can think of a $q\times p$ monomial $A\in\mathbf{M}%
_{\mathbf{x}}^{\mathbf{y}}$ as an operator on the positive integer lattice
$\mathbb{N}^{2},$ pictured as an arrow $\left(  \left\vert \mathbf{x}%
\right\vert ,q\right)  \mapsto\left(  p,\left\vert \mathbf{y}\right\vert
\right)  $ (see Figure 4 and Example \ref{prop}). While this representation is
helpful conceptually, it is unfortunately not faithful.
\pagebreak

\vspace{0.4in}
\hspace*{1.6in}
\setlength{\unitlength}{0.001in} 
\begin{picture}
(152,-140)(202,-100) \thicklines
\put(10,40){\line( 0,-1){100}} \put(310,40){\line( 0,-1){100}} \put(10,-60){\line(
-1,-1){140}} \put(10,-60){\line( 1,-1){288}} \put(310,-60){\line( -1,-1){288}}
\put(310,-60){\line( 1,-1){142}} \put(17,-355){\line( 0,-1){100}} \put(303,-355){\line(
0,-1){100}} \put(303,-355){\line( 1,1){147}} \put(17,-355){\line( -1,1){150}}
\put(155,-750){\line( -1,1){150}} \put(155,-750){\line( 1,1){150}}
\put(155,-750){\line( 1,3){50}} \put(155,-750){\line( -1,3){50}} \put(155,-750){\line(
0,-1){100}} \put(155,-850){\line( 1,-1){150}} \put(155,-850){\line( -1,-1){150}}
\put(700,40){\line( 0,-1){245}} \put(1000,40){\line( 0,-1){245}} \put(700,-60){\line(
-1,-1){140}} \put(700,-60){\line( 1,-1){288}} \put(700,-205){\line( 1,-1){147}}
\put(1000,-60){\line( -1,-1){288}} \put(1000,-60){\line( 1,-1){142}}
\put(1000,-205){\line( -1,-1){147}} \put(707,-355){\line( 0,-1){100}}
\put(850,-355){\line( 0,-1){100}} \put(993,-355){\line( 0,-1){100}}
\put(993,-355){\line( 1,1){147}} \put(707,-355){\line( -1,1){150}}
\put(850,-800){\line( -1,1){150}} \put(850,-800){\line( 1,1){150}}
\put(850,-800){\line( -1,-1){150}} \put(850,-800){\line( 1,-1){150}}
\put(850,-800){\line( 0,-1){165}} \put(850,-800){\line( 1,3){50}} \put(850,-800){\line(
-1,3){50}}
\put(-300,150){\line( 0,-1){1300}} \put(-300,150){\line( 1,0){100}}
\put(-300,-1150){\line( 1,0){100}}
\put(1300,150){\line( 0,-1){1300}} \put(1300,150){\line( -1,0){100}}
\put(1300,-1150){\line( -1,0){100}}
\put(-900,-500){\makebox(0,0) {$A \in \mathbf{M}_{\substack{ \  \\ 23}}^{\substack{ 2
\\ 4 \\ \ }} \hspace*{0.1in} \leftrightarrow  \hspace*{0.1in} $}}
\put(1600,-550){\makebox(0,0){$ \rightarrow \hspace*{0.1in} $}}
\put(2000,-1150){\line(1,0){1100}} \put(2000,-1150){\line(0,1){1300}}
\put(2220,-1183){\line(0,1){75}} \put(2760,-1183){\line(0,1){75}}
\put(1963,-930){\line(1,0){75}} \put(1963,-150){\line(1,0){75}}
\put(2220,-150){\makebox(0,0){$\bullet$}} \put(2760,-930){\makebox(0,0){$\bullet$}}
\put(2690,-830){\vector(-2,3){400}} \put(2700,-500){\makebox(0,0){$A$}}
\put(1850,-930){\makebox(0,0){$2$}} \put(1850,-150){\makebox(0,0){$6$}}
\put(2220,-1300){\makebox(0,0){$2$}} \put(2760,-1300){\makebox(0,0){$5$}}
\end{picture}
\vspace*{1.2in}
\begin{center}
Figure 4. A $2\times2$ monomial in $\mathbf{M}$ and its arrow representation.
\end{center}

Unless explicitly indicated otherwise, we henceforth assume that
$\mathbf{x\times y}=$\linebreak$\left(  x_{1},\ldots,x_{p}\right)
\times\left(  y_{1},\ldots,y_{q}\right)  ^{T}\in\mathbb{N}^{1\times p}%
\times\mathbb{N}^{q\times1}$ with $p,q\geq1$. When the context is clear, we
will often write $\mathbf{y}$ as a row vector. The \emph{bisequence vector
submodule} is the intersection\emph{\ }%
\[
\mathbf{V}=\overline{\mathbf{V}}\cap\mathbf{M}=\bigoplus\limits_{s,t\in
\mathbb{N}}\mathbf{M}_{s}^{\mathbf{y}}\oplus\mathbf{M}_{\mathbf{x}}^{t}.
\]
\begin{example}
\label{bisequence_vect}Let $M=M_{1,1}\oplus M_{2,2}=\left\langle \theta
_{1}^{1}\right\rangle \oplus\left\langle \theta_{2}^{2}\right\rangle .$ Then
the bisequence vector submodule of $TTM$ is
\begin{align*}
\mathbf{V}  &  =\mathbf{M}_{1}^{1}\oplus\mathbf{M}_{11}^{1}\oplus
\mathbf{M}_{111}^{1}\oplus\cdots\oplus\mathbf{M}_{1}^{11}\oplus\mathbf{M}%
_{1}^{111}\oplus\cdots\oplus\\
&  \mathbf{M}_{2}^{2}\oplus\mathbf{M}_{22}^{2}\oplus\mathbf{M}_{222}^{2}%
\oplus\cdots\oplus\mathbf{M}_{2}^{22}\oplus\mathbf{M}_{2}^{222}\oplus\cdots\\
&  \approx T^{+}\left(  M_{1,1}\right)  \oplus T^{+}\left(  M_{1,1}\right)
\oplus T^{+}\left(  M_{2,2}\right)  \oplus T^{+}\left(  M_{2,2}\right)  .
\end{align*}

\end{example}

A submodule
\[
\mathbf{W}=M\oplus\bigoplus\limits_{\mathbf{x},\mathbf{y}\notin\mathbb{N}%
;\text{ }s,t\in\mathbb{N}}\mathbf{W}_{s}^{\mathbf{y}}\oplus\mathbf{W}%
_{\mathbf{x}}^{t}\subseteq\mathbf{V}%
\]
is \emph{telescoping }if for all $\mathbf{x},\mathbf{y},s,t$

\begin{enumerate}
\item[\textit{(i)}] $\mathbf{W}_{s}^{\mathbf{y}}\subseteq\mathbf{M}%
_{s}^{\mathbf{y}}$ and $\mathbf{W}_{\mathbf{x}}^{t}\subseteq\mathbf{M}%
_{\mathbf{x}}^{t};$

\item[\textit{(ii)}] $\alpha_{s}^{y_{1}}\otimes\cdots\otimes\alpha_{s}^{y_{q}%
}\in\mathbf{W}_{s}^{\mathbf{y}}$ implies $\alpha_{s}^{y_{1}}\otimes
\cdots\otimes\alpha_{s}^{y_{j}}\in\mathbf{W}_{s}^{y_{1}\cdots y_{j}}$ for all
$j<q;$

\item[\textit{(iii)}] $\beta_{x_{1}}^{t}\otimes\cdots\otimes\beta_{x_{p}}%
^{t}\in\mathbf{W}_{\mathbf{x}}^{t}$ implies $\beta_{x_{1}}^{t}\otimes
\cdots\otimes\beta_{x_{i}}^{t}\in\mathbf{W}_{x_{1}\cdots x_{i}}^{t}$ for all
$i<p.$
\end{enumerate}

\noindent Thus the truncation\emph{\ }maps $\tau:\mathbf{W}_{s}^{y_{1}\cdots
y_{j}}\rightarrow\mathbf{W}_{s}^{y_{1}\cdots y_{j-1}}$ and $\tau
:\mathbf{W}_{x_{1}\cdots x_{i}}^{t}\rightarrow\mathbf{W}_{x_{1}\cdots x_{i-1}%
}^{t}$ determine the following \textquotedblleft telescoping\textquotedblright%
\ sequences of submodules:%
\begin{align*}
\tau\left(  \mathbf{W}_{s}^{\mathbf{y}}\right)   & \subseteq\tau^{2}\left(
\mathbf{W}_{s}^{\mathbf{y}}\right)  \subseteq\cdots\subseteq\tau^{q-1}\left(
\mathbf{W}_{s}^{\mathbf{y}}\right)  =\mathbf{W}_{s}^{y_{1}}\\
\tau\left(  \mathbf{W}_{\mathbf{x}}^{t}\right)   & \subseteq\tau^{2}\left(
\mathbf{W}_{\mathbf{x}}^{t}\right)  \subseteq\cdots\subseteq\tau^{p-1}\left(
\mathbf{W}_{\mathbf{x}}^{t}\right)  =\mathbf{W}_{x_{1}}^{t}.
\end{align*}
In general, $\mathbf{W}_{\mathbf{x}}^{t}$ is an \emph{additive }submodule of
$\mathbf{M}_{x_{1}}^{t}\otimes\cdots\otimes\mathbf{M}_{x_{p}}^{t}$ and does
\emph{not} necessarily decompose as $B_{1}\otimes\cdots\otimes B_{p}$ with
$B_{i}\subseteq\mathbf{M}_{x_{i}}^{t}$.

The \emph{telescopic extension}\ of a telescoping submodule
$\mathbf{W\subseteq V}$ is the submodule of matrices $\mathcal{W\subseteq
}\overline{\mathbf{M}}$ with the following properties: If $A=\left[
\alpha_{x_{i,j}}^{y_{i,j}}\right]  \ $is a $q\times p$ monomial in
$\mathcal{W}$ and

\begin{enumerate}
\item[(\textit{i)}] $\left[  \alpha_{x_{i,j}}^{y_{i,j}}\text{ }\cdots\text{
}\alpha_{x_{i,j+m}}^{y_{i,j+m}}\right]  $ is a string in the $i^{th}$ row of
$A$ such that $y_{i,j}=\cdots=y_{i,j+m}=t,$ then $\alpha_{x_{i,j}}^{t}%
\otimes\cdots\otimes\alpha_{x_{i,j+m}}^{t}\in\mathbf{W}_{x_{i,j}%
,...,x_{i,j+m}}^{t}.$

\item[\textit{(ii)}] $\left[  \alpha_{x_{i,j}}^{y_{i,j}}\text{ }\cdots\text{
}\alpha_{x_{i+r,j}}^{y_{i+r,j}}\right]  ^{T}$ is a string in the $j^{th}$
column of $A$ such that $x_{i,j}=\cdots=x_{i+r,j}=s,$ then $\alpha
_{s}^{y_{i,j}}\otimes\cdots\otimes\alpha_{s}^{y_{i+r,j}}\in\mathbf{W}%
_{s}^{y_{i,j},...,y_{i+r,j}}.$
\end{enumerate}

\noindent Thus if $A\in\mathbf{M}_{\mathbf{x}}^{\mathbf{y}}\cap\mathcal{W}$,
the $i^{th}$ row of $A$ lies in $\mathbf{W}_{\mathbf{x}}^{y_{i}}$ and $j^{th}$
column of $A$ lies in $\mathbf{W}_{x_{j}}^{\mathbf{y}}$.

\section{Prematrads}

\subsection{$\Upsilon$-products on $\overline{\mathbf{M}}$}

Given a family of maps $\bar{\gamma}=\{M^{\otimes q}\otimes M^{\otimes
p}\rightarrow M\}_{p,q\geq1},$ there is a canonical extension of the component
$\gamma=\{\gamma_{\mathbf{x}}^{\mathbf{y}}:\mathbf{M}_{p}^{\mathbf{y}}%
\otimes\mathbf{M}_{\mathbf{x}}^{q}\rightarrow\mathbf{M}_{\left\vert
\mathbf{x}\right\vert }^{\left\vert \mathbf{y}\right\vert }\}\ $to a global
product $\Upsilon:\overline{\mathbf{M}}\otimes\overline{\mathbf{M}%
}\mathbf{\ \rightarrow}\overline{\mathbf{M}}$. Pairs of bisequence matrices in
$\mathbf{M}_{p}^{\mathbf{y}}\otimes\mathbf{M}_{\mathbf{x}}^{q}$ are called
\textquotedblleft transverse pairs.\textquotedblright

\begin{definition}
A pair $A\otimes B=\left[  \alpha_{v_{kl}}^{y_{kl}}\right]  \otimes\left[
\beta_{x_{ij}}^{u_{ij}}\right]  $ of $\left(  q\times s,t\times p\right)  $
monomials in $\overline{\mathbf{M}}\mathbf{\otimes}\overline{\mathbf{M}}$ is a

\begin{enumerate}
\item[\textit{(i)}] \textbf{Transverse Pair} (TP) if $s=t=1,$ $u_{1,j}=q,$ and
$v_{k,1}=p$ for all $j$ and $k,$ i.e., setting $x_{j}=x_{1,j}$ and
$y_{k}=y_{k,1}$ gives
\[
A\otimes B=\left[
\begin{array}
[c]{c}%
\alpha_{p}^{y_{1}}\\
\vdots\\
\alpha_{p}^{y_{q}}%
\end{array}
\right]  \otimes\left[
\begin{array}
[c]{lll}%
\beta_{x_{1}}^{q} & \cdots & \beta_{x_{p}}^{q}%
\end{array}
\right]  \in\mathbf{M}_{p}^{\mathbf{y}}\otimes\mathbf{M}_{\mathbf{x}}^{q}.
\]

\item[\textit{(ii)}] \textbf{Block Transverse Pair} (BTP) if there exist
$t\times s$ block decompositions $A=\left[  A_{k^{\prime},l}^{\prime}\right]
$ and $B=\left[  B_{i,j^{\prime}}^{\prime}\right]  $ such that $A_{il}%
^{\prime}\otimes B_{il}^{\prime}$ is a TP for\textit{\ all }$i$ and $l$.
\end{enumerate}
\end{definition}

The block sizes in a BTP decomposition are uniquely determined. Unlike the
blocks in a standard block matrix, the blocks $A_{il}^{\prime}$ (or
$B_{il}^{\prime}$) of a BTP $A\otimes B\in\overline{\mathbf{M}}\otimes
\overline{\mathbf{M}}$ may vary in length within a given row (or column).
However, if $A\otimes B\in\mathbf{M}_{p_{1}\cdots p_{s}}^{\mathbf{y}_{1}%
\cdots\mathbf{y}_{t}}\otimes\mathbf{M}_{\mathbf{x}_{1}\cdots\mathbf{x}_{s}%
}^{q_{1}\cdots q_{t}}$ is a BTP, each TP $A_{il}^{\prime}\otimes
B_{il}^{\prime}\in\mathbf{M}_{p_{l}}^{\mathbf{y}_{i}}\otimes\mathbf{M}%
_{\mathbf{x}_{l}}^{q_{i}}$ so that for fixed $i\,$(or $l$) the blocks
$A_{il}^{\prime}$ (or $B_{il}^{\prime}$) have constant length $q_{i}$ (or
$p_{l}$). Furthermore, $A\otimes B\in\mathbf{M}_{\mathbf{v}}^{\mathbf{y}%
}\otimes\mathbf{M}_{\mathbf{x}}^{\mathbf{u}}$ is a BTP if and only if
$\mathbf{x\times y}\in\mathbb{N}^{1\times|\mathbf{v}|}\times\mathbb{N}%
^{|\mathbf{u}|\times1}$ if and only if the initial point of arrow $A$
coincides with the terminal point of arrow $B$ in $\mathbb{N}^{2}$.

\begin{example}
\label{example1}A $\left(  4\times2,2\times3\right)  $ monomial pair $A\otimes
B\in\mathbf{M}_{21}^{1543}\otimes\mathbf{M}_{123}^{31}$ is a $2\times2$ BTP
per the block decompositions\vspace{0.1in}\newline\hspace*{0.9in}%
\setlength{\unitlength}{.06in}\linethickness{0.4pt}
\begin{picture}(118.66,29.34)
\put(-5,16.77){\makebox(0,0)[cc]{$A=$}} \
\put(5.23,25.00){\makebox(0,0)[cc]{$\alpha^{1}_{2}$}}
\put(5.23,18.00){\makebox(0,0)[cc]{$\alpha^{5}_{2}$}}
\put(5.23,11.67){\makebox(0,0)[cc]{$\alpha^{4}_{2}$}}
\put(5.23,4.67){\makebox(0,0)[cc]{$\alpha^{3}_{2}$}}
\put(1.66,8.67){\dashbox{0.67}(6.33,20.00)[cc]{}}
\put(1.66,2.34){\dashbox{0.67}(6.33,4.67)[cc]{}}
\put(14.00,25.00){\makebox(0,0)[cc]{$\alpha^{1}_{1}$}}
\put(14.00,18.00){\makebox(0,0)[cc]{$\alpha^{5}_{1}$}}
\put(14.00,11.67){\makebox(0,0)[cc]{$\alpha^{4}_{1}$}}
\put(14.00,4.67){\makebox(0,0)[cc]{$\alpha^{3}_{1}$}}
\put(10.33,8.67){\dashbox{0.67}(6.33,20.00)[cc]{}}
\put(10.33,2.34){\dashbox{0.67}(6.33,4.67)[cc]{}} \put(-1.0,1.0){\line(0,1){29}}
\put(-1.0,1.0){\line(1,0){2}} \put(-1.0,30.0){\line(1,0){2}}
\put(19.0,1.0){\line(0,1){29}} \put(19.0,1.0){\line(-1,0){2}}
\put(19.0,30.0){\line(-1,0){2}}
\put(24.83,16.77){\makebox(0,0)[cc]{\hspace*{0.2in}$B=$}} \
\put(37.00,20.01){\makebox(0,0)[cc]{$\beta^{3}_{1}$}}
\put(43.58,20.01){\makebox(0,0)[cc]{$\beta^{3}_{2}$}}
\put(51.53,20.01){\makebox(0,0)[cc]{$\beta^{3}_{3}$}}
\put(33.67,17.67){\dashbox{0.67}(12.67,5.67)[cc]{}}
\put(48.33,17.67){\dashbox{0.67}(6.67,5.67)[cc]{}}
\put(37.00,12.34){\makebox(0,0)[cc]{$\beta^{1}_{1}$}}
\put(43.58,12.34){\makebox(0,0)[cc]{$\beta^{1}_{2}$}}
\put(51.53,12.34){\makebox(0,0)[cc]{$\beta^{1}_{3}$}}
\put(33.67,10.01){\dashbox{0.67}(12.67,5.67)[cc]{}}
\put(48.53,10.01){\dashbox{0.67}(6.67,5.67)[cc]{}} \put(31.0,7.0){\line(0,1){19}}
\put(31.0,7.0){\line(1,0){2}} \put(31.0,26.0){\line(1,0){2}}
\put(57.5,7.0){\line(0,1){19}} \put(57.5,7.0){\line(-1,0){2}}
\put(57.5,26.0){\line(-1,0){2}} \put(60.00,16.00){\makebox(0,0)[cc]{$.$}}
\end{picture}

\end{example}

Given a family of maps $\bar{\gamma}=\{M^{\otimes q}\otimes M^{\otimes
p}\rightarrow M\}_{p,q\geq1},$ extend the component $\gamma=\{\gamma
_{\mathbf{x}}^{\mathbf{y}}:\mathbf{M}_{p}^{\mathbf{y}}\otimes\mathbf{M}%
_{\mathbf{x}}^{q}\rightarrow\mathbf{M}_{\left\vert \mathbf{x}\right\vert
}^{\left\vert \mathbf{y}\right\vert }\}\ $to a global product $\Upsilon
:\overline{\mathbf{M}}\otimes\overline{\mathbf{M}}\rightarrow\overline
{\mathbf{M}}$ by defining%
\begin{equation}
\Upsilon\left(  A\otimes B\right)  =\left\{
\begin{array}
[c]{ll}%
\left[  \gamma\left(  A_{ij}^{\prime}\otimes B_{ij}^{\prime}\right)  \right]
, & A\otimes B\ \text{is a\ BTP}\\
0, & \text{otherwise,}%
\end{array}
\right. \label{upsilon}%
\end{equation}
where $A_{ij}^{\prime}\otimes B_{ij}^{\prime}$ is the $\left(  i,j\right)
^{th}$ TP in the BTP decomposition of $A\otimes B.$

We denote the $\Upsilon$-product by \textquotedblleft$\cdot$\textquotedblright%
\ or juxtaposition. When $A\otimes B=\left[  \alpha_{p}^{y_{j}}\right]
^{T}\otimes\left[  \beta_{x_{i}}^{q}\right]  $ is a TP, we write
\[
AB=\gamma(\alpha_{p}^{y_{1}},\ldots,\alpha_{p}^{y_{q}};\beta_{x_{1}}%
^{q},\ldots,\beta_{x_{p}}^{q}).
\]
As an arrow, $AB$ \textquotedblleft transgresses\textquotedblright\ from the
$x$-axis to the $y$-axis $\mathbb{N}^{2}$. When $A^{q\times s}\otimes
B^{t\times p}\in\mathbf{M}_{\mathbf{v}}^{\mathbf{y}}\otimes\mathbf{M}%
_{\mathbf{x}}^{\mathbf{u}}$ is a BTP$\ $and $A_{ij}^{\prime}\otimes
B_{ij}^{\prime}\in\mathbf{M}_{p_{j}}^{\mathbf{y}_{i}}\otimes\mathbf{M}%
_{\mathbf{x}_{j}}^{q_{i}}$ is its BTP decomposition, $AB$ is a $t\times s$
matrix in $\mathbf{M}_{\left\vert \mathbf{x}_{1}\right\vert \cdots\left\vert
\mathbf{x}_{s}\right\vert }^{\left\vert \mathbf{y}_{1}\right\vert
\cdots\left\vert \mathbf{y}_{t}\right\vert }.$ As an arrow, $AB$ runs from the
initial point of $B$ to the terminal point of $A.$

Note that $\Upsilon$-products always restrict to the submodules $\overline
{\mathbf{M}}_{{\operatorname*{row}}}\ $and $\overline{\mathbf{M}%
}^{{\operatorname*{col}}}, $ and consequently to $\mathbf{M}.$ To see this,
consider a BTP $A\otimes B\in\overline{\mathbf{M}}^{{\operatorname*{col}}%
}\otimes\overline{\mathbf{M}}^{{\operatorname*{col}}}$ with block
decomposition $\left[  A_{ij}^{\prime}\right]  \otimes\left[  B_{ij}^{\prime
}\right]  .$ Since each entry along the $i^{th}$ row of $B$ has $q$ outputs,
each block $A_{ij}^{\prime}$ is a column of length $q.$ Since all entries
along a row of $A$ have the same number of outputs, the total number of
outputs from each block $A_{ij}^{\prime}$ is the same for all $j.$ Thus
$AB\in\overline{\mathbf{M}}^{{\operatorname*{col}}}$ and $\Upsilon$ is closed
in $\overline{\mathbf{M}}^{{\operatorname*{col}}}.$ Dually, $\Upsilon$ is
closed in $\overline{\mathbf{M}}_{{\operatorname*{row}}},$ and consequently,
$\Upsilon$ is closed in $\mathbf{M.}$

\begin{example}
\label{upsilonex} Continuing Example \ref{example1}, the action of $\Upsilon$
on the $\left(  4\times2,2\times3\right)  $ monomial pair $A\otimes
B\in\mathbf{M}_{21}^{1543}\otimes\mathbf{M}_{123}^{31}$ produces a $2\times2$
monomial in $\mathbf{M}_{33}^{10,3}:$\vspace{0.1in}\newline\hspace
*{0.1in}%
\setlength{\unitlength}{1.1mm}\linethickness{0.4pt}\begin{picture}(118.66,29.34)
\put(7.20,23.00){\makebox(0,0)[cc]{$\alpha^1_2$}}
\put(7.20,16.00){\makebox(0,0)[cc]{$\alpha^5_2$}}
\put(7.20,9.67){\makebox(0,0)[cc]{$\alpha^4_2$}}
\put(7.20,2.67){\makebox(0,0)[cc]{$\alpha^3_2$}}
\put(3.66,6.67){\dashbox{0.8}(7.0,20.00)[cc]{}}
\put(3.66,0.0){\dashbox{0.8}(7.0,6.0)[cc]{}}
\put(15.85,23.00){\makebox(0,0)[cc]{$\alpha^1_1$}}
\put(15.85,16.00){\makebox(0,0)[cc]{$\alpha^5_1$}}
\put(15.85,9.67){\makebox(0,0)[cc]{$\alpha^4_1$}}
\put(15.85,2.67){\makebox(0,0)[cc]{$\alpha^3_1$}}
\put(12.33,6.67){\dashbox{0.8}(7.0,20.00)[cc]{}}
\put(12.33,0.0){\dashbox{0.8}(7.0,6.0)[cc]{}}
\put(1.0,-1.0){\line(0,1){29}}
\put(1.0,-1.0){\line(1,0){2}}
\put(1.0,28.0){\line(1,0){2}}
\put(22.0,-1.0){\line(0,1){29}}
\put(22.0,-1.0){\line(-1,0){2}}
\put(22.0,28.0){\line(-1,0){2}}
\put(29.00,17.01){\makebox(0,0)[cc]{$\beta^3_1$}}
\put(35.58,17.01){\makebox(0,0)[cc]{$\beta^3_2$}}
\put(43.53,17.01){\makebox(0,0)[cc]{$\beta^3_3$}}
\put(25.67,14.0){\dashbox{0.8}(13.30,6.0)[cc]{}}
\put(40.33,14.0){\dashbox{0.8}(7.0,6.0)[cc]{}}
\put(29.00,9.34){\makebox(0,0)[cc]{$\beta^1_1$}}
\put(35.58,9.34){\makebox(0,0)[cc]{$\beta^1_2$}}
\put(43.53,9.34){\makebox(0,0)[cc]{$\beta^1_3$}}
\put(25.67,6.90){\dashbox{0.8}(13.3,6.0)[cc]{}}
\put(40.53,6.90){\dashbox{0.8}(7.0,6.0)[cc]{}}
\put(24.0,3.5){\line(0,1){19}}
\put(24.0,3.5){\line(1,0){2}}
\put(24.0,22.5){\line(1,0){2}}
\put(49.5,3.5){\line(0,1){19}}
\put(49.5,3.5){\line(-1,0){2}}
\put(49.5,22.5){\line(-1,0){2}}
\put(54.00,13.00){\makebox(0,0)[cc]{$=$}}
\put(65.20,24.50){\makebox(0,0)[cc]{$\alpha^1_2$}}
\put(65.20,17.50){\makebox(0,0)[cc]{$\alpha^5_2$}}
\put(65.20,11.83){\makebox(0,0)[cc]{$\alpha^4_2$}}
\put(61.66,8.83){\dashbox{0.8}(7.0,20.00)[cc]{}}
\put(74.00,17.50){\makebox(0,0)[cc]{$\beta^3_1$}}
\put(80.58,17.50){\makebox(0,0)[cc]{$\beta^3_2$}}
\put(70.67,15.00){\dashbox{0.8}(14.0,6.0)[cc]{}}
\put(93.85,24.50){\makebox(0,0)[cc]{$\alpha^1_1$}}
\put(93.85,17.50){\makebox(0,0)[cc]{$\alpha^5_1$}}
\put(93.85,11.83){\makebox(0,0)[cc]{$\alpha^4_1$}}
\put(90.33,8.83){\dashbox{0.8}(7.0,20.00)[cc]{}}
\put(102.53,17.50){\makebox(0,0)[cc]{$\beta^3_3$}}
\put(99.00,15.0){\dashbox{0.8}(7.0,6.0)[cc]{}}
\put(65.20,0.83){\makebox(0,0)[cc]{$\alpha^3_2$}}
\put(61.66,-2.00){\dashbox{0.8}(7.0,6.0)[cc]{}}
\put(74.00,0.83){\makebox(0,0)[cc]{$\beta^1_1$}}
\put(80.58,0.83){\makebox(0,0)[cc]{$\beta^1_2$}}
\put(70.67,-2.00){\dashbox{0.8}(13.30,6.0)[cc]{}}
\put(93.85,0.33){\makebox(0,0)[cc]{$\alpha^3_1$}}
\put(90.33,-2.0){\dashbox{0.8}(7.0,6.0)[cc]{}}
\put(102.53,0.33){\makebox(0,0)[cc]{$\beta^1_3$}}
\put(99.00,-2.0){\dashbox{0.8}(7.0,6.0)[cc]{}}
\put(60.0,-3.5){\line(0,1){34}}
\put(60.0,-3.5){\line(1,0){2}}
\put(60.0,30.5){\line(1,0){2}}
\put(109.5,-3.5){\line(0,1){34}}
\put(109.5,-3.5){\line(-1,0){2}}
\put(109.5,30.5){\line(-1,0){2}}
\end{picture}\vspace{0.2in}\newline In the target, $\left(  \left\vert
\mathbf{x}_{1}\right\vert ,\left\vert \mathbf{x}_{2}\right\vert \right)
=\left(  1+2,3\right)  $ since $\left(  p_{1},p_{2}\right)  =\left(
2,1\right)  ;$ and $\left(  \left\vert \mathbf{y}_{1}\right\vert ,\left\vert
\mathbf{y}_{2}\right\vert \right)  =
\linebreak \left(  1+5+4,3\right)  $ since $\left(
q_{1},q_{2}\right)  =\left(  3,1\right)  .$ As an arrow in $\mathbb{N}^{2}$,
$AB$ initializes at $\left(  6,2\right)  $ and terminates at $\left(
2,13\right)  .$
\end{example}

\subsection{Prematrads Defined}

Let $1^{1\times p}=\left(  1,\ldots,1\right)  \in\mathbb{N}^{1\times p}$ and
$1^{q\times1}=\left(  1,\ldots,1\right)  ^{T}\in\mathbb{N}^{q\times1}$; we
often suppress the exponents when the context is clear.

\begin{definition}
\label{prematrad}A \textbf{prematrad} $(M,\gamma,\eta)$ is a bigraded
$R$-module $M=\left\{  M_{n,m}\right\}  _{m,n\geq1}$ together with a family of
structure maps $\gamma=\{\gamma_{\mathbf{x}}^{\mathbf{y}}:\mathbf{M}%
_{p}^{\mathbf{y}}\otimes\mathbf{M}_{\mathbf{x}}^{q}\rightarrow\mathbf{M}%
_{\left\vert \mathbf{x}\right\vert }^{\left\vert \mathbf{y}\right\vert }\}$
and a unit $\eta:R\rightarrow\mathbf{M}_{1}^{1}$ such that

\begin{enumerate}
\item[\textit{(i)}] $\Upsilon\left(  \Upsilon\left(  A;B\right)  ;C\right)
=\Upsilon\left(  A;\Upsilon\left(  B;C\right)  \right)  $ whenever $A\otimes
B$ and $B\otimes C$ are BTPs in $\overline{\mathbf{M}}\otimes\overline
{\mathbf{M}};$

\item[\textit{(ii)}] the following compositions are the canonical
isomorphisms:%
\begin{align*}
&  R^{\otimes b}\otimes\mathbf{M}_{a}^{b}\overset{\eta^{\otimes b}%
\otimes{\mathbf{1}}}{\longrightarrow}\mathbf{M}_{1}^{\mathbf{1}^{b\times1}%
}\otimes\mathbf{M}_{a}^{b}\overset{\gamma_{a}^{1^{b\times1}}}{\longrightarrow
}\mathbf{M}_{a}^{b};\\
&  \mathbf{M}_{a}^{b}\otimes R^{\otimes a}\overset{{\mathbf{1}}\otimes
\eta^{\otimes a}}{\longrightarrow}\mathbf{M}_{a}^{b}\otimes\mathbf{M}%
_{\mathbf{1}^{1\times a}}^{1}\overset{\gamma_{1^{1\times a}}^{b}%
}{\longrightarrow}\mathbf{M}_{a}^{b}.
\end{align*}
We denote the element $\eta(1_{R})$ by $\mathbf{1}_{\mathbf{M}}.$
\end{enumerate}

\noindent A \textbf{morphism of prematrads} $(M,\gamma)$ and $(M^{\prime
},\gamma^{\prime})$ is a map $f:{M}\rightarrow{M}^{\prime}$ such that
$f\gamma_{\mathbf{x}}^{\mathbf{y}}=\gamma{^{\prime}}_{\mathbf{x}}^{\mathbf{y}%
}(f^{\otimes q}\otimes f^{\otimes p})$ for all $\mathbf{x}\times\mathbf{y} $.
\end{definition}

Although $\Upsilon$ fails to act associatively on $\overline{\mathbf{M}},$
Axiom (i) implies that $\left(  AB\right)  C=A\left(  BC\right)  $ whenever
$A,B,C\in\overline{\mathbf{M}},$ $AB\neq0,$ and $BC\neq0.$ On the other hand,
$\Upsilon$ acts associatively on $\mathbf{M,}$ which is the content of
Proposition \ref{Upsilon-assoc} below. Given a bisequence matrix $A^{q\times
p}\in\mathbf{M}_{\mathbf{x}}^{\mathbf{y}},$ let $\mathbf{1}^{q\times\left\vert
\mathbf{x}\right\vert }$ and $\mathbf{1}^{\left\vert \mathbf{y}\right\vert
\times p}$ denote the (bisequence) matrices whose entries are constantly
$\mathbf{1}_{\mathbf{M}}.$ Then Axiom (ii) implies $\Upsilon\left(
\mathbf{1}^{\left\vert \mathbf{y}\right\vert \times p};A\right)
=A=\Upsilon\left(  A;\mathbf{1}^{q\times\left\vert \mathbf{x}\right\vert
}\right)  .$

In the discussion that follows, we think of a string of matrices as a
composition of operators and index the matrices in the order applied (from right-to-left).

\begin{definition}
\label{basic}Let $(M,\gamma,\eta)$ be a prematrad, and let $m,n \geq 1$. 
A string of matrices $A_{s}\cdots A_{1}$ is a \textbf{basic string of bidegree} $(m,n)$ if
\begin{enumerate}
\item[\textit{(i)}] $A_{1}\in\mathbf{M}_{\mathbf{x}}^{b},\ \left\vert
\mathbf{x}\right\vert =m,$
\item[\textit{(ii)}] $A_{i}\in\overline{\mathbf{M}}\smallsetminus\left\{
\mathbf{1}^{q\times p}\mid p,q\in\mathbb{N}\right\}  $ for all $i,$
\item[\textit{(iii)}] $A_{s}\in\mathbf{M}_{a}^{\mathbf{y}},$ $\left\vert
\mathbf{y}\right\vert =n,$ and
\item[\textit{(iv)}] some association of $A_{s}\cdots A_{1}$ defines a
sequence of BTPs.
\end{enumerate}
\end{definition}
\noindent In particular, if $A_{s}\cdots A_{1}$ is a basic string of bidegree $(m,n)$ and 
each BTP in Axiom (iv) defines a non-zero $\Upsilon$-product, then $A_{s}\cdots A_{1}$ 
defines a non-zero element of $\mathbf{M}_m^n=M_{n,m}.$ Indeed, this is exactly the situation 
when a prematrad $(M,\gamma,\eta)$ is ``free'' (see Definition \ref{freeprematrad} below). 
\begin{lemma}
\label{basic lemma}Let $(M,\gamma,\eta)$ be a prematrad.  Then
$\Upsilon$ acts associatively on a basic string $A_{s}\cdots A_{1}$ if and only if $A_{i}%
\in\mathbf{M}$ for all $i$.
\end{lemma}
\begin{proof}
Suppose $\Upsilon$ acts associatively on a basic string $A_{s}\cdots A_{1}$.
If $s=1,2$ there is nothing to prove, so assume that $s>2$ and
$1<i<s.$ Since $\Upsilon$ acts associatively, every association of
$A_{s}\cdots A_{1}$ defines a sequence of BTPs. Hence 
\[
BA_{i}C=\left(  A_{s}\cdots A_{i+1}\right)  A_{i}\left(
A_{i-1}\cdots A_{1}\right)
\] 
is a basic string and $B\otimes A_{i}$
and $A_{i}\otimes C$ are BTPs. So write $B=\left[  b_{p}^{y_{1}}\cdots
b_{p}^{y_{q}}\right]  ^{T}$ and $C=\left[  c_{x_{1}}^{r}\cdots c_{x_{t}}%
^{r}\right]  ;$ then the $i^{th}$ row of $A_{i}$ has the form $\left[
a_{u_{1}}^{w_{i}}\cdots a_{u_{p}}^{w_{i}}\right]  $ and the $j^{th}$ column of
$A_{i}$ has the form $\left[  a_{z_{j}}^{v_{1}}\cdots a_{z_{j}}^{v_{r}%
}\right]  ^{T}.$ Thus $A_{i}\in \overline{\mathbf{M}}_{\operatorname*{row}}%
\cap\overline{\mathbf{M}}^{\operatorname*{col}}=\mathbf{M}.$

Conversely, we proceed by induction on string length $s$. Consider a basic
string $ABC$ with $B\in\mathbf{M},$ and suppose that $A\otimes
B$ and $AB\otimes C$ are BTPs. Write%
\[
A\otimes B=\left[
\begin{array}
[c]{c}%
a_{p}^{y_{1}}\\
\vdots\medskip\\
a_{p}^{y_{r}}%
\end{array}
\right]  \otimes\left[
\begin{array}
[c]{ccc}%
b_{u_{1}}^{v_{1}} & \cdots & b_{u_{p}}^{v_{1}}\\
\vdots\medskip &  & \vdots\\
b_{u_{1}}^{v_{q}} & \cdots & b_{u_{p}}^{v_{q}}%
\end{array}
\right]  ;
\]
let $B^{i}$ and $B_{j}$ denote the $i^{th}$ row and $j^{th}$ column of $B.$
Set $v_{0}=0$ and express $A$ as the block matrix $\left[  A^{1}\cdots\text{
}A^{q}\right]  ^{T}$ where $A^{i}=\left[  a_{p}^{y_{v_{1}+\cdots+v_{i-1}+1}%
}\cdots\text{ }a_{p}^{y_{v_{1}+\cdots+v_{i}}}\right]  ^{T}.$ Then%
\[
AB=\left[
\begin{array}
[c]{c}%
A^{1}B^{1}\\
\vdots\medskip\\
A^{q}B^{q}%
\end{array}
\right]  ,\text{ \ where }A^{i}B^{i}=\left[
\begin{array}
[c]{c}%
a_{p}^{y_{v_{1}+\cdots+v_{i-1}+1}}\\
\vdots\medskip\\
a_{p}^{y_{v_{1}+\cdots+v_{i}}}%
\end{array}
\right]  \left[  b_{u_{1}}^{v_{i}}\ \cdots\ b_{u_{p}}^{v_{i}}\right]  .
\]
Since $AB\otimes C$ is a BTP, $C$ has the form $\left[  c_{x_{1}}^{q}\cdots
c_{x_{s}}^{q}\right]  ,$ where $s=u_{1}+\cdots+u_{p}.$ Set $u_{0}=0$ and
express $C$ as the block matrix $\left[  C_{1}\cdots C_{p}\right]  ,$ where
\[
C_{j}=\left[  c_{x_{u_{1}+\cdots+u_{j-1}+1}}^{q}\cdots\right.  \left.
c_{x_{u_{1}+\cdots+u_{j}}}^{q}\right].
\]
Then 
\[
B_{j}\otimes C_{j}=\left[
\begin{array}
[c]{c}%
b_{u_{j}}^{v_{1}}\\
\vdots\\
b_{u_{j}}^{v_{q}}%
\end{array}
\right]\ \otimes\ \left[  c_{x_{u_{1}+\cdots+u_{j-1}+1}}^{q}\cdots\ c_{x_{u_{1}+\cdots+u_{j}}}^{q}\right] 
\]
is a TP and $B\otimes C$ is a BTP.
Therefore $\left(  AB\right)  C=A\left(  BC\right)  $ by Definition
\ref{prematrad}, Axiom (i). Similarly, if $B\otimes C$ and $A\otimes BC$ are
BTPs, then $A\otimes B$ is a BTP.

Next consider a basic string $ABCD$ with $B,C\in\mathbf{M}$.
If $B\otimes C,$ $A\otimes BC,$ and $A\left(  BC\right)  \otimes D$ are BTPs,
$A\left(  BC\right)  $ is a column matrix whose entries are basic strings of
the form $A_{i}\left(  \left[  B_{i1}\cdots B_{ip}\right]  \left[
C_{i1}\cdots C_{ip}\right]  \right)  $. Hence $A\left(  BC\right)  =\left(
AB\right)  C$ by the calculations above; and dually, $\left(  BC\right)
D=B\left(  CD\right)  .$ Furthermore, the equalities $\left(  A\left(
BC\right)  \right)  D=\left(  \left(  AB\right)  C\right)  D$ and $A\left(
\left(  BC\right)  D\right)  =A\left(  B\left(  CD\right)  \right)  \ $imply
$\left(  \left(  AB\right)  C\right)  D=\left(  AB\right)  \left(  CD\right)
=$ 
\linebreak $A\left(  B\left(  CD\right)  \right)  .$

Inductively, let $k\geq4,$ and assume that $\Upsilon$ acts associatively on
basic strings $A_{s}\cdots A_{1}$ of length $s\leq k$ with
$A_{i}\in\mathbf{M}$ for all $i,$ and consider a basic string
$A_{k+1}\cdots A_{1}$ with $A_{j}\in\mathbf{M}$ for all $j.$
Since some association of $A_{k+1}\cdots A_{1}$ defines a sequence of BTPs,
there is an innermost BTP $A_{j+1}\otimes A_{j}.$ Let $B=A_{j+1}A_{j}$; then
$B$ is a bisequence matrix since $\Upsilon$ is closed in $\mathbf{M},\,$and
$\Upsilon$ acts associatively on $A_{k+1}\cdots A_{j+2}BA_{j-1}\cdots A_{1}. $
If $1<j<k,$ let $C=A_{k+1}\cdots A_{j+2}$ and $D=A_{j-1}\cdots A_{1};$ then
$\Upsilon$ acts associatively on $CA_{j+1}A_{j}D$, completing the proof.
\end{proof}

\begin{proposition}
\label{Upsilon-assoc}Let $(M,\gamma,\eta)$ be a prematrad; then $\Upsilon$ acts associatively on $\mathbf{M}.$
\end{proposition}

\begin{proof}
If $A,B,C\in\mathbf{M}$ such that $A\otimes B$ and $AB\otimes C$ are BTPs, the
entries of $\left(  AB\right)  C$ are basic strings on which
$\Upsilon$ acts associatively by Lemma \ref{basic lemma}. Hence $B\otimes C$
is a BTP and $\left(  AB\right)  C=A\left(  BC\right)  $. If $A\otimes B$ is a
BTP and $B\otimes C$ is not, neither is $AB\otimes C.$ Dually, if $B\otimes C$
is a BTP and $A\otimes B$ is not, neither is $A\otimes BC.$ In either case,
$\left(  AB\right)  C=A\left(  BC\right)  =0.$
\end{proof}

Some examples of prematrads now follow.

\begin{remark}
\label{2operads}If $\left(  M,\gamma,\eta\right)  $ is a prematrad, the
restrictions $\left(  M_{1,\ast},\text{ }\gamma_{\mathbf{\ast}}^{1}%
,\eta\right)  $ and $\left(  M_{\ast,1},\text{ }\gamma_{1}^{\mathbf{\ast}%
},\eta\right)  $ are non-$\Sigma$ operads in the sense of May (see \cite{MSS}).
\end{remark}

\begin{example}
\label{K}A non-$\Sigma$ operad $\left(  K,\gamma_{\ast}\right)  $ with
$\gamma_{\mathbf{x}}:K(p)\otimes K(x_{1})\otimes\cdots\otimes K(x_{p}%
)\rightarrow K(|\mathbf{x}|)$ is a prematrad via
\[
M_{n,m}=\left\{
\begin{array}
[c]{ll}%
K(m), & \text{if}\ n=1\\
0, & \text{otherwise}%
\end{array}
\right.  \ \ \text{and}\ \ \gamma_{\mathbf{x}}^{\mathbf{y}}=\left\{
\begin{array}
[c]{ll}%
\gamma_{\mathbf{x}}, & \text{if}\ \mathbf{y}=1\\
0, & \text{otherwise}%
\end{array}
\right.
\]
(c.f. Remark \ref{2operads}). For a discussion of the differential in the
special case $K=\mathcal{A}_{\infty},$ see Example \ref{freematrad}.
\end{example}

\begin{example}
\label{L}Let $(K=\bigoplus\nolimits_{n\geq1}K(n),\gamma_{\mathbf{\ast}})$ and
$(L=\bigoplus\nolimits_{m\geq1}L(m),\gamma^{\ast})$ be non-$\Sigma$ operads
with $K(1)=L(1)$ and the same unit $\eta.$ Set
\[
M_{n,m}=\left\{
\begin{array}
[c]{ll}%
K(n), & \text{if}\ m=1\\
L(m), & \text{if}\ n=1\\
0, & \text{otherwise}%
\end{array}
\right.  \ \ \text{and}\ \ \gamma_{\mathbf{x}}^{\mathbf{y}}=\left\{
\begin{array}
[c]{ll}%
\gamma_{\mathbf{x}}, & \text{if}\ \mathbf{y}=1\\
\gamma^{\mathbf{y}}, & \text{if}\ \mathbf{x}=1\\
0, & \text{otherwise,}%
\end{array}
\right.
\]
then $\left(  M,\gamma\right)  $ is a prematrad.
\end{example}

\begin{example}
\label{prop}\underline{The Prematrad PROP\textit{\ }$M.$} The free
PROP\textit{\ }$M,$ with its horizontal and vertical products $\times
:M_{n,m}\otimes M_{n^{\prime},m^{\prime}}\rightarrow M_{n+n^{\prime
},m+m^{\prime}}$ and $\circ:M_{r,q}\otimes M_{q,p}\rightarrow M_{r,p}$ (c.f.
\cite{Adams}, \cite{MSS}), is endowed with \textit{a canonical prematrad}
structure $\left(  M^{\operatorname*{pre}},\gamma,\eta\right)  ,$ with $\eta$
determined by $\eta\left(  1_{R}\right)  =\{$unit of the PROP $M\}.$ To define
the structure map $\gamma,$ define $\times^{0}=${$\mathbf{1}$} and iterate
$\times$ to obtain
\[
\times^{q-1}\otimes\times^{p-1}:\mathbf{M}_{p}^{\mathbf{y}}\otimes
\mathbf{M}_{\mathbf{x}}^{q}\longrightarrow\mathbf{M}_{pq}^{|\mathbf{y}%
|}\otimes\mathbf{M}_{|\mathbf{x}|}^{qp}.
\]
View $\alpha_{\left\vert \mathbf{x}\right\vert }^{qp}\in\mathbf{M}%
_{|\mathbf{x}|}^{qp}$ as a graph with $p$ groups of $q$ outputs $\left(
y_{1,1}\cdots y_{1,q}\right)  $ $\cdots$ $\left(  y_{p,1}\cdots y_{p,q}%
\right)  $ labeled from left-to-right. The leaf permutation
\[
\sigma_{q,p}:\left(  y_{1,1}\cdots y_{1,q}\right)  \cdots\left(  y_{p,1}\cdots
y_{p,q}\right)  \mapsto\left(  y_{1,1}\cdots y_{p,1}\right)  \cdots\left(
y_{1,q}\cdots y_{p,q}\right)
\]
induces a map $\sigma_{q,p}^{\ast}:\mathbf{M}_{|\mathbf{x}|}^{qp}%
\rightarrow\mathbf{M}_{|\mathbf{x}|}^{pq}.$ Then $\gamma$ is the sum of the
compositions
\begin{equation}
\gamma_{\mathbf{x}}^{\mathbf{y}}:\mathbf{M}_{p}^{\mathbf{y}}\otimes
\mathbf{M}_{\mathbf{x}}^{q}\overset{\times^{q-1}\otimes\times^{p-1}%
}{\longrightarrow}\mathbf{M}_{pq}^{|\mathbf{y}|}\otimes\mathbf{M}%
_{|\mathbf{x}|}^{qp}\overset{{\mathbf{1}}\otimes\sigma_{q,p}^{\ast
}}{\longrightarrow}\mathbf{M}_{pq}^{|\mathbf{y}|}\otimes\mathbf{M}%
_{|\mathbf{x}|}^{pq}\overset{\circ}{\longrightarrow}\mathbf{M}_{|\mathbf{x}%
|}^{|\mathbf{y}|}.\label{gamma}%
\end{equation}
The induced associative product $\Upsilon$ on $\mathbf{M}$ takes values in
matrices of (typically) non-planar graphs as in Figure 4. In particular, let
$H$ be a free DG $R$-module of finite type and view the universal PROP
$U_{H}=End(TH)$ as the bigraded $R$-module $M=\left\{  {Hom}\left(  H^{\otimes
m},H^{\otimes n}\right)  \right\}  _{n,m\geq1}.$ Then a $q\times p$ monomial
$A\in\mathbf{M}_{\mathbf{x}}^{\mathbf{y}}$ admits a representation as an
operator on $\mathbb{N}^{2}$ via the identification $\left(  H^{\otimes
p}\right)  ^{\otimes q}\leftrightarrow\left(  p,q\right)  $ with the action of
$A$ given by the composition
\[
\left(  H^{\otimes\left\vert \mathbf{x}\right\vert }\right)  ^{\otimes
q}\approx\left(  H^{\otimes x_{1}}\otimes\cdots\otimes H^{\otimes x_{p}%
}\right)  ^{\otimes q}\rightarrow\left(  H^{\otimes y_{1}}\right)  ^{\otimes
p}\otimes\cdots\otimes\left(  H^{\otimes y_{q}}\right)  ^{\otimes p}%
\]%
\[
\overset{\sigma_{y_{1},p}\otimes\cdots\otimes\sigma_{y_{q},p}}{\longrightarrow
}\left(  H^{\otimes p}\right)  ^{\otimes y_{1}}\otimes\cdots\otimes\left(
H^{\otimes p}\right)  ^{\otimes y_{q}}\approx\left(  H^{\otimes p}\right)
^{\otimes\left\vert \mathbf{y}\right\vert }.
\]
This motivates the representation of a general $A$ as an arrow $\left(
\left\vert \mathbf{x}\right\vert ,q\right)  \mapsto\left(  p,\left\vert
\mathbf{y}\right\vert \right)  $ in $\mathbb{N}^{2}$ (see Figure 4). The map
$\times^{q-1}\otimes\times^{p-1}$ in \ref{gamma} is the canonical isomorphism
and $\gamma$ agrees with the composition product on the universal preCROC
\cite{borya}.
\end{example}

\begin{example}
\label{UEP}\underline{The Universal Enveloping PROP $U$.} Recall that the
structure map $\gamma_{FP}$ in the free PROP\ $FP\left(  M\right)  $ is the
sum of all possible (iterated) \textquotedblleft horizontal\textquotedblright%
\ and \textquotedblleft vertical\textquotedblright\ products $\times
:M_{b,a}\otimes M_{b^{\prime},a^{\prime}}\rightarrow M_{b+b^{\prime
},a+a^{\prime}}$ and $\circ:M_{c,b}\otimes M_{b,a}\rightarrow M_{c,a}.$
Furthermore, the tensor product induces left $S_{n}$- and right $S_{m}%
$-actions
\[
M_{b,b}\otimes M_{b,a}\overset{\sigma\otimes\mathbf{1}}{\rightarrow}%
M_{b,b}\otimes M_{b,a}\overset{\gamma_{FP}}{\rightarrow}M_{b,a}\text{ \ and}%
\]%
\[
M_{b,a}\otimes M_{a,a}\overset{\mathbf{1}\otimes\sigma}{\rightarrow}%
M_{b,a}\otimes M_{a,a}\overset{\gamma_{FP}}{\rightarrow}M_{b,a}.
\]
Note that $FP$ is functorial: Given bigraded modules $M$ and $N,$ a map
$f=\{f_{b,a}:M_{b,a}\rightarrow N_{b,a}\}$ extends to a map
$FP(f):FP(M)\rightarrow FP(N)$ preserving horizontal and vertical products,
i.e., $FP(f)(M_{b,a}\times M_{b^{\prime},a^{\prime}})=f(M_{b,a})\times
f(M_{b^{\prime},a^{\prime}})$ and $FP(f)(M_{c,b}\circ M_{b,a}):=f(M_{c,b}%
)\circ f(M_{b,a}).$ Now if $\left(  M,\gamma_{M}\right)  $ is a prematrad,
$\gamma_{M}$ is a component of $\gamma_{_{FP}}$ on $\mathbf{M}_{p}%
^{\mathbf{y}}\otimes\mathbf{M}_{\mathbf{x}}^{q}$. Let $J$ be the two-sided
ideal generated by $\bigoplus\nolimits_{\mathbf{x\times y}}\left(
\gamma_{_{FP}}-\gamma_{_{M}}\right)  \left(  \mathbf{M}_{p}^{\mathbf{y}%
}\otimes\mathbf{M}_{\mathbf{x}}^{q}\right)  .$ The \textbf{universal
enveloping PROP of}\textit{\ }$M$ is the quotient%
\[
U(M)=FP(M)\diagup J.
\]
Note that the restriction of $U$ to operads is the standard functor from
operads to PROPs \cite{Adams}.
\end{example}

\subsection{Free Prematrads\label{freepre}}

\textquotedblleft Free prematrads\textquotedblright\ are fundamentally
important. Our definition of a free prematrad (below) involves an inductive
definition of the intermediate set $G^{\operatorname{pre}}=G_{\ast,\ast
}^{\operatorname*{pre}}$ in which $G_{n,m}^{\operatorname*{pre}}$ is defined
in terms of the set
\[
G_{[n,m]} =\bigcup\limits_{\substack{i\leq m,\text{ }j\leq n, \\i+j<m+n}%
}G_{j,i}^{\operatorname*{pre}}.
\]
We think of $G_{[n,m]}$ as an $n\times m$ array whose $\left(  j,i\right)
^{th}$ cell contains $G_{j,i}^{\operatorname*{pre}}$ when $i+j<m+n$ and whose
$\left(  n,m\right)  ^{th}$ cell is empty. For example, we picture $G_{[3,4]}$
as:
\begin{center}
\includegraphics[
height=1.0819in,
width=1.3482in
]%
{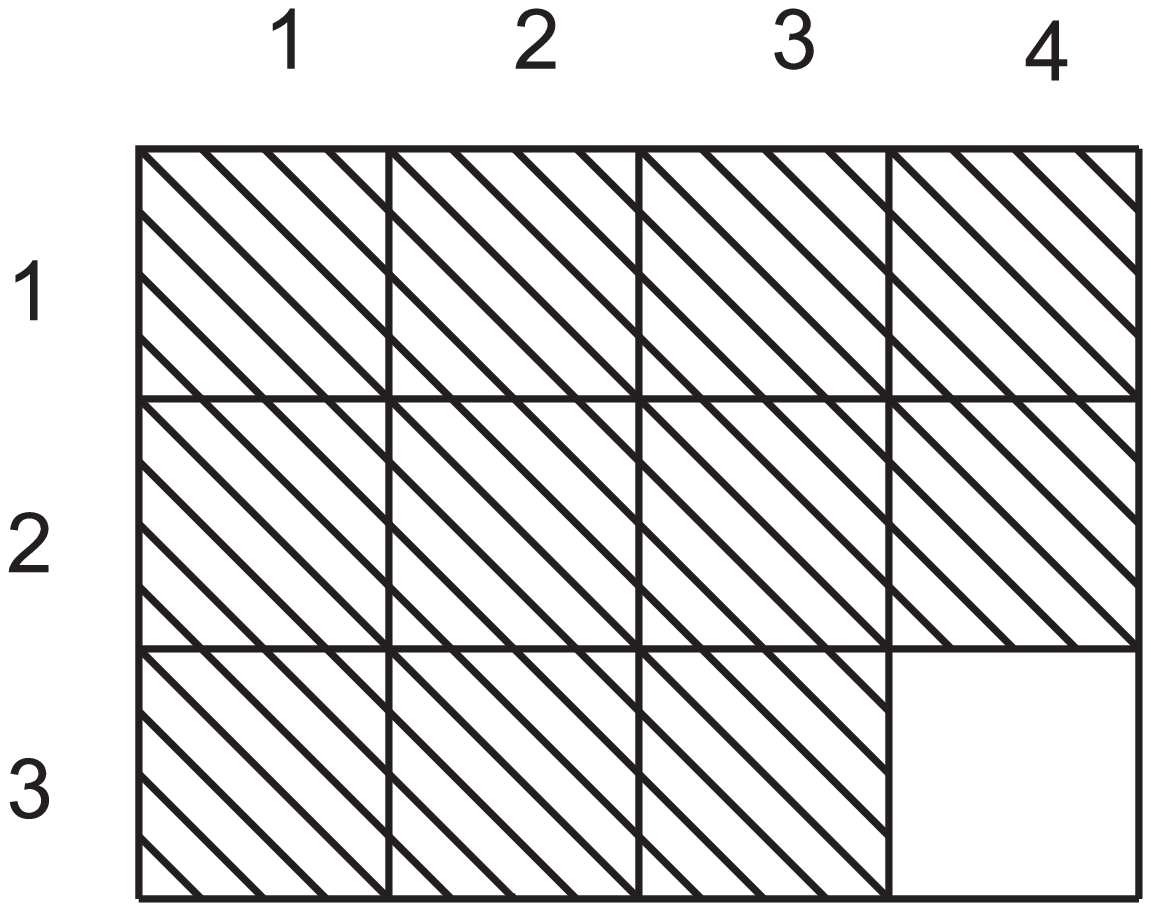}
\\
{}%
\end{center}
Borrowing our notation for the matrix and bisequence submodules of $TTM$, we denote 
the set of matrices over $G_{[n,m]}$ by $\overline{\mathbf{G}}_{[n,m]}$, the set of 
matrices over $G^{\operatorname{pre}}$ by $\overline{\mathbf{G}}$, and 
the subset of bisequence matrices in $\overline{\mathbf{G}}$ by $\mathbf{G}$.

\begin{definition}
\label{freeprematrad}Let $\Theta=\left\langle \theta_{m}^{n}\mid\theta_{1}%
^{1}=\mathbf{1}\neq0\right\rangle _{m,n\geq1}$ be a free bigraded $R$-module
generated by singletons $\theta_{m}^{n}$ and set $G_{1,1}^{\operatorname{pre}%
}=\mathbf{1}.$ Inductively, if $m+n\geq3$ and $\overline{\mathbf{G}}_{[n,m]}$
has been constructed, define%
\[
G_{n,m}^{\operatorname{pre}}=\theta_{m}^{n}\cup\left\{  \text{basic strings
}A_{s}\cdots A_{1}\text{ of bidegree }(m,n)\text{, }A_i\in \overline{\mathbf{G}}_{[n,m]},\ s\geq2\right\}  .
\]
Let $\sim$ be the equivalence relation on $G^{\operatorname{pre}}$ generated by $\left[  A_{ij}B_{ij}\right]
\sim\left[  A_{ij}\right]  \left[  B_{ij}\right]  $ if and only if $\left[
A_{ij}\right]  \times\left[  B_{ij}\right]  \in\overline{\mathbf{G}}%
\times\overline{\mathbf{G}}$ is a BTP, and let $F^{{\operatorname{pre}}%
}\left(  \Theta\right)  =\langle G^{\operatorname{pre}}\diagup\sim\rangle.$
The \textbf{free prematrad generated by} $\Theta$ is the prematrad
$(F^{{\operatorname{pre}}}(\Theta),\gamma,\eta),$ where $\gamma$ is
juxtaposition and $\eta:R\rightarrow F_{1,1}^{\operatorname{pre}}\left(
\Theta\right)  $ is given by $\eta\left(  1_{R}\right)  =\mathbf{1}$.
\end{definition}

\begin{example}
\label{non-sigma}\underline{The $\mathcal{A}_{\infty}$ operad}. Let
$\theta_{m}=\theta_{m}^{1}\neq0$ and $\theta^{n}=\theta_{1}^{n}\neq0$ for all
$m,n\geq1.$ The non-sigma operads $K=F^{^{{\operatorname{pre}}}}\left(
\theta_{\ast}\right)  $ and $L=F^{^{{\operatorname{pre}}}}(\theta^{\ast})$ are
isomorphic to the $\mathcal{A}_{\infty}$ operad and encode the combinatorial
structure of an $A_{\infty}$-algebra and an $A_{\infty}$-coalgebra,
respectively. Let $\mathbf{\theta}_{m,i}^{p}$ denote the $1\times p$ matrix
$\left[  \theta_{1}\cdots\theta_{m}\cdots\theta_{1}\right]  $ with $\theta
_{m}$ in the $i^{th}$ position, and let $\mathbf{\theta}_{q}^{n,j}$ denote the
$q\times1$ matrix $\left[  \theta^{1}\cdots\theta^{n}\cdots\theta^{1}\right]
^{T}$ with $\theta^{n}$ in the $j^{th}$ position. Then modulo prematrad axioms
(i) and (ii), the bases for $K$ and $L$ given by Definition
\ref{freeprematrad} are%
\[
\left\{  \theta_{p_{k}}\mathbf{\theta}_{m_{k},i_{k}}^{p_{k}}\cdots
\mathbf{\theta}_{m_{1},i_{1}}^{p_{1}}\in K\left(  m_{1}+p_{1}-1\right)  \mid
m_{r}=p_{r-1}-p_{r}+1\right\}
\]
and%
\[
\left\{  \mathbf{\theta}_{q_{l}}^{n_{l},j_{l}}\cdots\mathbf{\theta}_{q_{1}%
}^{n_{1},j_{1}}\theta^{q_{1}}\in L\left(  n_{l}+q_{l}-1\right)  \mid
n_{r}=q_{r+1}-q_{r}+1\right\}  .
\]

\end{example}

Given $A\in G_{n,m}^{\operatorname*{pre}}\diagup\sim$, choose a representative
$A_{s}\cdots A_{1}\in G_{n,m}^{\operatorname*{pre}}.$ In view of Definition
\ref{basic}, Axiom (iv), some association of $A_{s}\cdots A_{1}$ defines a
sequence of BTPs; thus $s-1$ successive applications of $\Upsilon$ produces a
$1\times1$ (bisequence) representative $B.$

\begin{definition}
Let $A\in G_{n,m}^{\operatorname*{pre}}\diagup\sim.$ A \textbf{factorization
of} $A$ is a representative 
\linebreak $A_{s}\cdots A_{1}\in A.$ A factorization
$A_{s}\cdots A_{1}\in A$ is a $\Theta$\text{-}\textbf{factorization} if the
entries of $A_{i}$ are elements of $\Theta$ for all $i.$ A factorization
$B_{t}\cdots B_{1}\in A$ is a \textbf{bisequence factorization }if $B_{i}%
\in\mathbf{G}$ for all $i$. A bisequence factorization $A_{s}\cdots A_{1}\in
A$ is \textbf{bisequence decomposable} if there is a bisequence factorization
$B_{t}\cdots B_{1}\in A$ such that $t>s$.
\end{definition}

Since bisequence factorizations are basic strings of matrices,
bisequence factorizations are characterized by Lemma \ref{basic lemma}:
$A_{s}\cdots A_{1}\in A$ is a bisequence factorization if and only if
$\Upsilon$ acts associatively on $A_{s}\cdots A_{1}$.

Given $A\in G^{\operatorname*{pre}}\diagup\sim,$ consider a factorization
$A_{s}\cdots A_{1}\in A.$ Each association of $A_{s}\cdots A_{1}$ determines a
fraction, any two of which will look quite different. For example, the two
associations of the bisequence $\Theta$-factorization
\begin{equation}
A_{3}A_{2}A_{1}=\left[
\begin{array}
[c]{r}%
\theta_{2}^{2}\\
\theta_{2}^{2}\\
\theta_{2}^{1}%
\end{array}
\right]  \left[
\begin{array}
[c]{cc}%
\theta_{2}^{2} & \theta_{1}^{2}\\
\theta_{2}^{1} & \mathbf{1}%
\end{array}
\right]  \left[
\begin{array}
[c]{rrr}%
\theta_{2}^{2} & \theta_{2}^{2} & \theta_{2}^{2}%
\end{array}
\right]  \in\mathbf{G}_{6}^{6},\label{assoc}%
\end{equation}
which are%
\[
\left(  A_{3}A_{2}\right)  A_{1}=\left[
\begin{array}
[c]{r}%
\left[
\begin{array}
[c]{c}%
\theta_{2}^{2}\\
\theta_{2}^{2}%
\end{array}
\right]  \left[
\begin{array}
[c]{rr}%
\theta_{2}^{2} & \theta_{1}^{2}%
\end{array}
\right] \\
\\
\left[
\begin{array}
[c]{c}%
\theta_{2}^{2}%
\end{array}
\right]  \left[
\begin{array}
[c]{cc}%
\theta_{2}^{1} & \mathbf{1}%
\end{array}
\right]
\end{array}
\right]  \left[
\begin{array}
[c]{rrr}%
\theta_{2}^{2} & \theta_{2}^{2} & \theta_{2}^{2}%
\end{array}
\right]
\]
and
\[
A_{3}\left(  A_{2}A_{1}\right)  =\left[
\begin{array}
[c]{r}%
\theta_{2}^{2}\\
\theta_{2}^{2}\\
\theta_{2}^{2}%
\end{array}
\right]  \left[
\begin{array}
[c]{ccc}%
\left[
\begin{array}
[c]{c}%
\theta_{2}^{2}\\
\theta_{2}^{1}%
\end{array}
\right]  \left[
\begin{array}
[c]{rr}%
\theta_{2}^{2} & \theta_{2}^{2}%
\end{array}
\right]  &  & \left[
\begin{array}
[c]{c}%
\theta_{1}^{2}\\
\mathbf{1}%
\end{array}
\right]  \left[
\begin{array}
[c]{c}%
\theta_{2}^{2}%
\end{array}
\right]
\end{array}
\right]  ,
\]
respectively determine the fractions
\[
\raisebox{-0.2707in}{\includegraphics[
height=0.6495in,
width=0.6746in
]%
{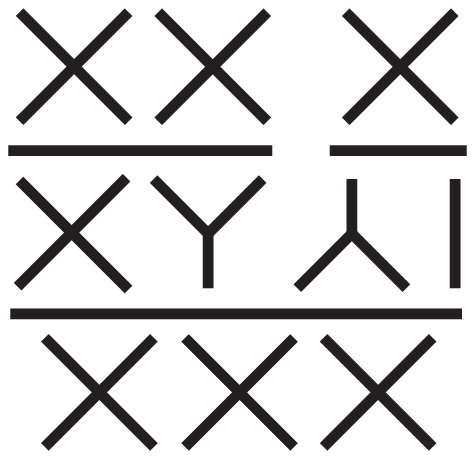}%
}
\text{ \ \ and \ }
\raisebox{-0.2707in}{\includegraphics[
height=0.6495in,
width=0.6547in
]%
{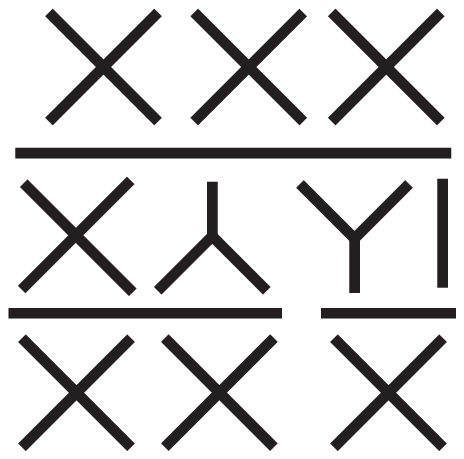}%
}
\text{ \ }.
\]
Define the \emph{graph of} $A$ to be the graph of the $1\times1$
representative $A^{\prime}\in A.$ If a fraction $f$ represents an association
of $A_{s}\cdots A_{1}\in A$, the graph of $A$ can be obtained from $f$ by
removing fraction bars and making the prescribed connections.

The arrows representing the factors of $A_{s}\cdots A_{1}\in A$ form a
polygonal path in $\mathbb{N}^{2}$ from the $x$-axis to the $y$-axis, and
evaluating subproducts in an association changes the path. For example, the
$3$-step path in Figure 5 represents the factorization $A_{3}A_{2}A_{1}$ in
(\ref{assoc}) thought of as the composition%
\[
H^{\otimes6}\overset{\mathsf{
\raisebox{-0.1946in}{\includegraphics[
trim=0.000000in 0.033928in 0.000000in 0.000000in,
height=0.1635in,
width=0.2776in
]%
{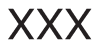}%
}
}}{\overrightarrow{\longrightarrow\overset{}{\left(  H^{\otimes2}\right)
^{\otimes3}}\overset{\sigma_{2,3}}{\longrightarrow}}}\left(  H^{\otimes
3}\right)  ^{\otimes2}\overset{\mathsf{
\raisebox{-0.0069in}{\includegraphics[
trim=-0.066614in 0.201281in 0.561893in 0.000000in,
height=0.224in,
width=0.2041in
]%
{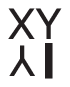}%
}%
}}{\overrightarrow{\longrightarrow\overset{}{\left(  H^{\otimes2}\right)
^{\otimes3}}\overset{\sigma_{2,2}\otimes\mathbf{1}}{\longrightarrow}}}\left(
H^{\otimes2}\right)  ^{\otimes3}\overset{\mathsf{
\raisebox{-0.0069in}{\includegraphics[
trim=0.023743in 0.078031in 0.023743in 0.000000in,
height=0.3425in,
width=0.1021in
]%
{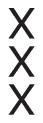}%
}
}}{\longrightarrow}H^{\otimes6},
\]
and the $2$-step paths represent the associations $\left(  A_{3}A_{2}\right)
A_{1}$ and $A_{3}\left(  A_{2}A_{1}\right)  ,$ where the subproducts in
parentheses have been evaluated and the remaining unevaluated products are
thought of as the compositions%
\[
H^{\otimes6}\overset{\mathsf{
\raisebox{-0.1946in}{\includegraphics[
trim=0.000000in 0.033928in 0.000000in 0.000000in,
height=0.1635in,
width=0.2776in
]%
{label1.eps}%
}
}}{\overrightarrow{\longrightarrow\overset{}{\left(  H^{\otimes2}\right)
^{\otimes3}}\overset{\sigma_{2,3}}{\longrightarrow}}}\left(  H^{\otimes
3}\right)  ^{\otimes2}\overset{\mathsf{
\raisebox{-0.0069in}{\includegraphics[
height=0.5431in,
width=0.1583in
]%
{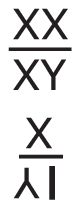}%
}
}}{\longrightarrow}H^{\otimes6}%
\]
and%
\[
H^{\otimes6}\overset{
{\includegraphics[
trim=0.000000in 0.060062in 0.000000in 0.000000in,
height=0.2318in,
width=0.3563in
]%
{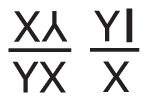}%
}
}{\overrightarrow{\longrightarrow\overset{}{\left(  H^{\otimes3}\right)
^{\otimes2}}\overset{\sigma_{3,2}}{\longrightarrow}}}\left(  H^{\otimes
2}\right)  ^{\otimes3}\overset{\mathsf{
\raisebox{-0.0069in}{\includegraphics[
trim=0.023743in 0.078031in 0.023743in 0.000000in,
height=0.3425in,
width=0.1021in
]%
{label4.eps}%
}
}}{\longrightarrow}H^{\otimes6}.
\]
\begin{center}
\includegraphics[
trim=0.000000in -0.141983in 0.215388in 0.215667in,
height=1.9778in,
width=1.8343in
]%
{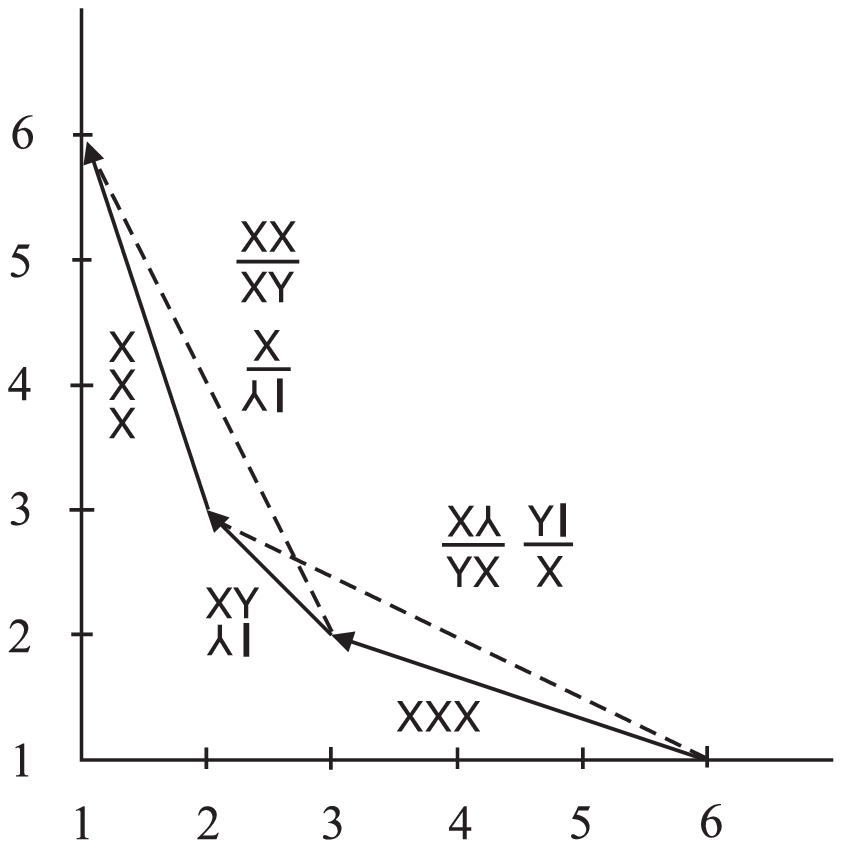}%
\\
Figure 5:\ Polygonal paths of $A_{3}A_{3}A_{1}.$%
\end{center}

In the special operadic situations of Example \ref{non-sigma}, every $\Theta
$-factorization is a bisequence factorization, but this is not true, in
general. For example, the $\Theta$-factorization%
\[
B_{3}B_{2}B_{1}=\left[
\begin{array}
[c]{r}%
\theta_{2}^{1}\\
\theta_{2}^{1}%
\end{array}
\right]  \left[
\begin{array}
[c]{cc}%
\theta_{2}^{1} & \mathbf{1}\\
\mathbf{1} & \theta_{2}^{1}%
\end{array}
\right]  \left[
\begin{array}
[c]{lll}%
\theta_{1}^{2} & \theta_{1}^{2} & \theta_{1}^{2}%
\end{array}
\right]  \in G_{2,3}^{{\operatorname{pre}}}\left(  \Theta\right)
\]
is not a bisequence factorization since $B_{2}$ is not a bisequence matrix.
Furthermore, $B_{3}B_{2}B_{1}$ only associates on the left since $B_{2}\otimes
B_{1}$ is\textit{\ not }a BTP. However, $C_{2}=B_{3}B_{2},$ $C_{1}=B_{1},$ and
$C=C_{2}C_{1}$ are bisequence matrices; hence $C_{2}C_{1}$ and $C$ are
bisequence factorizations of which%
\begin{equation}
C_{2}C_{1}=\left[
\begin{array}
[c]{c}%
\theta_{2}^{1}\left[
\begin{array}
[c]{cc}%
\theta_{2}^{1} & \mathbf{1}%
\end{array}
\right] \\
\\
\theta_{2}^{1}\left[
\begin{array}
[c]{cc}%
\mathbf{1} & \theta_{2}^{1}%
\end{array}
\right]
\end{array}
\right]  \left[
\begin{array}
[c]{ccc}%
\theta_{1}^{2} & \theta_{1}^{2} & \theta_{1}^{2}%
\end{array}
\right] \label{balanced}%
\end{equation}
is bisequence indecomposable (see Figure 6).\
\[
\raisebox{-0.32in}{\includegraphics[
trim=-0.026987in -0.054057in 0.026987in 0.054057in,
height=0.6746in,
width=0.8536in
]%
{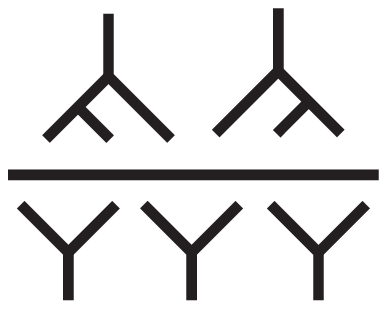}%
}
\text{ \ \ ;\ \ \ }
\raisebox{-0.3355in}{\parbox[b]{0.8838in}{\begin{center}
\includegraphics[
height=0.659in,
width=0.8838in
]%
{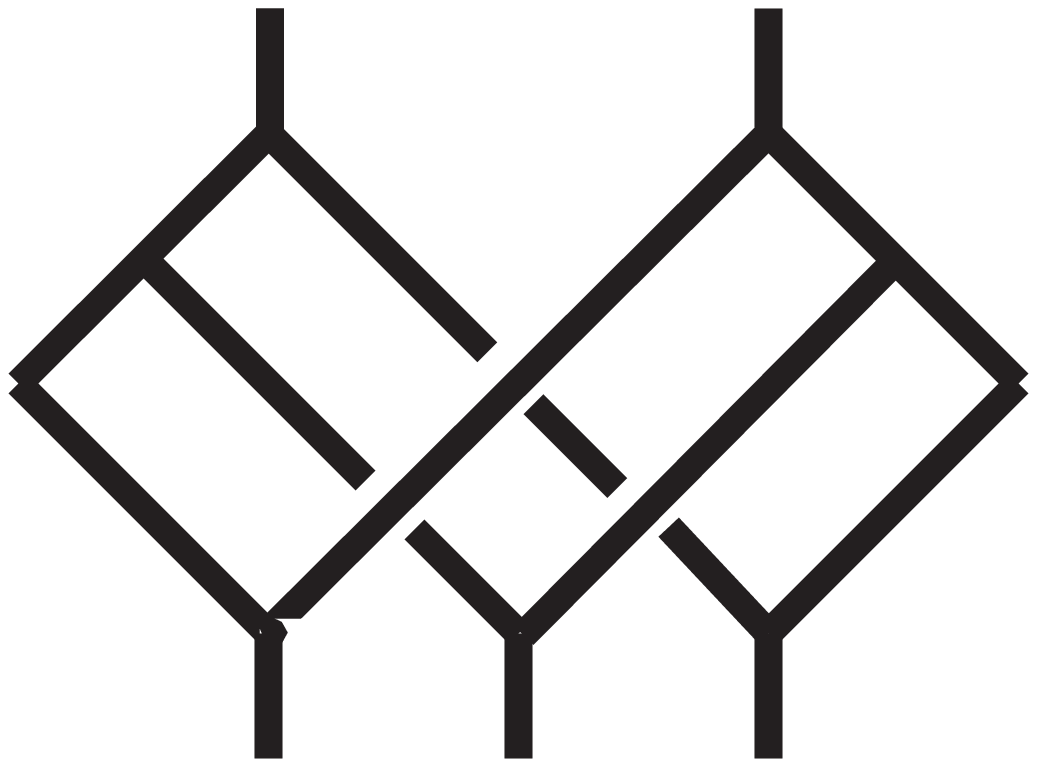}%
\\
{}%
\end{center}}}
\]

\begin{center}
Figure 6. Graphical representations of $C_{2}C_{1}$ and $C$.\vspace*{0.1in}
\end{center}

Bisequence indecomposables are especially important, and we wish to identify a
canonical bisequence indecomposable representative of a given class $A\in
G_{n,m}^{\operatorname*{pre}}\diagup\sim.$ First consider a bisequence
indecomposable $A_{k}\cdots A_{1}\in A$ with some $A_{i}\in\left\{
\theta_{\ast}^{1},\theta_{1}^{\ast}\right\}  $. If $A_{i}=\theta_{s}^{1}$ for
some $s,$ then $A_{i}\cdots A_{1}$ is identified with an up-rooted planar
rooted leveled trees (PLT) with $i$ levels and $m$ leaves; dually, if
$A_{i}=\theta_{1}^{t}$ for some $t,$ then $A_{k}\cdots A_{i}$ is identified
with a down-rooted PLT with $k-i+1$ levels and $n$ leaves. Note that both
situations occur in factorizations of the form $A_{k}\cdots\theta_{1}%
^{t}\theta_{s}^{1}\cdots A_{1}$. But in either case, the indicated PLT
represents an isomorphism class of PLTs whose elements determine distinct
bisequence indecomposable factorizations of $A$. For example, the isomorphic
PLTs%
\[
\raisebox{-0.1323in}{\includegraphics[
height=0.3347in,
width=0.3442in
]%
{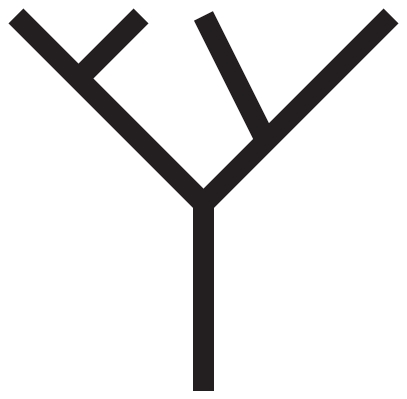}%
}
\text{ }:\text{ }\bullet\bullet\bullet\bullet\rightarrow(\bullet
\bullet)\bullet\bullet\rightarrow(\bullet\bullet)(\bullet\bullet)
\]
and%
\[
\raisebox{-0.1254in}{\includegraphics[
height=0.3347in,
width=0.3442in
]%
{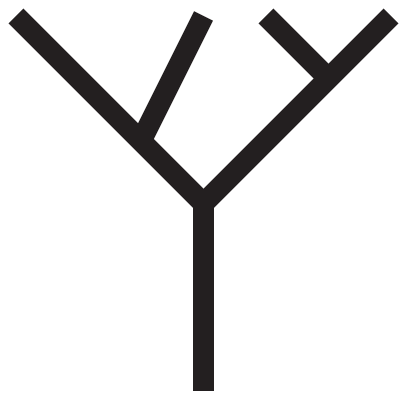}%
}
\text{ }:\text{ }\bullet\bullet\bullet\bullet\rightarrow\bullet\bullet\left(
\bullet\bullet\right)  \rightarrow(\bullet\bullet)(\bullet\bullet)
\]
respectively determine the distinct equivalent bisequence indecomposables%
\[
A_{4}A_{3}A_{2}A_{1}=\left[  \theta_{1}^{2}\text{ }\mathbf{1}\text{
}\mathbf{1}\right]  ^{T}\left[  \mathbf{1}\text{ }\theta_{1}^{2}\right]
^{T}\left[  \theta_{1}^{2}\right]  \left[  \theta_{3}^{1}\right]
\]
and
\[
A_{4}^{\prime}A_{3}^{\prime}A_{2}^{\prime}A_{1}^{\prime}=\left[
\mathbf{1}\text{ }\mathbf{1}\text{ }\theta_{1}^{2}\right]  ^{T}\left[
\theta_{1}^{2}\text{ }\mathbf{1}\right]  ^{T}\left[  \theta_{1}^{2}\right]
\left[  \theta_{3}^{1}\right]
\]
in $G_{4,3}^{\operatorname{pre}}.$

Now recall that a PLT with $s$ levels and $t$ leaves specifies the order in
which $s$ pairs of parentheses are inserted into a string of $t$
indeterminants. When faced with the situations described above, \emph{choose
the factorization indexed by the PLT that successively inserts parentheses as
far to the left as possible. }Such factorizations have \emph{preferred
operadic substrings}.

\begin{definition}
A \textbf{preferred factorization}\emph{\ }is a\emph{\ }factorization with
preferred oper- adic substrings (if any) and factors whose entries have this
property. A bisequence matrix $B\in G^{\operatorname{pre}}$ is
\textbf{balanced} if each entry of $B$ is a preferred bisequence
indecomposable factorization. A \textbf{balanced factorization} of $A\in
G^{\operatorname{pre}}\diagup\sim$ is a preferred bisequence indecomposable
factorization $A_{s}\cdots A_{1}$ in which $A_{i}$ is balanced for all $i$.
\end{definition}

In the example above, $A_{4}A_{3}A_{2}A_{1}$ is the (obviously unique)
balanced factorization of $A\in G_{4,3}^{\operatorname{pre}}\diagup\sim$. In
view of Definition \ref{prematrad}, Axiom (i), two associations of a basic
string $A_{s}\cdots A_{1}$ that define sequences of BTPs are equal. Hence each
class $A\in G_{n,m}^{\operatorname{pre}}\diagup\sim$ has a unique preferred
factorization of maximal length. In fact:

\begin{proposition}
Each class $A\in G_{n,m}^{\operatorname{pre}}\diagup\sim$ has a unique
balanced factorization.
\end{proposition}

\begin{proof}
If $m+n\leq4,$ balanced factorizations of $A$ are $\Theta$-factorizations, and
obviously unique. Inductively, assume the statement holds for all $B\in
G_{l,k}^{\operatorname{pre}}\diagup\sim$ with $k\leq m,$ $l\leq n$ and
$k+l<m+n,$ and consider a class $A\in G_{n,m}^{\operatorname{pre}}\diagup
\sim.$ If $\theta_{m}^{n}\in A,$ then $A$ is a balanced singleton class.
Otherwise, consider the unique preferred class representative $A_{s}\cdots
A_{1}$ of maximal length. If $A_{s}\cdots A_{1}$ is a bisequence
factorization, it is balanced. If not, evaluate $\Upsilon$-products and obtain
a bisequence factorization $B_{1}B_{2}$ of $A.$ Either $B_{1}=A_{s}$ or
$B_{1}=A_{s}B_{1}^{\prime}.$ In either case, there is a bisequence
factorization $A_{s}B$ of $A$. If $B$ is bisequence indecomposable, $A_{s}B$
is the unique balanced factorization of $A.$ If not, there is a bisequence
factorization $C_{1}C_{2}$ of $B.$ If $C_{1}$ is bisequence indecomposable,
consider $C_{2};$ otherwise, decompose $C_{1}.$ Continue in this manner until
the process terminates.
\end{proof}

We now construct a set that conveniently indexes the module generators
$G^{\operatorname{pre}}\diagup\sim$ for $F^{\operatorname{pre}}\left(
\Theta\right)  .$ Define a map $\phi$ that splits the projection
$G^{\operatorname{pre}}\rightarrow G^{\operatorname{pre}}\diagup\sim$ as
follows: For each pair $\left(  m,n\right)  $ with $m+n\leq3, $ set
$\phi(\operatorname{cls}\theta_{m}^{n})=\theta_{m}^{n}.$ Inductively, for each
pair $\left(  m,n\right)  $ with $m+n\geq4,$ assume that $\phi$ has been
defined on $G_{j,i}^{\operatorname{pre}}\diagup\sim$ for $i+j<m+n.$ Then given
a class $g\in G_{n,m}^{\operatorname{pre}}\diagup\sim,$ define $\phi_{n,m}(g)$
to be the balanced factorization of $g.$ Let $\mathcal{B}_{n,m}%
^{\operatorname{pre}}=\operatorname*{Im}\phi_{n,m}$ and $\mathcal{B}%
^{\operatorname{pre}}=%
{\textstyle\bigcup}
\mathcal{B}_{n,m}^{\operatorname{pre}}.$

Although a balanced factorization $\beta\in\mathcal{B}^{\operatorname{pre}}$
is not a $\Theta$-factorization, there is a related PLT $\Psi\left(
\beta\right)  $ whose leaves are $\Theta$-factorizations. Let $\beta
=B_{s}\cdots B_{1}\in\mathcal{B}^{\operatorname{pre}}$ and set $\Psi
_{1}\left(  \beta\right)  =\beta.$ If $\beta$ is a $\Theta$-factorization, set
$\Psi\left(  \beta\right)  =\Psi_{1}\left(  \beta\right)  .$ Otherwise, let
$\left(  \beta_{k}\right)  $ denote the tuple of entries of $B_{k}$ listed in
row order and replace each entry $\beta_{k}$ of $\left(  \beta_{k}\right)  $
with its balanced factorization $\beta_{k}^{\prime}\in\mathcal{B}%
^{\operatorname{pre}}$. Form the $2$-level tree $\Psi_{2}\left(  \beta\right)
$ with root $\beta$ and leaves labeled by the entries of $\left(  \beta
_{k}^{\prime}\right)  $ (see Example \ref{3-level}). If each $\beta
_{k}^{\prime}$ is a $\Theta$-factorization, set $\Psi\left(  \beta\right)
=\Psi_{2}\left(  \beta\right)  $. Otherwise, repeat this process for each
$\beta_{k}^{\prime},$ i.e., let $\Psi_{1}\left(  \beta_{k}^{\prime}\right)
=\beta_{k}^{\prime}; $ if $\beta_{k}^{\prime}$ that is not a $\Theta
$-factorization, construct $\Psi_{2}\left(  \beta_{k}^{\prime}\right)  .$ Now
construct $3$-level tree $\Psi_{3}\left(  \beta\right)  $ either by appending
the tree $\Psi_{2}\left(  \beta_{k}^{\prime}\right)  $ to the leaf $\beta
_{k}^{\prime}$ or by extending the branch of $\beta_{k}^{\prime}$ otherwise.
If each level $3$ entry of $\Psi_{3}\left(  \beta\right)  $ is a $\Theta
$-factorization, set $\Psi\left(  \beta\right)  =\Psi_{3}\left(  \beta\right)
$; otherwise continue inductively. This process terminates after $r$ steps and
uniquely determines an $r$-level tree $\Psi\left(  \beta\right)  $ whose
leaves are balanced $\Theta\,$-factorizations.

\begin{example}
\label{3-level}The balanced factorization
\[
\beta=\left[
\begin{array}
[c]{c}%
\theta_{2}^{1}\left[
\begin{array}
[c]{cc}%
\theta_{2}^{1} & \theta_{1}^{1}%
\end{array}
\right] \\
\\
\theta_{2}^{1}\left[
\begin{array}
[c]{cc}%
\theta_{1}^{1} & \theta_{2}^{1}%
\end{array}
\right]
\end{array}
\right]  \left[
\begin{array}
[c]{ccc}%
\theta_{1}^{2} & \theta_{1}^{2} & \theta_{1}^{2}%
\end{array}
\right]
\]
is associated with the $2$-level tree%
\[
\Psi\left(  \beta\right)  =%
\begin{array}
[c]{ccccc}
&  & \beta &  & \\
&  &  &  & \\
&
\raisebox{-0.083in}{\includegraphics[
height=0.2534in,
width=0.2534in
]%
{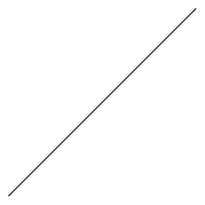}%
}
\hspace*{-0.2in} & \text{{\LARGE
$\vert$
}} & \hspace*{-0.8in}
\raisebox{-0.0346in}{\includegraphics[
height=0.2093in,
width=0.8778in
]%
{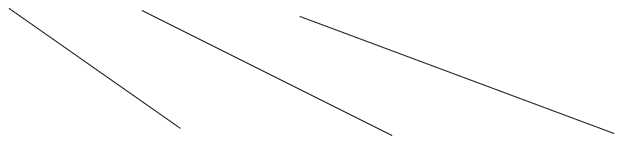}%
}
& \\
&  &  &  & \\
\theta_{2}^{1}\left[
\begin{array}
[c]{cc}%
\theta_{2}^{1} & \theta_{1}^{1}%
\end{array}
\right]  \hspace*{-0.2in} &  & \theta_{2}^{1}\left[
\begin{array}
[c]{cc}%
\theta_{1}^{1} & \theta_{2}^{1}%
\end{array}
\right]  & \text{ \ }%
\begin{array}
[c]{ccc}%
\theta_{1}^{2} & \theta_{1}^{2} & \theta_{1}^{2}%
\end{array}
& \text{.}%
\end{array}
\]
\medskip
\end{example}

Finally, let $\mathcal{C}_{n,m}^{\operatorname{pre}}=\left\{  \Psi\left(
\beta\right)  \mid\beta\in\mathcal{B}_{n,m}^{\operatorname{pre}}\right\}  $;
then $\mathcal{C}^{\operatorname{pre}}=\bigcup_{m,n}\mathcal{C}_{n,m}%
^{\operatorname{pre}}$ indexes the set of module generators for
$F^{\operatorname{pre}}\left(  \Theta\right)  .$ In Subsection 6.2 we identify
a subset $\mathcal{C\subset C}^{\operatorname{pre}}$ whose elements
simultaneously index module generators of the \textquotedblleft free
matrad\textquotedblright\ $F\left(  \Theta\right)  $ and cells of the
biassociahedra $KK$. This identification relates the module structure of
$F\left(  \Theta\right)  $ to the combinatorics of permutahedra.

\subsection{The Bialgebra Prematrad}

As is the case for operads and PROPs, prematrads can be described by
generators and relations.

\begin{definition}
Let $\mathit{{\mathcal{H}}}^{{\operatorname{pre}}}=F^{{\operatorname{pre}}%
}\left(  \mathbf{\Theta}\right)  /\sim,$ where $\mathbf{\Theta=}\left\langle
\theta_{1}^{1}=\mathbf{1,}\theta_{2}^{1},\theta_{1}^{2}\right\rangle ,$
$\theta_{i}^{j}\neq0,$ and $A\sim B$ if and only if ${\operatorname{bideg}%
}(A)={\operatorname{bideg}}(B).$ Let $\eta\left(  1\right)  =\mathbf{1}$ and
let $\gamma$ be the structure map induced by projection. Then $\left(
\mathit{{\mathcal{H}}}^{{\operatorname{pre}}},\eta,\gamma\right)  $ is the
\textbf{bialgebra prematrad} if the following axioms hold:

\begin{enumerate}
\item[\textit{(i)}] Associativity: $\gamma(\theta_{2}^{1};\theta_{2}%
^{1},\mathbf{1})=\gamma(\theta_{2}^{1};\mathbf{1},\theta_{2}^{1})$.

\item[\textit{(ii)}] Coassociativity: $\gamma(\theta_{1}^{2},\mathbf{1}%
;\theta_{1}^{2})=\gamma(\mathbf{1},\theta_{1}^{2};\theta_{1}^{2})$.

\item[\textit{(iii)}] Hopf compatibility: $\gamma(\theta_{1}^{2};\theta
_{2}^{1})=\gamma(\theta_{2}^{1},\theta_{2}^{1};\theta_{1}^{2},\theta_{1}^{2})$.
\end{enumerate}
\end{definition}

Each bigraded component ${\mathcal{H}}_{n,m}^{{\operatorname{pre}}}$ is
generated by a singleton $c_{n,m};$ for example, $c_{1,2}=\theta_{2}^{1},$
$c_{2,1}=\theta_{1}^{2},$ $c_{1,3}=\gamma(\theta_{2}^{1};\theta_{2}%
^{1},\mathbf{1}),$ $c_{3,1}=\gamma(\theta_{1}^{2},\mathbf{1};\theta_{1}^{2}),$
$c_{2,2}=\gamma(\theta_{1}^{2};\theta_{2}^{1}),$ and so on. Note that
${\mathcal{H}}_{1,\ast}^{{\operatorname{pre}}}$ and ${\mathcal{H}}_{\ast
,1}^{{\operatorname{pre}}}$ are operads; the first is generated by $\theta
_{2}^{1}$ subject to (\textit{i}) while the second is generated by $\theta
_{1}^{2}$ subject to (\textit{ii}). Both are isomorphic to the associativity
operad $\underline{{\mathcal{A}}ss}$ (\cite{MSS}). And furthermore, for the
bialgebra PROP $\mathcal{B}$ we have $\mathcal{B}^{\operatorname*{pre}%
}=\mathcal{H}^{\operatorname*{pre}} $ and $U(\mathcal{H}^{\operatorname*{pre}%
})=\mathcal{B}.$

Given a graded $R$-module $H,$ a map of prematrads ${\mathcal{H}%
}^{{\operatorname{pre}}}\rightarrow U_{H}$ defines a bialgebra structure on
$H$ and vise versa. Since each path of arrows from $(m,1)$ to $(1,n)$ in
$\mathbb{N}^{2}$ represents some $\Upsilon$-factorization of $c_{n,m}$ (see
Figure 5), we think of all such paths as equal. Therefore ${\mathcal{H}%
}^{{\operatorname{pre}}}$ is the smallest module among existing general
constructions that describe bialgebras (c.f. \cite{Markl2}, \cite{borya}).
Although the symmetric groups do not act on prematrads, the permutation
$\sigma_{n,m}$ built into the associativity axiom minimizes the modules involved.

\subsection{Local Prematrads}

Let $M=\left\{  M_{n,m}\right\}  _{m,n\geq1}$ be a bigraded $R$-module, let
$\mathbf{W}$ be a telescoping submodule of $TTM,$ let $\mathcal{W}$ be its
telescopic extension, and let $\gamma_{_{\mathbf{W}}}=\left\{  \gamma
_{\mathbf{x}}^{\mathbf{y}}:\mathbf{W}_{p}^{\mathbf{y}}\otimes\mathbf{W}%
_{\mathbf{x}}^{q}\rightarrow\mathbf{W}_{\left\vert \mathbf{x}\right\vert
}^{\left\vert \mathbf{y}\right\vert }\right\}  $ be a structure map. If
$A\otimes B$ is a BTP in $\mathcal{W}\otimes\mathcal{W}$, each TP $A^{\prime
}\otimes B^{\prime}$ in $A\otimes B$ lies in $\mathbf{W}_{p}^{\mathbf{y}%
}\otimes\mathbf{W}_{\mathbf{x}}^{q}$ for some $\mathbf{x},\mathbf{y},p,q.$
Consequently, $\gamma_{_{\mathbf{W}}}$ extends to a global product
$\Upsilon:\mathcal{W}\otimes\mathcal{W}\rightarrow\mathcal{W}$ as in
(\ref{upsilon}). In fact, $\mathcal{W}$ is the smallest matrix submodule
containing $\mathbf{W}$ on which $\Upsilon$ is well-defined.

\begin{definition}
Let $\mathbf{W}$ be a telescoping submodule of $TTM$, let $\gamma
_{_{\mathbf{W}}}=\left\{  \gamma_{\mathbf{x}}^{\mathbf{y}}:\right.
$\linebreak$\left.  \mathbf{W}_{p}^{\mathbf{y}}\otimes\mathbf{W}_{\mathbf{x}%
}^{q}\rightarrow\mathbf{W}_{\left\vert \mathbf{x}\right\vert }^{\left\vert
\mathbf{y}\right\vert }\right\}  $ be a structure map, and let $\eta
:R\rightarrow\mathbf{M}$. The triple $\left(  M,\gamma_{_{\mathbf{W}}}%
,\eta\right)  $ is a \textbf{local prematrad (with domain} $\mathbf{W}$) if
the following axioms are satisfied:

\begin{enumerate}
\item[\textit{(i)}] $\mathbf{W}_{\mathbf{x}}^{1}=\mathbf{M}_{\mathbf{x}}^{1}$
and $\mathbf{W}_{1}^{\mathbf{y}}=\mathbf{M}_{1}^{\mathbf{y}}$ for all
$\mathbf{x},\mathbf{y};$

\item[\textit{(ii)}] $\Upsilon$\ is associative on $\mathcal{W}\cap\mathbf{M}$;

\item[\textit{(iii)}] the prematrad unit axiom holds for $\gamma
_{_{\mathbf{W}}}.$
\end{enumerate}
\end{definition}

\begin{example}
\label{extremes}Let $\mathbf{V}$ be the bisequence vector submodule of $TTM$.
If $\left(  M,\gamma\right)  $ is a prematrad, then $\gamma=\gamma
_{_{\mathbf{V}}}$ and $\left(  M,\gamma\right)  $ is a local prematrad with
domain $\mathbf{V}$. Local subprematrads $\left(  M,\gamma_{\mathbf{W}%
}\right)  \subset\left(  M,\gamma\right)  $ are obtained by restricting the
domain to submodules such as $\mathbf{M}_{1}^{\ast},$ $\mathbf{M}_{\ast}^{1},$
and $\mathbf{M}_{\ast}^{1}\cup\mathbf{M}_{1}^{\ast}.$ Note that $\left(
M,\gamma_{_{\mathbf{M}_{1}^{\ast}}}\right)  $ and $\left(  M,\gamma
_{_{\mathbf{M}_{\ast}^{1}}}\right)  $ are operads.
\end{example}

\section{Diagonal Approximations}

A diagonal approximation $\Delta_{X}$ on a cellular complex $X$ determines a
\textquotedblleft$k$-subdivi- sion\textquotedblright\ $X^{\left(  k\right)  }$
of $X$ and a cellular inclusion $\Delta^{\left(  k\right)  }:X^{\left(
k\right)  }\hookrightarrow X^{\times k+1}$ whose image is the subcomplex of
$X^{\times k+1}$ we denote by $\Delta^{\left(  k\right)  }\left(  X\right)  .$
In particular, the subcomplex $\Delta^{\left(  k\right)  }\left(
P_{n}\right)  \subset P_{n}^{\times k+1}$ determined by the S-U diagonal
$\Delta_{P}$ defines the selection rule in Subsection \ref{low-products} and,
more generally, in the next section.

Recall that a map $f:X\rightarrow Y$ of $CW$-complexes is homotopic to a
cellular map $g:X\rightarrow Y,$ which in turn induces a chain map $g:C_{\ast
}\left(  X\right)  \rightarrow C_{\ast}\left(  Y\right)  .$ Given a geometric
diagonal $\Delta:X\rightarrow X\times X$, a cellular map $\Delta
_{X}:X\rightarrow X\times X$ homotopic to $\Delta$ is called a \emph{diagonal
approximation.} A diagonal approximation $\Delta_{X}$ induces a chain map
$\Delta_{X}:C_{\ast}\left(  X\right)  \rightarrow C_{\ast}\left(  X\right)
\otimes C_{\ast}\left(  X\right)  ,$ called a \emph{diagonal}. A brief review
of the S-U diagonals $\Delta_{P}$ and $\Delta_{K}$ on cellular chains of
permutahedra $P=\sqcup_{n\geq1}P_{n}$ and associahedra $K=\sqcup_{n\geq2}%
K_{n}$ (up to sign) now follows (see \cite{SU2} for details).

\subsection{The S-U Diagonals $\Delta_{P}$ and $\Delta_{K}$\label{S-U}}

Let $\underline{n}=\{1,2,\dots,n\},$ $n\geq1.$ A matrix $E$ with entries from
$\left\{  0\right\}  \cup\underline{n}$ is a \emph{step matrix} if the
following conditions hold:

\begin{enumerate}
\item[\textit{(i)}] Each element of $\underline{n}$ appears as an entry of $E
$ exactly once.

\item[\textit{(ii)}] Elements of $\underline{n}$ in each row and column of $E
$ form an increasing contiguous block.

\item[\textit{(iii)}] Each diagonal parallel to the main diagonal of $E$
contains exactly one element of $\underline{n}.$
\end{enumerate}

\noindent The non-zero entries in a step matrix form a continuous staircase
connecting the lower-left and upper-right most entries. There is a bijective
correspondence between step matrices with non-zero entries in $\underline{n}$
and permutations of $\underline{n}$.

Given a $q\times p$ integer matrix $M=\left(  m_{ij}\right)  ,$ choose proper
subsets $S_{i}\subset\left\{  \text{non-zero }\right.  $ $\left.
\text{entries in }\operatorname*{row}\left(  i\right)  \right\}  $and
$T_{j}\subset\left\{  \text{non-zero entries in }\operatorname*{col}\left(
j\right)  \right\}  ,$ and define \emph{down-shift} and \emph{right-shift}
operations\emph{\ }$D_{S_{i}} $ and $R_{T_{j}}$ on $M$ as follows:

\begin{enumerate}
\item[\textit{(i)}] If $S_{i}\neq\varnothing,$ $\max\operatorname*{row}%
(i+1)<\min S_{i}=m_{ij},$ and $m_{i+1,k}=0$ for all $k\geq j$, then $D_{S_{i}%
}M$ is the matrix obtained from $M$ by interchanging each $m_{ik}\in S_{i}$
with $m_{i+1,k};$ otherwise $D_{S_{i}}M=M.$

\item[\textit{(ii)}] If $T_{j}\neq\varnothing,$ $\max\operatorname*{col}%
(j+1)<\min T_{j}=m_{ij},$ and $m_{k,j+1}=0$ for all $k\geq i,$ then\emph{\ }%
$R_{T_{j}}M$ is the matrix obtained from $M$ by interchanging each $m_{k,j}\in
T_{j}$ with $m_{k,j+1};$ otherwise $R_{Tj}M=M.$
\end{enumerate}

\noindent Given a $q\times p$ step matrix $E\ $together with subsets
$S_{1},\ldots,S_{q}$ and $T_{1},\ldots,T_{p}$ as above, there is the
\emph{derived matrix}%
\[
R_{T_{p}}\cdots R_{T_{2}}R_{T_{1}}D_{S_{q}}\cdots D_{S_{2}}D_{S_{1}}E.
\]
In particular, step matrices are derived matrices under the trivial action
with $S_{i}=T_{j}=\varnothing$ for all $i,j$.

Let $a=A_{1}|A_{2}|\cdots|A_{p}$ and $b=B_{q}|B_{q-1}|\cdots|B_{1}$ be
partitions of $\underline{n}.$ Then $a\times b$ is a $(p,q)$%
\emph{-complementary pair} (CP) if there is a $q\times p$ derived matrix
$M=\left(  m_{ij}\right)  $ such that $A_{j}=\{m_{ij}\neq0\mid1\leq i\leq q\}$
and $B_{i}=\{m_{ij}\neq0\mid1\leq j\leq p\}.$ Thus $\left(  p,q\right)  $-CPs,
which are in one-to-one correspondence with derived matrices, identify a
particular set of product cells in $P_{n}\times P_{n}$.

\begin{definition}
Define $\Delta_{P}:C_{0}\left(  P_{1}\right)  \rightarrow C_{0}\left(
P_{1}\right)  \otimes C_{0}\left(  P_{1}\right)  $ by $\Delta_{P}%
(\underline{1})=\underline{1}\otimes\underline{1}$. Inductively, having
defined $\Delta_{P}:C_{\ast}\left(  P_{k}\right)  \rightarrow C_{\ast}\left(
P_{k}\right)  \otimes C_{\ast}\left(  P_{k}\right)  $ for all $k\leq n$,
define $\Delta_{P} $ on $\underline{n+1}\in C_{n}(P_{n+1})$ by
\[
\Delta_{P}(\underline{n+1})=\sum_{\substack{(p,q)\text{-CPs }a\times b
\\p+q=n+2}}\pm\text{ }a\otimes b
\]
and extend $\Delta_{P}$ multiplicatively, i.e., $\Delta_{P}\left(
u_{1}|\cdots|u_{r}\right)  =\Delta_{P}\left(  u_{1}\right)  |\cdots|\Delta
_{P}\left(  u_{r}\right)  .$
\end{definition}

The diagonal $\Delta_{p}$ induces a diagonal $\Delta_{K}$ on $C_{\ast}\left(
K\right)  $. Recall that faces of $P_{n}$ in codimension $k$ are indexed by
planar rooted trees with $n+1$ leaves and $k+1$ levels (PLTs), and forgetting
levels defines the cellular projection $\vartheta_{0}:P_{n}\rightarrow
K_{n+1}$ given by A. Tonks \cite{Tonks}. Thus faces of $P_{n}$ indexed by PLTs
with multiple nodes in the same level degenerate under $\vartheta_{0}$, and
corresponding generators span the kernel of the induced map $\vartheta
_{0}:C_{\ast}\left(  P_{n}\right)  \rightarrow C_{\ast}\left(  K_{n+1}\right)
$. The diagonal $\Delta_{K}$ is given by
\[
\Delta_{K}\vartheta_{0}=(\vartheta_{0}\otimes\vartheta_{0})\Delta_{P}.
\]

\subsection{$k$-Subdivisions and $k$-Approximations\label{topological}}

When $X$ is a polytope one can choose a diagonal approximation $\Delta
_{X}:X\rightarrow X\times X$ such that

\begin{enumerate}
\item[\textit{(i)}] $\Delta_{X}$ acts on each face $e\subset X$ as a
(topological) inclusion $\Delta_{X}(e)\subset e\times e,$ and

\item[\textit{(ii)}] there is an induced (cellular) $1$\emph{-subdivision
}$X^{(1)}$ of $X$ that converts $\Delta_{X}$ into a cellular inclusion
$\Delta_{X}^{(1)}:X^{(1)}\rightarrow X\times X.$
\end{enumerate}

\noindent If $e=\bigcup_{i=1}^{m}e_{i}$ and $e_{i}\subseteq X^{\left(
1\right)  },$ there are faces $u_{i},v_{i}\subseteq e_{i}$ such that
${\Delta_{X}^{(1)}}(e_{i})=u_{i}\times v_{i}.$ Thus $\Delta_{X}(e)=\bigcup
_{i=1}^{m}{\Delta_{X}^{(1)}}(e_{i})=\bigcup_{i=1}^{m}u_{i}\times v_{i},$ and
in particular, $\Delta_{X}^{(1)}$ agrees with the geometric diagonal $\Delta$
only on vertices of $X.$

The $1$-subdivision $X^{\left(  1\right)  }$ arising from an explicit diagonal
approximation $\Delta_{X}$ can be thought of as the cell complex obtained by
gluing the cells in $\Delta_{X}\left(  X\right)  $ together along their common
boundaries in the only way possible. For example, the A-W diagonal on the
simplex $\Delta^{n}$, the Serre diagonal on the cube $I^{n},$ and the S-U
diagonals on the permutahedron $P_{n+1}$ and the associahedron $K_{n+2}$ (see
\cite{SU2}) induce explicit 1-subdivisions $\left(  \Delta^{n}\right)
^{\left(  1\right)  },$\ $\left(  I^{n}\right)  ^{\left(  1\right)  },$
$P_{n+1}^{(1)}$ and $K_{n+2}^{(1)}$ (see Figures 7-10), and it is a good
exercise to determine how the vertices of $\left(  I^{n}\right)  ^{\left(
1\right)  }$ resolve in $P_{n+1}^{(1)}$ and $K_{n+2}^{(1)}$.

Algebraically, the assignment
\begin{equation}
e\mapsto\{u_{i}\times v_{i}\}_{1\leq i\leq m}\label{celldecomp}%
\end{equation}
determines a DG\ diagonal $\Delta_{X}:C_{\ast}\left(  X\right)  \rightarrow
C_{\ast}\left(  X\right)  \otimes C_{\ast}\left(  X\right)  $ on cellular
chains such that $\Delta_{X}\left(  C_{\ast}(X)\right)  \subseteq C_{\ast
}\left(  \Delta_{X}^{(1)}(X^{(1)})\right)  ,$ where equality holds in the
(unique) case $X=\ast$. Conversely, if a DG diagonal approximation $\Delta
_{X}$ on $C_{\ast}\left(  X\right)  $ is determined by a cellular
decomposition as in (\ref{celldecomp}), there is a corresponding
$1$-subdivision $X^{\left(  1\right)  }$ and a $1$\emph{-approximation}
$\Delta_{X}^{(1)}:X^{(1)}\rightarrow X\times X.$

Furthermore, there is a $2$-\emph{subdivision} $X^{\left(  2\right)  }$ of $X$
(see Figure 11) and a corresponding $2$-\emph{approximation} $\tilde{\Delta
}_{X}^{(2)}:X^{\left(  2\right)  }\rightarrow X^{\left(  1\right)  }\times X$
that sends each cell of $X^{\left(  2\right)  }$ onto a single cell of
$X^{\left(  1\right)  }\times X;$ consequently, the composition $\Delta
_{X}^{(2)}=\left(  \Delta_{X}^{(1)}\times\mathbf{1}\right)  \tilde{\Delta}%
_{X}^{(2)}$ sends each cell of $X^{(2)}$ onto a single cell of $X^{\times3}. $
Inductively, for each $k,$ there is a $k$-subdivision $X^{\left(  k\right)  }
$ and a $k$-approximation $\tilde{\Delta}_{X}^{(k)}:X^{\left(  k\right)
}\rightarrow X^{\left(  k-1\right)  }\times X$ such that $\Delta_{X}%
^{(k)}=\left(  \Delta_{X}^{(k-1)}\times\mathbf{1}\right)  \tilde{\Delta}%
_{X}^{(k)}$ sends each cell of $X^{(k)}$ onto a single cell of $X^{\times
k+1}.$ Thus, for each $k\geq0,$ a diagonal approximation $\Delta_{X}$
$\emph{fixes}$ the subcomplex $\Delta^{(k)}\left(  X\right)  :=\Delta
_{X}^{(k)}(X^{(k)})\subset X^{\times k+1},$ which says that $\Delta_{X}$ acts
on $\Delta^{(k)}\left(  X\right)  $ as an inclusion.
\vspace{0.2in}

\begin{center}
\includegraphics[
height=0.563in,
width=1.8862in
]%
{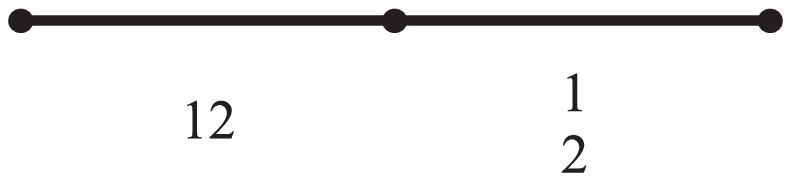}%
\\
Figure 7. The 1-subdivision of $P_{2}=K_{3}=I.$%
\end{center}

\begin{center}
\includegraphics[
trim=0.000000in -0.745955in 0.000000in 0.000000in,
height=2.2393in,
width=4.783in
]%
{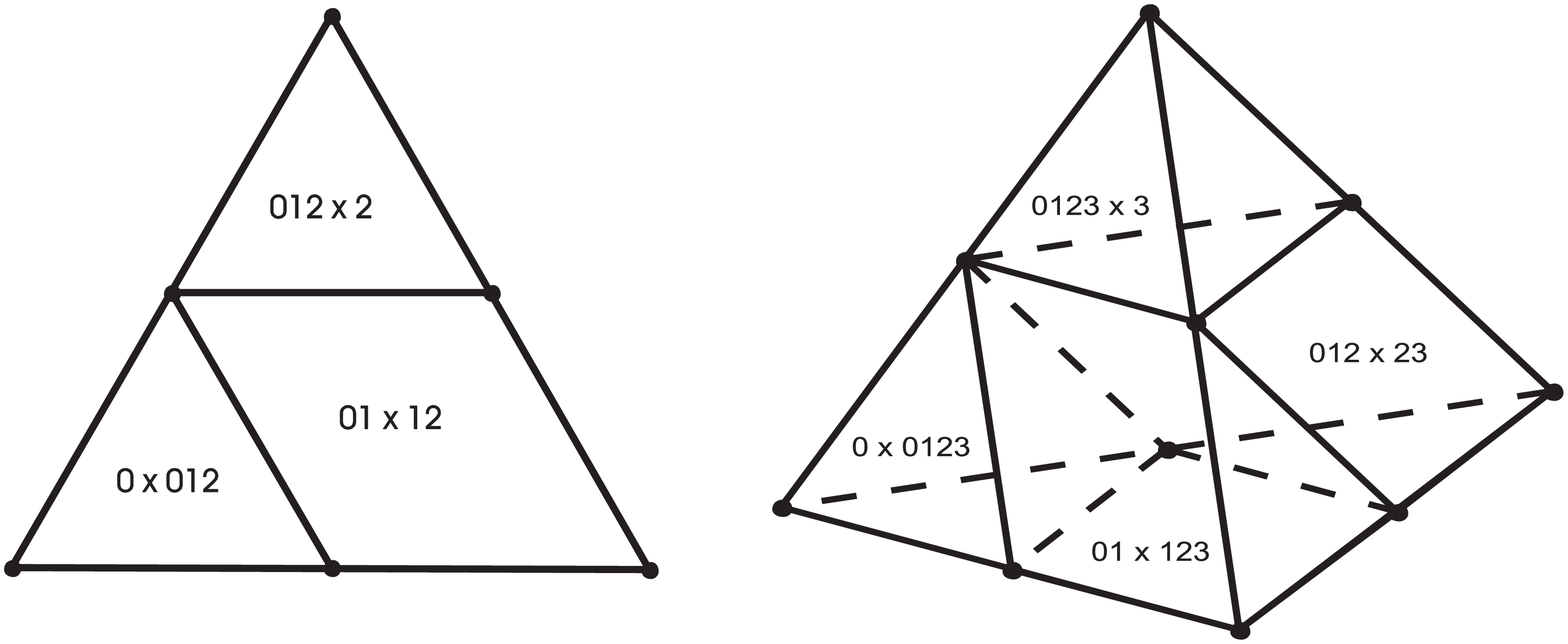}%
\\
Figure 8. 1-subdivisions of $\Delta^{2}$ and $\Delta^{3}$ via the A-W
diagonal.
\end{center}
\begin{center}
\includegraphics[
trim=0.000000in -0.867743in 0.000000in -0.638888in,
height=2.4967in,
width=4.5723in
]%
{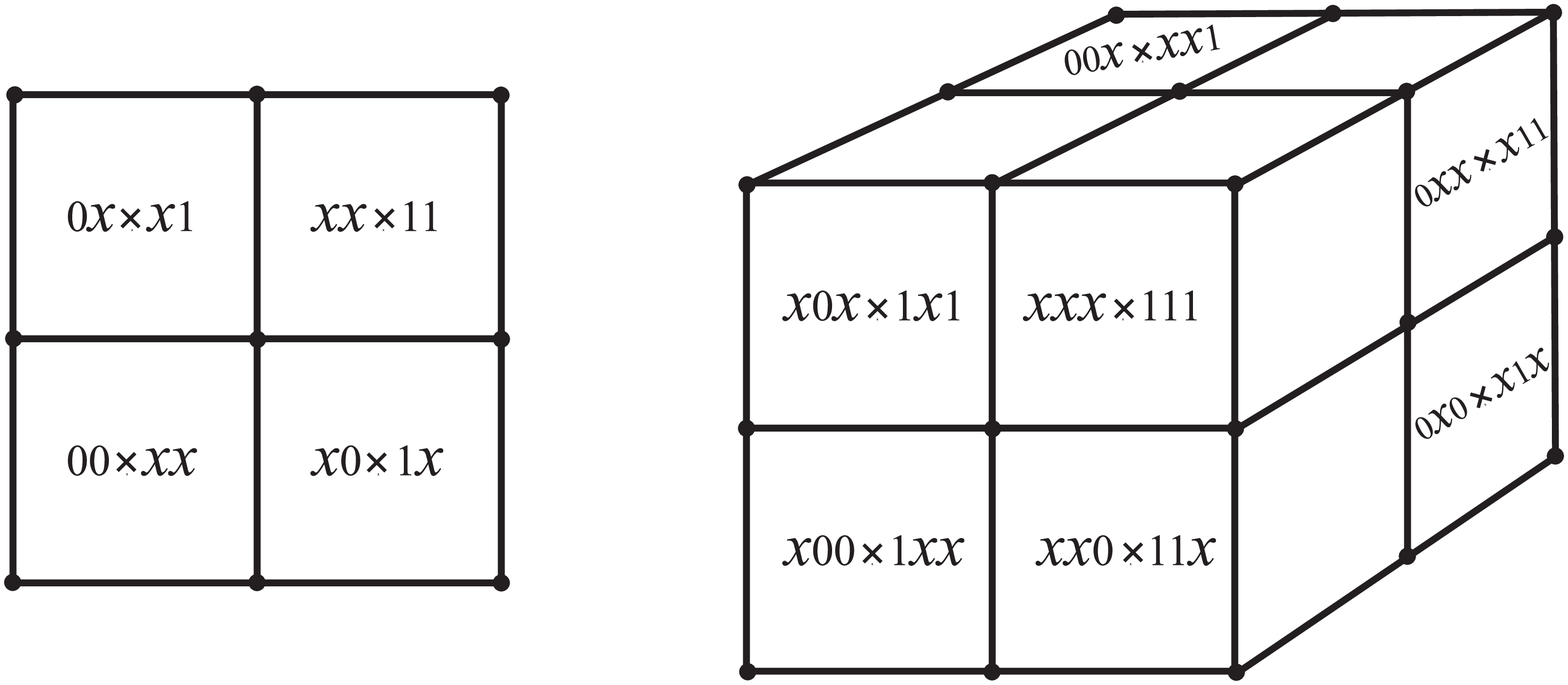}%
\\
Figure 9. 1-subdivisions of $I^{2}$ and $I^{3}$ via the Serre diagonal.
\end{center}
\begin{center}
\includegraphics[
trim=0.000000in -0.286110in 0.000000in -0.572219in,
height=2.1958in,
width=4.9822in
]%
{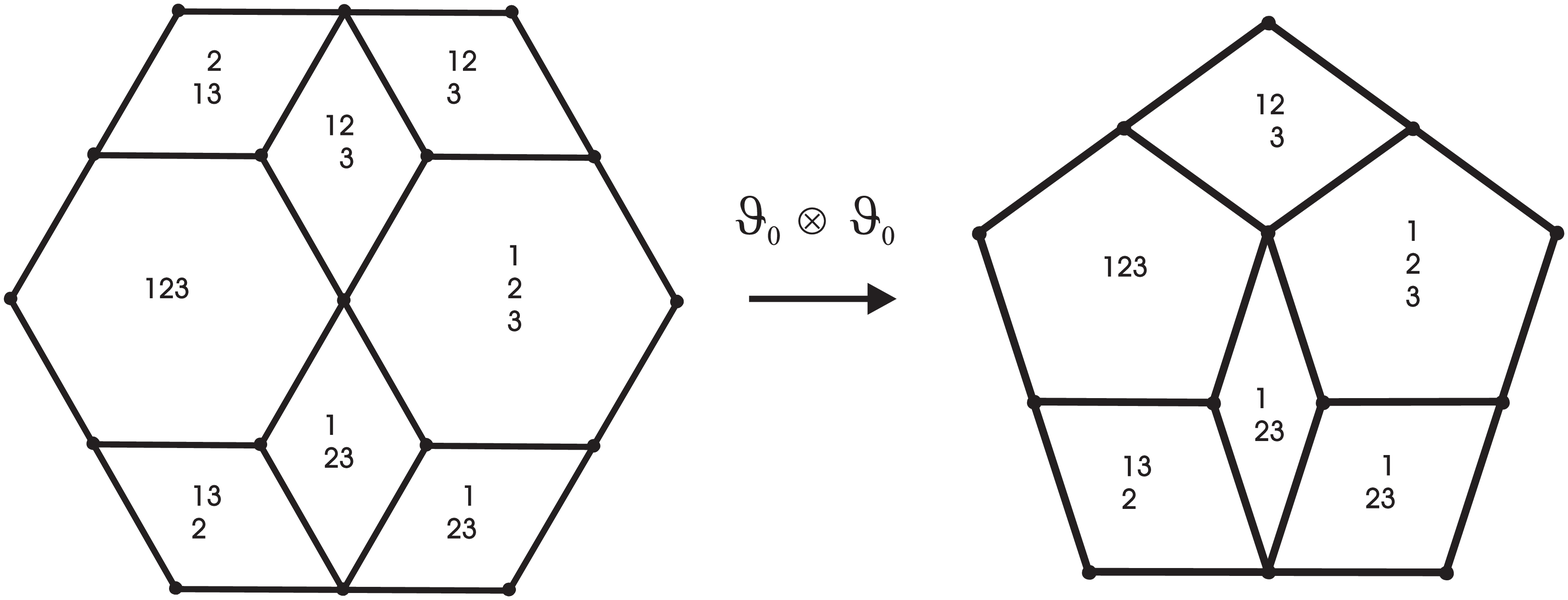}%
\\
Figure 10. 1-subdivisions of $P_{3}$ and $K_{4}$ via S-U diagonals $\Delta
_{P}$ and $\Delta_{K}.$%
\end{center}
\begin{center}
\includegraphics[
trim=0.000000in -0.243413in 0.000000in -0.486826in,
height=2.3229in,
width=4.8352in
]%
{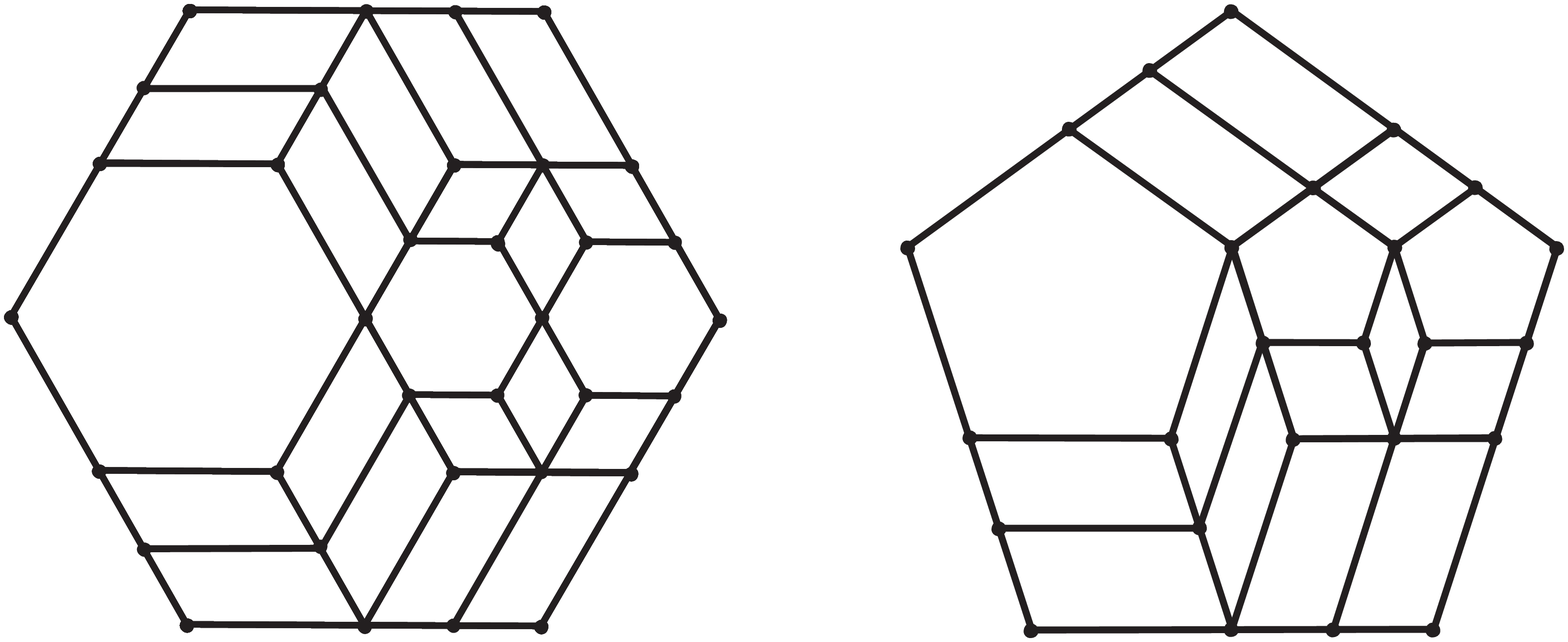}%
\\
Figure 11. 2-subdivisions of $P_{3}$ and $K_{4}.$%
\end{center}

The subcomplex $\Delta^{(k)}\left(  P_{n}\right)  \subset P_{n}^{\times k+1}$
defines the \textquotedblleft configuration module\textquotedblright\ of a
local prematrad in the next section.

\section{Matrads\label{smatrad}}

In this section we introduce the objects in the category of matrads; morphisms
require a relative theory that will be constructed in the sequel. As
motivation, let $\Theta=\left\langle \theta_{m}^{n}\neq0\mid\theta_{1}%
^{1}=\mathbf{1}\right\rangle _{m,n\geq1},$ and consider the canonical
projection $\rho^{{\operatorname{pre}}}:F^{^{{\operatorname{pre}}}}%
(\Theta)\rightarrow{\mathcal{H}}^{^{{\operatorname{pre}}}};$ then%

\[
\rho^{^{{\operatorname{pre}}}}(\theta_{m}^{n})=\left\{
\begin{array}
[c]{ll}%
\theta_{m}^{n}, & m+n\leq3\\
0, & \text{otherwise.}%
\end{array}
\right.
\]
Now consider a differential $\partial^{^{{\operatorname{pre}}}}$ on
$F^{\operatorname{pre}}(\Theta)$ such that $\rho^{^{{\operatorname{pre}}}}$ is
a free resolution in the category of prematrads. Then the induced isomorphism
$\varrho^{{\operatorname{pre}}}:H_{\ast}\left(  F^{^{{\operatorname{pre}}}%
}(\Theta),\partial^{^{{\operatorname{pre}}}}\right)  \approx{\mathcal{H}%
}^{^{{\operatorname{pre}}}}$ implies
\[%
\begin{array}
[c]{c}%
\partial^{^{{\operatorname{pre}}}}(\theta_{2}^{1})=\partial
^{^{{\operatorname{pre}}}}(\theta_{1}^{2})=0\\
\\
\partial^{^{{\operatorname{pre}}}}(\theta_{3}^{1})=\gamma(\theta_{2}%
^{1};\mathbf{1},\theta_{2}^{1})-\gamma(\theta_{2}^{1};\theta_{2}%
^{1},\mathbf{1})\\
\\
\partial^{^{{\operatorname{pre}}}}(\theta_{2}^{2})=\gamma(\theta_{1}%
^{2};\theta_{2}^{1})-\gamma(\theta_{2}^{1}\theta_{2}^{1};\theta_{1}^{2}%
\theta_{1}^{2})\\
\\
\partial^{^{{\operatorname{pre}}}}(\theta_{1}^{3})=\gamma(\mathbf{1}%
,\theta_{1}^{2};\theta_{1}^{2})-\gamma(\theta_{1}^{2},\mathbf{1};\theta
_{1}^{2}).
\end{array}
\]
However, defining $\partial^{^{{\operatorname{pre}}}}$ on all of $\Theta$ is
quite subtle, and while it is possible to canonically extend $\partial
^{^{{\operatorname{pre}}}}$ to all of $\Theta$, acyclicity is difficult to
verify. Instead, there is a canonical proper submodule $\mathcal{H}_{\infty
}\subset{F^{^{{\operatorname{pre}}}}(\Theta)}$ and a differential $\partial$
on $\mathcal{H}_{\infty}$ such that the canonical projection $\varrho
:\mathcal{H}_{\infty}\rightarrow{\mathcal{H}}^{^{{\operatorname{pre}}}}$ is a
free resolution in the category of \textquotedblleft
matrads.\textquotedblright\ Furthermore, we conjecture that the minimal
resolution of $\mathcal{H}$ in the category of PROPs (and consequently, in the
category of prematrads) is recovered by the universal enveloping functor $U$
discussed in Example \ref{universal} below, i.e., the minimal resolution of
the bialgebra PROP $\mathcal{B}$ is $U(\rho):U(\mathcal{H}_{\infty
})\rightarrow U(\mathcal{H}^{\operatorname*{pre}})=\mathcal{B}$, in which case
$\mathcal{H}_{\infty}$ is the smallest extension of ${\mathcal{H}%
}^{\operatorname*{pre}}$ in the category of modules.

The precise definition of $\mathcal{H}_{\infty}$ requires more machinery.

\subsection{Matrads Defined}

Consider a family of pairs $(\mathbf{W}_{\alpha},\gamma_{\alpha}),$ where
$\mathbf{W}_{\alpha}\subset TTM$ is a telescoping submodule, and the
corresponding family of telescopic extensions $\left(  \mathcal{W}_{\alpha
},\Upsilon_{\alpha}\right)  $. To each pair $\left(  \mathcal{W}_{\alpha
},\Upsilon_{\alpha}\right)  $ the $\Upsilon$-factorizations given by via
Definition \ref{asc-desc} below determine a unique \textquotedblleft
configuration module\textquotedblright\ $\Gamma\left(  \mathbf{W}_{\alpha
}\right)  \subseteq$ $\mathcal{W}_{\alpha}$ with the following property: If
$\mathbf{W}_{\alpha}\subseteq\mathbf{W}_{\beta}$ and $\gamma_{_{\mathbf{W}%
_{\alpha}}}=\left.  \gamma_{_{\mathbf{W}_{\beta}}}\right\vert _{\mathbf{W}%
_{\alpha}},$ then $\Gamma\left(  \mathbf{W}_{\alpha}\right)  \subseteq
\Gamma\left(  \mathbf{W}_{\beta}\right)  .$ The local prematrad\ $\left(
M,\gamma_{_{\mathbf{W}_{\alpha}}}\right)  $ is a \textquotedblleft
matrad\textquotedblright\ if $\mathbf{W}_{\alpha}$ is \textquotedblleft%
$\Gamma$-stable,\textquotedblright\ i.e., $\mathbf{W}_{\alpha}$ is the
smallest telescoping submodule such that $\Gamma\left(  \mathbf{W}_{\alpha
}\right)  =\Gamma\left(  \mathbf{W}_{\beta}\right)  $ whenever $\mathbf{W}%
_{\alpha}\subset\mathbf{W}_{\beta}$ and $\gamma_{_{\mathbf{W}_{\alpha}}%
}=\left.  \gamma_{_{\mathbf{W}_{\beta}}}\right\vert _{\mathbf{W}_{\alpha}}$.

Matrads are intimately related to the permutahedra. Recall that each
codimension $k-1$ face of $P_{m-1}$ is identified with two PLTs--an up-rooted
PLT and its down-rooted mirror image--with $m$ leaves and $k\geq2$ levels (see
\cite{Loday}, \cite{SU2}). Define the $\left(  m,1\right)  $\emph{-row descent
sequence} of $\curlywedge_{m}$ to be $\mathbf{m}=\left(  m\right)  .$ Given an
up-rooted PLT $T=T^{k}$ with $k$ levels, express $T^{i}=T^{i-1}\diagup
\curlywedge_{m_{i,1}}\cdots\curlywedge_{m_{i,r_{i}}}$ for $k\geq i>1$ and
$T^{1}=\curlywedge_{m_{1,1}}.$ Define the $i^{th}$ \emph{leaf sequence} of $T$
to be the row matrix $\mathbf{m}_{i}=(m_{i,1},\ldots,m_{i,r_{i}})$ and the
$\left(  m,k\right)  $\emph{-row descent sequence of} $T$ to be the $k$-tuple
of row matrices $\left(  \mathbf{m}_{1},\ldots,\mathbf{m}_{k}\right)  .$ Note
that the vertices of $P_{m-1}$ are identified with $\left(  m,m-1\right)
$-row descent sequences $\left(  \mathbf{m}_{1},\ldots,\mathbf{m}%
_{m-1}\right)  ,$ where $\mathbf{m}_{i}=\left(  1,\ldots,2,\ldots,1\right)
^{T}\in\mathbb{N}^{r_{i}}$ with exactly one $2$ in position $j$ for some
$1\leq j\leq r_{i}$, in which case $\mathbf{m}_{1}=\left(  2\right)  .$
Dually, define the $\left(  n,1\right)  $\emph{-column descent sequence} of
$\curlyvee^{n}$ to be $\mathbf{n}=\left(  n\right)  .$ Given an down-rooted
PLT $T=T^{l}$ with $l$ levels, express $T^{i}=\curlyvee^{n_{i,1}}%
\cdots\curlyvee^{n_{i,s_{i}}}\diagup T^{i-1}$ for $l\geq i>1$ and
$T^{1}=\curlyvee^{n_{1,1}}.$ Define the $i^{th}$\emph{\ leaf sequence} of $T$
to be the column matrix $\mathbf{n}_{i}=(n_{i,1},\ldots,n_{i,s_{i}})^{T}$ and
the $\left(  n,l\right)  $-\emph{column} \emph{descent sequence of }$T$ to be
the $l$-tuple of column matrices $\left(  \mathbf{n}_{l},\ldots,\mathbf{n}%
_{1}\right)  $.

Given a telescoping submodule $\mathbf{W}$ and its telescopic extension
$\mathcal{W},$ let $\mathcal{W}_{{\operatorname{row}}}=\mathcal{W}%
\cap\overline{\mathbf{M}}_{{\operatorname{row}}}$ and $\mathcal{W}%
^{\operatorname{col}}=\mathcal{W}\cap\overline{\mathbf{M}}^{\operatorname{col}%
}.$

\begin{definition}
\label{asc-desc}Given a local prematrad $(M,\gamma_{\mathbf{W}})$ with domain
$\mathbf{W},$ let $\zeta\in{M}_{\ast,m}$ and $\xi\in M_{n,\ast}$ be elements
with $m,n\geq2.$

\begin{enumerate}
\item[\textit{(i)}] A \textbf{row factorization of} $\zeta$ \textbf{with
respect to} $\mathbf{W}$ is a $\Upsilon\!$-factorization $A_{1}\!\cdots
\!A_{k}\linebreak=\zeta$ such that $A_{j}\in\mathcal{W}_{\operatorname{row}}$
and ${\operatorname{rls}}(A_{j})\neq\mathbf{1}$ for all $j$. The sequence
$\left({\operatorname{rls}}(A_{1}),\ldots,\right.$
\linebreak $\left. {\operatorname{rls}}(A_{k}) \right )  $ is the related $\left(  m,k\right)  $-row \textbf{descent
sequence of} $\zeta$.

\item[(\textit{ii)}] A \textbf{column factorization of }$\xi$ \textbf{with
respect to} $\mathbf{W}$ is a $\Upsilon$-factorization\linebreak\ $B_{l}\cdots
B_{1}=\xi$ such that $B_{i}\in\mathcal{W}^{{\operatorname{col}}}$\thinspace
and ${\operatorname{cls}}(B_{i})\neq\mathbf{1}$ for all $i.$ The sequence
\linebreak $\left(  {\operatorname{cls}}(B_{l}),\ldots,{\operatorname{cls}}%
(B_{1})\right)  $ is the related $\left(  n,l\right)  $-\textbf{column}
\textbf{descent sequence of} $\xi$.
\end{enumerate}
\end{definition}

\noindent Column and row factorizations are not unique. Note that an element
$A\in{M}_{n,\ast}$ always has a trivial column factorization as the $1\times1$
matrix $\left[  A\right]  $. When matrix entries in a row factorization are
pictured as graphs, terms of the row descent sequence are \textquotedblleft
lower (input) leaf sequences\textquotedblright\ of the graphs along any row,
and dually for column factorizations.

\begin{example}
\label{colrow}An $\Upsilon$-product of bisequence matrices is simultaneously a
row and column factorization. For example, consider the $\Upsilon$-product%
\[
C=C_{1}C_{2}C_{3}=\left[
\begin{array}
[c]{c}%
\theta_{2}^{1}\\
\alpha_{2}^{2}\\
\beta_{2}^{2}%
\end{array}
\right]  \left[
\begin{array}
[c]{cc}%
\theta_{1}^{2} & \zeta_{2}^{2}\\
\theta_{1}^{1} & \theta_{2}^{1}%
\end{array}
\right]  \left[  \theta_{1}^{2}\text{ }\xi_{2}^{2}\text{ }\theta_{1}%
^{2}\right]  \in{M}_{5,4}.
\]
As a row factorization, the $\left(  4,3\right)  $-row descent sequence of $C$
is
\[
\left(  {\operatorname{rls}}\left(  C_{1}\right)  ,{\operatorname{rls}}\left(
C_{2}\right)  ,{\operatorname{rls}}\left(  C_{3}\right)  \right)  =\left(
\left(  2\right)  ,\left(  12\right)  ,\left(  121\right)  \right)  ,
\]
and as a column factorization, the $\left(  5,3\right)  $-column descent
sequence of $C$ is
\[
\left(  {\operatorname{cls}}\left(  C_{1}\right)  ,{\operatorname{cls}}\left(
C_{2}\right)  ,{\operatorname{cls}}\left(  C_{3}\right)  \right)  =\left(
\begin{pmatrix}
1\\
2\\
2
\end{pmatrix}
,%
\begin{pmatrix}
2\\
1
\end{pmatrix}
,\left(  2\right)  \right)  .
\]

\end{example}

\begin{center}
\includegraphics[
height=0.9513in,
width=2.3601in
]%
{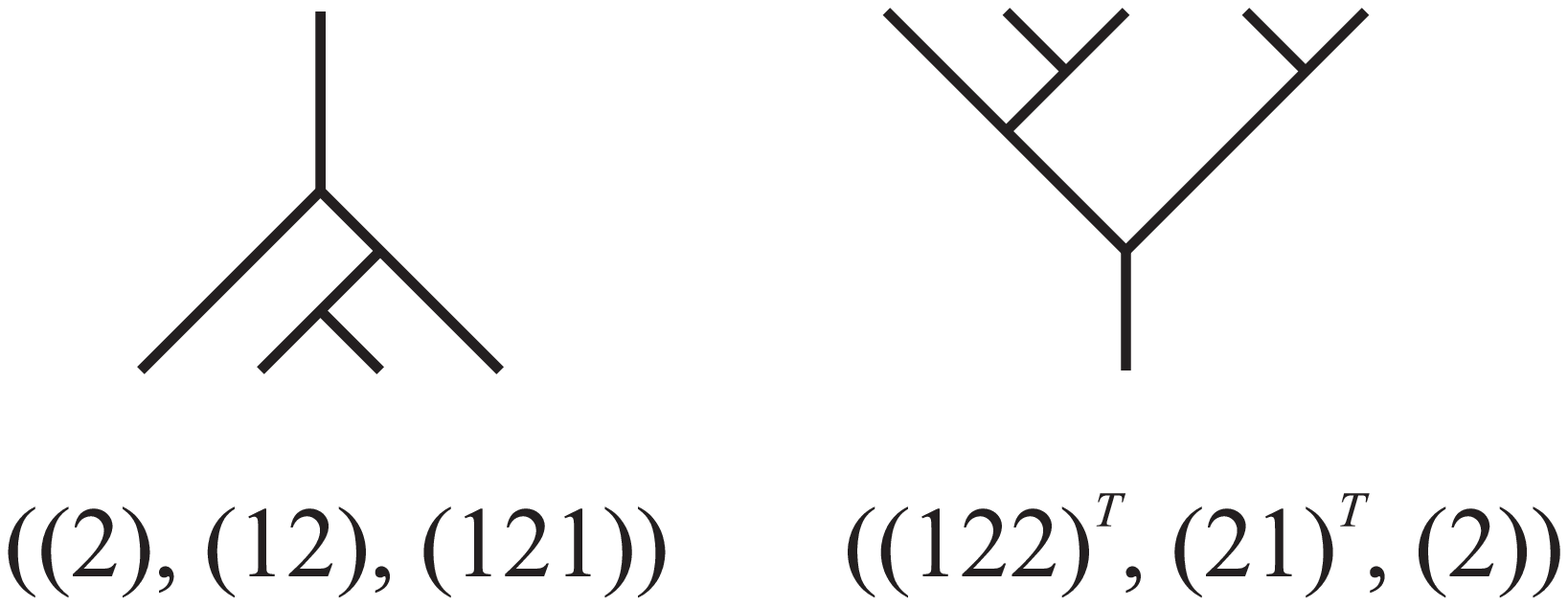}%
\end{center}
\vspace*{0.1in}

Given a local prematrad $(M,\gamma_{\mathbf{W}})$ and elements $A\in{M}%
_{\ast,s}$ and $B\in{M}_{t,\ast}$ with $s,t\geq2$, choose a row factorization
$A_{1}\cdots A_{k}$ of $A$ with respect to $\mathbf{W}$ and $\left(
s,k\right)  $-row descent sequence $\mathbf{\alpha,}$ and a column
factorization $B_{l}\cdots B_{1}$ of $B\ $with respect to $\mathbf{W}$ and
$\left(  t,l\right)  $-column descent sequence $\mathbf{\beta}$. Then
$\mathbf{\alpha}$ identifies $A$ with an up-rooted $s$-leaf, $k$-level PLT and
a codimension $k-1$ face $\overset{_{\wedge}}{e}_{A}$ of $P_{s-1}$, and
$\mathbf{\beta}$ identifies $B$ with a down-rooted $t$-leaf, $l$-level PLT and
a codimension $l-1$ face $\overset{_{\vee}}{e}_{B}$ of $P_{t-1}$. Extending to
Cartesian products, identify the monomials $A=A_{1}\otimes\cdots\otimes
A_{q}\in\left(  {M}_{\ast,s}\right)  ^{\otimes q}$ and $B=B_{1}\otimes
\cdots\otimes B_{p}\in\left(  {M}_{t,\ast}\right)  ^{\otimes p}$ with the
product cells
\[
\overset{\wedge}{e}_{A}=\overset{_{\wedge}}{e}_{A_{1}}\times\cdots
\times\overset{_{\wedge}}{e}_{A_{q}}\subset P_{s-1}^{\times q}\ \ \text{and}%
\ \ \overset{_{\vee}}{e}_{B}=\overset{_{\vee}}{e}_{B_{1}}\times\cdots
\times\overset{_{\vee}}{e}_{B_{p}}\subset P_{t-1}^{\times p}.
\]

Now consider the S-U diagonal $\Delta_{P}$ and recall from Section
\ref{topological} that there is a $k$-subdivision $P_{r}^{\left(  k\right)  }$
of $P_{r}$ and a cellular inclusion $P_{r}^{\left(  k\right)  }\hookrightarrow
\Delta^{\left(  k\right)  }\left(  P_{r}\right)  \subset P_{r}^{\times k+1}$
for each $k$ and $r$. Thus for each $q\geq2$, the product cell
$\overset{\wedge}{e}_{A}$ either\emph{\ is} or \emph{is not} a subcomplex of
${\Delta}^{(q-1)}(P_{s-1})\subset P_{s-1}^{\times q}$, and dually for
$\overset{\vee}{e}_{B}.$ This leads to the notion of \textquotedblleft
configuration module.\textquotedblright

Let $\mathbf{x}_{m,i}^{p}=\left(  1,\ldots,m,\ldots,1\right)  \in
\mathbb{N}^{1\times p}$ with $m$ in the $i^{th}$ position and let
$\mathbf{y}_{q}^{n,j}=\left(  1,\ldots,n,\ldots,1\right)  ^{T}\in
\mathbb{N}^{q\times1} $ with $n$ in the $j^{th}$ position.

\begin{definition}
The \textbf{(left) configuration module} of a local prematrad $\left(
M,\gamma_{_{\mathbf{W}}}\right)  $ is the $R$-module
\[
\Gamma(M,\gamma_{_{\mathbf{W}}})=M\oplus\bigoplus_{\mathbf{x,y}\notin%
\mathbb{N};\text{ }s,t\geq1}\Gamma_{s}^{\mathbf{y}}(M)\oplus\Gamma
_{\mathbf{x}}^{t}(M),\text{ where}%
\]%
\begin{align*}
\Gamma_{s}^{\mathbf{y}}(M) &  =\left\{
\begin{array}
[c]{cl}%
\mathbf{M}_{1}^{\mathbf{y}}, & s=1;\text{ }\mathbf{y}=\mathbf{y}_{q}%
^{n,j}\text{ for some }n,j,q\\
\left\langle A\in\mathbf{M}_{s}^{\mathbf{y}}\mid\overset{_{\wedge}}{e}%
_{A}\subset{\Delta}^{(q-1)}(P_{s-1})\right\rangle , & s\geq2\\
0, & \text{otherwise,}%
\end{array}
\right. \\
& \\
\Gamma_{\mathbf{x}}^{t}(M) &  =\left\{
\begin{array}
[c]{cl}%
\mathbf{M}_{\mathbf{x}}^{1}, & t=1;\text{ }\mathbf{x}=\mathbf{x}_{m,i}%
^{p}\text{ for some }m,i,p\\
\left\langle B\in\mathbf{M}_{\mathbf{x}}^{t}\mid\overset{_{\vee}}{e}%
_{B}\subset{\Delta}^{(p-1)}(P_{t-1})\right\rangle , & t\geq2\\
0, & \text{otherwise.}%
\end{array}
\right.
\end{align*}

\end{definition}

\noindent Thus $\Gamma_{\mathbf{x}}^{t}(M)$ is generated by those tensor
monomials $B=\beta_{x_{1}}^{t}\otimes\cdots\otimes\beta_{x_{p}}^{t}\in
{M}_{t,x_{1}}\otimes\cdots\otimes{M}_{t,x_{p}}$ whose tensor factor
$\beta_{x_{i}}^{t}$ is identified with some factor of a product cell in
${\Delta}^{\left(  p-1\right)  }(P_{t-1})$ corresponding to some column
factorization $\beta_{x_{i}}^{t}=B_{i,l}\cdots B_{i,1}$ with respect to
$\mathbf{W}$, and dually for $\Gamma_{s}^{\mathbf{y}}(M).$

\begin{example}
\label{vectormodule}Referring to Example \ref{bisequence_vect}, let
$M=M_{1,1}\oplus M_{2,2}=\left\langle \theta_{1}^{1}\right\rangle
\oplus\left\langle \theta_{2}^{2}\right\rangle ,$ and consider the local
prematrad $\left(  M,\gamma_{_{\mathbf{V}}}\right)  .$ Note that the action of
$\gamma_{_{\mathbf{V}}}$ is trivial modulo unit (e.g., $\mathbf{M}_{2}%
^{2}\cdot\mathbf{M}_{11}^{1}=\mathbf{M}_{2}^{2}$ and $\mathbf{M}_{2}^{22}%
\cdot\mathbf{M}_{22}^{2}\subseteq\mathbf{M}_{4}^{4}=0$). Then
\[
\Gamma_{\mathbf{\ast}}^{1}(M)\approx\Gamma_{1}^{\ast}(M)\approx T^{+}\left(
M_{1,1}\right)  .
\]
Since $A=\left(  \theta_{2}^{2}\right)  ^{\otimes q}$ can be thought of as an
element of either $\mathbf{M}_{2\cdots2}^{2}$ or $\mathbf{M}_{2}^{2\cdots2},$
we have either $\overset{\wedge}{e}_{A}={\Delta}^{\left(  q-1\right)  }%
(P_{1})=P_{1}^{\times q}$ or $\overset{\vee}{e}_{A}={\Delta}^{\left(
q-1\right)  }(P_{1})=P_{1}^{\times q}$ so that
\[
\Gamma_{\mathbf{\ast}}^{2}(M)\approx\Gamma_{2}^{\ast}(M)\approx T^{+}\left(
M_{2,2}\right)  .
\]
Thus
\[
\Gamma\left(  M,\gamma_{_{\mathbf{V}}}\right)  =\mathbf{V}.
\]

\end{example}

Since $r=1$ is the only case in which the equality ${\Delta}^{\left(
k\right)  }(P_{r})=P_{r}^{\times k+1}$ holds for each $k,$ it follows
immediately that if $A\in{M}_{n,m}$ is $\gamma_{_{\mathbf{W}}}$-indecomposable
(in which case its row and column factorizations with respect to $\mathbf{W}$
are trivial) and $\overset{\wedge}{e}_{A}\times\cdots\times\overset{\wedge
}{e}_{A}=e^{m-2}\times\cdots\times e^{m-2}$ is a subcomplex of ${\Delta
}^{\left(  n-1\right)  }(P_{m-1}),$ then either $m=2$ and $n$ is arbitrary or
$m>2$ and $n=1;$ dually, if $\overset{\vee}{e}_{A}\times\cdots\times
\overset{\vee}{e}_{A}=e^{n-2}\times\cdots\times e^{n-2}$ is a subcomplex of
${\Delta}^{\left(  m-1\right)  }(P_{n-1}),$ then either $m$ is arbitrary and
$n=2$ or $m=1$ and $n>2$. Consequently,
\[
\bigoplus\nolimits_{\mathbf{x}}{\Gamma}_{\mathbf{x}}^{2}(M)\approx
T^{+}(M_{2,\ast})\ \ \text{and}\ \ \bigoplus\nolimits_{\mathbf{y}}{\Gamma}%
_{2}^{\mathbf{y}}(M)\approx T^{+}(M_{\ast,2}).
\]
Furthermore, if $m+n\geq4,$ $(m,n)\neq(2,2),$ and $A^{\otimes r}\in
\Gamma(M,\gamma_{_{\mathbf{W}}}),$ then $r=1,$ and the inclusion $\Gamma
_{s}^{\mathbf{y}}(M)\subseteq\mathbf{M}_{s}^{\mathbf{y}}$ is proper whenever
$s\geq3,$ $\mathbf{y}\in\mathbb{N}^{q\times1}$ with $q\geq2,$ and ${M}%
_{y_{j},s}$ contains a $\gamma_{\mathbf{W}}$-indecomposable element for each
$1\leq j\leq q$, and dually for $\Gamma_{\mathbf{x}}^{t}(M)\subseteq
\mathbf{M}_{\mathbf{x}}^{t}$ .

We are ready to define the notion of a matrad.

\begin{definition}
A local prematrad $(M,\gamma_{_{\mathbf{W}}},\eta)$ is a
\textbf{(left)\thinspace matrad} if
\[
\Gamma_{p}^{\mathbf{y}}(M,\gamma_{_{\mathbf{W}}})\otimes\Gamma_{\mathbf{x}%
}^{q}(M,\gamma_{_{\mathbf{W}}})=\mathbf{W}_{p}^{\mathbf{y}}\otimes
\mathbf{W}_{\mathbf{x}}^{q}%
\]
for all $p,q\geq2.$ A \textbf{morphism of matrads} is a map of underlying
local prematrads.
\end{definition}

\begin{example}
\underline{The Bialgebra Matrad ${\mathcal{H}}.$} The bialgebra prematrad
${\mathcal{H}}^{^{{\operatorname{pre}}}}$satisfies\linebreak$\Gamma
_{p}^{\mathbf{y}}(M)\otimes\Gamma_{\mathbf{x}}^{q}(M)=\mathbf{M}%
_{p}^{\mathbf{y}}\otimes\mathbf{M}_{\mathbf{x}}^{q}$ for $p,q\geq2.$ Hence
${\mathcal{H}}^{^{{\operatorname{pre}}}}$ is also a matrad, called the
\emph{bialgebra matrad} and henceforth denoted by $\mathcal{H}$.
\end{example}

\begin{example}
Continuing Example \ref{extremes}, the inclusions $(M_{1,\ast},\gamma
_{_{\mathbf{M}_{\ast}^{1}}})\subset(M,\gamma_{_{\mathbf{M}_{\ast}^{1}}})$ and
$(M_{\ast,1},\gamma_{_{\mathbf{M}_{1}^{\ast}}})\subset(M,\gamma_{_{\mathbf{M}%
_{1}^{\ast}}})$ are inclusions of matrads since $\Gamma_{p}^{\mathbf{y}%
}(M,\gamma_{_{\mathbf{M}_{1}^{\ast}}})\otimes\Gamma_{\mathbf{x}}^{q}%
(M,\gamma_{_{\mathbf{M}_{\ast}^{1}}})=\mathbf{M}_{p}^{\mathbf{y}}%
\otimes\mathbf{M}_{\mathbf{x}}^{q}=0$ for $p,q\geq2.$ In particular, when
$M=F^{^{{\operatorname{pre}}}}\left(  \Theta\right)  $ and $\theta_{m}^{n}%
\neq0$ for all $m,n\geq1,$ the free operad $(\mathcal{A}_{\infty}%
,\gamma_{_{\mathbf{M}_{\ast}^{1}}})$ embeds in $(F^{^{{\operatorname{pre}}}%
}\left(  \Theta\right)  ,\gamma_{_{\mathbf{M}_{\ast}^{1}}})$ as a submatrad
(c.f. Example \ref{non-sigma} and Definition \ref{defnfreematrad} below).
\end{example}

\begin{example}
\label{universal} \underline{The Universal Enveloping Functor $U.$} The
universal enveloping PROP\ $U$ discussed in Example \ref{UEP} induces the
\textbf{universal enveloping functor }$U$ from the category of matrads to the
category of PROPs. Given a matrad $(M,\gamma_{_{\mathbf{W}}}),$ let $FP(M)$ be
the free PROP generated by $M$ and let $J$ be the two-sided ideal generated by
the elements $\bigoplus\nolimits_{\mathbf{x\times y}}\left(  \gamma
_{_{\mathbf{W}}}-\gamma_{_{FP}}\right)  \left(  \Gamma_{p}^{\mathbf{y}%
}(M)\otimes\Gamma_{\mathbf{x}}^{q}(M)\right)  .$ Then $U(M)=FP(M)\diagup J$.
\end{example}

\subsection{Free Matrads}

Recall that the domain of the free prematrad $(M=F^{^{{\operatorname{pre}}}%
}(\Theta),\gamma,$ $\eta)$ generated by $\Theta=\left\langle \theta_{m}%
^{n}\right\rangle _{m,n\geq1}$ is $\mathbf{V}=M\oplus\bigoplus
\limits_{\mathbf{x},\mathbf{y}\notin\mathbb{N};\text{ }s,t\in\mathbb{N}%
}\mathbf{M}_{s}^{\mathbf{y}}\oplus\mathbf{M}_{\mathbf{x}}^{t}\ $whose
submodules $M,$ $\mathbf{M}_{\mathbf{x}}^{1},$ and $\mathbf{M}_{1}%
^{\mathbf{y}}$ are contained in the configuration module $\Gamma\left(
M\right)  $. As above, the symbol \textquotedblleft$\cdot$\textquotedblright%
\ denotes the $\gamma$ product.

\begin{definition}
\label{defnfreematrad}Let $(M=F^{^{{\operatorname{pre}}}}(\Theta),\gamma
,\eta)$ be the free prematrad generated by $\Theta=\left\langle \theta_{m}%
^{n}\right\rangle _{m,n\geq1},$ let $F(\Theta)=\Gamma(M)\cdot\Gamma(M),$ and
let $\gamma_{_{F\left(  \Theta\right)  }}=\gamma|_{_{\Gamma(M)\otimes
\Gamma(M)}}.$ The \textbf{free matrad generated by} $\Theta$ is the triple
$\left(  F(\Theta),\gamma_{_{F\left(  \Theta\right)  }},\eta\right)  .$
\end{definition}

\begin{remark}
Let $\omega=\sum\theta_{m}^{n}$ and consider the biderivative $d_{\omega}.$
Then $F(\Theta)$ is generated by the components of $d_{\omega}\circledcirc
d_{\omega}$ in $F^{\operatorname*{pre}}(\Theta)$ (the admissible fractions).
Thus $\Gamma(F(\Theta),\gamma_{F(\Theta)})\setminus F(\Theta)=\Gamma
(F^{\operatorname*{pre}}(\Theta),\gamma_{\mathbf{V}})\setminus
F^{\operatorname*{pre}}(\Theta)$ (c.f. Proposition \ref{bases}).
\end{remark}

Let $\beta\in\mathcal{B}^{\operatorname{pre}}.$ In Subsection \ref{freepre} we
constructed the $r$-level tree $\Psi(\beta)$ whose leaves are balanced
$\Theta$-factorizations; the set $\mathcal{C}^{\operatorname{pre}}=\left\{
\Psi(\beta)\mid\beta\in\mathcal{B}^{\operatorname{pre}}\right\}  $ indexes the
set 
\linebreak $G^{\operatorname{pre}}\diagup\sim$ of module generators of the free
prematrad $F^{\operatorname{pre}}(\Theta).$ Let $\mathfrak{G}=F(\Theta)\cap
G^{\operatorname{pre}}\diagup\sim$ and\ let $\mathcal{B}=\phi\left(
\mathfrak{G}\right)  .$ Then
\[
\mathcal{C=}\left\{  \Psi(\beta))\mid\beta\in\mathcal{B}\right\}
\]
indexes the set $\mathfrak{G}$ of module generators of $F(\Theta)$.

To establish the relationship between elements of $\mathfrak{G}$ and cells of
$KK_{n,m},$ let $\beta=A_{s}\cdots A_{1}\in\mathcal{B}_{n,m}%
^{\operatorname{pre}}$ and observe that $\beta\in\mathcal{B}_{n,m}$ if and
only if the tensor monomials along each row and column of $A_{k}$ lie in
$\Gamma(F^{\operatorname{pre}}\left(  \Theta\right)  )$ for all $k$ (see
(\ref{iota}) below). Let
\[
\mathcal{C}_{n,m}^{\prime\prime}=\left\{  \Psi(\beta)\mid\beta=A_{s}\cdots
A_{1}\in\mathcal{B}_{n,m}\text{ and either }A_{1}\text{ or }A_{s}\text{ is a
}1\times1\right\}  ;
\]
then$\;$in particular, $\mathcal{C}_{1,m}^{^{\prime\prime}}=\mathcal{C}_{1,m}
$ and $\mathcal{C}_{n,1}^{^{\prime\prime}}=\mathcal{C}_{n,1}.$ Let
$\mathcal{C}_{n,m}^{\prime}=\mathcal{C}_{n,m}\setminus\mathcal{C}%
_{n,m}^{\prime\prime}; $ then $\mathcal{C}_{n,m}=\mathcal{C}_{n,m}^{\prime
}\sqcup\mathcal{C}_{n,m}^{\prime\prime}$ for each $m,n\geq1$. Elements of
$\mathcal{C}^{\prime} $ are defined in terms of $\Delta_{P};$ elements of
$\mathcal{C}^{\prime\prime}$ are independent of $\Delta_{P}$.

Define the dimension of $\theta_{1}^{1}=\mathbf{1}$ to be zero and the
dimension of $\theta_{m}^{n}$ to be $m+n-3;$ if $A\in\mathfrak{G}$, then the
dimension of $A,$ denoted by $|A|$, is the sum of the dimensions of the matrix
entries in any representative monomial in $G,$ and in particular, in its
balanced factorization in $\mathcal{B}$. Clearly, given $\mathbf{x\times y}%
\in{\mathbb{N}^{1\times p}}\times\mathbb{N}^{q\times1}$, the set $\{\left\vert
A\right\vert \mid A\in\mathbf{G}_{\mathbf{x}}^{\mathbf{y}}\}$ is bounded and,
consequently, has a maximal element. For example, if $A$ is a monomial in
$\Gamma_{s}^{\mathbf{y}}(F\left(  \Theta\right)  )$ with $s\geq2$ and
$\overset{_{\wedge}}{e}_{A}$ is the corresponding subcomplex of $\Delta
^{\left(  q-1\right)  }\left(  P_{s-1}\right)  $, then $\left\vert
A\right\vert =\left\vert \overset{\wedge}{e}_{A}\right\vert +\left\vert
\mathbf{y}\right\vert -q;$ and dually, if $B$ is a monomial in $\Gamma
_{\mathbf{x}}^{t}(F\left(  \Theta\right)  )$ with $t\geq2$ and
$\overset{_{\vee}}{e}_{B}$ is the corresponding subcomplex of $\Delta^{\left(
p-1\right)  }\left(  P_{t-1}\right)  ,$ then $\left\vert B\right\vert
=\left\vert \overset{\vee}{e}_{B}\right\vert +\left\vert \mathbf{x}\right\vert
-p.$ Consequently, $\max\{\left\vert A\right\vert \mid A\in\Gamma
_{s}^{\mathbf{y}}(F\left(  \Theta\right)  )\}=|\mathbf{y}|+s-q-2$ and
$\max\{\left\vert B\right\vert \mid B\in\Gamma_{\mathbf{x}}^{t}(F\left(
\Theta\right)  )\}=|\mathbf{x}|+t-p-2$. In particular, if $A\in\mathfrak{G}%
_{n+1,m+1}$ has balanced factorization $\beta=A_{s}\cdots A_{1},$ $1\leq s\leq
m+n,$ then $\operatorname{codim}A\geq s-1$ and $\operatorname{codim}A=s-1$ if
and only if the dimensions of each bisequence matrix $A_{k}$ is maximal (see
(\ref{codim}) below).

Given $m+n\geq3$ and $\mathbf{x\times y}\in\mathbb{N}^{1\times p}%
\times\mathbb{N}^{q\times1}$ such that $\left\vert \mathbf{x}\right\vert =m$
and $\left\vert \mathbf{y}\right\vert =n,$ define the codimension $1$ face
$e_{\left(  \mathbf{y,x}\right)  }\subset P_{m+n-2}$ as follows: If
$\left\vert \mathbf{x}\right\vert =m>p\geq2,$ let $A_{\mathbf{x}%
}|B_{\mathbf{x}}$ be the codimension $1$ face of $P_{m-1}$ with leaf sequence
$\mathbf{x};$ dually, if $\left\vert \mathbf{y}\right\vert =n>q\geq2,$ let
$A_{\mathbf{y}}|B_{\mathbf{y}}$ be the codimension $1$ face of $P_{n-1}$ with
leaf sequence $\mathbf{y}$. If $A=\left\{  a_{1},\ldots,a_{r}\right\}
\subset\mathbb{Z}$ and $z\in\mathbb{Z}$, define $-A=\left\{  -a_{1}%
,\ldots,-a_{n}\right\}  $ and $A+z=\left\{  a_{1}+z,\ldots,a_{r}+z\right\}  ;$
then set%
\[%
\begin{array}
[c]{ll}%
A_{1}=\left\{
\begin{array}
[c]{cl}%
\underline{m-1}, & \text{if }\mathbf{x}=\mathbf{1}^{m},\,m\geq1\medskip\\
\varnothing, & \text{if }\mathbf{x}=m\geq2\medskip\\
-A_{\mathbf{x}}+m, & \text{otherwise,}%
\end{array}
\right.  & A_{2}=\left\{
\begin{array}
[c]{cl}%
\varnothing, & \text{if }\mathbf{y}=\mathbf{1}^{n},\,n\geq1\medskip\\
\underline{n-1}, & \text{if }\mathbf{y}=n\geq2\medskip\\
A_{\mathbf{y}}, & \text{otherwise,}%
\end{array}
\right. \\
& \\
B_{1}=\left\{
\begin{array}
[c]{cl}%
\varnothing, & \text{if }\mathbf{x}=\mathbf{1}^{m},\,m\geq1\medskip\\
\underline{m-1}, & \text{if }\mathbf{x}=m\geq2\medskip\\
-B_{\mathbf{x}}+m, & \text{otherwise,}%
\end{array}
\right.  & B_{2}=\left\{
\begin{array}
[c]{cl}%
\underline{n-1}, & \text{if }\mathbf{y}=\mathbf{1}^{n},\,n\geq1\medskip\\
\varnothing, & \text{if }\mathbf{y}=n\geq2\medskip\\
B_{\mathbf{y}}, & \text{otherwise,}%
\end{array}
\right.
\end{array}
\]
and define
\begin{equation}
e_{(\mathbf{y},\mathbf{x})}={A}_{1}\cup(A_{2}+m-1)|{B}_{1}\cup({B}%
_{2}+m-1).\label{special}%
\end{equation}
For example, $e_{(\mathbf{1}^{n},\mathbf{1}^{m})}=\underline{m-1}\,|\left(
\underline{n-1}+m-1\right)  ,$ \ $e_{(n,m)}=\left(  \underline{n-1}%
+m-1\right)  |\,\underline{m-1},$ \ and $\ e_{((21),(21))}=13|24.$

\begin{example}
\label{freematrad} \underline{The $A_{\infty}$-bialgebra Matrad $\mathcal{H}%
_{\infty}.$} Let $\Theta=\left\langle \theta_{m}^{n}\neq0\mid\theta_{1}%
^{1}=\mathbf{1}\right\rangle _{m,n\geq1}$. We say that $\beta\in
\mathcal{B}_{n,m}$ has \textbf{word length} $2$ if $\beta=C_{2}C_{1}$ for some
$C_{2}\times C_{1}\in\mathbf{G}_{p}^{\mathbf{y}}\times\mathbf{G}_{\mathbf{x}%
}^{q},$ where $\mathbf{x\times y}\in{\mathbb{N}^{1\times p}}\times
\mathbb{N}^{q\times1},$ $\left\vert \mathbf{x}\right\vert =m,$ and $\left\vert
\mathbf{y}\right\vert =n.$ Let%
\[
\mathcal{A}_{p}^{\mathbf{y}}\times\mathcal{B}_{\mathbf{x}}^{q}=\left\{
\Psi(\beta)\mid\beta\in\mathcal{B}_{n,m}\text{ has word length }2\text{
}\right\}  .
\]
Denote the corresponding bases of ${\Gamma}_{\mathbf{x}}^{q}(F\left(
\Theta\right)  )$ and ${\Gamma}_{p}^{\mathbf{y}}(F\left(  \Theta\right)  )$ by
$\{\left(  {B}_{\mathbf{x}}^{q}\right)  _{\beta}\}_{\beta\in\mathcal{B}%
_{\mathbf{x}}^{q}}$ and $\{\left(  {A}_{p}^{\mathbf{y}}\right)  _{\alpha
}\}_{\alpha\in\mathcal{A}_{p}^{\mathbf{y}}},$ respectively. Then $A_{i}%
^{j}=B_{i}^{j}=\theta_{i}^{j}$ with $\left\vert \theta_{i}^{j}\right\vert
=i+j-3;$ $B_{\mathbf{x}}^{1}=\mathbf{\theta}_{m}^{p,i}$ with $\mathbf{x}%
=\mathbf{x}_{m}^{p,i}$ and $|B_{\mathbf{x}}^{1}|=m-2,$ and $A_{1}^{\mathbf{y}%
}=\mathbf{\theta}_{q,j}^{n}$ with $\mathbf{y}=\mathbf{y}_{q,j}^{n}$ and
$|A_{1}^{\mathbf{y}}|=n-2$ (c.f. Example \ref{non-sigma}). In general, for
$p,q\geq2,$ $|(B_{\mathbf{x}}^{q})_{\beta}|=|\mathbf{x}|+q-p-2$ and
$|(A_{p}^{\mathbf{y}})_{\alpha}|=|\mathbf{y}|+p-q-2.$ Then each
$\overset{_{\wedge}}{e}_{A_{\alpha}}$ is a subcomplex of $\Delta^{\left(
q-1\right)  }\left(  P_{p-1}\right)  $ with the associated sign $\left(
-1\right)  ^{\epsilon_{\alpha}}$ and each $\overset{_{\vee}}{e}_{B_{\beta}}$
is a subcomplex of $\Delta^{\left(  p-1\right)  }\left(  P_{q-1}\right)  $
with the associated sign $\left(  -1\right)  ^{\epsilon_{\beta}}$. Define a
differential $\partial:F(\Theta)\rightarrow F(\Theta)$ of degree $-1$ as
follows: Define $\partial$ on generators by
\begin{equation}
\partial(\theta_{m}^{n})=\sum_{(\alpha,\beta)\in\mathcal{A}\mathcal{B}_{m,n}%
}(-1)^{\epsilon+\epsilon_{\alpha}+\epsilon_{\beta}}\gamma\left[  ({A}%
_{p}^{\mathbf{y}})_{\alpha};({B}_{\mathbf{x}}^{q})_{\beta}\right]
,\label{diff}%
\end{equation}
where $\left(  -1\right)  ^{\epsilon}$ is the standard sign of $e_{(\mathbf{y}%
,\mathbf{x})}\subset P_{m+n-2}$. Extend $\partial$ as a derivation of $\gamma
$; then $\partial^{2}=0$ follows from the associativity of $\gamma$. The
DG\ matrad $\left(  F(\Theta),\partial\right)  ,$ denoted by $\mathcal{H}%
_{\infty}$ and called the $A_{\infty}$\emph{-bialgebra matrad}, is realized by
the biassociahedra $\left\{  KK_{n,m}\leftrightarrow\theta_{m}^{n}\right\}  $
(see Theorem \ref{kkh} below). One recovers $\mathcal{A}_{\infty}$ by
restricting $\partial$ to $(\mathcal{H}_{\infty})_{1,\ast}$ or $(\mathcal{H}%
_{\infty})_{\ast,1}.$ Note that $\epsilon=i(m-1)$ in (\ref{diff}) gives the
sign of the cell $e_{(\mathbf{y}_{p,i}^{m},\mathbf{1})}=e_{(\mathbf{1}%
,\mathbf{x}_{m}^{p,i})}\subset P_{m+p-2}$ (see \cite{SU2}). This simplifies
the standard sign in the differential on $\mathcal{A}_{\infty}$ \cite{MSS}.
\end{example}

\begin{example}
\label{index1.234}For $\mathbf{x}=(2,1)$ and $\mathbf{y}=(1,1,1)^{T},$ we have
$\mathcal{A}_{2}^{111}=\{\alpha\}$ and $\mathcal{B}_{21}^{3}= 
\linebreak \{\beta_{1}%
,\beta_{2},\beta_{3}\}.$ The corresponding bases are%
\[%
\begin{array}
[c]{l}%
A_{\alpha}=\left[
\begin{array}
[c]{r}%
\theta_{2}^{1}\\
\theta_{2}^{1}\\
\theta_{2}^{1}%
\end{array}
\right]  \text{ \ and \ }B_{\beta_{1}}=\left[
\begin{array}
[c]{r}%
\theta_{2}^{3}\,\ \ \left[
\begin{array}
[c]{c}%
\theta_{1}^{1}\\
\theta_{1}^{2}%
\end{array}
\right]  \theta_{1}^{2}%
\end{array}
\right]  ,\ \\
\\
B_{\beta_{2}}=\left[
\begin{array}
[c]{r}%
\left[
\begin{array}
[c]{r}%
\theta_{2}^{2}\\
\theta_{2}^{1}%
\end{array}
\right]  [\theta_{1}^{2}\,\theta_{1}^{2}]\ \ \ \theta_{1}^{3}%
\end{array}
\right]  ,\text{ \ }B_{\beta_{3}}=\left[
\begin{array}
[c]{r}%
\left[
\begin{array}
[c]{c}%
\theta_{1}^{2}\\
\theta_{1}^{1}%
\end{array}
\right]  \theta_{2}^{2}\ \ \ \theta_{1}^{3}%
\end{array}
\right]  .
\end{array}
\]
Thus, $\partial(\theta_{3}^{3})=-\gamma(A_{\alpha}\,;\,B_{\beta_{1}}%
+B_{\beta_{2}}+B_{\beta_{3}})+\cdots$ (see Example \ref{faces1.234}).
\end{example}

\subsection{The Biderivative}

\label{sbiderivative}In \cite{SU3} we used the canonical prematrad structure
$\gamma$ on the universal PROP $U_{A}=End(TA)$ to define the biderivative
operator. By replacing $U_{A}$ with an arbitrary prematrad $(M,\gamma)$ we
obtain the general biderivative operator $Bd_{\gamma}:\mathbf{M}%
\rightarrow\mathbf{M}$ having the property $Bd_{\gamma}\circ Bd_{\gamma
}=Bd_{\gamma}. $ An element $A\in\mathbf{M}$ is a \emph{$\gamma$-biderivative}
if $A=Bd_{\gamma}(A).$ Note that $Bd_{\gamma}\left(  \mathbf{M}\right)
\subseteq(\Gamma(M),\gamma);$ when $M$ is generated by singletons in
each\ bidegree, the image $Bd_{\gamma}(\mathbf{M})$ is the module of fixed
points of $Bd_{\gamma}$ and gives rise to an algorithmic construction of an
additive basis for the $A_{\infty}$-bialgebra matrad $\mathcal{H}_{\infty}$.
More precisely:

\begin{proposition}
Let $(M,\gamma)$ be a prematrad generated by $\Theta=\{\theta_{m}%
^{n}\}_{m,n\geq1},$ and let $Bd_{\gamma}:\mathbf{M}\rightarrow\mathbf{M}$
denote the associated biderivative operator. Then

\begin{enumerate}
\item[\textit{(i)}] $Bd_{\gamma}(\mathbf{M})\subseteq\Gamma(M,\gamma)$ and
$Bd_{\gamma}\circ Bd_{\gamma}=Bd_{\gamma}.$

\item[\textit{(ii)}] Each element $\theta\in\Theta$ has a unique $\gamma
$-biderivative $d_{\theta}^{\gamma}\in\mathbf{M}.$

\item[\textit{(iii)}] $\Gamma(M,\gamma)=\left\langle d_{\theta}^{\gamma
}\right\rangle $.
\end{enumerate}
\end{proposition}

\noindent Thus the $\gamma$-biderivative can be viewed as a non-linear map
$d_{\bullet}^{\gamma}:M\rightarrow\mathbf{M}.$ When $M=End\left(  TA\right)  $
we omit the symbol $\gamma$ and denote the biderivative of $\theta$ by
$d_{\theta}$ as in \cite{SU3}.

In particular, the modules $\langle Bd_{\gamma}(M_{\ast,2})\rangle\subseteq
T^{+}(M_{\ast,2})$ and $\langle Bd_{\gamma}(M_{2,\ast})\rangle\subseteq
T^{+}(M_{2,\ast})$ are spanned by symmetric tensors (compare Example
\ref{vectormodule}); furthermore, $Bd_{\gamma}(M_{n,m})=M_{n,m}$ for $m,n>2.$

Finally, the algorithm that produces $d_{\theta}^{\gamma}$ for $(M,\gamma
)=(F^{\operatorname*{pre}}\left(  \Theta\right)  ,\gamma)$ simultaneously
produces an additive basis for $\mathcal{H}_{\infty}.$

\begin{proposition}
\label{bases}Let $\Theta=\left\langle \theta_{m}^{n}\neq0\right\rangle
_{m,n\geq1}$ as in Example \ref{freematrad}. Elements of the bases $\{({A}%
_{p}^{\mathbf{y}})_{\alpha}\}_{\alpha\in\mathcal{A}_{p}^{\mathbf{y}}}$ and
$\{({B}_{\mathbf{x}}^{q})_{\beta}\}_{\beta\in\mathcal{B}_{\mathbf{x}}^{q}}$
are exactly the components of $d_{\theta}^{\gamma}$ in ${\Gamma}%
_{p}^{\mathbf{y}}(F\left(  \Theta\right)  ,\gamma)$ and ${\Gamma}_{\mathbf{x}%
}^{q}(F\left(  \Theta\right)  ,\gamma)$ with degrees $|\mathbf{y}|+p-q-2$ and
$|\mathbf{x}|+q-p-2,$ respectively. Thus
\[
\partial\left(  \theta_{m}^{n}\right)  =\sum_{\substack{\left\vert
\mathbf{x}\right\vert =m;\text{ }\left\vert \mathbf{y}\right\vert =n \\A\times
B\in(d_{\theta}^{\gamma})_{p}^{\mathbf{y}}\times(d_{\theta}^{\gamma
})_{\mathbf{x}}^{q}}}\gamma\left(  A;B\right)  .
\]

\end{proposition}

\begin{proof}
The proof follows from the definition of $d^{\gamma}_{\theta}$ and is straightforward.
\end{proof}

\subsubsection{The $\circledcirc$-Product}

Given a prematrad $\left(  M,\gamma\right)  ,$ define a (non-bilinear)
operation
\begin{equation}
\circledcirc:M\times M\overset{d_{\bullet}^{\gamma}\times d_{\bullet}^{\gamma
}}{\longrightarrow}\mathbf{M}\times\mathbf{M}\overset{\Upsilon
}{\longrightarrow}\mathbf{M}\overset{{proj}}{\rightarrow}M,\label{operation1}%
\end{equation}
where ${proj}$ is the canonical projection. The following facts are now obvious:

\begin{proposition}
The $\circledcirc$ operation acts bilinearly on $M_{\ast,1}$ and
${\ M_{1,\ast}.}$ In fact, when $M=End\left(  TH\right)  ,$ the $\circledcirc$
operation coincides with Gerstenhaber's $\circ$-operation on $M_{1.\ast}$ (see
\cite{Gersten}) and dually on $M_{\ast,1}$.
\end{proposition}

\begin{remark}
The bilinear part of the $\circledcirc$ operation, i.e., its restriction to
either $M_{\ast,1}$ or $M_{1,\ast},$ is completely determined by the
associahedra $K=\sqcup K_{n}$ (rather than permutahedra) and induces the
cellular projection $\vartheta_{0}:P_{n}\rightarrow K_{n+1}$ due to A. Tonks
\cite{Tonks}.
\end{remark}

\section{The Posets $\mathcal{PP}$ and $\mathcal{KK}$}

In this section we construct a poset $\mathcal{PP}$ and an appropriate
quotient poset $\mathcal{KK}$. The elements of $\mathcal{KK}$ correspond with
the 0-dimensional module generators of the free matrad ${\mathcal{H}}_{\infty
}.$ The geometric realization of $\mathcal{KK}$, constructed in Section 8, is
the disjoint union of biassociahedra $KK=\left\{  KK_{n,m}\right\}
_{m,n\geq1}\ $whose cellular chains are identified with ${\mathcal{H}}%
_{\infty}$.

Let $\mathcal{V}_{n}$ denote the set of vertices of $P_{n}$ and identify
$\mathcal{V}_{n}$ with the set $S_{n}$ of permutations of $\underline{n}%
=\left\{  1,2,\ldots,n\right\}  $ via the standard bijection $\mathcal{V}%
_{n}\leftrightarrow S_{n}.$ The Bruhat partial ordering on $S_{n}$ generated
by the relation $a_{1}|\cdots|a_{n}<a_{1}|\cdots|a_{i+1}|a_{i}|\cdots|a_{n}$
if and only if $a_{i}<a_{i+1}$ imposes a poset structure on $\mathcal{V}_{n}.
$ For $n\geq1,$ set ${\mathcal{PP}}_{n,0}={\mathcal{PP}}_{0,n}=\mathcal{V}%
_{n}$ and define the \emph{geometric realization} $PP_{n,0}=\left\vert
\mathcal{PP}_{n,0}\right\vert =\left\vert \mathcal{PP}_{0,n}\right\vert
=PP_{0,n}$ to be the permutahedron $P_{n}.$ Then $KK_{n+1,1}=\left\vert
\mathcal{KK}_{n+1,1}\right\vert =\left\vert \mathcal{KK}_{1,n+1}\right\vert
=KK_{1,n+1}$ is the Stasheff associahedron $K_{n+1}$ (see \cite{Stasheff},
\cite{Stasheff2}, \cite{SU2}). In the discussion that follows, we construct
the posets ${\mathcal{PP}}_{n,m}$ and $\mathcal{KK}_{n+1,m+1}$ and their
geometric realizations $PP_{n,m}$ and $KK_{n+1,m+1}$ for all $m,n\geq1$.

Denote the sets of up-rooted and down-rooted binary trees with $n+1$ leaves
and $n$ levels by $\wedge_{n}$ and $\vee_{n},$ respectively; then each vertex
of $P_{n}$ is indexed by two trees, one the reflection of the other. These
indexing sets have a poset structure induced by the standard bijections
$\hat{\ell}:\wedge_{n}\rightarrow\mathcal{V}_{n}$ and $\check{\ell}:\vee
_{n}\rightarrow\mathcal{V}_{n}$, and the products $\wedge_{m}^{\times n},$
$\vee_{n}^{\times m},$ and ${\wedge}_{m}^{\times n}\times{\vee_{n}}^{\times
m}$ are posets with respect to lexicographic ordering. Now consider the
subcomplexes $\Delta^{\left(  n\right)  }\left(  P_{m}\right)  \subseteq
P_{m}^{\times n+1}$ and $\Delta^{\left(  m\right)  }\left(  P_{n}\right)
\subseteq P_{n}^{\times m+1}$ with faces of $P_{m}$ and $P_{n}$ indexed by
up-rooted and down-rooted PLTs, respectively. Then the $0$-skeletons
$X_{m}^{n+1}\subseteq\Delta^{\left(  n\right)  }\left(  P_{m}\right)  $ and
$Y_{n}^{m+1}\subseteq\Delta^{\left(  m\right)  }\left(  P_{n}\right)  $ are
subposets of $\wedge_{m}^{\times n+1}$ and $\vee_{n}^{\times m+1}$ and there
is the inclusion of posets
\[
X_{m}^{n+1}\times Y_{n}^{m+1}\hookrightarrow{\wedge}_{m}^{\times n+1}%
\times{\vee_{n}^{\times m+1}}.
\]

Express $x\in\wedge_{m}^{\times n}$ as an $n\times1$ column matrix of
up-rooted binary trees and replace $x$ with its (unique) BTP\ $\Upsilon
$-factorization $x_{1}\cdots x_{m}\in\overline{\mathbf{M}},$ where $x_{i}$ is
an $n\times i$ matrix over $\left\{  \mathbf{1,\curlywedge}\right\}  $ with
$\mathbf{\curlywedge}$ appearing in each row exactly once. Dually, express
$y\in\vee_{n}^{\times m}$ as a $1\times m$ row matrix of down-rooted binary
trees and replace $y$ with its (unique) BTP\ factorization as a $\Upsilon
$-product $y_{n}\cdots y_{1}\in\overline{\mathbf{M}},$ where $y_{j}$ is an
$j\times m$ matrix over $\left\{  \mathbf{1,}\curlyvee\right\}  $ with
$\curlyvee$ appearing in each column exactly once.

\begin{example}
\label{Ex1}Whereas the product $\Delta^{\left(  1\right)  }\left(
P_{1}\right)  =P_{1}^{\times2}$ can be thought of as either
$\mathbf{\curlywedge}\times\mathbf{\curlywedge}$ or $\curlyvee\times
\curlyvee,$ we have%
\[
X_{1}^{2}=\left[
\begin{array}
[c]{c}%
\curlywedge\\
\curlywedge
\end{array}
\right]  \text{ \ and \ }Y_{1}^{2}=\left[  \curlyvee\curlyvee\right]  \text{
\ so that }X_{1}^{2}\times Y_{1}^{2}=\left[
\begin{array}
[c]{c}%
\curlywedge\\
\curlywedge
\end{array}
\right]  \left[  \curlyvee\curlyvee\right]  .
\]
The poset of vertices in $\Delta^{\left(  1\right)  }\left(  P_{2}\right)
\subset P_{2}^{\times2}$ expresses the following products of permutations and
matrix sequences:%
\[%
\begin{tabular}
[c]{ccccccc}%
$a|b\times c|d$ & $:$ & $1|2\times1|2$ & $<$ & $1|2\times2|1$ & $<$ &
$2|1\times2|1$\\
&  &  &  &  &  & \\
$X_{2}^{2}$ & $\emph{:}$ & $\left[
\begin{array}
[c]{c}%
\curlywedge\\
\curlywedge
\end{array}
\right]  \left[
\begin{array}
[c]{c}%
\curlywedge\text{ }\mathbf{1}\\
\curlywedge\text{ }\mathbf{1}%
\end{array}
\right]  $ & $<$ & $\left[
\begin{array}
[c]{c}%
\curlywedge\\
\curlywedge
\end{array}
\right]  \left[
\begin{array}
[c]{c}%
\curlywedge\text{ }\mathbf{1}\\
\mathbf{1}\text{ }\curlywedge
\end{array}
\right]  $ & $<$ & $\left[
\begin{array}
[c]{c}%
\curlywedge\\
\curlywedge
\end{array}
\right]  \left[
\begin{array}
[c]{c}%
\mathbf{1}\text{ }\curlywedge\\
\mathbf{1}\text{ }\curlywedge
\end{array}
\right]  $\\
&  &  &  &  &  & \\
$Y_{2}^{2}$ & $:$ & $\left[
\begin{array}
[c]{c}%
\curlyvee\text{ }\curlyvee\\
\mathbf{1}\text{ \ }\mathbf{1}%
\end{array}
\right]  \left[  \curlyvee\curlyvee\right]  $ & $<$ & $\left[
\begin{array}
[c]{c}%
\curlyvee\text{ }\mathbf{1}\\
\mathbf{1}\text{ }\curlyvee
\end{array}
\right]  \left[  \curlyvee\curlyvee\right]  $ & $<$ & $\left[
\begin{array}
[c]{c}%
\mathbf{1}\text{ \ }\mathbf{1}\\
\curlyvee\text{ }\curlyvee
\end{array}
\right]  \left[  \curlyvee\curlyvee\right]  $%
\end{tabular}
\]
Furthermore, thinking of the product $\Delta^{\left(  2\right)  }\left(
P_{1}\right)  =P_{1}^{\times3}$ as $\curlyvee\times\curlyvee\times\curlyvee,$
we have $Y_{1}^{3}=\left[  \curlyvee\curlyvee\curlyvee\right]  ;$ consequently
$X_{2}^{2}\times Y_{1}^{3}=$
\[
\left\{  \left[
\begin{array}
[c]{c}%
\curlywedge\\
\curlywedge
\end{array}
\right]  \left[
\begin{array}
[c]{c}%
\curlywedge\text{ }\mathbf{1}\\
\curlywedge\text{ }\mathbf{1}%
\end{array}
\right]  \left[  \curlyvee\curlyvee\curlyvee\right]  <\left[
\begin{array}
[c]{c}%
\curlywedge\\
\curlywedge
\end{array}
\right]  \left[
\begin{array}
[c]{c}%
\curlywedge\text{ }\mathbf{1}\\
\mathbf{1}\text{ }\curlywedge
\end{array}
\right]  \left[  \curlyvee\curlyvee\curlyvee\right]  <\left[
\begin{array}
[c]{c}%
\curlywedge\\
\curlywedge
\end{array}
\right]  \left[
\begin{array}
[c]{c}%
\mathbf{1}\text{ }\curlywedge\\
\mathbf{1}\text{ }\curlywedge
\end{array}
\right]  \left[  \curlyvee\curlyvee\curlyvee\right]  \right\}  .
\]

\end{example}

\begin{definition}
Let $A=\left[  a_{ij}\right]  $ be an $(n+1)\times m$ matrix over
$\{\mathbf{1},\curlywedge\},$ each row of which contains the entry
$\curlywedge$ exactly once. Let $B=\left[  b_{ij}\right]  $ be an
$n\times(m+1)$ matrix over $\{\mathbf{1},\curlyvee\},$ each column of which
contains the entry $\curlyvee$ exactly once. Then $(A,B)$ is an $\left(
i,j\right)  $\textbf{-edge pair} if

\begin{enumerate}
\item[\textit{(i)}] $A\otimes B$ is a BTP,

\item[\textit{(ii)}] $a_{ij}=a_{i+1,j}=\curlywedge$ and $b_{ij}=b_{i,j+1}%
=\curlyvee$.
\end{enumerate}
\end{definition}

\noindent For $A_{1}\cdots A_{m}B_{n}\cdots B_{1}\in X_{m}^{n+1}\times
Y_{n}^{m+1},$ the only possible edge pair in is $\left(  A_{m},B_{n}\right)  .
$ In $X_{2}^{3}\times Y_{2}^{3},$ for example, the respective matrix
sequences
\[
\left[
\begin{array}
[c]{c}%
\curlywedge\\
\curlywedge\\
\curlywedge
\end{array}
\right]  \left[
\begin{array}
[c]{c}%
\curlywedge\text{ }\mathbf{1}\\
\curlywedge\text{ }\mathbf{1}\\
\curlywedge\text{ }\mathbf{1}%
\end{array}
\right]  \left[
\begin{array}
[c]{c}%
\curlyvee\curlyvee\text{ }\mathbf{1}\\
\mathbf{1}\text{ }\mathbf{1}\text{ }\curlyvee
\end{array}
\right]  \left[  \curlyvee\curlyvee\curlyvee\right]  \text{ and }\left[
\begin{array}
[c]{c}%
\curlywedge\\
\curlywedge\\
\curlywedge
\end{array}
\right]  \left[
\begin{array}
[c]{c}%
\curlywedge\text{ }\mathbf{1}\\
\curlywedge\text{ }\mathbf{1}\\
\curlywedge\text{ }\mathbf{1}%
\end{array}
\right]  \left[
\begin{array}
[c]{c}%
\curlyvee\text{ }\mathbf{1}\text{ }\mathbf{1}\\
\mathbf{1}\text{ }\curlyvee\curlyvee
\end{array}
\right]  \left[  \curlyvee\curlyvee\curlyvee\right]
\]
\emph{do} and \emph{do not} contain an edge pair.

\begin{definition}
Let $Q$ be a poset and let $x_{1}\leq x_{2}\in Q$. The pair $\left(
x_{1},x_{2}\right)  $ is an \textbf{edge of} $Q$ if $x\in Q$ and $x_{1}\leq
x\leq x_{2}$ implies $x=x_{1}$ or $x=x_{2}.$
\end{definition}

\noindent Edges of $X_{m}^{n+1}\times Y_{n}^{m+1}$ correspond to 1-dimensional
elements of $\mathcal{H}_{\infty}$ generated by $\left\{  \mathbf{1}%
,\theta_{2}^{1},\theta_{1}^{2},\right.  $ $\left.  \theta_{3}^{1},\theta
_{1}^{3}\right\}  $; 1-dimensional elements of $\mathcal{H}_{\infty}$
generated by $\left\{  \mathbf{1},\theta_{2}^{1},\theta_{1}^{2},\theta_{2}%
^{2}\right\}  $ correspond to edges of a poset $Z_{n,m}$ related to but
disjoint from $X_{m}^{n+1}\times Y_{n}^{m+1},$ which we now define.

Let $A^{i\ast}$ and $B^{\ast j}$ denote the matrices obtained by deleting the
$i^{th}$ row of $A$ and the $j^{th}$ column of $B$.

\begin{definition}
\label{transposition}Let $c=C_{1}\cdots C_{r}$ be a sequence of matrices such
that $(C_{k},C_{k+1})$ is an $\left(  i,j\right)  $-edge pair for some $k\leq
r-1$ and some association of $C_{1}\cdots\left(  C_{k}C_{k+1}\right)  \cdots
C_{r}$ defines a sequence of BTPs. The $(i,j)$\textbf{-transposition of} $c$
\textbf{in position} $k$ is the sequence
\[
\mathcal{T}_{ij}^{k}(c)=C_{1}\cdots C_{k+1}^{\ast j}C_{k}^{i\ast}\cdots C_{r}.
\]
The symbol $\mathcal{T}_{ij}^{k}\left(  c\right)  $ implies that the action of
$\mathcal{T}_{ij}^{k}$ on $c$ is defined.
\end{definition}

\noindent Note that if $\mathcal{T}_{ij}^{k}$ acts on $u=A_{1}\cdots
A_{m}B_{n}\cdots B_{1}\in X_{m}^{n+1}\times Y_{n}^{m+1},$ then $k=m$ and the
potential edge pairs of consecutive matrices in
\[
\mathcal{T}_{ij}^{m}\left(  u\right)  =A_{1}\cdots A_{m-1}B_{n}^{\ast j}%
A_{m}^{i\ast}B_{n-1}\cdots B_{1}%
\]
are $(A_{m-1},B_{n}^{\ast j})$ and $(A_{m}^{i\ast},B_{n-1})$. If
$(A_{m-1},B_{n}^{\ast j})$ is an edge pair and $\mathcal{T}_{kl}^{m-1}$ is
defined on $\mathcal{T}_{ij}^{m}\left(  u\right)  $, then
\[
\mathcal{T}_{kl}^{m-1}\mathcal{T}_{ij}^{m}\left(  u\right)  =A_{1}\cdots
A_{m-2}B_{n}^{\ast j\ast l}A_{m-1}^{k\ast}A_{m}^{i\ast}B_{n-1}\cdots B_{1},
\]
and so on. In this manner, iterate $\mathcal{T}$ on each element $u\in
X_{m}^{n+1}\times Y_{n}^{m+1}$ in all possible ways and obtain%
\[
Z_{n,m}=\left\{  \left.  \mathcal{T}_{i_{t}j_{t}}^{k_{t}}\cdots\mathcal{T}%
_{i_{1}j_{1}}^{k_{1}}\left(  u\right)  \text{ }\right\vert \text{ }u\in
X_{m}^{n+1}\times Y_{n}^{m+1},\text{ }t\geq1\right\}  .
\]
Then
\[
{\mathcal{PP}}_{n,m}=X_{m}^{n+1}\times Y_{n}^{m+1}\cup Z_{n,m}.
\]

To extend the partial ordering to $Z_{n,m},$ first define $c<\mathcal{T}%
_{ij}^{k}\left(  c\right)  $ for $c\in\mathcal{PP}_{n,m}.$ To define a
generating relation on $Z_{n,m},$ note that each composition $\mathcal{T}%
_{i_{t}j_{t}}^{k_{t}}\cdots\mathcal{T}_{i_{1}j_{1}}^{k_{1}}$ defined on $u\in
X_{m}^{n+1}\times Y_{n}^{m+1}$ uniquely determines an $\left(  m,n\right)
$-shuffle $\sigma,$ in which case we denote%
\[
\mathcal{T}_{\sigma}\left(  u\right)  =\mathcal{T}_{i_{t}j_{t}}^{k_{t}}%
\cdots\mathcal{T}_{i_{1}j_{1}}^{k_{1}}\left(  u\right)
\]
and define $\mathcal{T}_{\operatorname{Id}}=\operatorname{Id}.$ When
$\mathcal{T}_{\sigma}\left(  u\right)  $ is defined, multiple compositions of
$\left(  i,j\right)  $-transposi- tions on $u$ may determine the same $\sigma$;
thus $\mathcal{T}_{\sigma}\left(  u\right)  $ is a set, in general. For
$u_{1}\leq u_{2}\in X_{m}^{n+1}\times Y_{n}^{m+1},$ define $\mathcal{T}%
_{\sigma}\left(  u_{1}\right)  \leq\mathcal{T}_{\sigma}\left(  u_{2}\right)  $
if $\left(  u_{1},u_{2}\right)  $ is an edge of $X_{m}^{n+1}\times Y_{n}%
^{m+1}$ or $u_{2}$ is \textquotedblleft$\sigma$ -compatible\textquotedblright%
\ with $u_{1}$ in the following sense: Let $a=a_{1}|\cdots|a_{m}\in S_{m}$ and
$b=b_{1}|\cdots|b_{n}\in S_{n}$. The action of $\sigma$ on $\left(
a;b\right)  $ decomposes $a$ and $b$ into subsequences $\mathbf{m}_{1}%
,\ldots,\mathbf{m}_{k}$ and $\mathbf{n}_{1},\ldots,\mathbf{n}_{l}$ in one of
the following four ways:%
\[
\sigma\left(  a;b\right)  =\left\{
\begin{array}
[c]{ll}%
\mathbf{m}_{1},\mathbf{n}_{1},\mathbf{m}_{2},\mathbf{n}_{2},\ldots
,\mathbf{n}_{k-1},\mathbf{m}_{k}, & \sigma\left(  a_{1}\right)  =a_{1},\text{
}\sigma\left(  b_{n}\right)  \neq b_{n}\\
\mathbf{m}_{1},\mathbf{n}_{1},\mathbf{m}_{2},\mathbf{n}_{2},\ldots
,\mathbf{m}_{k},\mathbf{n}_{k}, & \sigma\left(  a_{1}\right)  =a_{1},\text{
}\sigma\left(  b_{n}\right)  =b_{n}\\
\mathbf{n}_{1},\mathbf{m}_{1},\mathbf{n}_{2},\mathbf{m}_{2},\ldots
,\mathbf{n}_{k},\mathbf{m}_{k}, & \sigma\left(  a_{1}\right)  \neq
a_{1},\text{ }\sigma\left(  b_{n}\right)  \neq b_{n}\\
\mathbf{n}_{1},\mathbf{m}_{1},\mathbf{n}_{2},\mathbf{m}_{2},\ldots
,\mathbf{m}_{k},\mathbf{n}_{k+1}, & \sigma\left(  a_{1}\right)  \neq
a_{1},\text{ }\sigma\left(  b_{n}\right)  =b_{n}.
\end{array}
\right.
\]
Define $\mathbf{I}_{\sigma}=\{\left(  \alpha_{1},...,\alpha_{n+1}\right)  \in
S_{m}^{\times n+1}\mid\alpha_{i}\in S_{\#\mathbf{m}_{1}}\times\cdots\times
S_{\#\mathbf{m}_{k}}\subset S_{m}\}\ $for all $i$ and $\mathbf{J}_{\sigma
}=\{\left(  \beta_{1},...,\beta_{m+1}\right)  \in S_{n}^{\times m+1}\mid
\beta_{j}\in S_{\#\mathbf{n}_{1}}\times\cdots\times S_{\#\mathbf{n}_{l}%
}\subset S_{n}\}$ for all $j.$ Let $\chi:\mathcal{V}_{m}\rightarrow
\mathcal{V}_{m}$ be the involutory bijection defined by
\[
\chi(a_{1}|\cdots|a_{m})=\left(  m+1-a_{m}\right)  |\cdots|\left(
m+1-a_{1}\right)
\]
and fix the inclusion of posets
\begin{equation}
X_{m}^{n+1}\times Y_{n}^{m+1}\overset{\kappa}{\hookrightarrow}\mathcal{V}%
_{m}^{\times n+1}\times\mathcal{V}_{n}^{\times m+1}\leftrightarrow
S_{m}^{\times n+1}\times S_{n}^{\times m+1},\label{chia}%
\end{equation}
where $\kappa=\left(  \chi\circ\hat{\ell}\right)  ^{\times n+1}\times\left(
{\check{\ell}}\right)  ^{\times m+1}.$ Then $u_{2}$ is $\sigma$%
\emph{-compatible} with $u_{1}$ if $u_{2}=(\alpha\times\beta)(u_{1})$ for some
$\alpha\times\beta\in\mathbf{I}_{\sigma}\times\mathbf{J}_{\sigma}.$

To view this geometrically, suppose $u_{2}=A_{1}^{\prime}\cdots A_{m}^{\prime
}B_{n}^{\prime}\cdots B_{1}^{\prime}$ is $\sigma$-compatible with $u_{1}%
=A_{1}\cdots A_{m}B_{n}\cdots B_{1}$ in $X_{m}^{n+1}\times Y_{n}^{m+1}.$ For
each $i,$ let $a_{i}=i_{1}|\cdots|i_{m}$ and $a_{i}^{\prime}=i_{1}^{\prime
}|\cdots|i_{m}^{\prime}$ be the permutations of $\underline{m}$ corresponding
with the up-rooted trees given by $\gamma$-products $A_{i,1}\cdots A_{i,m}$
and $A_{i,1}^{\prime}\cdots A_{i,m}^{\prime}$of $i^{th}$ rows, respectively;
dually, for each $j,$ let $b_{j}=j_{1}|\cdots|j_{n}$ and $b_{j}^{\prime}%
=j_{1}^{\prime}|\cdots|j_{n}^{\prime}$ be the permutations of $\underline{n}$
corresponding with the down-rooted trees given by the $\gamma$-products
$B_{n,j}\cdots B_{1,j}$ and $B_{n,j}^{\prime}\cdots B_{1,j}^{\prime}$ of
$j^{th}$ columns, respectively. Then for each $\left(  i,j\right)  ,$ the
$\sigma$-partition of $\left(  a_{i};b_{j}\right)  $ determines a product face
$\mathbf{m}_{1}|\cdots|\mathbf{m}_{k}\times\mathbf{n}_{1}|\cdots
|\mathbf{n}_{l}\subset P_{m}\times P_{n}$ containing the vertices $a_{i}\times
b_{j}$ and $a_{i}^{\prime}\times b_{j}^{\prime}$ and an oriented path of edges
from $a_{i}\times b_{j}$ to $a_{i}^{\prime}\times b_{j}^{\prime}.$

\begin{remark}
\label{edge}\ \vspace*{-0.15in}\newline

\begin{enumerate}
\item[\textit{(i)}] The map $\chi$ used to define the poset structure of
$\mathcal{PP}_{n,m}$ is evoked to induce the correct orientation of the
quotient poset $\mathcal{KK}_{n+1,m+1}$ (see below), and is necessary to
establish the bijection in Theorem 1 (see also item (iii) below). For
geometric realizations of $\mathcal{KK}_{n+1,m+1}$ and $\mathcal{KK}%
_{m+1,n+1}$ compare Figures 21 and 22.\smallskip

\item[\textit{(ii)}] Note that if $\left(  u_{1},u_{2}\right)  $ is an edge of
$X_{m}^{n+1}\times Y_{n}^{m+1},$ the partial ordering in $\mathcal{PP}_{n,m}$
implies that $\left(  \mathcal{T}(u_{1}),\mathcal{T}(u_{2})\right)  $ is an
edge of $\mathcal{PP}_{n,m}.$\smallskip

\item[\textit{(iii)}] The transpose map $X_{m}^{n+1}\times Y_{n}%
^{m+1}\rightarrow X_{n}^{m+1}\times Y_{m}^{n+1}$ given by
\[
A_{1}\cdots A_{m}B_{n}\cdots B_{1}\mapsto B_{1}^{T}\cdots B_{n}^{T}A_{m}%
^{T}\cdots A_{1}^{T}%
\]
induces a canonical order-preserving bijection ${\mathcal{PP}}_{n,m}%
\leftrightarrow{\mathcal{PP}}_{m,n}.$
\end{enumerate}
\end{remark}

\begin{example}
\label{compatible}Using the notation of Example \ref{Ex1}, let us determine
those elements $u_{i}\in X_{2}^{2}\times Y_{1}^{3}$ that are $\sigma
$-compatible with $u_{1}$. Since all matrices in $u_{1}$ have constant columns
or rows, $a_{i}\times b_{j}=1|2\times1\subset P_{2}\times P_{1}$ for all
$i,j$. The $\left(  2,1\right)  $-shuffles of $\left(  1,2;3\right)  $ are
$\sigma_{0}=1|2|3,$ $\sigma_{1}=1|3|2$ and $\sigma_{2}=3|1|2.$ The $\sigma
_{1}$-partition of $\left(  1|2;1\right)  $ determines the face $1|2\times
1\subset P_{2}\times P_{1}$ whose only vertex is $u_{1}.$ Hence the only
element of $X_{2}^{2}\times Y_{1}^{3}$ that is $\sigma_{1}$-compatible with
$u_{1}$ is itself. If $\sigma\in\left\{  \sigma_{0},\sigma_{2}\right\}  ,$ the
$\sigma$-partition of $\left(  1|2;1\right)  $ determines the face
$12\times1\subset P_{2}\times P_{1}$ with vertices are $1|2\times1$ and
$2|1\times1.$ Since all matrices in $u_{3}$ have constant columns or rows,
$a_{i}^{\prime}\times b_{j}^{\prime}=2|1\times1$ for all $i,j$ implies that
$u_{3}$ is $\sigma$-compatible with $u_{1}.$ Furthermore, $a_{1}^{\prime
}\times b_{j}^{\prime}=1|2\times1$ and $a_{2}^{\prime}\times b_{j}^{\prime
}=2|1\times1$ for all $j$ implies that $u_{2}$ is also $\sigma$-compatible
with $u_{1}.$ Since $u_{1}<u_{3}$ we have $\mathcal{T}_{\sigma_{2}}\left(
u_{1}\right)  <\mathcal{T}_{\sigma_{2}}\left(  u_{3}\right)  .$
\end{example}

\begin{example}
\label{Ex2}Since $Z_{1,1}=\left[  \curlyvee\right]  \left[  \curlywedge
\right]  $ we have%
\[
{\mathcal{PP}}_{1,1}=\left\{  \left[
\begin{array}
[c]{c}%
\curlywedge\\
\curlywedge
\end{array}
\right]  [\curlyvee\curlyvee]<[\curlyvee][\curlywedge]\right\}  .
\]

\noindent Using the notation of Example \ref{compatible}, the action of
$\mathcal{T}$ on%
\[
u_{1}=\left[
\begin{array}
[c]{c}%
\curlywedge\\
\curlywedge
\end{array}
\right]  \left[
\begin{array}
[c]{c}%
\curlywedge\text{ }\mathbf{1}\\
\curlywedge\text{ }\mathbf{1}%
\end{array}
\right]  \left[  \curlyvee\curlyvee\curlyvee\right]  \text{ and \ }%
u_{3}=\left[
\begin{array}
[c]{c}%
\curlywedge\\
\curlywedge
\end{array}
\right]  \left[
\begin{array}
[c]{c}%
\mathbf{1}\text{ }\curlywedge\\
\mathbf{1}\text{ }\curlywedge
\end{array}
\right]  \left[  \curlyvee\curlyvee\curlyvee\right]
\]
produces the following four elements of $Z_{1,2}:$
\[
u_{1}\overset{\mathcal{T}_{\sigma_{1}}}{\longmapsto}\left[
\begin{array}
[c]{c}%
\curlywedge\\
\curlywedge
\end{array}
\right]  [\curlyvee\curlyvee]\left[  \curlywedge\text{ }\mathbf{1}\right]
\overset{\mathcal{T}_{\sigma_{2}}}{\longmapsto}[\curlyvee][\curlywedge]\left[
\curlywedge\text{ }\mathbf{1}\right]  ;
\]%
\[
u_{3}\overset{\mathcal{T}_{\sigma_{1}}}{\longmapsto}\left[
\begin{array}
[c]{c}%
\curlywedge\\
\curlywedge
\end{array}
\right]  [\curlyvee\curlyvee]\left[  \mathbf{1}\text{ }\curlywedge\right]
\overset{\mathcal{T}_{\sigma_{2}}}{\longmapsto}[\curlyvee][\curlywedge]\left[
\mathbf{1}\text{ }\curlywedge\right]  .
\]
Thus ${\mathcal{PP}}_{1,2}=\left\{  u_{1}<\mathcal{T}_{\sigma_{1}}\left(
u_{1}\right)  <\mathcal{T}_{\sigma_{2}}\left(  u_{1}\right)  ,\text{ }%
u_{2},\text{\ }u_{3}<\mathcal{T}_{\sigma_{1}}\left(  u_{3}\right)
<\mathcal{T}_{\sigma_{2}}\left(  u_{3}\right)  \right\}  .$ Recall\linebreak
that the action of $\mathcal{T}$ on $u_{2}$ is undefined, and as mentioned in
Examples \ref{Ex1} and \ref{compatible}, $u_{1}<u_{2}<u_{3}$ and
$\mathcal{T}_{\sigma_{2}}\left(  u_{1}\right)  <\mathcal{T}_{\sigma_{2}%
}\left(  u_{3}\right)  $ (see Figure 12).
\end{example}

\begin{center}
\includegraphics[
trim=0.000000in -0.208613in 0.000000in 0.000000in,
height=2.7856in,
width=2.7925in
]%
{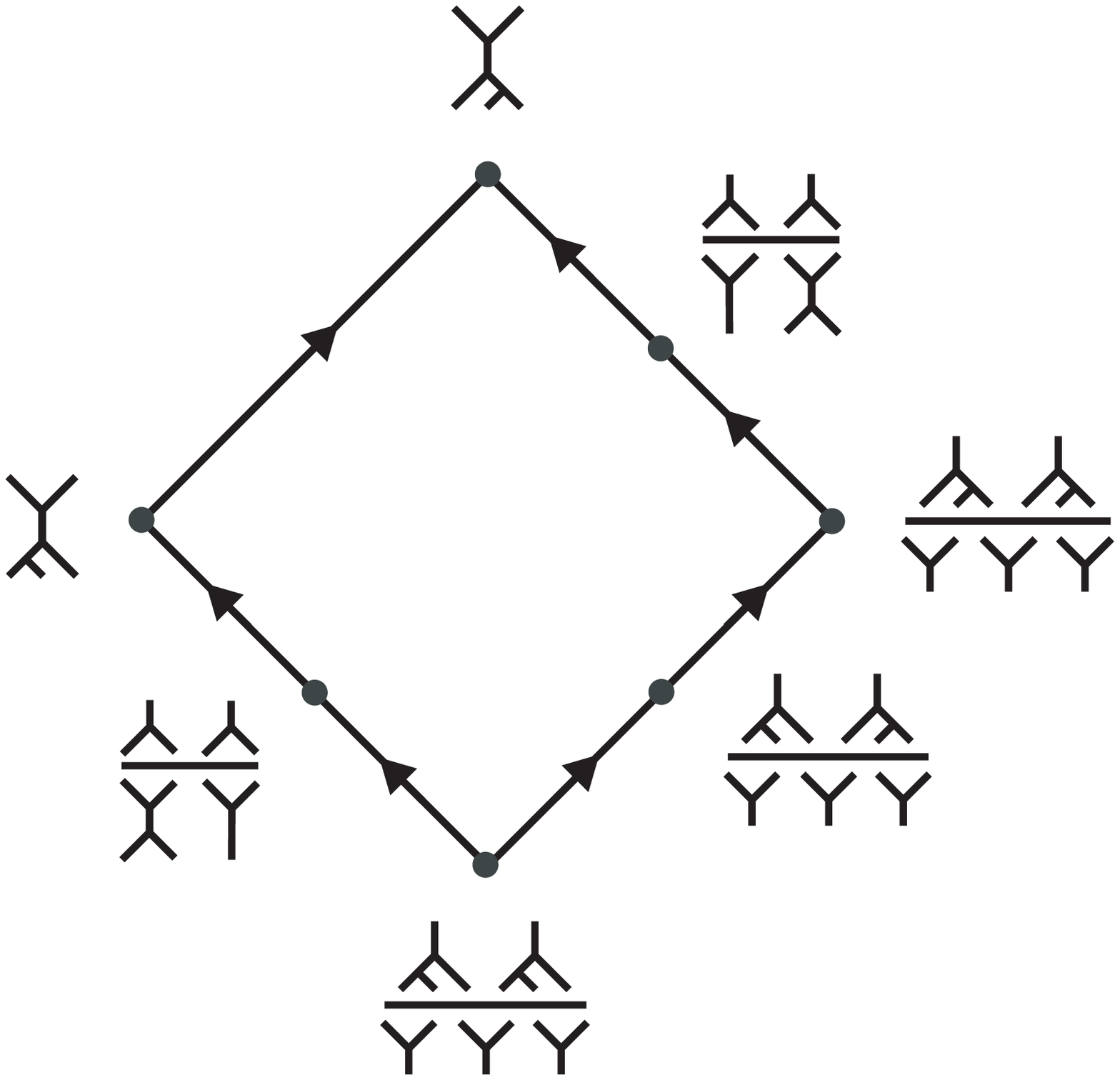}%
\\
Figure 12. The digraph of $\mathcal{PP}_{1,2}$.
\end{center}

One can represent $u=A_{1}\cdots A_{m}B_{n}\cdots B_{1}\in X_{m}^{n+1}\times
Y_{n}^{m+1}$ and $z=\mathcal{T}_{\sigma}\left(  u\right)  \in Z_{n,m}\ $ as
piecewise linear paths of from $(m+1,1)$ to $(1,n+1)$ in the integer lattice
$\mathbb{N}^{2}$ with $m+n$ horizontal and vertical directed components. The
arrow $(i+1,n+1)\rightarrow(i,n+1)$ represents $A_{i},$ while the arrow
$(m+1,j)\rightarrow(m+1,j+1)$ represents $B_{j}.$ Consequently, $u$ is
represented by the path $(m+1,1)\rightarrow\cdots\rightarrow
(m+1,n+1)\rightarrow\cdots\rightarrow(1,n+1)$, and $z$ by some other path. In
general, if the path $\left(  r+1,s-1\right)  \rightarrow\left(  r+1,s\right)
\linebreak \rightarrow\left(  r,s\right)  $ represents the edge pair $\left(
A_{k}^{\prime},B_{l}^{\prime}\right)  $ in $z,$ the path $\left(
r+1,s-1\right)  \rightarrow\left(  r,s-1\right)  \rightarrow\left(
r,s\right)  $ represents its transposition $\left(  B_{l}^{\prime\prime}%
,A_{k}^{\prime\prime}\right)  $ in $\mathcal{T}\left(  z\right)  $ (see Figure 13).

\hspace*{-0.7in}%
\unitlength=0.055in\special{em:linewidth 0.4pt}\linethickness{0.4pt}\begin{picture}(80.67,38)
\put(40,10){\line(1,0){40}}
\put(40,10){\line(0,1){25}}
\put(40,10){\circle*{1}}
\put(70,10){\circle*{1}}
\put(60,10){\circle*{1}}
\put(50,10){\circle*{1}}
\put(40,20){\circle*{1}}
\put(40,30){\circle*{1}}
\put(50,20){\circle*{1}}
\put(60,20){\circle*{1}}
\put(70,20){\circle*{1}}
\put(50,30){\circle*{1}}
\put(60,30){\circle*{1}}
\put(70,30){\circle*{1}}
\put(70,10){\vector(0,1){10}}
\put(70,20){\vector(-1,0){10}}
\put(50,30){\vector(-1,0){10}}
\put(55,17){\makebox(0,0)[cc]{$A^{\prime}_2$}}
\put(65,17){\makebox(0,0)[cc]{$A^{\prime}_3$}}
\put(53,25){\makebox(0,0)[cc]{$B^{\prime\prime}_2$}}
\put(63,25){\makebox(0,0)[cc]{$B^{\prime}_2$}}
\put(73,25){\makebox(0,0)[cc]{$B_2$}}
\put(73,15){\makebox(0,0)[cc]{$B_1$}}
\put(40,6){\makebox(0,0)[cc]{$1$}}
\put(50,6){\makebox(0,0)[cc]{$2$}}
\put(60,6){\makebox(0,0)[cc]{$3$}}
\put(70,6){\makebox(0,0)[cc]{$4$}}
\put(36.33,10.0){\makebox(0,0)[cc]{$1$}}
\put(36.33,20){\makebox(0,0)[cc]{$2$}}
\put(36.33,30){\makebox(0,0)[cc]{$3$}}
\put(60,20){\vector(0,1){10.33}}
\put(60,30.00){\vector(-1,0){9.33}}
\put(65.00,33.){\makebox(0,0)[cc]{$A_3$}}
\put(55,33.){\makebox(0,0)[cc]{$A_2$}}
\put(45,33.){\makebox(0,0)[cc]{$A_1$}}
\put(70,30.00){\vector(-1,0){10}}
\put(70,20){\vector(0,1){10}}
\put(70,10){\line(0,1){10}}
\put(60,20){\vector(-1,0){10}}
\put(50,20){\vector(0,1){10}}
\end{picture}\vspace{-0.1in}

\begin{center}
Figure 13. $A_{1}A_{2}A_{3}B_{2}B_{1}<A_{1}A_{2}B_{2}^{\prime}A_{3}^{\prime
}B_{1}<A_{1}B_{2}^{\prime\prime}A_{2}^{\prime}A_{3}^{\prime}B_{1}.$%
\vspace*{0.1in}
\end{center}

The poset $\mathcal{KK}$ is a quotient of $\mathcal{PP}$, which we now
describe. Recall Tonks' projection $\vartheta_{0}:P_{m}\rightarrow K_{m+1}$
\cite{Tonks}: If $a$ and $b$ are faces of $P_{m},$ then $\vartheta_{0}\left(
a\right)  =\vartheta_{0}\left(  b\right)  $ if and only if corresponding PLTs
are isomorphic as planar rooted trees (forgetting levels). Define $a\sim b$ if
$\vartheta_{0}\left(  a\right)  =\vartheta_{0}\left(  b\right)  .$ Then
$\mathcal{\tilde{V}}_{m+1}=\vartheta_{0}\left(  \mathcal{V}_{m}\right)  $ is
the set of vertices of $K_{m+1}$. For example, $3|1|2=1|3|2\in\mathcal{\tilde
{V}}_{4},$ since $3|1|2$ and $1|3|2$ are end points of the degenerate edge
$13|2\subset P_{3},$ and in terms matrix sequences we have $\left[
\curlywedge\right]  \left[  \curlywedge\text{ }\mathbf{1}\right]  \left[
\mathbf{1}\text{ }\mathbf{1}\text{ }\curlywedge\right]  =\left[
\curlywedge\right]  \left[  \mathbf{1}\text{ }\curlywedge\right]  \left[
\curlywedge\text{ }\mathbf{1}\text{ }\mathbf{1}\right]  $ (and dually $\left[
\mathbf{1}\text{ }\mathbf{1}\text{ }\curlyvee\right]  ^{T}\left[
\curlyvee\text{ }\mathbf{1}\right]  ^{T}\left[  \curlyvee\right]  =\left[
\curlyvee\text{ }\mathbf{1}\text{ }\mathbf{1}\right]  ^{T}\left[
\mathbf{1}\text{ }\curlyvee\right]  ^{T}\left[  \curlyvee\right]  $ )$.$ Of
course, $\mathcal{V}_{1}=\mathcal{\tilde{V}}_{2}=\left[  \curlywedge\right]  $
and $\mathcal{V}_{2}=\mathcal{\tilde{V}}_{3}=\left\{  \left[  \curlywedge
\right]  \left[  \curlywedge\text{ }\mathbf{1}\right]  ,\left[  \curlywedge
\right]  \left[  \mathbf{1}\text{ }\curlywedge\right]  \right\}  .$

For matrix sequences in $X_{m}^{n+1},\ $define $x_{1}^{\prime}\cdots
x_{m}^{\prime}\sim_{X}x_{1}\cdots x_{m}$ if the trees produced by $\gamma
$-products of $i^{th}$ rows are equivalent for each $i$. Dually, for matrix
sequences in $Y_{n}^{m+1},\ $define $y_{n}^{\prime}\cdots y_{1}^{\prime}%
\sim_{Y}y_{n}\cdots y_{1}$ if the trees produced by $\gamma$-products of
$j^{th}$ columns are equivalent for each $j$. Define $a\times b\sim c\times d$
in $X_{m}^{n+1}\times Y_{n}^{m+1}$ if $a\sim_{X}c$ and $b\sim_{Y}d.$ Finally,
for $u_{1}\leq u_{2}\in X_{m}^{n+1}\times Y_{n}^{m+1}$ and $z_{1}%
=\mathcal{T}_{\sigma}(u_{1})\leq z_{2}=\mathcal{T}_{\sigma}(u_{2}),$ define
$z_{1}\sim z_{2}$ if $u_{1}\sim u_{2}.$ Then
\[
{\mathcal{K}}{\mathcal{K}}_{n+1,m+1}=\left.  \mathcal{PP}_{n,m}\right/
\sim\text{ };
\]
and $\vartheta\vartheta:\mathcal{PP}_{n,m}\rightarrow{\mathcal{K}}%
{\mathcal{K}}_{n+1,m+1}$ denotes the projection.

\section{The Combinatorial Join of Permutahedra}

The combinatorial join of permutahedra, which resembles the standard join of
spaces, plays an important role in our construction of the biassociahedra to
follow. The \emph{combinatorial join }$P_{m}\ast_{c}P_{n}$ of permutahedra
$P_{m}$ and $P_{n}$ is the permutahedron $P_{m+n}$ constructed as follows:
Given faces $A_{1}|\cdots|A_{k}\subseteq P_{m}$ and $B_{1}|\cdots
|B_{l}\subseteq P_{n},$ let $s$ be an integer such that $\max\left\{
k,l\right\}  \leq s\leq k+l,$ and let $\left(  \mathbf{i};\mathbf{j}\right)
=(i_{1}<\cdots<i_{k};j_{1}<\cdots<j_{l}),$ where $\mathbf{i}\cup
\mathbf{j}=\underline{s}.$ Obtain $A_{1}^{\prime}|\cdots|A_{s}^{\prime}$ and
$B_{1}^{\prime}|\cdots|B_{s}^{\prime}$ by setting $A_{i_{r}}^{\prime}=A_{r},$
$B_{j_{t}}^{\prime}=B_{t},$ and $A_{i}^{\prime}=B_{j}^{\prime}=\varnothing$
otherwise. Note that $(A_{r}^{\prime},B_{r}^{\prime})\neq\left(
\varnothing,\varnothing\right)  $ for all $r.$ Given a set $B=\left\{
b_{1},\ldots,b_{k}\right\}  \subset\mathbb{N}$ and $m\in\mathbb{N}$, define
$B+m=\left\{  b_{1}+m,\ldots,b_{k}+m\right\}  $ and consider the codimension
$s-1$ face $A_{1}|\cdots|A_{k}\ast_{\left(  \mathbf{i};\mathbf{j}\right)
}B_{1}|\cdots|B_{l}=A_{1}^{\prime}\cup(B_{1}^{\prime}+m)\mid\cdots\mid
A_{s}^{\prime}\cup(B_{s}^{\prime}+m)\subset P_{m+n}.$ When $s=m+n,$ each pair
of vertices $A_{1}|\cdots|A_{m}\times B_{1}|\cdots|B_{n}\subset P_{m}\times
P_{n}$ generates $\tbinom{m+n}{m}$ vertices $A_{1}|\cdots|A_{m}\ast_{\left(
\mathbf{i};\mathbf{j}\right)  }B_{1}|\cdots|B_{n}$ of $P_{m+n}$ as $\left(
\mathbf{i};\mathbf{j}\right)  $ ranges over all $\left(  m,n\right)
$-shuffles of $\left(  A_{1},\ldots,A_{m};B_{1}+m,\ldots,B_{n}+m\right)  $.
Define%
\[
P_{m}\ast_{c}P_{n}=\bigcup\limits_{\substack{A_{1}|\cdots|A_{k}\times
B_{1}|\cdots|B_{l}\subseteq P_{m}\times P_{m} \\\mathbf{i}\cup\mathbf{j}%
=\underline{s};\text{ }\max\left\{  k,l\right\}  \leq s\leq k+l}}A_{1}%
|\cdots|A_{k}\ast_{\left(  \mathbf{i};\mathbf{j}\right)  }B_{1}|\cdots|B_{l}.
\]
Thus, given $m,n\geq1$ and a cell $e\subseteq P_{m+n},$ there is a unique
decomposition $e=A_{1}|\cdots|A_{k}\ast_{\left(  \mathbf{i};\mathbf{j}\right)
}B_{1}|\cdots|B_{l}$ with $A_{1}|\cdots|A_{k}\subset P_{m}$ and $B_{1}%
|\cdots|B_{l}\subset P_{n}.$

\begin{example}
Setting $s=2$ produces the $14$ codimension $1$ faces of $P_{2}\ast_{c}%
P_{2}=P_{4}:$
\[%
\begin{tabular}
[c]{|c|c|c|c|c|c|c|c|c|}\cline{1-4}\cline{6-9}%
$\left(  \mathbf{i};\mathbf{j}\right)  $ & $A$ & $B$ & $A\ast_{\left(
\mathbf{i};\mathbf{j}\right)  }B$ & $%
\begin{array}
[c]{c}%
\text{ }\\
\text{ }%
\end{array}
$ & $\left(  \mathbf{i};\mathbf{j}\right)  $ & $A$ & $B$ & $A\ast_{\left(
\mathbf{i};\mathbf{j}\right)  }B$\\\cline{1-4}\cline{6-9}%
$\left(  1,2;1,2\right)  $ & $1|2$ & $1|2$ & $13|24$ &  & $\left(
2;1,2\right)  $ & $12$ & $1|2$ & $3|124$\\\cline{1-4}\cline{6-9}
& $1|2$ & $2|1$ & $14|23$ &  &  & $12$ & $2|1$ & $4|123$\\\cline{1-4}%
\cline{6-9}
& $2|1$ & $1|2$ & $23|14$ &  & $\left(  1,2;1\right)  $ & $1|2$ & $12$ &
$134|2 $\\\cline{1-4}\cline{6-9}
& $2|1$ & $2|1$ & $24|13$ &  &  & $2|1$ & $12$ & $234|1$\\\cline{1-4}%
\cline{6-9}%
$\left(  1;1,2\right)  $ & $12$ & $1|2$ & $123|4$ &  & $\left(  1,2;2\right)
$ & $1|2$ & $12$ & $1|234$\\\cline{1-4}\cline{6-9}
& $12$ & $2|1$ & $124|3$ &  &  & $2|1$ & $12$ & $2|134$\\\cline{1-4}%
\cline{6-9}%
$\left(  1;2\right)  $ & $12$ & $12$ & $12|34$ &  & $\left(  2;1\right)  $ &
$12$ & $12$ & $34|12$\\\cline{1-4}\cline{6-9}%
\end{tabular}
\]

\end{example}

Fraction products $a/b$ reappear as combinatorial joins $a\ast_{c}b$ in Step 2
of the construction that follows in the next section.

\section{Constructions of $PP$ and $KK$}

We conclude the paper with constructions of the geometric realizations
$PP=\left\vert \mathcal{PP}\right\vert $ and $KK=\left\vert \mathcal{KK}%
\right\vert .$ While the edges of $PP$ and $KK$ realize the edges of
$\mathcal{PP}$ and $\mathcal{KK}$, it is difficult to imagine their higher
dimensional faces. Fortunately, $PP_{n,m}$ is a subdivision of the
permutahedron $P_{m+n}$, which is a subdivision of $I^{n+m-1}$. Thus the
higher dimensional combinatorics of $PP_{n,m}$ are determined by the
orientation on the faces of $I^{m+n-1}$. 

Our construction of $PP_{n,m}$ has two steps: (1) Perform an
\textquotedblleft$\left(  m,n\right)  $-subdivision\textquotedblright\ of the
codimension 1 cell $\underline{m}\,|\left(  \underline{n}+m\right)  \subset
P_{m+n}$ and (2) use the $\left(  m,n\right)  $-subdivision to subdivide
certain other cells of $P_{m+n}$. We emphasize that $\Delta_P$ is used only 
in step (1) and only in terms of its geometrical definition. Thus the 
non-coassociativity and non-cocommutativity of $\Delta_P$ are not in play here 
(see also Remark \ref{independent} below).  We begin with some preliminaries.

\subsection{Matrices with constant rows or columns}

Given a set $Q$ of matrix sequences, let%
\[
\operatorname{con}Q=\left\{  C_{1}\cdots C_{s}\in Q\text{ }|\text{ }%
C_{k}\text{ has constant rows or constant columns}\right\}  .
\]
Note that if $A_{1}\cdots A_{m}\in\operatorname{con}X_{m}^{n+1},$ each $A_{i}$
has constant columns; dually, if $B_{n}\cdots B_{1}\in\operatorname{con}%
Y_{n}^{m+1},$ each $B_{j}$ has constant rows. Consequently, the inclusion of
posets $\kappa:X_{m}^{n+1}\times Y_{n}^{m+1}\hookrightarrow\mathcal{V}%
_{m}^{\times n+1}\times\mathcal{V}_{n}^{\times m+1}$ given in (\ref{chia})
restricts to an order-preserving bijection%
\[
\left(  \chi\circ\hat{\ell}\right)  \times{\check{\ell}}:{\operatorname{con}%
}\left(  X_{m}^{n+1}\times Y_{n}^{m+1}\right)  \leftrightarrow
{\operatorname{con}}X_{m}^{n+1}\times{\operatorname{con}}Y_{n}^{m+1}%
\hspace*{1in}%
\]%
\begin{equation}
\hspace*{1in}\leftrightarrow\Delta\left(  \mathcal{V}_{m}^{\times n+1}\right)
\times\Delta\left(  \mathcal{V}_{n}^{\times m+1}\right)  \leftrightarrow
\mathcal{V}_{m}\times\mathcal{V}_{n},\label{bijection}%
\end{equation}
where $\Delta\left(  \mathcal{V}_{m}^{\times n+1}\right)  \leftrightarrow
\mathcal{V}_{m}$ is given by the embedding $\mathcal{V}_{m}\hookrightarrow
\mathcal{V}_{m}^{\times n+1}$ along the diagonal subposet $\Delta\left(
\mathcal{V}_{m}^{\times n+1}\right)  =\left\{  \left(  v,\ldots,v\right)
\text{ }|\text{ }v\in\mathcal{V}_{m}\right\}  .$ Thus elements of
$\mathcal{V}_{m}\times\mathcal{V}_{n}$ may be represented as matrix sequences
in ${\operatorname{con}}\left(  X_{m}^{n+1}\times Y_{n}^{m+1}\right)  $.

Note that $\left(  i,j\right)  $-transpositions preserve constant rows and
columns, i.e., $u\in$ 
\linebreak ${\operatorname{con}}\mathcal{PP}_{n,m}$ if and only if
$\mathcal{T}_{ij}^{k}\left(  u\right)  \in{\operatorname{con}}${$\mathcal{PP}%
$}$_{n.m}$. And furthermore, if $u=A_{1}\cdots A_{m}$ $B_{n}\cdots B_{1}%
\in{\operatorname{con}}\left(  X_{m}^{n+1}\times Y_{n}^{m+1}\right)  ,$ and
$\sigma$ is an $\left(  m,n\right)  $-shuffle, ${\mathcal{T}}_{\sigma}\left(
u\right)  $ is defined since each $A_{i}$ has a constant column of
$\curlywedge$'s and each $B_{j}$ has a constant row of $\curlyvee$'s. Thus%
\begin{equation}
{\operatorname{con}}\mathcal{PP}_{n,m}=\bigcup_{\left(  m,n\right)
\text{-shuffles }\sigma}{\mathcal{T}}_{\sigma}\left(  {\operatorname{con}%
}\left(  X_{m}^{n+1}\times Y_{n}^{m+1}\right)  \right)  .\label{conPP}%
\end{equation}

The order-preserving bijection%
\[
{\operatorname{con}}\mathcal{PP}_{1,2}={\mathcal{PP}}_{1,2}\setminus\left\{
u_{2}\right\}  \leftrightarrow\mathcal{V}_{3}%
\]
discussed in Example \ref{Ex2} illustrates the following remarkable fact:

\begin{proposition}
\label{Prop1}The bijection $\left(  \chi\circ\hat{\ell}\right)  \times
{\check{\ell}}:{\operatorname*{con}}\left(  X_{m}^{n+1}\times Y_{n}%
^{m+1}\right)  \rightarrow\mathcal{V}_{m}\times\mathcal{V}_{n}$ extends to a
canonical order-preserving bijection
\[
\kappa_{_{\#}}:{\operatorname*{con}}\mathcal{PP}_{n,m}\rightarrow
\mathcal{V}_{m+n}.
\]
Thus $\left\vert \kappa_{_{\#}}\right\vert :\left\vert {\operatorname*{con}%
}\mathcal{PP}_{n,m}\right\vert \overset{\approx}{\longrightarrow}P_{m+n}.$
\end{proposition}

\begin{proof}
There is the order-preserving bijection $\left(  \chi\circ\hat{\ell}\right)
\times{\check{\ell}}:{\operatorname{con}}\left(  X_{m}^{n+1}\times Y_{n}%
^{m+1}\right)  \leftrightarrow S_{m}\times S_{n}$ via the identification
$\mathcal{V}_{m}\times\mathcal{V}_{n}\leftrightarrow S_{m}\times S_{n}.$ Thus
\[
{\operatorname{con}}\mathcal{PP}_{n,m}\leftrightarrow\{\sigma\circ\left(
\sigma_{m}\times\sigma_{n}\right)  \mid\sigma\text{ is an }\left(  m,n\right)
\text{-shuffle; }\sigma_{m}\times\sigma_{n}\in S_{m}\times S_{n}\}
\]
by formula (\ref{conPP}). But each permutation in $S_{m+n}$ factors as
$\sigma\circ\left(  \sigma_{m}\times\sigma_{n}\right)  $ for some $\left(
m,n\right)  $-shuffle $\sigma$ and some $\sigma_{m}\times\sigma_{n}\in
S_{m}\times S_{n}.$ Therefore $\left(  \chi\circ\hat{\ell}\right)
\times{\check{\ell}}$ extends to $\kappa_{_{\#}}:{\operatorname{con}%
}\mathcal{PP}_{n,m}\leftrightarrow S_{m+n}\leftrightarrow\mathcal{V}_{m+n}.$
\end{proof}

\begin{corollary}
\label{cor1}For all $m,n\geq1,$there is the commutative diagram
\[%
\begin{array}
[c]{ccc}%
\left\vert {{\operatorname{con}}}\left(  X_{m}^{n+1}\times Y_{n}^{m+1}\right)
\right\vert  & \hookrightarrow & \left\vert \operatorname{con}\mathcal{PP}%
_{n,m}\right\vert \vspace*{0.1in}\\
\!|\kappa|\text{ }\downarrow\,\approx\text{ \ \ \ } &  & \ \approx
\,\downarrow\text{ }|\kappa_{_{\#}}|\vspace*{0.1in}\\
\underline{m}|\left(  \underline{n}+m\right)  & \hookrightarrow & P_{m+n}.
\end{array}
\]

\end{corollary}

Note that if $\chi_{_{\#}}:\mathcal{V}_{m+n}\rightarrow\mathcal{V}_{m+n}$ were
induced by the composition $\mathcal{V}_{m}\times\mathcal{V}_{n}%
\overset{\chi\times\mathbf{1}}{\longrightarrow}\mathcal{V}_{m}\times
\mathcal{V}_{n}\hookrightarrow\mathcal{V}_{m+n}$ in the same manner as
$\kappa_{_{\#}},$ then $\chi_{_{\#}}$ and $\chi$ would differ on
$\mathcal{V}_{m+n}$ even with $n=1.$

\begin{example}
\label{Ex3}Continuing Example \ref{Ex2}, the identification
${\operatorname{con}}\mathcal{PP}_{1,1}\leftrightarrow\mathcal{V}_{2}$ is
given by%
\[%
\begin{array}
[c]{ccc}%
\left[
\begin{array}
[c]{c}%
\curlywedge\\
\curlywedge
\end{array}
\right]  [\curlyvee\curlyvee] & \overset{\mathcal{T}}{\longrightarrow} &
\left[  \curlyvee\right]  \left[  \curlywedge\right] \\
&  & \\
\updownarrow &  & \updownarrow\\
&  & \\
1|\left(  1+1\right)  & \underset{\left(  1,1\right)  \text{-shuffle}%
}{\longrightarrow} & 2|1.
\end{array}
\]
The identification ${\operatorname{con}}\mathcal{PP}_{1,2}\leftrightarrow
\mathcal{V}_{3}:$%
\[
\left[
\begin{array}
[c]{c}%
\curlywedge\\
\curlywedge
\end{array}
\right]  \left[
\begin{array}
[c]{c}%
\curlywedge\text{ }\mathbf{1}\\
\curlywedge\text{ }\mathbf{1}%
\end{array}
\right]  \leftrightarrow1|2\in S_{2}\text{ and }\left[  \curlyvee
\curlyvee\curlyvee\right]  \leftrightarrow1\in S_{1}%
\]
so that
\[
u_{1}=\left[
\begin{array}
[c]{c}%
\curlywedge\\
\curlywedge
\end{array}
\right]  \left[
\begin{array}
[c]{c}%
\curlywedge\text{ }\mathbf{1}\\
\curlywedge\text{ }\mathbf{1}%
\end{array}
\right]  \left[  \curlyvee\curlyvee\curlyvee\right]  \leftrightarrow
1|2|\left(  1+2\right)  =1|2|3\in S_{2}\times S_{1}.
\]
Similarly,
\[
u_{3}=\left[
\begin{array}
[c]{c}%
\curlywedge\\
\curlywedge
\end{array}
\right]  \left[
\begin{array}
[c]{c}%
\mathbf{1}\text{ }\curlywedge\\
\mathbf{1}\text{ }\curlywedge
\end{array}
\right]  \left[  \curlyvee\curlyvee\curlyvee\right]  \leftrightarrow
2|1|\left(  1+2\right)  =2|1|3\in S_{2}\times S_{1}%
\]
and we have%
\[%
\begin{array}
[c]{c}%
\begin{array}
[c]{cc}%
\left[
\begin{array}
[c]{c}%
\curlywedge\\
\curlywedge
\end{array}
\right]  \left[
\begin{array}
[c]{c}%
\curlywedge\text{ }\mathbf{1}\\
\curlywedge\text{ }\mathbf{1}%
\end{array}
\right]  \left[  \curlyvee\curlyvee\curlyvee\right]  & \text{ }%
\overset{\mathcal{T}_{\sigma_{1}}}{\longrightarrow}%
\end{array}
\vspace*{0.1in}\\%
\begin{array}
[c]{c}%
\updownarrow\vspace*{0.1in}\\
1|2|3\vspace*{0.45in}%
\end{array}
\end{array}%
\begin{array}
[c]{c}%
\begin{array}
[c]{ccc}%
\left[
\begin{array}
[c]{c}%
\curlywedge\\
\curlywedge
\end{array}
\right]  [\curlyvee\curlyvee]\left[  \curlywedge\text{ }\mathbf{1}\right]  &
\overset{\mathcal{T}_{\sigma_{2}}}{\longrightarrow} & [\curlyvee
][\curlywedge]\left[  \curlywedge\text{ }\mathbf{1}\right]
\end{array}
\vspace*{0.1in}\\
\underset{}{\underbrace{%
\begin{array}
[c]{ccc}%
\updownarrow & \text{ \ \ \ \ \ \ \ \ \ \ \ \ \ \ \ \ \ \ } & \updownarrow
\vspace*{0.1in}\\
1|3|2 &  & 3|1|2\vspace*{0.1in}%
\end{array}
}}\\
\left(  2,1\right)  \text{-shuffles of }1|2|3
\end{array}
\vspace*{0.1in}%
\]%
\[%
\begin{array}
[c]{c}%
\begin{array}
[c]{cc}%
\left[
\begin{array}
[c]{c}%
\curlywedge\\
\curlywedge
\end{array}
\right]  \left[
\begin{array}
[c]{c}%
\mathbf{1}\text{ }\curlywedge\\
\mathbf{1}\text{ }\curlywedge
\end{array}
\right]  \left[  \curlyvee\curlyvee\curlyvee\right]  & \text{ }%
\overset{\mathcal{T}_{\sigma_{1}}}{\longrightarrow}%
\end{array}
\vspace*{0.1in}\\%
\begin{array}
[c]{c}%
\updownarrow\vspace*{0.1in}\\
2|1|3\vspace*{0.45in}%
\end{array}
\end{array}%
\begin{array}
[c]{c}%
\begin{array}
[c]{ccc}%
\left[
\begin{array}
[c]{c}%
\curlywedge\\
\curlywedge
\end{array}
\right]  [\curlyvee\curlyvee]\left[  \mathbf{1}\text{ }\curlywedge\right]  &
\overset{\mathcal{T}_{\sigma_{2}}}{\longrightarrow} & [\curlyvee
][\curlywedge]\left[  \mathbf{1}\text{ }\curlywedge\right]
\end{array}
\vspace*{0.1in}\\
\underset{}{\underbrace{%
\begin{array}
[c]{ccc}%
\updownarrow & \text{ \ \ \ \ \ \ \ \ \ \ \ \ \ \ \ \ \ \ } & \updownarrow
\vspace*{0.1in}\\
2|3|1 &  & 3|2|1\vspace*{0.1in}%
\end{array}
}}\\
\left(  2,1\right)  \text{-shuffles of }2|1|3.
\end{array}
\]

\end{example}

The projection $\vartheta:\operatorname{con}\mathcal{PP}_{n,m}\rightarrow
\operatorname{con}\mathcal{PP}_{n,m}\diagup\sim$ has the following simple
geometrical interpretation: An element of $\operatorname{con}\left(
X_{m}^{n+1}\times Y_{n}^{m+1}\right)  $ is represented as a fraction with
multiple copies of the \emph{same }leveled binary tree in the numerator and
likewise in the denominator. Two such elements are equivalent if and only if
the trees in their numerators or denominators (possibly both) are isomorphic
as PRTs. So equivalence in $\operatorname{con}\left(  X_{m}^{n+1}\times
Y_{n}^{m+1}\right)  $ amounts to forgetting levels as in Tonks' projection.
The poset structure then propagates this equivalence to general elements of
$\operatorname{con}\mathcal{PP}_{n,m}$.

\begin{remark}
\label{independent}
Our constructions are independent of the various choices involved here. If
$\tilde{\Delta}_{P}^{(k)}$ iterates $\Delta_{P}$ on factors other than the
those on the extreme left, let $\widetilde{X}_{m}^{n+1}\times\widetilde{Y}%
_{n}^{m+1}$ be the poset defined in terms of $\tilde{\Delta}_{P}^{(k)}$ and
let $\widetilde{{\mathcal{PP}}}_{n,m}$ be the poset produced by our
construction. Then there is a canonical bijection ${{\mathcal{PP}}}%
_{n,m}\leftrightarrow\widetilde{{\mathcal{PP}}}_{n,m}$ and the corresponding
geometric realizations are canonically homeomorphic. When $\Delta_{P}$ acts on
the extreme right, for example, a (combinatorial) isomorphism $|\mathcal{PP}%
_{n,m}|\cong|\widetilde{\mathcal{PP}}_{n,m}|$ is evident pictorially: The
picture of $|\mathcal{PP}_{n,m}|$ uses the standard orientation of the
interval $P_{2}$, while the picture of $|\widetilde{\mathcal{PP}}_{n,m}|$ uses
the opposite orientation, but nevertheless, these pictures are identical.
\end{remark}

\subsection{Step 1: The $(m,n)$-subdivision of $\protect\underline{m}%
\,|\left(  \protect\underline{n}+m\right)  $}

The first step in our construction of $PP_{n,m}$ performs an \textquotedblleft%
$(m,n)$-subdivision\textquotedblright\ of the codimension $1$ cell
$\underline{m}\,|\left(  \underline{n}+m\right)  \subset P_{n+m}$. In
Subsection \ref{topological} we applied the left-iterated diagonal $\Delta
_{P}^{\left(  n\right)  }$ to construct the $n$-subdivision $P_{m}^{\left(
n\right)  }$ of $P_{m}$. Since the poset $X_{m}^{n+1}$ is the $0$-skeleton of
$\Delta^{\left(  n\right)  }\left(  P_{m}\right)  ,$ the geometric realization
$\left\vert X_{m}^{n+1}\right\vert =P_{m}^{\left(  n\right)  },$ and dually
$\left\vert Y_{n}^{m+1}\right\vert =P_{n}^{(m)}.$ The cellular subdivision
$\left\vert X_{m}^{n+1}\times Y_{n}^{m+1}\right\vert =P_{m}^{(n)}\times
P_{n}^{(m)}$ of $\underline{m}\,|\left(  \underline{n}+m\right)  =P_{m}\times
P_{n}$ is called the $(m,n)$-\emph{subdivision} of $\underline{m}\,|\left(
\underline{n}+m\right)  ;$ thus each cell in this subdivision has a canonical
Cartesian product decomposition. The \emph{basic subdivision vertices} of
$PP_{n,m}$ are elements of%
\[
\mathcal{BS}_{n,m}=\left(  X_{m}^{n+1}\times Y_{n}^{m+1}\right)
\setminus\mathcal{V}_{m+n}.
\]
Each subdivision cell of $PP_{n,m}$ is a proper subset of some cell of
$P_{m+n}$ and is the geometric realization of its poset of vertices.

\begin{example}
\label{subdivisions}The $1$-subdivision $P_{2}^{\left(  1\right)  }$ consists
of two $1$-cells obtained by subdividing the interval $P_{2}$ at its
midpoint$\ $(see Figure 7). Thus the $\left(  2,1\right)  $-subdivision
$P_{2}^{\left(  1\right)  }\times P_{1}^{\left(  2\right)  }$ of the edge
$12|3\subset P_{3}$ contains one basic subdivision vertex represented by the
midpoint
\[
u_{2}=\left[
\begin{array}
[c]{c}%
\curlywedge\\
\curlywedge
\end{array}
\right]  \left[
\begin{array}
[c]{c}%
\curlywedge\text{ }\mathbf{1}\\
\mathbf{1}\text{ }\curlywedge
\end{array}
\right]  [\curlyvee\curlyvee\curlyvee]\in X_{2}^{2}\times Y_{1}^{3}%
\setminus{\operatorname{con}}\left(  X_{2}^{2}\times Y_{1}^{3}\right)  ,
\]
and two $1$-cells of $PP_{1,2}$. In fact, $PP_{1,2}$ is exactly the heptagon
obtained by subdividing $P_{3}$ in this way. The $2$-subdivision
$P_{2}^{\left(  2\right)  }$ consists of three $1$-cells obtained by
subdividing $P_{2}$ at its midpoint and again at its three-quarter point. Thus
the $(2,2)$-subdivision $P_{2}^{\left(  2\right)  }\times P_{2}^{\left(
2\right)  }$ of the square $12|34\subset P_{4}$ contains twelve basic
subdivision vertices and nine $2$-cells of $PP_{2,2}$ as pictured in Figure
14\vspace{0.1in}.
\begin{center}
\includegraphics[
height=1.1713in,
width=1.1652in
]%
{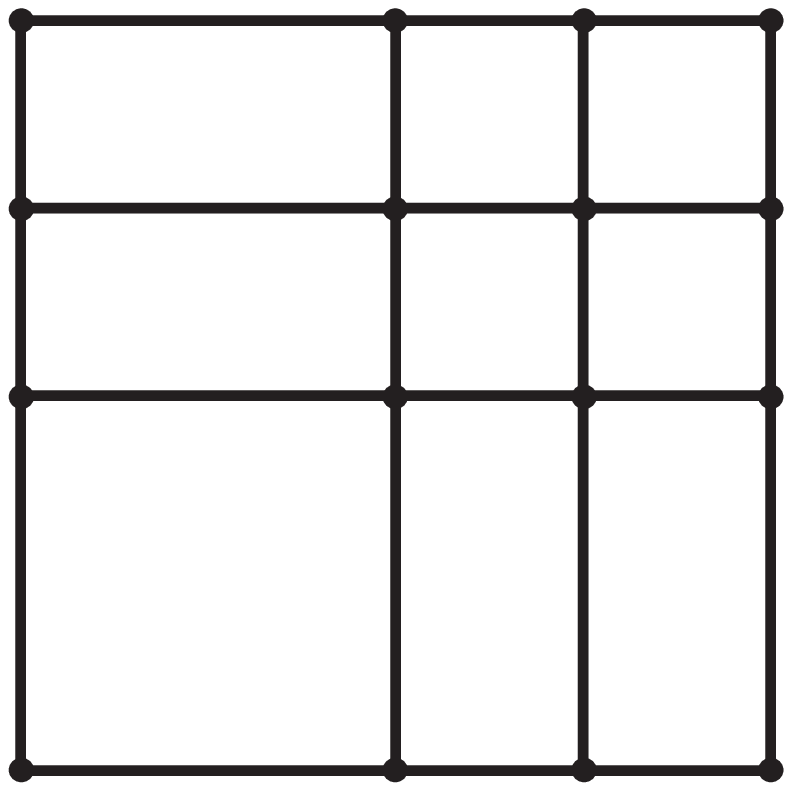}%
\\
Figure 14. The $\left(  2,2\right)  $-subdivision $P_{2}^{\left(  2\right)
}\times P_{2}^{\left(  2\right)  }$.
\end{center}
The $(3,1)$-subdivision $P_{3}^{\left(  1\right)  }\times P_{1}^{\left(
3\right)  }$ of the hexagon $123|4\subset P_{4}$ is identified with the
$1$-subdivision $P_{3}^{\left(  1\right)  }$ and contains eleven basic
subdivision vertices and eight $2$-cells of $PP_{3,1}$ as pictured in Figure
15\vspace{0.1in}.
\begin{center}
\includegraphics[
height=1.5368in,
width=1.7599in
]%
{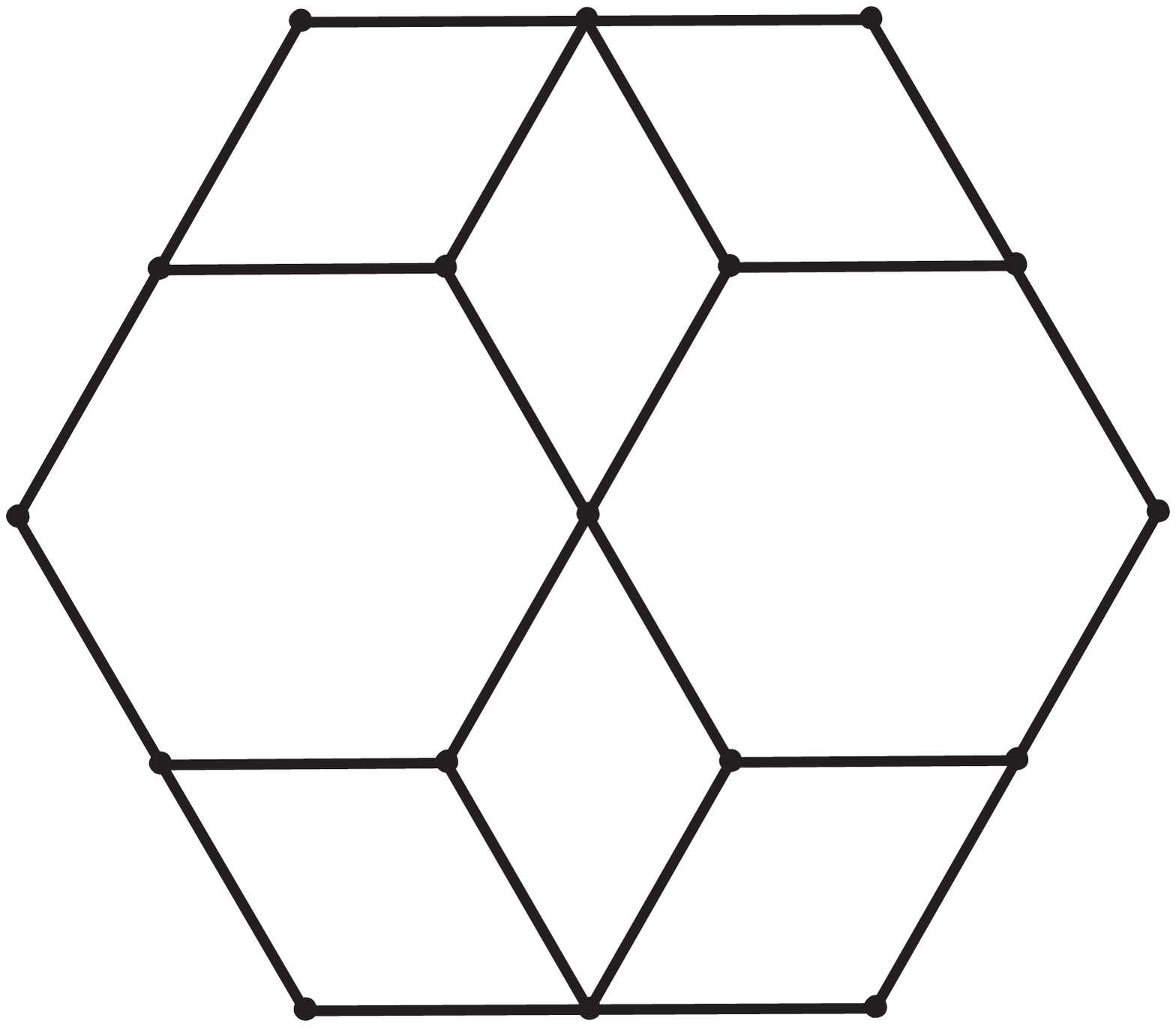}%
\\
Figure 15. The $\left(  3,1\right)  $-subdivision $P_{3}^{\left(  1\right)
}\times P_{1}^{\left(  3\right)  }.$%
\end{center}
The $\left(  3,2\right)  $-subdivision $P_{3}^{\left(  2\right)  }\times
P_{2}^{\left(  3\right)  }$ of the cylinder $123|45\subset P_{5},$ obtained
from $P_{3}^{\left(  2\right)  }\times I$ by subdividing along the horizontal
cross-sections $P_{3}^{\left(  2\right)  }\times\frac{1}{2},$ $P_{3}^{\left(
2\right)  }\times\frac{3}{4},$ and $P_{3}^{\left(  2\right)  }\times\frac
{7}{8},$ contains $140$ basic subdivision vertices and eighty-four $3$-cells
of $PP_{3,2}$ ($P_{3}^{\left(  2\right)  }$ is pictured in Figure 11).
\end{example}

\subsection{Step 2: Subdividing cells of $P_{m+n}\setminus
\protect\underline{m}|\left(  \protect\underline{n}+m\right)  $}

Recall that elements of $Z_{n,m}$ arise from the non-trivial action of
$\mathcal{T}_{\sigma}$ on $X_{m}^{n+1}\times Y_{n}^{m+1}$. When $\sigma$
ranges over all $\left(  m,n\right)  $-shuffles (including the identity), we
obtain the poset%
\[
\mathcal{S}_{n,m}=\bigcup_{\left(  m,n\right)  \text{-shuffles }\sigma
}\mathcal{T}_{\sigma}\left(  \mathcal{BS}_{n,m}\right)
\]
of \emph{subdivision vertices }of $\mathcal{PP}_{n,m}.$ Thus as sets,
$\mathcal{PP}_{n,m}=\operatorname{con}\mathcal{PP}_{n,m}\sqcup\mathcal{S}%
_{n,m}.$

The second step of our construction detects those cells of $P_{m+n}%
\setminus\underline{m}|\left(  \underline{n}+m\right)  $ that contain
subdivision vertices. We refer to such cells as Type I cells; all other cells
have Type II. We use the poset structure of subdivision vertices to subdivide
Type I cells, and having done so, our construction of $PP_{n,m}$ will be complete.

To begin, let us characterize those Type I cells of minimal dimension that
contain \emph{non-basic} subdivision vertices. If $e$ is a cell of some
polytope, denote the set of vertices of $e$ by $\mathcal{V}_{e}.$ Given an
$(m,n)$-shuffle $\sigma$ and a cell $e\subset\underline{m}\,|\left(
\underline{n}+m\right)  ,$ let $\mathcal{T}\left(  \sigma,e\right)  $ denote
the cell of $P_{m+n}$ of minimal dimension such that $\mathcal{T}_{\sigma
}\left(  \mathcal{V}_{e}\right)  \subseteq\mathcal{V}_{\mathcal{T}\left(
\sigma,e\right)  }.$ This defines a map%
\[
\mathcal{T}:\left\{  \left(  m,n\right)  \text{-shuffles}\right\}
\times\{\text{partitions of}\ \underline{m}|\left(  \underline{n}+m\right)
\}\rightarrow\{\text{partitions of}\ \underline{m+n}\},
\]
which extends the map $\left(  \sigma,\sigma_{m}\times\sigma_{n}\right)
\mapsto\sigma\circ\left(  \sigma_{m}\times\sigma_{n}\right)  $ in the proof of
Proposition \ref{Prop1}. To define $\mathcal{T}$ at a particular shuffle
$\sigma$ and partition $e=A_{1}|\cdots|A_{k}|B_{1}|\cdots|B_{l}\subseteq
\underline{m}|\left(  \underline{n}+m\right)  ,$ remove all block delimiters
of $e$ and think of $e$ as a permutation of $\underline{m+n} $ in which
$A_{i}$ and $B_{j}$ are contiguous subsequences. Consider the set $\left\{
D_{1},\ldots,D_{r}\right\}  $ of all contiguous subsequences $\sigma(A_{i})$
and $\sigma(B_{j})$ of $\sigma\left(  e\right)  $ that preserve the contiguity
of the $B_{j^{\prime}}$'s and $A_{i^{\prime}}$'s, respectively, then reinsert
block delimiters so that
\[
\mathcal{T}(\sigma,e)=C_{1}|D_{i_{1}}|\cdots|C_{r}|D_{i_{r}}|C_{r+1}.
\]
Since each cell of $P_{m+n}$ can be expressed uniquely as a component of the
combinatorial join $P_{m}\ast_{c}P_{n}$, we have
\begin{equation}
\mathcal{T}(\sigma,e)=E\ast_{\left(  \mathbf{i}{;}\mathbf{j}\right)  }%
F=E_{1}^{\prime}\cup(F_{1}^{\prime}+m)\mid\cdots\mid E_{s}^{\prime}\cup
(F_{s}^{\prime}+m),\label{specialjoin}%
\end{equation}
where $E_{i}$ and $F_{j}$ are the unions of consecutive blocks $A_{i^{\prime}%
}|\cdots|A_{i^{\prime}+i^{\prime\prime}}$ and $B_{j^{\prime}}|\cdots$
$|B_{j^{\prime}+j^{\prime\prime}}$ of $e,$ respectively. Thus $\sigma$ acts on
the blocks of $e$ as a $(k,l)$-shuffle if and only if $C_{i}=\varnothing$ for
all $i$ if and only if $\mathcal{T}(\sigma,e)=A_{1}|\cdots|A_{k}\ast_{\left(
\mathbf{i}{;}\mathbf{j}\right)  }B_{1}|\cdots|B_{l}$ for some $(k,l)$%
-unshuffle $(\mathbf{i}{;}\mathbf{j})=(i_{1}<\cdots<i_{k}\,;\,j_{1}%
<\cdots<j_{l})$ of $\underline{k+l}.$ Clearly, a cell $a\subset P_{n+m}$
contains a non-basic subdivision vertex $\mathcal{T}_{\sigma}\left(  u\right)
$ if and only if $a=\mathcal{T}(\sigma,e) $ for some cell $e\subset
\underline{m}\,|\left(  \underline{n}+m\right)  $ containing a basic
subdivision vertex $u$ on which $\mathcal{T}_{\sigma}$ acts non-trivially. In
fact, $a$ contains at most one non-basic\emph{\ }subdivision vertex when
$m+n\leq4$.

The following proposition incorporates the property of $\mathcal{T}$ described
in Remark \ref{edge} and will be applied in our subsequent examination of the
poset structure of $\mathcal{PP}_{n,m}$.

\begin{proposition}
\label{Prop5}If a cell $e\subset\underline{m}\,|\left(  \underline{n}%
+m\right)  $ contains a subdivision cell $a\subset|X_{m}^{n+1}\times
Y_{n}^{m+1}|$ and $\mathcal{T}_{\sigma}(\mathcal{V}_{a})\subset\left(
\mathcal{S}_{n,m}\cap\mathcal{T}(\sigma,e)\right)  \cup\mathcal{V}%
_{\mathcal{T}(\sigma,e)},$ then $|\mathcal{T}_{\sigma}(\mathcal{V}_{a})|$ is a
subdivision cell of $\mathcal{T}(\sigma,e)$ (combinatorially) isomorphic to
$a;$ in particular, if $a=a_{1}\times a_{2},$ then
\[
|\mathcal{T}_{\sigma}(\mathcal{V}_{a_{1}\times a_{2}})|=|\mathcal{T}_{\sigma
}(\mathcal{V}_{a_{1}})|\times|\mathcal{T}_{\sigma}(\mathcal{V}_{a_{2}})|.
\]

\end{proposition}

\begin{proof}
Since $\mathcal{T}_{\sigma}\left(  u\right)  $ is defined for all
$u\in\mathcal{V}_{a}$ and $\mathcal{T}_{\sigma}$ preserves the poset structure
of $\mathcal{V}_{a}$, the cells $a=\left\vert \mathcal{V}_{a}\right\vert $ and
$|\mathcal{T}_{\sigma}(\mathcal{V}_{a})|$ are combinatorially isomorphic.
\end{proof}

\begin{example}
\label{1.234}The action of $\mathcal{T}$ on the four vertices of $12|34$
partitions the $24$ vertices of $P_{4}$ into four mutually disjoint sets of
six vertices each. The vertices $v_{1}=1|2|3|4$ and $v_{2}=1|2|4|3$ of edge
$e=1|2|34$ correspond respectively to
\[
\left[
\begin{array}
[c]{c}%
\curlywedge\\
\curlywedge\\
\curlywedge
\end{array}
\right]  \left[
\begin{array}
[c]{c}%
\curlywedge\text{ }\mathbf{1}\\
\curlywedge\text{ }\mathbf{1}\\
\curlywedge\text{ }\mathbf{1}%
\end{array}
\right]  \left[
\begin{array}
[c]{c}%
\curlyvee\curlyvee\curlyvee\\
\mathbf{1}\text{ }\mathbf{1}\text{ }\mathbf{1}%
\end{array}
\right]  \left[  \curlyvee\curlyvee\curlyvee\right]  \text{ and }\left[
\begin{array}
[c]{c}%
\curlywedge\\
\curlywedge\\
\curlywedge
\end{array}
\right]  \left[
\begin{array}
[c]{c}%
\curlywedge\text{ }\mathbf{1}\\
\curlywedge\text{ }\mathbf{1}\\
\curlywedge\text{ }\mathbf{1}%
\end{array}
\right]  \left[
\begin{array}
[c]{c}%
\mathbf{1}\text{ }\mathbf{1}\text{ }\mathbf{1}\\
\curlyvee\curlyvee\curlyvee
\end{array}
\right]  \left[  \curlyvee\curlyvee\curlyvee\right]  .
\]
There are two basic subdivision vertices $u_{1}$ and $u_{2}$ along $e$,
exactly one of which admits a non-trivial action of $\mathcal{T},$ namely,
\[
u_{1}=\left[
\begin{array}
[c]{c}%
\curlywedge\\
\curlywedge\\
\curlywedge
\end{array}
\right]  \left[
\begin{array}
[c]{c}%
\curlywedge\text{ }\mathbf{1}\\
\curlywedge\text{ }\mathbf{1}\\
\curlywedge\text{ }\mathbf{1}%
\end{array}
\right]  \left[
\begin{array}
[c]{c}%
\curlyvee\curlyvee\text{ }\mathbf{1}\\
\mathbf{1}\text{ }\mathbf{1}\text{ }\curlyvee
\end{array}
\right]  \left[  \curlyvee\curlyvee\curlyvee\right]
\]%
\[%
\begin{array}
[c]{lll}%
\begin{array}
[c]{cc}
& \mathcal{T}_{\sigma_{1}}=\mathcal{T}_{11}^{2}\\
u_{1}\text{\ } & \longrightarrow\\
&
\end{array}
& z_{1}= & \medskip\left[
\begin{array}
[c]{c}%
\curlywedge\\
\curlywedge\\
\curlywedge
\end{array}
\right]  \left[
\begin{array}
[c]{c}%
\curlyvee\text{ }\mathbf{1}\\
\mathbf{1}\text{ }\curlyvee
\end{array}
\right]  \left[
\begin{array}
[c]{c}%
\curlywedge\text{ }\mathbf{1}\\
\curlywedge\text{ }\mathbf{1}%
\end{array}
\right]  \left[  \curlyvee\curlyvee\curlyvee\right] \\
\text{ \ \ \ \ \ \ }\mathcal{T}_{\sigma_{2}}\searrow\medskip &  &
\hspace*{0.5in}\downarrow\mathcal{T}_{11}^{3}\\
& z_{2}= & \left[
\begin{array}
[c]{c}%
\curlywedge\\
\curlywedge\\
\curlywedge
\end{array}
\right]  \left[
\begin{array}
[c]{c}%
\curlyvee\text{ }\mathbf{1}\\
\mathbf{1}\text{ }\curlyvee
\end{array}
\right]  \left[  \curlyvee\curlyvee\right]  \left[  \curlywedge\text{
}\mathbf{1}\right]  .
\end{array}
\]
\vspace*{0.5in}\newline\hspace*{0.5in}\unitlength=0.3mm
\thicklines\begin{picture}(100,100)(-100,100)
\put(-50.00,0){\line(1,0){150}}
\put(-50.00,200){\line(1,0){150}}
\put(0,0){\line(0,1){200}}
\put(-50,100){\line(1,0){50}}
\put(0,150){\line(1,0){100}}
\put(0,0){\line(3,5){51}}
\put(0,125){\line(4,-3){100}}
\put(100,100){\line(2,1){40}}
\put(100,200){\line(2,1){40}}
\put(100,75){\line(2,1){24}}
\put(100,50){\line(2,1){12}}
\put(100.00,0){\line(0,1){200}}
\put(100.00,0){\line(2,1){40}}
\put(0,0){\makebox(0,0){$\bullet$}}
\put(0,100){\makebox(0,0){$\bullet$}}
\put(0,125){\makebox(0,0){$\bullet$}}
\put(0,150){\makebox(0,0){$\bullet$}}
\put(0,200){\makebox(0,0){$\bullet$}}
\put(52,86){\makebox(0,0){$\bullet$}}
\put(40.00,86){\makebox(0,0){$z_1$}}
\put(-10.00,125){\makebox(0,0){$z_2$}}
\put(90,49){\makebox(0,0){$u_1$}}
\put(90,75){\makebox(0,0){$u_2$}}
\put(100,0){\makebox(0,0){$\bullet$}}
\put(100,75){\makebox(0,0){$\bullet$}}
\put(100,50){\makebox(0,0){$\bullet$}}
\put(100,100){\makebox(0,0){$\bullet$}}
\put(100,150){\makebox(0,0){$\bullet$}}
\put(100,200){\makebox(0,0){$\bullet$}}
\put(100,-15){\makebox(0,0){$1|2|3|4$}}
\put(0.00,-15){\makebox(0,0){$1|3|2|4$}}
\put(-25.00,85){\makebox(0,0){$1|3|4|2$}}
\put(-25.00,150){\makebox(0,0){$1|4|3|2$}}
\put(75,105){\makebox(0,0){$1|2|4|3$}}
\put(125.00,150){\makebox(0,0){$1|4|2|3$}}
\put(55.00,40){\makebox(0,0){$1|234$}}
\put(135.00,60){\makebox(0,0){$12|34$}}
\put(-40,50){\makebox(0,0){$13|24$}} \put(-40,175){\makebox(0,0){$134|2$}}
\put(135,180){\makebox(0,0){$124|3$}} \put(50,175){\makebox(0,0){$14|23$}}
\put(50,-50){\makebox(0,0){\text{Figure 16. The subdivision of} $1|234$ \text{in} $PP_{3,3}.$}}
\end{picture}$\vspace*{2in}$\newline To physically position $z_{1}$ and
$z_{2},$ first note that $\mathcal{T}_{\sigma_{1}}\left(  \mathcal{V}%
_{e}\right)  = \left\{  \mathcal{T}_{\sigma_{1}}\left(  v_{1}\right)
=1|3|2|4 \right. ,$
\linebreak  $\left. \mathcal{T}_{\sigma_{1}}\left(  v_{2}\right)=1|4|2|3\right\}  $ 
and $\mathcal{T}_{\sigma_{2}}\left(  \mathcal{V}%
_{e}\right)  =\left\{  \mathcal{T}_{\sigma_{2}}\left(  v_{1}\right)
=1|3|4|2,\text{ }\mathcal{T}_{\sigma_{2}}\left(  v_{2}\right)
=1|4|3|2\right\}  .$ Now $e=A_{1}|A_{2}|B_{1}=1|2|34$ and $\sigma_{1}\left(
e\right)  =1324;$ thus $\sigma_{1}\left(  B_{1}\right)  $ is not contiguous in
$\sigma_{1}\left(  e\right)  $ and $\sigma_{1}\left(  A_{2}\right)  $ breaks
the contiguity of $B_{1}$ in $\sigma_{1}\left(  e\right)  .$ Thus
$D_{1}=\sigma_{1}\left(  A_{1}\right)  $ and $\mathcal{T}(\sigma_{1},e)=1|234.
$ On the other hand, $\sigma_{2}\left(  e\right)  =1342;$ in this case
$\mathcal{T}(\sigma_{2},e)=1|34|2$ since $\sigma_{2}$ acts on the blocks of
$e$ as a $\left(  2,1\right)  $-shuffle$.$ Consequently, we represent the
vertices $z_{1}$ and $z_{2}$ as interior points of the faces $1|234$ and
$1|34|2,$ respectively. To complete the subdivision of $1|234,$ use the poset
structure to construct new edges from $u_{1}$ to $z_{1}$ and from $z_{1}$ to
$z_{2},$ and apply Proposition \ref{Prop5} to the subdivision cell
$a=(v_{1},u_{1})\subset e$ to construct the edge $|\mathcal{T}_{\sigma_{1}%
}(\mathcal{V}_{a})|$ from $1|3|2|4$ to $z_{1}.\ $Then $1|234=d_{1}\cup
d_{2}\cup d_{3}$ in which $\mathcal{V}_{d_{1}}=\{u_{1},u_{2},v_{2}%
,\mathcal{T}_{\sigma_{1}}(v_{2}),\mathcal{T}_{\sigma_{2}}(v_{2}),z_{2}%
,z_{1}\},$ $\mathcal{V}_{d_{2}}=\{u_{1},v_{1},\mathcal{T}_{\sigma_{1}}%
(v_{1}),z_{1}\},$ and $\mathcal{V}_{d_{3}}=\{\mathcal{T}_{\sigma_{1}}%
(v_{1}),\mathcal{T}_{\sigma_{2}}(v_{1}),z_{2},z_{1}\}$ (see Figure 16 and
Example \ref{1.234c}). An algebraic interpretation of these cells appears in
the discussion of $KK_{3,3}$ following Theorem 1.\ 
\end{example}

\subsection{$PP$-factorization of proper cells}

Recall that an element of $\mathcal{PP}_{n,m}\ $ is assigned to a unique
directed piece-wise linear path from $(m+1,1)$ to $(1,n+1)$ in $\mathbb{N}%
^{2}$ with $m+n$ components of unit length (see Figure 13). Let $\Pi_{n,m}$
denote the set of all such paths and consider the map $\pi:\mathcal{PP}%
_{n,m}\rightarrow\Pi_{n,m}$. If $u\in{\operatorname*{con}}\left(  X_{m}%
^{n+1}\times Y_{n}^{m+1}\right)  ,$ i.e., $u$ is a vertex of $\underline{m}%
|\left(  \underline{n}+m\right)  ,$ then $\pi$ restricts to a bijection
$\left\{  \mathcal{T}_{\sigma}\left(  u\right)  \mid\left(  m,n\right)
\text{-shuffles }\sigma\right\}  \leftrightarrow\Pi_{n,m},$ and in view of
Proposition \ref{Prop1}, $\pi$ assigns each vertex of $P_{m+n}$ to a path in
$\Pi_{n,m}$ albeit non-injectively.

Now consider a proper cell $c=C_{1}|\cdots|C_{s}\subset|{\operatorname{con}%
}\mathcal{PP}_{n,m}|\leftrightarrow P_{m+n}.$ Each factor $C_{t}$ is a
permutahedron $P_{m_{t}+n_{t}}$ whose vertices are assigned to connected
subpaths of paths in $\Pi_{n,m}.$ Assign $c$ to a directed piece-wise linear
path $\varepsilon_{c}=\cup\varepsilon_{t}$ in the following way: Write
$c=E_{1}|\cdots|E_{f}\ast_{\left(  \mathbf{i}{;}\mathbf{j}\right)  }%
F_{1}|\cdots|F_{g}$ and obtain sequences
\begin{equation}
\mathbf{\gamma}=\left\{  m+1=\gamma_{0}>\cdots>\gamma_{f}=1\right\}
\ \ \text{and}\ \ \mathbf{\delta}=\left\{  1=\delta_{0}<\cdots<\delta
_{g}=n+1\right\}  ,\label{subsequences}%
\end{equation}
where $\gamma_{i+1}=\gamma_{i}-\#E_{f-i}$ and $\delta_{j+1}=\delta
_{j}+\#F_{g-j},$ and assign $C_{t}=E_{t}^{\prime}\cup\left(  F_{t}^{\prime
}+m\right)  $ to the path

\begin{enumerate}
\item $\varepsilon_{t}:(\gamma_{t^{\prime}-1}\,,n+1-j)\rightarrow
(\gamma_{t^{\prime}}\,,n+1-j),$ if $C_{t}=E_{\gamma_{t^{\prime}}}$ for some
$t^{\prime}$ and maximal $j$ such that $F_{s_{1}}^{\prime},...,F_{s_{j}%
}^{\prime}\neq\varnothing$ and $s_{1}<\cdots<s_{j}<t;$

\item $\varepsilon_{t}:(i,\delta_{t^{\prime}})\rightarrow(i,\delta_{t^{\prime
}+1}),$ if $C_{t}=F_{\delta_{t^{\prime}}}+m$ for some $t^{\prime}$ and maximal
$i$ such that $E_{s_{1}}^{\prime},...,E_{s_{i}}^{\prime}\neq\varnothing$ and
$s_{1}<\cdots<s_{i}<t;$

\item $\varepsilon_{t}:(\gamma_{i},\delta_{j})\rightarrow(\gamma_{i+1}%
,\delta_{j+1}),$ if $C_{t}=E_{\gamma_{i}}\cup(F_{\delta_{j}}+m)$ with
$E_{\gamma_{i}},F_{\delta_{j}}\neq\varnothing$ for some $i,j.$
\end{enumerate}

\noindent In particular, a cell $a=A_{1}\cdots A_{k}B_{l}\cdots B_{1}%
\subset\left\vert {\operatorname{con}}\left(  X_{m}^{n+1}\times Y_{n}%
^{m+1}\right)  \right\vert \leftrightarrow\underline{m}|\left(  \underline{n}%
+m\right)  $ is assigned to the path
\[
\varepsilon_{a}:(m+1,\beta_{0})\overset{B_{1}}{\rightarrow}\cdots
\overset{B_{l}}{\rightarrow}\left(  m+1,\beta_{l}\right)  \overset{A_{k}%
}{\rightarrow}\left(  \alpha_{1},n+1\right)  \overset{A_{k-1}}{\rightarrow
}\cdots\overset{A_{1}}{\rightarrow}\left(  \alpha_{k},n+1\right)  ,
\]
where $\mathbf{\alpha}=\left\{  m+1=\alpha_{0}>\cdots>\alpha_{k}=1\right\}  $
and $\mathbf{\beta}=\left\{  1=\beta_{0}<\cdots<\beta_{l}=n+1\right\}  $ (case
(3) does not occur). Thus if $c=\mathcal{T}\left(  \sigma,a\right)  ,$ the
observation in (\ref{specialjoin}) implies that $\mathbf{\gamma}%
\subseteq\mathbf{\alpha}$ and $\mathbf{\delta}\subseteq\mathbf{\beta.}$

Given a subdivision cell $d\subseteq\mathcal{T}(\sigma,a),$ there is a
subdivision subcomplex $u\subset a$ such that $d=|\mathcal{T}_{\sigma
}(\mathcal{V}_{u})|.$ Representing $a$ as a partition $U_{1}|\cdots|U_{s}$ of
$\underline{m}|\left(  \underline{n}+m\right)  ,$ there is a Cartesian
product decomposition $d=D_{1}\times\cdots\times D_{s}$ in which $D_{t}$ is a
subdivision cell of $\mathcal{T}(\sigma,U_{t})$. The representation
$U_{1}|\cdots|U_{s}=E\ast_{\left(  \mathbf{i}{;}\mathbf{j}\right)  }F$ relates
the paths associated with the vertices of $\mathcal{T}(\sigma,U_{t})$ to the
vertices of $D_{t}$, and in view of case (3) above, the vertices of $D_{t}$
are assigned to paths related to those $z\in\mathcal{T}_{\sigma}%
(\mathcal{V}_{u})$ given by the action of $\mathcal{T}_{\sigma}$ on the matrix
sequences $x_{\gamma_{i+1}}\cdots x_{\gamma_{i}-1}y_{\delta_{j+1}-1}\cdots
y_{\delta_{j}}$ associated with the vertices of $u$ as a $(\gamma_{i}%
-\gamma_{i+1},\,\delta_{j+1}-\delta_{j})$-shuffle. But in every case, there is
the Cartesian product decomposition
\[
D_{t}=(e_{y_{1}^{t},x_{1}^{t}}\times\cdots\times e_{y_{1}^{t},x_{p_{t}}^{t}%
})\times\cdots\times(e_{y_{q_{t}}^{t},x_{1}^{t}}\times\cdots\times
e_{y_{q_{t}}^{t},x_{p_{t}}^{t}}),
\]
where $e_{y_{j}^{t},x_{i}^{t}}$ is some cell of $PP_{y_{j}^{t},x_{i}^{t}}$ and
$(p_{t},q_{t})\in\{(\gamma_{t^{\prime}},n+1-j),(i,\delta_{t^{\prime}}),$
$(\gamma_{i+1},\delta_{j})\},$ $p_{1}+\cdots+p_{t}\in\{\gamma_{t^{\prime}%
-1},i,\gamma_{i}\},$ and $q_{1}+\cdots+q_{t}\in\{n+1-j,\delta_{t^{\prime}%
+1},\delta_{j+1}\}$ (Cartesian 
\linebreak product decomposition of $D_{t}$ is trivial
whenever $m=1$ or $n=1$). Therefore every proper cell $e_{n,m}\subset
PP_{n,m}$ has a \emph{Cartesian matrix factorization}%
\begin{equation}
e_{n,m}=[(e_{y_{1}^{1},x_{1}^{1}}\times\cdots\times e_{y_{1}^{1},x_{p_{k}}%
^{1}})\times\cdots\times(e_{y_{q_{k}}^{1},x_{1}^{1}}\times\cdots\times
e_{y_{q_{k}}^{1},x_{p_{k}}^{1}})]\times\cdots\hspace*{0.7in}\label{proper}%
\end{equation}%
\[
\hspace*{1in}\times\lbrack(e_{y_{1}^{s},x_{1}^{s}}\times\cdots\times
e_{y_{1}^{s},x_{p_{k}}^{s}})\times\cdots\times(e_{y_{q_{k}}^{s},x_{1}^{s}%
}\times\cdots\times e_{y_{q_{k}}^{s},x_{p_{k}}^{s}})],
\]
where $p_{1}=q_{s}=1$ and $s\geq2.$ The decomposition in (\ref{proper}) is a
\emph{$PP$-factorization} if each factor $e_{y_{j}^{k},x_{i}^{k}}$ lies in the
family $PP.$

Indeed, each factor $e_{y_{j}^{k},x_{i}^{k}}$ of $e_{n,m}$ has a Cartesian
matrix factorization with $x_{i}^{k}+y_{j}^{k}<m+n,$ and we may inductively
apply the decomposition in (\ref{proper}) to obtain a decomposition of
$e_{n,m}$ as a Cartesian product of polytopes in the family $PP.\,\ $This
decomposition involves Cartesian products in two settings:\ Those within
bracketed quantities correspond to tensor products of entries in a bisequence
monomial (controlled by $\Delta_{P}^{\left(  k\right)  }$) and those between
bracketed quantities correspond to $\Upsilon$-products of bisequence
monomials. And indeed, this decomposition is encoded by a leveled tree
$\Psi(e_{n,m})$ constructed in the same way we constructed $\Psi(\phi(\xi))$
for $\xi\in\mathfrak{G}$. Whereas the levels and the leaves of $\Psi(\phi
(\xi))$ are bisequence and $\Theta$-factorizations, the levels and leaves of
$\Psi(e_{n,m})$ are Cartesian matrix and $PP$-factorizations.

\begin{example}
\label{1.234c}Refer to Example \ref{1.234} and consider the codimension $1$
cell $c=C_{1}|C_{2}=1|234\subset P_{2+2}.$ Write $c=E_{1}|E_{2}\ast
F_{1}=1|2\ast12$ and obtain $\gamma_{0}=3,$ $\gamma_{1}=2,$ $\gamma_{2}=1;$
and $\delta_{0}=1,$ $\delta_{1}=3.$ Then $C_{1}=E_{1}^{\prime}\cup
F_{1}^{\prime}=E_{1}\cup\varnothing$ and $C_{2}=E_{2}^{\prime}\cup\left(
F_{2}^{\prime}+2\right)  =E_{2}\cup\left(  F_{1}+2\right)  .$ The path $C_{1}$
is assigned to the path component $\varepsilon_{1}:\left(  2,3\right)
\rightarrow\left(  1,3\right)  \ $and $C_{2}$ is assigned to $\varepsilon
_{2}:\left(  3,1\right)  \rightarrow\left(  2,3\right)  ;$ in this case there
is the action of a $(\gamma_{0}-\gamma_{1},\delta_{1}-\delta_{0}%
)=(1,2)$-shuffle on $x_{\gamma_{0}-1}y_{\delta_{1}-1}y_{\delta_{0}}=x_{2}%
y_{2}y_{1},$ which generates (classes of ) vertices of $C_{2}.$ Let $u$ be the
subdivision subcomplex of $1|2|34$ consisting of the two edges $\left(
u_{1},u_{2}\right)  $ and $\left(  u_{2},v_{2}\right)  $ (see Figure 16). Then
for $i=1,2,3,$ the subdivision cell $d_{i}=D_{1}^{i}\times D_{2}^{i},$ where
$D_{1}^{i}=C_{1}=P_{1}$ is a vertex and $D_{2}^{i}\subset C_{2}$ has the form
$D_{2}^{i}=e_{3,2}^{i}\times e_{3,1}^{i},$ where $(\dim e_{3,2}^{1},\dim
e_{3,1}^{1})=(2,0)$ and $(\dim e_{3,2}^{i},\dim e_{3,1}^{i})=(1,1)$
($e_{3,2}^{1}=PP_{2,1}$ is a heptagon and $e_{3,1}^{1}$is a vertex of
$PP_{2,0};$ $e_{3,2}^{i}$ is an edge of $PP_{2,1}$ and $e_{3,1}^{i}=PP_{2,0}$
for $i=2,3$). Thus up to homeomorphism we have%
\begin{align*}
d_{1} &  =[PP_{0,1}\times PP_{0,1}\times PP_{0,1}]\times\left[  PP_{2,1}%
\times(PP_{1,0}\times PP_{1,0})\right] \\
d_{2} &  =[PP_{0,1}\times PP_{0,1}\times PP_{0,1}]\times\left[  \left(
PP_{1,1}\times PP_{1,0}\right)  \times PP_{2,0}\right] \\
d_{3} &  =[PP_{0,1}\times PP_{0,1}\times PP_{0,1}]\times\left[  \left(
PP_{1,0}\times PP_{1,1}\right)  \times PP_{2,0}\right]  .
\end{align*}

\end{example}

\subsection{The projection $\vartheta\vartheta:PP\rightarrow KK$}

The final piece of our construction establishes a geometric interpretation of
the projection $\vartheta\vartheta:PP_{m,n}\rightarrow KK_{n+1,m+1}$ induced
by the quotient map $\mathcal{PP}_{n,m}\rightarrow
\linebreak \mathcal{KK}_{n+1,m+1}.$ Let
$\mathcal{P}_{n,m}={con}\mathcal{PP}_{n,m},$ $P_{n,m}=\left\vert
\mathcal{P}_{n,m}\right\vert ,$ and $K_{n+1,m+1}=|\mathcal{P}_{n,m}\diagup
\sim|;$ we obtain
 $KK_{n+1,m+1}$ as the subdivision of $K_{n+1,m+1}$ that
commutes the following diagram:
\[%
\begin{array}
[c]{ccc}%
\ \ \ \ PP_{n,m} & \overset{\approx}{\longrightarrow} & P_{m+n}\vspace
*{0.1in}\\
\vartheta\vartheta\downarrow\text{ \ \ \ } &  & \ \ \ \ \ \downarrow
\vartheta\vspace*{0.1in}\\
KK_{n+1,m+1} & \underset{\approx}{\longrightarrow} & K_{n+1,m+1}%
\end{array}
\]
(the horizontal maps are non-cellular homeomorphisms induced by the
subdivision process). We identify the cellular chains $C_{\ast}\left(
KK\right)  $ with the free matrad ${\mathcal{H}}_{\infty}$ and prove that the
restriction of the free resolution of prematrads $\rho^{^{{\operatorname*{pre}%
}}}:F^{^{{\operatorname*{pre}}}}(\Theta)\rightarrow{\mathcal{H}}$ to
${\mathcal{H}}_{\infty}$ is a free resolution in the category of matrads.

To simplify notation, we suppress the subscripts of $\vartheta_{n,m}%
:P_{n,m}\rightarrow K_{n+1,m+1}$ when $m$ and $n$ are clear from context.
Since $\left\vert \mathcal{P}_{n,m}\right\vert =P_{n,m}=P_{m+n},$ a proper
face $e\subset P_{n,m}$ is a product of permutahedra
\[
e=P_{n_{1},m_{1}}\times\cdots\times P_{n_{s},m_{s}}%
\]
and projects to a product
\[
\tilde{e}=\vartheta\left(  e\right)  =\vartheta(P_{n_{1},m_{1}})\times
\cdots\times\vartheta(P_{n_{s},m_{s}})=K_{n_{1}+1,m_{1}+1}\times\cdots\times
K_{n_{s}+1,m_{s}+1}.
\]
The fact that $\vartheta_{n,m}={Id}$ when $1\leq m,n\leq2$ implies
$K_{n+1,m+1}=P_{m+n}$; also, $K_{n,2}\cong K_{2,n}$ is the {multiplihedron }
$J_{n}$ for all $n$ (see \cite{Stasheff2}, \cite{Iwase}, \cite{SU2},
\cite{SU4}). The faces $24|13$ and $1|24|3$ of $P_{3,1}$ are degenerate in
$K_{4,2}$ since $\vartheta_{3,1}\left(  24|13\right)  =24|1|3$ and
$\vartheta_{3,1}\left(  1|24|3\right)  =1|2|4|3;$ and dually, the faces
$24|13$ and $2|13|4$ of $P_{1,3}$ are degenerate in $KK_{2,4}$ since
$\vartheta_{1,3}\left(  24|13\right)  =2|4|13$ and $\vartheta_{1,3}\left(
2|13|4\right)  =2|1|3|4.$ Observe that the product cell $K_{m+1}\times
K_{n+1}=\vartheta(P_{m}\times P_{n})\subset\vartheta\left(  P_{m+n}\right)
=\vartheta\left(  P_{n,m}\right)  =K_{n+1,m+1}\ $admits the $(m,n)$%
-subdivision $K_{m+1}^{(n)}\times K_{n+1}^{(m)}=\vartheta\vartheta(P_{m}%
^{(n)}\times P_{n}^{(m)})\subset\vartheta\vartheta\left(  PP_{n,m}\right)
=KK_{n+1,m+1}. $

\begin{example}
The $(2,2)$-subdivision $K_{3}^{(2)}\times K_{3}^{(2)}$ of the face
$\vartheta_{2,2}(12|34)=K_{3}\times K_{3}\subset K_{3,3}$ produces 9 cells of
$KK_{3,3}$ (see Figures 14 and 20); the $(3,1)$-subdivision $K_{2}^{(3)}\times
K_{4}^{\left(  1\right)  }$ of the face $\vartheta_{1,3}(1|234)=K_{2}\times
K_{4}\subset K_{4,2}$ produces 6 cells of $KK_{4,2} $ (see Figures 3 and 21).
\end{example}

Define
\[
\widetilde{\mathcal{BS}}_{n+1,m+1}=(\vartheta_{\ast}\times\vartheta_{\ast
})\left(  \mathcal{BS}_{n,m}\right)  .
\]

\begin{example}
The biassociahedron $KK_{3,3}=PP_{2,2}$ and has $44$ vertices, $16$ of which
lie in $12|34.$ Of these $16$ vertices, $4$ lie in $\mathcal{P}_{2,2}$ and
generate the other $19$ vertices of $K_{3,3}=P_{4};$ another $4$ lie
$\widetilde{\mathcal{BS}}_{3,3}$ and generate the $8$ remaining vertices of
$KK_{3,3}$ (see Figure 20). By contrast, $KK_{4,2}$ is a non-trivial quotient
of $PP_{3,1}.$ As in Tonks' projection $\vartheta_{0}:P_{n}\rightarrow
K_{n+1},$ we identify faces of $PP_{3,1}$ indexed by isomorphic graphs
(forgetting levels) as pictured in Figure 17. Here an equivalence class of
graphs, which labels a face of the target interval, contains the three graphs
horizontally to its left.
\end{example}

\begin{center}
\includegraphics[
height=1.97in,
width=2.9793in
]%
{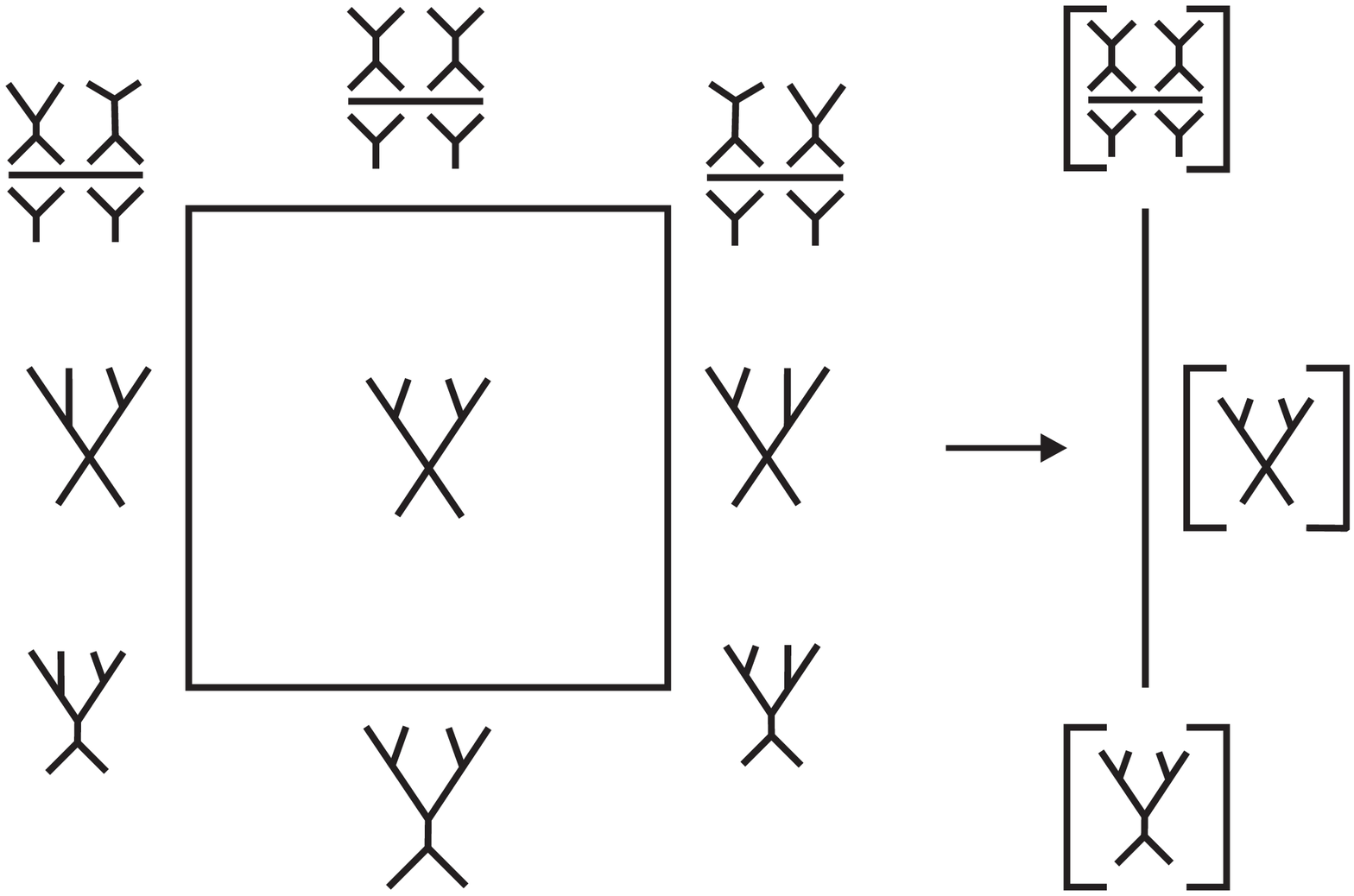}%
\\
Figure 17: Projection of a degenerate square in $PP_{3,1}$ to $KK_{4,2}.$%
\end{center}

The biassociahedra $KK_{1,1},$ $KK_{2,1}$ and $KK_{1,2}$ are isolated vertices
and correspond to the free matrad generators $\mathbf{1},$ $\theta_{1}^{2}$
and $\theta_{2}^{1},$ respectively. The biassociahedra $KK_{n,m}$ with $4\leq
m+n\leq6$ are pictured in Figures 18 through 22 below and labelled by
partitions and (co)derivation leaf sequences. Note that $KK_{n,m}\cong
KK_{m,n}$ for all $m,n\geq1$ and $KK_{2,m}$ is a subdivision of $J_{m}$ when
$m\geq3.$

If $\tilde{e}_{n,m}$ is a proper face of $KK_{n,m},$ the decomposition in
(\ref{proper}) induces a product decomposition of the form%
\begin{equation}
\tilde{e}_{n,m}=[(\tilde{e}_{y_{1}^{1},x_{1}^{1}}\times\cdots\times\tilde
{e}_{y_{1}^{1},x_{p_{k}}^{1}})\times\cdots\times(\tilde{e}_{y_{q_{k}}%
^{1},x_{1}^{1}}\times\cdots\times\tilde{e}_{y_{q_{k}}^{1},x_{p_{k}}^{1}%
})]\hspace*{0.7in}\label{proper2}%
\end{equation}%
\[
\hspace*{1.3in}\times\cdots\times\lbrack(\tilde{e}_{y_{1}^{s},x_{1}^{s}}%
\times\cdots\times\tilde{e}_{y_{1}^{s},x_{p_{k}}^{s}})\times\cdots
\times(\tilde{e}_{y_{q_{k}}^{s},x_{1}^{s}}\times\cdots\times\tilde
{e}_{y_{q_{k}}^{s},x_{p_{k}}^{s}})],
\]
where $p_{1}=q_{s}=1$ and $s\geq2.$ Here \textquotedblleft$\times
$\textquotedblright\ within a bracketed quantity corresponds to the tensor
product in a bisequence monomial (controlled by certain iterations of
$\Delta_{P}$ and the product cell within $k^{th}$ bracket is thought of as a
subdivision cell of $K_{n_{k}+1,m_{k}+1}$ in the decomposition of $e_{n,m}$ in
(\ref{proper})) and \textquotedblleft$\times$\textquotedblright\ between
bracketed quantities corresponds to a $\Upsilon$-product, and each $\tilde
{e}_{y_{j}^{k},x_{i}^{k}}$ has the form given by (\ref{proper2}) with
$x_{i}^{k}+y_{j}^{k}<m+n.$

We distinguish between two kinds of faces in (\ref{proper2}). A Type I face is
detected by the diagonal $\Delta_{P}\ $and its representation in
(\ref{proper2}) has $(p_{k},q_{k})>(1,1)$ for all $k$; thus $\Delta_{P}$ is
only involved in forming the Cartesian products in parentheses. A Type II face
$\tilde{e}_{n,m}$ is independent of $\Delta_{P}$ and its representation in
(\ref{proper2}) satisfies $(p_{k},q_{k})=(1,1)$ for all $k;$ thus $\tilde
{e}_{n,m}$ has the form $KK_{n,i_{2}}\times K_{i_{2}-i_{1}+1}\times
\cdots\times K_{i_{s}-i_{s-1}+1},\,1\leq i_{2}<\cdots<i_{s}=m,$ or
$K_{j_{0}-j_{1}+1}\times\cdots\times K_{j_{s-2}-_{s-1}+1}\times KK_{j_{s-1}%
,m},\,1\leq j_{s-1}<\cdots<j_{0}=n.$ In particular, a codimension 1 face (when
$s=2$) has the form $KK_{n,i}\times K_{m+1-i}$ or $K_{n+1-j}\times KK_{j,m}.$
Consequently, each cell $\tilde{e}_{n,m}\subset KK_{n,m}$ is associated with a
levelled tree $\Psi(\tilde{e}_{n,m}),$ whose levels are representations given
by (\ref{proper2}) and whose leaves are \emph{$KK$-factorizations}.

The assignment $\iota:\theta_{p}^{q}\mapsto KK_{q,p},$ which preserves levels,
induces a set map
\begin{equation}
\iota:\mathfrak{G}_{n,m}\rightarrow\left\{  \text{faces of}\,KK_{n,m}\right\}
\label{iota}%
\end{equation}
that sends balanced factorizations to Cartesian matrix factorizations and has
the following properties:\smallskip

\begin{enumerate}
\item[\textit{(i)}] The restriction of $\iota$ to $0$-dimensional module
generators of $F_{n,m}(\Theta)$ establishes a bijection with vertices of
$KK_{n,m}$ by replacing $\theta_{2}^{1}$ with $\curlywedge$ and $\theta
_{1}^{2}$ with $\curlyvee$ in each entry of $\Psi(\beta).$\smallskip

\item[\textit{(ii)}] There is a \emph{location map}
\begin{equation}
\tilde{\iota}:\mathfrak{G}_{q+1,p+1}\rightarrow\left\{  \text{faces
of}\,P_{p+q}\right\} \label{tiota}%
\end{equation}
that commutes the following diagram of set maps
\[%
\begin{array}
[c]{cccc}%
\mathfrak{G}_{q+1,p+1} & \overset{\tilde{\iota}}{\longrightarrow} & \left\{
\text{ faces of}\,P_{p+q}\right\}  & \\
{\iota}\downarrow &  & \downarrow\vartheta_{q,p} & \\
\left\{  \text{faces of}\,KK_{q+1,p+1}\right\}  & \overset{\nu
}{\longrightarrow} & \left\{  \text{faces of}\,K_{q+1,p+1}\right\}  , &
\end{array}
\]
where $\nu$ sends a cell of $KK$ to the cell of $K$ of minimal dimension
containing it. Indeed, if $\beta=C_{s}\cdots C_{1}\in\mathcal{B}$ is balanced
representative of $\theta$ with $C_{k}\in\mathbf{G}_{\mathbf{x}^{k}%
}^{\mathbf{y}^{k}},$ consider the (co)derivation leaf sequences $((\mathbf{x}%
^{1},\mathbf{y}^{1}),...,(\mathbf{x}^{s},\mathbf{y}^{s})),$ and let
$\mathbf{x}^{i_{1}},...,\mathbf{x}^{i_{k}}$ and $\mathbf{y}^{j_{1}%
},...,\mathbf{y}^{j_{l}}$ be the subsequences obtained by removing all
$\mathbf{x}^{i},\mathbf{y}^{j}=\mathbf{1}.$ Thinking of these subsequences as
row and column descent sequences, consider the corresponding faces
$A=A_{1}|\cdots|A_{k}\subset P_{m-1}$ and $B=B_{1}|\cdots|B_{l}\subset
P_{n-1}$ and set%
\[
e_{\theta}=\chi(A)\ast_{\left(  \mathbf{i,j}\right)  }B\subset P_{m+n-2},
\]
where $\chi:P_{m-1}\rightarrow P_{m-1}$ is the cellular involution defined by
$\chi(A_{1}|\cdots|A_{k})$\linebreak$=(m-A_{k})|\cdots|(m-A_{1})$ and
$(\mathbf{i},\mathbf{j})=(i_{1}<\cdots<i_{k};j_{1}<\cdots<j_{l})$. Then
$e_{\theta}$ is the unique cell of minimal dimension $\geq k$ such that
$\iota(\theta)\subset\vartheta_{n-1,m-1}(e_{\theta});$ in particular, when
$s=2,$ $\mathbf{x}^{1}=\mathbf{x}$ and $\mathbf{y}^{2}=\mathbf{y},$ and we
recover the special cell $e_{\theta}=e_{(\mathbf{y},\mathbf{x})}$ defined in
(\ref{special}). Thus, the term \textquotedblleft location
map\textquotedblright\ suggests the fact that $\tilde{\iota}$ points out the
position of the image cells $\iota(\theta)$ with respect to cells of the
permutahedron $P_{m+n}.$ Under $\tilde{\iota},$ the associativity of the
$\Upsilon$-product on $\mathfrak{G}$ is compatible with the associativity of
the partitioning procedure in $\underline{p+q}$ by which $a_{1}|\cdots|a_{k} $
is obtained from the (ordered) set $a_{1}\cdots a_{k}$ by inserting bars:
Given $\theta,\zeta\in\mathfrak{G},$ let $\xi\in\mathfrak{G}$ be the component
shared by $\partial(\theta)$ and $\partial(\zeta)$ in $\mathcal{H}_{\infty}.$
Then $\vartheta\left(  \tilde{\iota}(\xi)\right)  \subset\vartheta\left(
\tilde{\iota}(\theta)\right)  \cap\vartheta\left(  \tilde{\iota}%
(\zeta)\right)  $ in $K_{q+1,p+1}=\vartheta(P_{p+q})$ (see Example
\ref{twoasso}).\smallskip

\item[\textit{(iii)}] Let $\dim\theta=k$ and let $\sigma_{\theta}$ be the set
of all $0$-dimensional elements of $\mathcal{H}_{\infty}$ obtained by all
possible compositions $\partial_{i_{k}}\cdots\partial_{i_{1}}(\theta)$ where
$\partial_{i}$ is a component of $\partial=\sum_{i}\partial_{i}.$ Then
$\iota(\theta)$ is the $k$-face of $KK_{n,m}$ spanned on the set $\iota
(\sigma_{\theta}).$\smallskip

\item[\textit{(iv)}] If $(m,n)\in\{(2,0),(1,1),(0,2)\}$ in item \textit{(ii)}
and $k=1$ in item \textit{(iii) }, then
\linebreak  $KK_{n+1,m+1}=K_{n+1,m+1}=P_{m+n}$ is
an interval and \textit{(ii)} agrees with \textit{(iii)} under the equality
$\iota=\tilde{\iota}|_{\mathfrak{G}_{n+1,m+1}}$ for $m+n\leq2.$
\end{enumerate}

\begin{remark}
Since $\tilde{\iota}$ is not surjective, the action of the (pre)matrad axioms
on Type II generators forces us to obtain $KK_{n+1,m+1}$ as a quotient of
$PP_{n,m}$ modulo combinatorial relations in $PP_{n,m}$ as indicated in Figure
16 above, and thereby extend the equality $K_{n+1}=P_{n}\diagup\sim$ induced
by Tonks' projection (see Theorem \ref{kkh} below).
\end{remark}

\begin{example}
\label{twoasso}The action of the map $\theta\rightarrow\tilde{\iota}(\theta)$
involving associativity is illustrated by the example in (\ref{assoc}):
\[
((221),(31))\rightarrow146|2357\ \ \text{and}\ \ ((41),(211))\rightarrow
12456|37
\]
while
\[
\left(  ((221),(2)),\,((21),(21)),\,((2),(211))\right)  \rightarrow146|25|37
\]
(on left-hand sides only (co)derivation leaf sequences of underlying matrad
module generator are shown). Also, from Example \ref{colrow} we have
\[
C_{3}C_{2}C_{1}\overset{\tilde{\iota}}{\rightarrow}357|14|26.
\]

\end{example}

The properties above imply that $\iota$ is a bijection so that $\mathcal{C}%
_{n,m}$ indexes the faces of $KK_{n,m}.$ Define the boundary map in the
cellular chain complex $C_{\ast}\left(  KK_{n,m}\right)  $ by
\begin{equation}
\partial(
\raisebox{-0.0553in}{\includegraphics[
trim=0.000000in 0.140808in 0.000000in 0.254377in,
height=0.1833in,
width=0.1842in
]%
{dblvee_nm.eps}%
}
\ )=\sum\limits_{(\alpha,\beta)\in\mathcal{AB}_{m,n}}(-1)^{\epsilon
+\epsilon_{\alpha}+\epsilon_{\beta}}e_{\alpha,\beta},\label{kkdiff}%
\end{equation}
where $\epsilon$ is the sign of the cell $e_{(\mathbf{y},\mathbf{x})}\subset
P_{n-1,m-1}$ defined by (\ref{special}). This sign reflects the fact that the
sign of a subdivision cell in the boundary inherits (as a component) the sign
of the boundary. Therefore, we immediately obtain:

\begin{theorem}
\label{kkh} For each $m,n\geq1,$ there is a canonical isomorphism of chain
complexes
\begin{equation}
\iota_{\ast}:({\mathcal{H}}_{\infty})_{n,m}\overset{\approx}{\longrightarrow
}C_{\ast}(KK_{n,m})\label{inftymatrad}%
\end{equation}
extending the standard isomorphisms%
\[
\mathcal{A}_{\infty}(n)=({\mathcal{H}}_{\infty})_{n,1}\overset{\approx
}{\longrightarrow}C_{\ast}(KK_{n,1})=C_{\ast}(K_{n})
\]
and%
\[
\mathcal{A}_{\infty}(m)=({\mathcal{H}}_{\infty})_{1,m}\overset{\approx
}{\longrightarrow}C_{\ast}(KK_{1,m})=C_{\ast}(K_{m}).
\]

\end{theorem}

In other words, the cellular chains of the biassociahedra $KK$ realize the
free matrad resolution ${\mathcal{H}}_{\infty}$ of the bialgebra matrad
$\mathcal{H}.$ In particular, consider the submodule $\tilde{\mathcal{H}%
}_{\infty}\subset{\mathcal{H}}_{\infty}$ spanned on the generating set
$\tilde{\Theta}$ fixed by summing of all (distinct) elements of $\mathfrak{G}$
in $\mathcal{H}_{\infty}$ that have the same leaf sequence form. Then
(\ref{tiota}) induces an isomorphism
\begin{equation}
(\vartheta\circ\tilde{\iota})_{\ast}:(\tilde{\mathcal{H}}_{\infty}%
)_{n,m}\overset{\approx}{\longrightarrow}C_{\ast}(K_{n,m})\label{tiotaiso}%
\end{equation}
by $(\vartheta\circ\tilde{\iota})_{\ast}(\tilde{\theta})=\vartheta\left(
\tilde{\iota}(\theta_{s})\right)  ,$ where $\theta_{s}\in\mathfrak{G}_{n,m} $
is any summand component of $\tilde{\theta}\in\tilde{\Theta},$ and the
following diagram commutes:
\[%
\begin{array}
[c]{cccc}%
\ \ \ \ \ \ \ \ \ (\tilde{\mathcal{H}}_{\infty})_{n,m} & {\longrightarrow} &
({\mathcal{H}}_{\infty})_{n,m} & \\
\vspace{1mm}(\vartheta\!\circ\!\tilde{\iota})_{\ast}\downarrow\,\approx &  &
\approx\,\downarrow{\iota}_{\ast} & \\
\ \ \ C_{\ast}(K_{n,m}) & \underset{\nu^{_{\#}}}{\longrightarrow} & C_{\ast
}(KK_{n,m}). &
\end{array}
\]

\begin{example}
\label{faces1.234} We have $\tilde{\iota}((111),(21))=1|234\subset P_{4},$ and
apply Example \ref{1.234c} for which bijection (\ref{inftymatrad})
implies\smallskip\newline\hspace*{1.2in}
$
\begin{array}
[c]{lll}%
\left[
\begin{array}
[c]{r}%
\theta_{2}^{1}\\
\theta_{2}^{1}\\
\theta_{2}^{1}%
\end{array}
\right]  \left[
\begin{array}
[c]{r}%
\theta_{2}^{3}\,\ \ \left[
\begin{array}
[c]{c}%
\mathbf{1}\\
\theta_{1}^{2}%
\end{array}
\right]  \theta_{1}^{2}%
\end{array}
\right]  & \leftrightarrow & d_{1}\\
&  & \\
\left[
\begin{array}
[c]{r}%
\theta_{2}^{1}\\
\theta_{2}^{1}\\
\theta_{2}^{1}%
\end{array}
\right]  \left[
\begin{array}
[c]{r}%
\left[
\begin{array}
[c]{r}%
\theta_{2}^{2}\\
\theta_{2}^{1}%
\end{array}
\right]  [\theta_{1}^{2}\text{ }\theta_{1}^{2}]\ \ \ \theta_{1}^{3}%
\end{array}
\right]  & \leftrightarrow & d_{2}\\
&  & \\
\left[
\begin{array}
[c]{r}%
\theta_{2}^{1}\\
\theta_{2}^{1}\\
\theta_{2}^{1}%
\end{array}
\right]  \left[
\begin{array}
[c]{r}%
\left[
\begin{array}
[c]{c}%
\theta_{1}^{2}\\
\mathbf{1}%
\end{array}
\right]  \theta_{2}^{2}\ \ \ \theta_{1}^{3}%
\end{array}
\right]  & \leftrightarrow & d_{3}.
\end{array}
$
\smallskip\newline 
The edge $\left(  1|3|2|4,z_{1}\right)  $ in Figure 16 is
the intersection $d_{2}\cap d_{3}$ corresponding to
\begin{equation}
A_{2}A_{1}=\left[
\begin{array}
[c]{r}%
\theta_{2}^{1}\\
\theta_{2}^{1}\\
\theta_{2}^{1}%
\end{array}
\right]  \left[
\begin{array}
[c]{r}%
\left(  \left[
\begin{array}
[c]{c}%
\theta_{1}^{2}\\
\mathbf{1}%
\end{array}
\right]  \left[
\begin{array}
[c]{r}%
\theta_{2}^{1}\\
\theta_{2}^{1}%
\end{array}
\right]  [\theta_{1}^{2}\text{ }\theta_{1}^{2}]\right)  \text{ }\theta_{1}^{3}%
\end{array}
\right]  \in\mathcal{B}_{3,3}.\label{codim}%
\end{equation}

\end{example}

The following proposition applies Proposition \ref{bases} to reformulate
Theorem \ref{kkh} in terms of the $\circledcirc$-operation defined in
(\ref{operation1}) for $(M,\gamma)=(F^{{\operatorname*{pre}}}(\Theta),\gamma).$

\begin{proposition}
\label{altkkf}Let $(F^{{\operatorname*{pre}}}(\Theta),\gamma)$ be the free
prematrad and $(\mathcal{H}_{\infty},\partial)$ be the $A_{\infty}$-matrad. If
$\xi=[(\theta\setminus\theta_{m}^{n})\circledcirc(\theta\setminus\theta_{m}^{n})]_{m}^{n}$ 
with $mn\geq3,$ the components of $\xi$ fit the boundary
of $KK_{n,m}$ and $\partial(\xi)=0.$
\end{proposition}

\noindent Thus, an $A_{\infty}$-bialgebra structure on a DGM $H$ is defined by
a morphism of matrads $\mathcal{H}_{\infty}\rightarrow U_{H}$ (compare
\cite{SU3}).

In our forthcoming paper \cite{SU4}, we construct the theory of relative
matrads and use it to define a morphism of $A_{\infty}$-bialgebras. Using
\emph{relative $A_{\infty}$-matrads,} we prove that over a field, the homology
of every biassociative DG bialgebra admits a canonical $A_{\infty}$-bialgebra
structure.
\newpage

\vspace{.2in}
\noindent\underline{For $KK_{2,2}$:}

$
\begin{array}
[c]{lllll}%
1|2 & \leftrightarrow & _{11}^{11} & = & \gamma(\theta_{2}^{1}\theta_{2}%
^{1}\,;\theta_{1}^{2}\theta_{1}^{2})\vspace{1mm}\newline\\
2|1 & \leftrightarrow & _{2}^{2} & = & \gamma(\theta_{1}^{2}\,;\theta_{2}%
^{1})\\
&  &  &  &
\end{array}
$

\vspace{0.2in}
\hspace{-0.6in} 
\unitlength=1.00mm\special{em:linewidth 0.4pt}
\linethickness{0.4pt} \begin{picture}(86.33,13.33)
\put(43.00,8.67){\line(1,0){43.00}} \put(43.33,8.67){\circle*{1.33}}
\put(85.67,8.67){\circle*{1.33}} \put(43.33,13.33){\makebox(0,0)[cc]{$1|2$}}
\put(43.33,03.33){\makebox(0,0)[cc]{$_{11}^{11}$}}
\put(85.67,13.00){\makebox(0,0)[cc]{$2|1$}}
\put(85.67,03.00){\makebox(0,0)[cc]{$_{2}^{2}$}}
\end{picture}
\begin{center}
Figure 18. The biassociahedron $KK_{2,2}$ (an interval)\vspace{0.1in}
\end{center}
\vspace{0.1in}

\noindent\underline{For $KK_{3,2}$:}

$
\begin{array}
[c]{lllll}%
1|23 & \leftrightarrow & _{11}^{111} & = & \gamma(\theta_{2}^{1}\theta_{2}%
^{1}\theta_{2}^{1}\,;\gamma(\theta_{1}^{2}{\theta_{1}^{1}};\theta_{1}%
^{2})\theta_{1}^{3}+\theta_{1}^{3}\gamma({\theta_{1}^{1}}\theta_{1}^{2}%
;\theta_{1}^{2}))\vspace{1mm}\newline\\
13|2 & \leftrightarrow & _{11}^{12} & = & \gamma(\theta_{2}^{1}\theta_{2}%
^{2}\,;\theta_{1}^{2}\theta_{1}^{2})\vspace{1mm}\newline\\
3|12 & \leftrightarrow & _{2}^{12} & = & \gamma(\theta_{1}^{1}\theta_{1}%
^{2}\,;\theta_{2}^{2})\vspace{1mm}\newline\\
12|3 & \leftrightarrow & _{11}^{21} & = & \gamma(\theta_{2}^{2}\theta_{2}%
^{1}\,;\theta_{1}^{2}\theta_{1}^{2})\vspace{1mm}\newline\\
2|13 & \leftrightarrow & _{2}^{21} & = & \gamma(\theta_{1}^{2}\theta_{1}%
^{1}\,;\theta_{2}^{2})\vspace{1mm}\newline\\
23|1 & \leftrightarrow & _{2}^{3} & = & \gamma(\theta_{1}^{3}\,;\theta_{2}%
^{1})\\
&  &  &  &
\end{array}
$

\noindent\underline{For $KK_{2,3}$:}

\hspace{-0.1in}%
$
\begin{array}
[c]{lllll}%
1|23 & \leftrightarrow & _{21}^{11} & = & \gamma(\theta_{2}^{1}\theta_{2}%
^{1}\,;\theta_{2}^{2}\theta_{1}^{2})\vspace{1mm}\newline\\
13|2 & \leftrightarrow & _{21}^{2} & = & \gamma(\theta_{2}^{2}\,;\theta
_{2}^{1}\theta_{1}^{1})\vspace{1mm}\newline\\
3|12 & \leftrightarrow & _{2}^{3} & = & \gamma(\theta_{1}^{2}\,;\theta_{3}%
^{1})\vspace{1mm}\newline\\
12|3 & \leftrightarrow & _{111}^{11} & = & \gamma(\gamma(\theta_{2}%
^{1}\,;\theta_{2}^{1}\theta_{1}^{1})\theta_{3}^{1}+\theta_{3}^{1}\gamma
(\theta_{2}^{1}\,;\theta_{1}^{1}\theta_{2}^{1})\,;\theta_{1}^{2}\theta_{1}%
^{2}\theta_{1}^{2})\vspace{1mm}\newline\\
2|13 & \leftrightarrow & _{12}^{11} & = & \gamma(\theta_{2}^{1}\theta_{2}%
^{1}\,;\theta_{1}^{2}\theta_{2}^{2})\vspace{1mm}\newline\\
23|1 & \leftrightarrow & _{12}^{2} & = & \gamma(\theta_{2}^{2}\,;\theta
_{1}^{1}\theta_{2}^{1})\\
&  &  &  &
\end{array}
$
\vspace{0.1in}

\hspace*{-0.5in}\unitlength=1.00mm \special{em:linewidth 0.4pt}
\linethickness{0.4pt} \begin{picture}(88.33,43.67)
\put(6.67,40.67){\line(1,0){48.33}} \put(55.00,40.67){\line(0,-1){31.67}}
\put(55.00,9.00){\line(-1,0){48.33}} \put(6.67,9.00){\line(0,1){31.67}}
\put(6.67,40.67){\circle*{1.33}} \put(6.67,25.00){\circle*{1.33}}
\put(6.67,9.00){\circle*{1.33}} \put(55.00,40.67){\circle*{1.33}}
\put(55.00,9.00){\circle*{1.33}} \put(55.00,25.00){\circle*{1.33}}
\put(6.67,17.00){\circle*{0.67}} \put(0.33,33.33){\makebox(0,0)[cc]{$13|2$}}
\put(10.33,33.33){\makebox(0,0)[cc]{$_{11}^{12}$}}
\put(0.00,17.00){\makebox(0,0)[cc]{$1|23$}}
\put(10.80,17.00){\makebox(0,0)[cc]{$_{11}^{111}$}}
\put(29.00,5.00){\makebox(0,0)[cc]{$12|3$}}
\put(29.00,13.00){\makebox(0,0)[cc]{$_{11}^{21}$}}
\put(30.00,43.67){\makebox(0,0)[cc]{$3|12$}}
\put(30.00,35.67){\makebox(0,0)[cc]{$_{2}^{12}$}}
\put(61.33,17.00){\makebox(0,0)[cc]{$2|13$}}
\put(52.33,17.00){\makebox(0,0)[cc]{$_{2}^{21}$}}
\put(61.33,32.67){\makebox(0,0)[cc]{$23|1$}}
\put(52.33,32.67){\makebox(0,0)[cc]{$_{2}^{3}$}} \put(80.67,40.67){\line(1,0){48.33}}
\put(129.00,40.67){\line(0,-1){31.67}} \put(129.00,9.00){\line(-1,0){48.33}}
\put(80.67,9.00){\line(0,1){31.67}} \put(80.67,40.67){\circle*{1.33}}
\put(80.67,25.00){\circle*{1.33}} \put(80.67,9.00){\circle*{1.33}}
\put(129.00,40.67){\circle*{1.33}} \put(129.00,9.00){\circle*{1.33}}
\put(129.00,25.00){\circle*{1.33}} \put(74.33,33.33){\makebox(0,0)[cc]{$13|2$}}
\put(84.33,33.33){\makebox(0,0)[cc]{$_{21}^{2}$}}
\put(74.00,17.00){\makebox(0,0)[cc]{$1|23$}}
\put(84.80,17.00){\makebox(0,0)[cc]{$_{21}^{11}$}}
\put(103.00,5.00){\makebox(0,0)[cc]{$12|3$}}
\put(103.00,13.00){\makebox(0,0)[cc]{$_{111}^{11}$}}
\put(104.00,43.67){\makebox(0,0)[cc]{$3|12$}}
\put(104.00,35.67){\makebox(0,0)[cc]{$_{3}^{2}$}}
\put(135.33,17.00){\makebox(0,0)[cc]{$2|13$}}
\put(126.33,17.00){\makebox(0,0)[cc]{$_{12}^{11}$}}
\put(135.33,32.67){\makebox(0,0)[cc]{$23|1$}}
\put(126.33,32.67){\makebox(0,0)[cc]{$_{12}^{2}$}} \put(104.33,9.00){\circle*{0.67}}
\end{picture}
\begin{center}
Figure 19. The biassociahedra $KK_{3,2}$ and $KK_{2,3}$ (heptagons)
\end{center}
\pagebreak

\noindent\underline{For $KK_{3,3}$:}

$
\begin{array}
[c]{lllll}%
1|234 & \leftrightarrow & _{21}^{111} & = & \gamma(\theta_{2}^{1}\theta
_{2}^{1}\theta_{2}^{1}\,;\theta_{2}^{3}\gamma(\theta_{1}^{1}\theta_{1}%
^{2}\,;\theta_{1}^{2})+a+b),\text{ where}\\
&  &  &  & \ \ \ \ a+b=\gamma(\theta_{2}^{2}\theta_{2}^{1}\,;\theta_{1}%
^{2}\theta_{1}^{2})\theta_{1}^{3}+\gamma(\theta_{1}^{2}\theta_{1}^{1}%
\,;\theta_{2}^{2})\theta_{1}^{3}\vspace{1mm}\newline\\
123|4 & \leftrightarrow & _{111}^{21} & = & \gamma(c+d+\theta_{3}^{2}%
\gamma(\theta_{2}^{1}\,;\theta_{1}^{1}\theta_{2}^{1})\,;\theta_{1}^{2}%
\theta_{1}^{2}\theta_{1}^{2}),\text{ where}\\
&  &  &  & \ \ \ \ c+d=\gamma(\theta_{2}^{1}\theta_{2}^{1}\,;\theta_{2}%
^{2}\theta_{1}^{2})\theta_{3}^{1}+\gamma(\theta_{2}^{2}\,;\theta_{2}^{1}%
\theta_{1}^{1})\theta_{3}^{1}\vspace{1mm}\newline\\
2|134 & \leftrightarrow & _{12}^{111} & = & \gamma(\theta_{2}^{1}\theta
_{2}^{1}\theta_{2}^{1}\,;\gamma(\theta_{1}^{2}\theta_{1}^{1}\,;\theta_{1}%
^{2})\theta_{2}^{3}+e+f),\text{ where}\\
&  &  &  & \ \ \ \ e+f=\theta_{1}^{3}\gamma(\theta_{2}^{1}\theta_{2}%
^{2}\,;\theta_{1}^{2}\theta_{1}^{2})+\theta_{1}^{3}\gamma(\theta_{1}^{1}%
\theta_{1}^{2}\,;\theta_{2}^{2})\vspace{1mm}\newline\\
124|3 & \leftrightarrow & _{111}^{12} & = & \gamma(g+h+\gamma(\theta_{2}%
^{1}\,;\theta_{2}^{1}\theta_{1}^{1})\theta_{3}^{2}\,;\theta_{1}^{2}\theta
_{1}^{2}\theta_{1}^{2}),\text{ where}\\
&  &  &  & \ \ \ \ g+h=\theta_{3}^{1}\gamma(\theta_{2}^{1}\theta_{2}%
^{1}\,;\theta_{1}^{2}\theta_{2}^{2})+\theta_{3}^{1}\gamma(\theta_{2}%
^{2}\,;\theta_{1}^{1}\theta_{2}^{1})\vspace{1mm}\newline\\
134|2 & \leftrightarrow & _{21}^{3} & = & \gamma(\theta_{2}^{3}\,;\theta
_{2}^{1}\theta_{1}^{1})\vspace{1mm}\newline\\
234|1 & \leftrightarrow & _{12}^{3} & = & \gamma(\theta_{2}^{3}\,;\theta
_{1}^{1}\theta_{2}^{1})\vspace{1mm}\newline\\
3|124 & \leftrightarrow & _{3}^{21} & = & \gamma(\theta_{1}^{2}\theta_{1}%
^{1}\,;\theta_{3}^{2})\vspace{1mm}\newline\\
4|123 & \leftrightarrow & _{3}^{12} & = & \gamma(\theta_{1}^{1}\theta_{1}%
^{2}\,;\theta_{3}^{2})\vspace{1mm}\newline\\
13|24 & \leftrightarrow & _{21}^{21} & = & \gamma(\theta_{2}^{2}\theta_{2}%
^{1}\,;\theta_{2}^{2}\theta_{1}^{2})\vspace{1mm}\newline\\
24|13 & \leftrightarrow & _{12}^{12} & = & \gamma(\theta_{2}^{1}\theta_{2}%
^{2}\,;\theta_{1}^{2}\theta_{2}^{2})\vspace{1mm}\newline\\
14|23 & \leftrightarrow & _{21}^{12} & = & \gamma(\theta_{2}^{1}\theta_{2}%
^{2}\,;\theta_{2}^{2}\theta_{1}^{2})\vspace{1mm}\newline\\
23|14 & \leftrightarrow & _{12}^{21} & = & \gamma(\theta_{2}^{2}\theta_{2}%
^{1}\,;\theta_{1}^{2}\theta_{2}^{2})\vspace{1mm}\newline\\
34|12 & \leftrightarrow & _{3}^{3} & = & \gamma(\theta_{1}^{3}\,;\theta
_{3}^{1})\\
12|34 & \leftrightarrow & _{111}^{111} & = & \gamma\lbrack\theta_{3}^{1}%
\gamma(\theta_{2}^{1}\,;\theta_{1}^{1}\theta_{2}^{1})\gamma(\theta_{2}%
^{1}\,;\theta_{1}^{1}\theta_{2}^{1})+\gamma(\theta_{2}^{1}\,;\theta_{2}%
^{1}\theta_{1}^{1})\theta_{3}^{1}\gamma(\theta_{2}^{1}\,;\theta_{1}^{1}%
\theta_{2}^{1})\\
&  &  &  & \ \ \ \ \ \ +\gamma(\theta_{2}^{1}\,;\theta_{2}^{1}\theta_{1}%
^{1})\gamma(\theta_{2}^{1}\,;\theta_{2}^{1}\theta_{1}^{1})\theta_{3}%
^{1}\,;\vspace{1mm}\newline\\
&  &  &  & \ \ \ \theta_{1}^{3}\gamma(\theta_{1}^{1}\theta_{1}^{2}%
\,;\theta_{1}^{2})\gamma(\theta_{1}^{1}\theta_{1}^{2}\,;\theta_{1}^{2}%
)+\gamma(\theta_{1}^{2}\theta_{1}^{1}\,;\theta_{1}^{2})\theta_{1}^{3}%
\gamma(\theta_{1}^{1}\theta_{1}^{2}\,;\theta_{1}^{2})+\\
&  &  &  & \ \ \ \ \ \ +\gamma(\theta_{1}^{2}\theta_{1}^{1}\,;\theta_{1}%
^{2})\gamma(\theta_{1}^{2} \theta_{1}^{1}\,;\theta_{1}^{2})\theta_{1}^{3}]\\
&  &  &  &
\end{array}
$
\vspace{0.2in}

\hspace*{-0.5in}\unitlength=1.00mm \special{em:linewidth 0.4pt}
\linethickness{0.4pt} \begin{picture}(98.33,76.99)
\put(67.00,5.33){\line(1,0){30.33}} \put(97.33,15.33){\line(-1,0){30.33}}
\put(97.33,25.33){\line(-1,0){30.33}} \put(97.33,35.33){\line(-1,0){30.33}}
\put(77.00,5.33){\line(0,1){30.00}} \put(67.33,29.33){\line(1,0){4.67}}
\put(72.00,35.00){\line(0,-1){5.67}} \put(97.67,55.33){\line(0,-1){50.00}}
\put(66.67,5.33){\line(0,1){49.67}} \put(66.67,55.00){\line(1,0){31.00}}
\put(87.00,5.33){\line(0,1){49.67}} \put(97.67,45.33){\line(-1,0){10.67}}
\put(92.00,40.00){\line(1,1){5.67}} \put(92.33,40.00){\line(0,1){0.00}}
\put(93.67,35.66){\line(-2,5){2.00}} \put(97.67,55.00){\line(-6,5){25.67}}
\put(66.67,5.33){\line(-6,5){28.67}} \put(38.33,29.33){\line(0,1){47.00}}
\put(66.67,55.00){\line(-4,3){28.33}} \put(38.33,76.33){\line(1,0){33.67}}
\put(72.00,35.66){\line(0,1){19.00}} \put(72.00,55.33){\line(0,1){21.00}}
\put(77.00,5.33){\line(-5,4){10.00}} \put(79.33,70.66){\line(0,-1){15.33}}
\put(79.33,54.66){\line(0,-1){19.00}} \put(79.33,35.00){\line(0,-1){9.33}}
\put(79.33,25.00){\line(0,-1){2.33}} \put(79.33,44.33){\line(-1,1){7.33}}
\put(66.67,45.00){\line(-6,5){20.67}} \put(46.00,45.33){\line(-1,1){7.67}}
\put(46.00,53.66){\line(1,-4){6.67}} \put(46.00,22.66){\line(3,2){6.67}}
\put(38.33,29.33){\line(1,0){7.33}} \put(46.33,29.33){\line(1,0){20.00}}
\put(66.67,15.33){\line(-6,5){14.00}} \put(66.00,14.00){\line(-6,5){13.67}}
\put(51.00,26.33){\line(-6,5){3.67}} \put(38.33,70.66){\line(1,0){6.33}}
\put(47.67,70.66){\line(1,0){24.33}} \put(46.00,22.66){\line(0,1){48.00}}
\put(46.00,70.66){\circle*{0.67}} \put(46.00,62.33){\circle*{1.33}}
\put(46.00,70.66){\circle*{1.33}} \put(38.33,76.33){\circle*{1.33}}
\put(38.33,70.66){\circle*{1.33}} \put(72.00,76.33){\circle*{1.33}}
\put(72.00,70.66){\circle*{1.33}} \put(79.33,70.33){\circle*{1.33}}
\put(79.33,61.00){\circle*{1.33}} \put(97.67,55.00){\circle*{1.33}}
\put(97.67,35.33){\circle*{1.33}} \put(66.67,55.00){\circle*{1.33}}
\put(66.67,45.33){\circle*{1.33}} \put(72.00,51.33){\circle*{1.33}}
\put(79.33,44.66){\circle*{1.33}} \put(72.00,29.33){\circle*{0.67}}
\put(72.00,29.33){\circle*{1.33}} \put(79.33,23.00){\circle*{1.33}}
\put(97.67,45.33){\circle*{1.33}} \put(66.67,35.33){\circle*{1.33}}
\put(66.67,5.33){\circle*{1.33}} \put(97.67,5.33){\circle*{1.33}}
\put(46.00,45.66){\circle*{1.33}} \put(38.33,52.66){\circle*{1.33}}
\put(38.33,29.33){\circle*{1.33}} \put(46.00,22.66){\circle*{1.33}}
\put(66.67,25.33){\circle*{0.67}} \put(66.67,15.33){\circle*{0.67}}
\put(77.00,35.33){\circle*{0.00}} \put(77.00,25.33){\circle*{0.67}}
\put(77.00,15.33){\circle*{0.00}} \put(77.00,5.33){\circle*{0.67}}
\put(87.00,35.33){\circle*{0.67}} \put(87.00,25.33){\circle*{0.67}}
\put(87.00,15.33){\circle*{0.67}} \put(87.00,5.33){\circle*{0.67}}
\put(97.67,25.33){\circle*{0.67}} \put(97.67,15.33){\circle*{0.67}}
\put(92.00,40.33){\circle*{0.67}} \put(87.00,45.33){\circle*{0.67}}
\put(79.33,52.00){\circle*{0.67}} \put(47.67,29.33){\circle*{0.67}}
\put(87.00,55.00){\circle*{0.67}} \put(46.00,22.33){\line(1,0){9.67}}
\put(55.67,22.33){\circle*{0.67}} \put(52.67,27.00){\circle*{0.67}}
\put(77.00,35.33){\circle*{0.67}} \put(77.00,15.33){\circle*{0.67}}
\put(58.00,18.66){\makebox(0,0)[cc]{$a$}} \put(48.33,33.00){\makebox(0,0)[cc]{$b$}}
\put(44.00,27.33){\makebox(0,0)[cc]{$d$}} \put(63.67,12.00){\makebox(0,0)[cc]{$c$}}
\put(95.33,39.33){\makebox(0,0)[cc]{$e$}} \put(84.00,52.66){\makebox(0,0)[cc]{$f$}}
\put(90.00,39.33){\makebox(0,0)[cc]{$g$}} \put(93.67,51.00){\makebox(0,0)[cc]{$h$}}
\put(87.67,44.66){\line(1,-1){4.33}} \put(79.33,52.00){\line(6,-5){7.00}}
\put(87.67,54.33){\line(6,-5){10.00}} \put(86.33,55.66){\line(-4,3){7.00}}
\put(46.00,53.00){\circle*{0.67}} \put(88.00,14.33){\line(1,-1){9.67}}
\put(86.33,16.00){\line(-1,1){7.00}} \put(79.33,23.00){\line(-6,5){7.33}}
\put(94.00,34.67){\line(2,-5){3.67}}
\end{picture}
\vspace{0.1in}
\begin{center}
Figure 20. The biassociahedron $KK_{3,3}$ (a subdivision of $P_{4}$).
\end{center}

\noindent\underline{For $KK_{4,2}$:}

$%
\begin{array}
[c]{lllll}%
123|4 & \leftrightarrow & _{11}^{31} & = & \gamma(\theta_{2}^{3}\theta_{2}%
^{1}\,;\theta_{1}^{2}\theta_{1}^{2})\vspace{1mm}\newline\\
2|134 & \leftrightarrow & _{2}^{211} & = & \gamma(\theta_{1}^{2}\theta_{1}%
^{1}\theta_{1}^{1}\,;\theta_{2}^{3})\vspace{1mm}\newline\\
124|3 & \leftrightarrow & _{11}^{22} & = & \gamma(\theta_{2}^{2}\theta_{2}%
^{2}\,;\theta_{1}^{2}\theta_{1}^{2})\vspace{1mm}\newline\\
134|2 & \leftrightarrow & _{11}^{13} & = & \gamma(\theta_{2}^{1}\theta_{2}%
^{3}\,;\theta_{1}^{2}\theta_{1}^{2})\vspace{1mm}\newline\\
234|1 & \leftrightarrow & _{2}^{4} & = & \gamma(\theta_{1}^{4}\,;\theta
_{2}^{1})\vspace{1mm}\newline\\
3|124 & \leftrightarrow & _{2}^{121} & = & \gamma(\theta_{1}^{1}\theta_{1}%
^{2}\theta_{1}^{1}\,;\theta_{2}^{3})\vspace{1mm}\newline\\
4|123 & \leftrightarrow & _{2}^{112} & = & \gamma(\theta_{1}^{1}\theta_{1}%
^{1}\theta_{1}^{2}\,;\theta_{2}^{3})\vspace{1mm}\newline\\
13|24 & \leftrightarrow & _{11}^{121} & = & \gamma\lbrack\theta_{2}^{1}%
\theta_{2}^{2}\theta_{2}^{1}\,;\theta_{1}^{3}\gamma(\theta_{1}^{1}\theta
_{1}^{2}\,;\theta_{1}^{2})+\gamma(\theta_{1}^{2}\theta_{1}^{1}\,;\theta
_{1}^{2})\theta_{1}^{3}]\vspace{1mm}\newline\\
14|23 & \leftrightarrow & _{11}^{112} & = & \gamma\lbrack\theta_{2}^{1}%
\theta_{2}^{1}\theta_{2}^{2}\,;\theta_{1}^{3}\gamma(\theta_{1}^{1}\theta
_{1}^{2}\,;\theta_{1}^{2})+\gamma(\theta_{1}^{2}\theta_{1}^{1}\,;\theta
_{1}^{2})\theta_{1}^{3}]\vspace{1mm}\newline\\
23|14 & \leftrightarrow & _{2}^{31} & = & \gamma(\theta_{1}^{3}\theta_{1}%
^{1}\,;\theta_{2}^{2})\vspace{1mm}\newline\\
34|12 & \leftrightarrow & _{2}^{13} & = & \gamma(\theta_{1}^{1}\theta_{1}%
^{3}\,;\theta_{2}^{2})\vspace{1mm}\newline\\
12|34 & \leftrightarrow & _{11}^{211} & = & \gamma\lbrack\theta_{2}^{2}%
\theta_{2}^{1}\theta_{2}^{1}\,;\theta_{1}^{3}\gamma(\theta_{1}^{1}\theta
_{1}^{2}\,;\theta_{1}^{2})+\gamma(\theta_{1}^{2}\theta_{1}^{1}\,;\theta
_{1}^{2})\theta_{1}^{3}]\vspace{1mm}\newline\\
1|234 & \leftrightarrow & _{11}^{1111} & = & \gamma(\theta_{2}^{1}\theta
_{2}^{1}\theta_{2}^{1}\theta_{2}^{1}\,;a+b-c+d+e+f),\text{ where}\\
&  &  &  & \ \ \ a=\theta_{1}^{4}\gamma(\theta_{1}^{1}\gamma(\theta_{1}%
^{1}\theta_{1}^{2}\,;\theta_{1}^{2})\,;\theta_{1}^{2})\\
&  &  &  & \ \ \ b=\gamma(\gamma(\theta_{1}^{2}\theta_{1}^{1}\,;\theta_{1}%
^{2})\theta_{1}^{1}\,;\theta_{1}^{2})\theta_{1}^{4}\\
&  &  &  & \ \ \ c=\gamma(\theta_{1}^{2}\theta_{1}^{1}\theta_{1}^{1}%
\,;\theta_{1}^{3})\,\gamma(\theta_{1}^{1}\theta_{1}^{1}\theta_{1}^{2}%
\,;\theta_{1}^{3})\\
&  &  &  & \ \ \ d=\gamma(\theta_{1}^{3}\theta_{1}^{1}\,;\theta_{1}%
^{2})\,\gamma(\theta_{1}^{1}\theta_{1}^{2}\theta_{1}^{1}\,;\theta_{1}^{3})\\
&  &  &  & \ \ \ e=\gamma(\theta_{1}^{3}\theta_{1}^{1}\,;\theta_{1}%
^{2})\,\gamma(\theta_{1}^{1}\theta_{1}^{3}\,;\theta_{1}^{2})\\
&  &  &  & \ \ \ f=\gamma(\theta_{1}^{1}\theta_{1}^{2}\theta_{1}^{1}%
\,;\theta_{1}^{3})\,\gamma(\theta_{1}^{1}\theta_{1}^{3}\,;\theta_{1}^{2})\\
&  &  &  &
\end{array}
$

\unitlength=0.75mm\special{em:linewidth 0.4pt} \linethickness{0.4pt}
\begin{picture}(98.33,76.99)
\put(67.00,5.33){\line(1,0){30.33}} \put(67.33,29.33){\line(1,0){4.67}}
\put(97.67,55.33){\line(0,-1){50.00}} \put(66.67,5.33){\line(0,1){49.67}}
\put(66.67,55.00){\line(1,0){31.00}} \put(97.67,55.00){\line(-6,5){25.67}}
\put(66.67,5.33){\line(-6,5){28.67}} \put(38.33,29.33){\line(0,1){47.00}}
\put(66.67,55.00){\line(-4,3){28.33}} \put(38.33,76.33){\line(1,0){33.67}}
\put(72.00,55.33){\line(0,1){21.00}} \put(79.33,70.66){\line(0,-1){15.33}}
\put(79.33,44.33){\line(-1,1){7.33}} \put(46.00,45.33){\line(-1,1){7.67}}
\put(38.33,29.33){\line(1,0){7.33}} \put(46.33,29.33){\line(1,0){20.00}}
\put(38.33,70.66){\line(1,0){6.33}} \put(47.67,70.66){\line(1,0){24.33}}
\put(46.00,22.66){\line(0,1){48.00}} \put(46.00,70.66){\circle*{0.67}}
\put(46.00,70.66){\circle*{1.33}} \put(38.33,76.33){\circle*{1.33}}
\put(38.33,70.66){\circle*{1.33}} \put(72.00,76.33){\circle*{1.33}}
\put(72.00,70.66){\circle*{1.33}} \put(79.33,70.33){\circle*{1.33}}
\put(97.67,55.00){\circle*{1.33}} \put(66.67,55.00){\circle*{1.33}}
\put(72.00,51.33){\circle*{1.33}} \put(79.33,44.66){\circle*{1.33}}
\put(72.00,29.33){\circle*{0.67}} \put(72.00,29.33){\circle*{1.33}}
\put(79.33,23.00){\circle*{1.33}} \put(66.67,5.33){\circle*{1.33}}
\put(97.67,5.33){\circle*{1.33}} \put(46.00,45.66){\circle*{1.33}}
\put(38.33,52.66){\circle*{1.33}} \put(38.33,29.33){\circle*{1.33}}
\put(46.00,22.66){\circle*{1.33}} \put(66.67,20.33){\line(1,0){31.00}}
\put(56.00,14.33){\line(0,1){48.67}} \put(66.67,20.33){\line(-4,3){20.67}}
\put(82.33,20.00){\line(1,-1){15.33}} \put(81.67,20.67){\line(-6,5){9.67}}
\put(56.00,14.00){\circle*{0.67}} \put(56.00,20.67){\circle*{0.67}}
\put(56.00,28.00){\circle*{0.67}} \put(56.00,62.67){\circle*{0.67}}
\put(46.00,35.67){\line(-6,5){7.67}} \put(46.00,54.33){\line(1,-4){5.67}}
\put(97.67,20.33){\circle*{0.67}} \put(66.67,20.33){\circle*{0.67}}
\put(38.33,41.67){\circle*{0.67}} \put(51.33,32.00){\circle*{0.67}}
\put(46.00,36.00){\line(2,-3){10.00}} \put(46.00,36.00){\circle*{0.67}}
\put(46.00,53.67){\circle*{0.67}} \put(66.67,42.67){\line(1,0){31.00}}
\put(66.67,42.67){\line(-6,5){20.67}} \put(46.00,59.67){\circle*{1.33}}
\put(66.67,42.67){\circle*{1.33}} \put(97.67,42.67){\circle*{1.33}}
\put(56.00,51.67){\circle*{0.67}} \put(79.33,43.00){\line(0,1){11.67}}
\put(79.33,22.67){\line(0,1){19.67}} \put(72.00,29.33){\line(0,1){13.00}}
\put(72.00,43.00){\line(0,1){11.67}}
\end{picture}\vspace*{-0.1in}

\begin{center}
Figure 21. The biassociahedron $KK_{4,2}$ (a subdivision of $J_{4}%
=K_{4,2}=\vartheta_{3,1}(P_{4})$).
\end{center}

\unitlength=0.80mm \special{em:linewidth 0.4pt} \linethickness{0.4pt}
\begin{picture}(98.33,76.99)
\put(67.00,5.33){\line(1,0){30.33}} \put(97.67,55.33){\line(0,-1){50.00}}
\put(66.67,5.33){\line(0,1){49.67}} \put(66.67,55.00){\line(1,0){31.00}}
\put(66.67,5.33){\line(-6,5){28.67}} \put(38.33,29.33){\line(0,1){47.00}}
\put(66.67,55.00){\line(-4,3){28.33}} \put(38.33,76.33){\line(1,0){33.67}}
\put(72.00,55.33){\line(0,1){21.00}} \put(38.33,76.33){\circle*{1.33}}
\put(72.00,76.33){\circle*{1.33}} \put(97.67,55.00){\circle*{1.33}}
\put(66.67,55.00){\circle*{1.33}} \put(66.67,5.33){\circle*{1.33}}
\put(97.67,5.33){\circle*{1.33}} \put(38.33,29.33){\circle*{1.33}}
\put(57.33,62.33){\line(0,-1){49.33}} \put(38.33,29.00){\line(1,0){18.33}}
\put(58.00,29.00){\line(1,0){7.67}} \put(67.33,29.00){\line(1,0){4.67}}
\put(72.00,29.00){\line(0,1){25.33}} \put(66.67,19.00){\line(1,0){31.00}}
\put(72.33,29.00){\line(5,-4){12.00}} \put(85.67,18.33){\line(1,-1){12.33}}
\put(57.33,52.67){\line(4,-3){9.33}} \put(59.67,61.00){\line(1,0){12.33}}
\put(56.67,61.00){\line(-1,0){18.33}} \put(57.33,62.00){\circle*{1.33}}
\put(57.33,52.33){\circle*{1.33}} \put(66.67,46.00){\circle*{1.33}}
\put(97.67,19.00){\circle*{1.33}} \put(66.67,19.00){\circle*{1.33}}
\put(57.33,13.33){\circle*{1.33}} \put(72.00,29.00){\circle*{1.33}}
\put(38.33,61.00){\circle*{0.67}} \put(72.00,61.00){\circle*{0.67}}
\put(72.00,61.00){\circle*{0.67}} \put(82.67,20.33){\line(-1,0){15.33}} \
\put(66.00,20.33){\line(-1,0){8.00}} \put(56.67,20.33){\line(-1,0){7.67}}
\put(59.33,29.00){\line(5,-4){6.67}} \put(67.33,22.33){\line(3,-2){4.67}}
\put(57.33,33.00){\line(-3,2){19.00}} \ \put(49.00,20.33){\line(0,1){18.33}}
\put(86.33,5.33){\line(0,1){13.67}} \ \put(62.00,9.00){\line(1,5){2.33}}
\put(49.00,20.33){\line(2,1){8.00}} \put(58.00,24.33){\line(2,1){4.67}}
\put(88.67,34.33){\line(-3,2){16.67}} \ \put(88.33,34.33){\line(1,-3){5.00}} \ \
\put(93.67,18.67){\line(1,-3){4.33}} \put(83.33,20.00){\line(0,1){18.00}} \ \
\put(57.33,33.67){\circle*{1.33}} \put(38.33,45.33){\circle*{1.33}}
\put(38.33,61.00){\circle*{1.33}} \put(72.00,61.00){\circle*{1.33}}
\put(72.00,45.00){\circle*{1.33}} \put(88.33,34.33){\circle*{1.33}}
\put(86.33,5.33){\circle*{0.67}} \put(49.00,20.33){\circle*{0.67}}
\put(49.00,38.67){\circle*{0.67}} \put(83.33,37.67){\circle*{0.67}} \
\put(83.33,20.33){\circle*{0.67}} \ \put(70.30,20.33){\circle*{0.67}}
\put(64.33,20.33){\circle*{0.67}} \put(62.00,9.33){\circle*{0.67}}
\put(62.33,26.67){\circle*{0.67}} \put(59.67,29.00){\circle*{0.67}}
\put(72.33,76.00){\line(6,-5){25.33}} \put(86.33,5.33){\line(-1,1){13.33}}
\put(86.33,19.00){\circle*{0.67}} \put(97.67,46.00){\circle*{1.33}}
\put(88.33,34.33){\line(4,5){9.33}}
\end{picture}\vspace*{-0.1in}

\begin{center}
Figure 22. The biassociahedron $KK_{2,4}$ (a subdivision of $J_{4}%
=K_{2,4}=\vartheta_{1,3}(P_{4})$).
\end{center}

\end{document}